\documentclass[11pt,twoside,bezier]{article}
\usepackage{amssymb,amsthm}
\newcommand{\pr}{\pageref}
\newcommand{\iv}{^{-1}}
\newcommand{\vk}{ van Kampen }
\newcommand{\ct}{ contiguity }
\newcommand{\topp}{{\bf top}}
\newcommand{\hra}{H_{ra}}
\newcommand{\nee}{\not\equiv}
\newcommand{\maa}{H_a}
\newcommand{\bott}{{\bf bot}}
\newcommand{\hnka}{H^*_{ka}}
\newcommand{\hkra}{H_{kra}}

\newcommand{\oo}{{\bf 0}}

\newcommand{\la}{\langle}
\newcommand{\ra}{\rangle}

\newcommand{\bc}{{\bf c}}
\newcommand{\bbc}[1]{\overline{{\bf c}_{#1}}}

\newcommand{\bk}{{\bf \rho}}
\newcommand{\bt}[1]{\overline{\tau(#1)}}
\newcommand{\tk}{{\bf \lambda}}
\newcommand{\rt}{\hbox{root}}
\newcommand{\kk}{\overline{\bf K}}

\newcommand{\lk}{\overleftarrow{{\bf \kappa}}}
\newcommand{\rk}{\overrightarrow{{\bf \kappa}}}
\newcommand{\tool}{\stackrel{\ell}{\too} }
\newcommand{\aaa}{{\bf A}}
\newcommand{\bkk}{{\bf K}}
\newcommand{\rr}{{\bf R}}
\newcommand{\ttt}{{\cal T}}

\newcommand{\Lab}{\phi}

\newcommand{\sss}{{\mathbb{S}}}
\newcommand{\too}{\to}

\newcommand{\xxx}{{\bf X}}
\newcommand{\bb}{{\cal B}}

\newcommand{\tz}{Z(\sss,ra)}
\newcommand{\tza}{Z(\sss,\Lambda,ra)}
\newcounter{ppp}
\newcounter{pdten}
\newcounter{pdeleven}
\newcounter{pdtwelve}

\begin{document}
\theoremstyle{plain}
\newtheorem{theo}{Theorem}[section]
\newtheorem{lm}[theo]{Lemma}
\newtheorem{cy}[theo]{Corollary}
\newtheorem{df}[theo]{Definition}
\newtheorem{remark}[theo]{Remark}
\newtheorem{prop}[theo]{Proposition}
\newtheorem{prob}[theo]{Problem}

\title{Non-amenable finitely presented torsion-by-cyclic groups}
 \author{A.Yu. Ol'shanskii, M.V. Sapir\thanks{Both authors were supported in part by
the NSF grant DMS 0072307. In addition, the research of the first
author was supported in part by the Russian Fund for Basic
Research 99-01-00894 and by the INTAS grant,  the research of the
second author was supported in part by the NSF grant DMS
9978802.}}
\date{}
\maketitle

\begin{abstract} We construct a finitely presented non-amenable
group without free non-cyclic subgroups thus providing a finitely
presented counterexample to von Neumann's problem. Our group is an
extension of a group of finite exponent $n>>1$ by a cyclic group,
so it satisfies the identity $[x,y]^n=1$.
\end{abstract}

\tableofcontents

\section{Introduction}

\subsection{Short history of the problem} Hausdorff \cite{Haus} proved in 1914 that one can
subdivide the 2-sphere minus a countable set of points into 3
parts $A$, $B$, $C$, such that each of these three parts can be
obtained from each of the other two parts by a rotation, and the
union of two of these parts can be obtained by rotating the third
part. This implied that one cannot define a finitely additive
measure on the 2-sphere which is invariant under the group
$SO(3)$. In 1924 Banach and Tarski \cite{BT} generalized
Hausdorff's result by proving, in particular, that in
$\mathbb{R}^3$, every two bounded sets $A, B$ with non-empty
interiors can be decomposed $A=\bigcup_{i=1}^n A_i$,
$B=\bigcup_{i=1}^n B_i$ such that $A_i$ can be rotated to $B_i$,
$i=1,...,n$ (the so called Banach-Tarski paradox). Von Neumann
\cite{vN} was first who noticed that the cause of the
Banach-Tarski paradox is not the geometry of $\mathbb{R}^3$ but an
algebraic property of the group $SO(3)$. He introduced the concept
of an amenable group (he called such groups ``measurable") as a
group $G$ which has a left invariant finitely additive
measure\footnote{Later amenable groups were characterized in many
different ways. There exist a geometric characterization by F\o
lner \cite{Folner}, a characterization by Kesten \cite{Kesten2} in
terms of random walks, and a combinatorial characterization by
Grigorchuk \cite{Grigorchuk} among others.} $\mu$, $\mu(G)=1$,
noticed that if a group is amenable then any set it acts upon
freely also has an invariant measure and proved that a group is
not amenable provided it contains a free non-abelian subgroup. He
also showed that groups like $PSL(2, \mathbb{Z})$,
$SL(2,\mathbb{Z})$ contain free non-abelian subgroups. So analogs
of Banach-Tarski paradox can be found in $\mathbb{R}^2$ and even
$\mathbb{R}$. Von Neumann showed that the class of amenable groups
contains abelian groups, finite groups and is closed under taking
subgroups, extensions, and infinite unions of increasing sequences
of groups. Day and Specht showed that this class is closed under
homomorphic images. The class of groups without free non-abelian
subgroups  is also closed under these operations and contains
abelian and finite groups.

The problem of existence of non-amenable groups without
non-abelian free subgroups probably goes back to von Neumann and
became known as the ``von Neumann problem" in the fifties. As far
as we know, the first paper where this problem was formulated was
the paper by Day \cite{day}. It is also mentioned in the
monograph by Greenleaf \cite{Greenleaf} based on his lectures
given in Berkeley in 1967. Tits \cite{Tits} proved that every
non-amenable matrix group over a field of characteristic $0$
contains a non-abelian free subgroup. In particular every
semisimple Lie group over a field of characteristic $0$ contains
such a subgroup (see also Vershik \cite{Vershik}).

First counterexamples to the von Neumann problem were constructed
by Ol'shanskii \cite{OlAmen}. He proved that the groups with all
proper subgroups cyclic constructed by him, both torsion-free
\cite{OlTar} and torsion \cite{OlTar1} (the so called ``Tarski
monsters"), are not amenable. Later Adian \cite{Adian} showed that
the free Burnside group of odd exponent $n\ge 665$ with at least
two generators (that is the group $B(m,n)$ given by $m$ generators
and all relations of the form $v^n=1$ where $v$ is a word in
generators) is not amenable. It is interesting that the
possibility of using torsion groups and, in particular, free
Burnside groups as potential counterexamples to von Neumann's
problem was mentioned by Day \cite{day}, conjectured by Kesten
\cite{Kesten1}, and by B.H. Neumann's in his review of Specht's
paper \cite{Specht} in Zentralblatt.

Both Ol'shanskii's and Adian's examples are not finitely
presented: in the modern terminology these groups are inductive
limits of word hyperbolic groups, but they are not hyperbolic
themselves. Since many mathematicians (especially topologists) are
mostly interested in groups acting ``nicely" on manifolds, it is
natural to ask if there exists a finitely presented non-amenable
group without non-abelian free subgroups. This question was
explicitly formulated, for example, by Grigorchuk in \cite{Kour}
and by Cohen in \cite{Cohen}. This question is one of a series of
similar questions about finding finitely presented ``monsters",
i.e. groups with unusual properties. Probably the most famous
problem in that series is the problem about finding a finitely
presented infinite torsion group. Other similar problems ask for
finitely presented divisible group (group where every element has
roots of every degree), finitely presented Tarski monster, etc. In
each case a finitely generated example can be constructed as a
limit of hyperbolic groups (see \cite{book}), and there is no hope
to construct finitely presented examples as such limits.

One difficulty in constructing a finitely presented non-amenable
group without free non-abelian subgroups is that there are ``very
few" known finitely presented groups without free non-abelian
subgroups. Most non-trivial examples are solvable or ``almost"
solvable (see \cite{KhSap}), and so they are amenable. The only
known examples of finitely presented groups without free
non-abelian subgroups for which the problem of amenability is
non-trivial, are R.Thompson's group $F$ and its close
``relatives". The fact that $F$ does not contain free non-abelian
subgroups was proved by Brin and Squier in \cite{BS}. A conjecture
that $F$ is not amenable was formulated first by Geoghegan
\cite{Geogh}. A considerable amount of work has been done to prove
this conjecture (see \cite{CFP}) but it is still open.

One approach to constructing a finitely presented counterexample
to the von Neumann problem would be in using the Higman embedding
theorem which states that every recursively presented group can be
embedded into a finitely presented group. So one can take a known
finitely generated non-amenable group without non-abelian free
subgroups and embed it into a finitely presented group. Of course,
the resulting group will be non-amenable since the class of
amenable groups is closed under taking subgroups. Unfortunately
all known constructions of Higman embeddings (see, for example,
\cite{talk}) use amalgamated products and non-ascending HNN
extensions, which immediately leads to non-abelian free subgroups.
Nevertheless Higman-like embeddings play an important role in our
construction.

Our main result is the following.

\begin{theo} \label{main} For every sufficiently
large odd $n$, there exists a finitely presented group \label{g1}
${\cal G}$ which satisfies the following conditions.
\begin{enumerate}
\item\label{2} ${\cal G}$ is an ascending HNN extension of a finitely generated
infinite group of exponent $n$.
\item\label{3} ${\cal G}$ is an extension of a non-locally finite group of
exponent $n$ by an infinite cyclic group.
\item\label{4} ${\cal G}$ contains a subgroup isomorphic to a free
Burnside group of exponent $n$ with $2$ generators.
\item\label{5} ${\cal G}$ is a non-amenable finitely presented group
without free non-cyclic subgroups.
\end{enumerate}
\end{theo}

Notice that part \ref{2} of Theorem \ref{main} immediately implies
part \ref{3}. By a theorem of Adian \cite{Adian}, part \ref{4}
implies that ${\cal G}$ is not amenable. Thus parts \ref{3} and
\ref{4} imply part \ref{5}.

Note that the first example of a finitely presented group which is
a cyclic extension of an infinite torsion group was constructed by
Grigorchuk \cite{Grig}. But the torsion subgroup in Grigorchuk's
group does not have a bounded exponent and his group is amenable
(it was the first example of a finitely presented amenable but not
elementary amenable group).

In the two subsequent subsections, we present the main ideas of
our construction (we simplify the notation for the sake of
clarity).

\subsection{Congruence extension embeddings of the free Burnside groups}
\label{intro}

We first formulate two theorems on embeddings of the free Burnside
group \label{bmn}$B(m,n)$ of sufficiently large\footnote{The
number $n$ in the paper is chosen after one chooses several
auxiliary parameters $\alpha>\beta>\gamma...$, related by a
system of inequalities. We first choose $\alpha$, then $\beta,
\gamma$, etc.  Although the values of parameters and $n$ can be
calculated precisely (for example $n > 10^{10}$), we are not
doing it here because the consistency of the system of
inequalities is obvious. Indeed, each inequality has the form
$f(x,y,...)>0$ where $x<y<...$, and for already chosen positive
values of $y,...$, this inequality is obviously true for any
sufficiently small positive value of the smallest parameter $x$
(this is the \label{lpp}Lowest Parameter Principle from
\cite{book}).} odd exponent $n$ with $m>1$ generators, and deduce
Theorem \ref{main} from theorems \ref{SQintr} and \ref{CEP}
which, in turn, will be proved in sections $2-10.$

We say that a subgroup $H$ of a group $G$ satisfies \label{CEPr}
the {\it Congruence Extension Property} (CEP) if any homomorphism
$H\to H_1$ of $H$ onto a group $H_1$ extends to a homomorphism
$G\to G_1$ of $G$ onto some group $G_1$ containing $H_1$ as a
subgroup. Equivalently, for any normal subgroup $L$ of $H$, there
exists a normal subgroup $M$ of $G$ such that $H\cap M=L$. There
is another obvious reformulation, which is more convenient when
dealing with defining relations. Namely, the subgroup $H$
satisfies CEP if, for any subset $S\subseteq H$, we have $H\cap
S^G=S^H$ (where the normal closure of a subset $T$ in a group $K$
is denoted as $T^K$).

We will write $H\le_{CEP} G$ if the subgroup $H$ of a group $G$
satisfies CEP. We also say that $H$ is CEP-embedded in $G$, or $H$
is a CEP-subgroup of $G$.

Clearly every retract $H$ of any group $G$ satisfies CEP. It is
also clear that if $K$ is a CEP-subgroup of $H$ and $H$ is a
CEP-subgroup of $G$ then $K$ is a CEP-subgroup of $G$.

But there exist less obvious examples. For instance, the free
group $F_\infty$ of infinite rank can be CEP-embedded into a
2-generated free group $F_2$. Moreover, $F_\infty$ can be
CEP-embedded into any non-elementary hyperbolic group \cite{Ol95}.
This fact was used in \cite{Ol95} to prove that every
non-elementary hyperbolic group is SQ-universal.

\begin{theo} \label{SQintr} For any sufficiently large odd $n,$ there
exists a natural number\footnote{By using our proof of theorem
\ref{SQintr} and some additional arguments, D.Sonkin recently
showed that one can set $s=2$.} $s=s(n)$ such that the free
Burnside group $B(\infty,n)$ of infinite countable rank is
CEP-embedded into the free Burnside group $B(s,n)$ of rank
$s$.\footnote{ Theorem \ref{SQintr} is of independent interest
because it immediately implies, in particular, that every
countable group of exponent $n$ is embedable into a finitely
generated group of exponent $n$. This was  first proved by
Obraztsov (see \cite{book})}.
\end{theo}

In the next theorem, ${\cal H}^n$ denotes the subgroup of a group
${\cal H}$ generated by all $n$-th powers $h^n$ of elements $h\in
{\cal H}$.

\begin{theo} \label{CEP}
Let $n$ be a sufficiently large odd number. Then for every
positive integer $s$, there is a finitely presented group  ${\cal
H}$ containing subgroup $B(s,n)$, such that

(1) ${\cal H}^n\cap B(s,n)=1,$ so there exists a canonical
embedding of $B(s,n)$ into ${\cal H}/{\cal H}^n$.

(2) the embedding of $B(s,n)$ into ${\cal H}/{\cal H}^n$ from (1)
satisfies CEP.
\end{theo}

To derive here Theorem \ref{main} from Theorems \ref{SQintr} and
\ref{CEP}, we assume that $n$ is a sufficiently large odd integer
and $s$ is the number given by Theorem \ref{SQintr}. Then let
\label{ash}${\cal H}$ be a group provided by Theorem \ref{CEP}. It
has a finite presentation $\langle c_1,\dots,c_m\ |{\cal
R}\rangle$.

Since retracts are CEP-subgroups, we can chose a subgroup chain
\begin{equation}
B(m,n)\le_{CEP}B(m+2,n)\le_{CEP} B(\infty,n)\le_{CEP} B_s(s,n)\le
{\cal H} \label{vved1} \end{equation}
  where the first and the second
embeddings are given by the embeddings of the generating sets, the
third and the fourth ones are given by theorems \ref{SQintr} and
\ref{CEP}, respectively. Since CEP is a transitive property, we
obtain from (\ref{vved1}) and Theorem \ref{CEP} the chain
\begin{equation}
B(m,n)\le_{CEP} B(m+2,n)\le_{CEP} {\cal H}(n) \label{vved2}
\end{equation}

The factor-group ${\cal H}(n)={\cal H}/{\cal H}^n$ admits
presentation
\begin{equation}
{\cal H}(n)=\langle c_1,\dots,c_m \ | {\cal R}\cup {\cal
V}\rangle, \label{vved3}
\end{equation}
where ${\cal V}$ is the set of all $n$-th powers $v^n$ of group
words $v=v(c_1,\dots,c_m)$.

Let $b_1,\dots,b_m$ and $b_1,\dots,b_m, b_{m+1}, b_{m+2}$ be free
generating sets of the subgroup $B(m,n)$ and $B(m+2,n)$,
respectively, in (\ref{vved1}) and in (\ref{vved2}). It will be
convenient to denote these subgroups of ${\cal H}$ (and of ${\cal
H}(n)$) by $B_b(m,n)$ and by $B_b(m+2,n)$, respectively. There are
words $w_1(c_1,\dots,c_m),\dots,w_m(c_1,\dots,c_m)$ such that
$$b_1=w_1(c_1,\dots,c_m),\dots,b_m=w_m(c_1,\dots,c_m)$$ in ${\cal
H}$. For every $r(c_1,\dots,c_m)\in {\cal R}$, we define {\em
derived} words
\begin{eqnarray}
r'(c_1,\dots,c_m)\equiv
r(w_1(c_1,\dots,c_m),\dots,w_m(c_1,\dots,c_m)),\nonumber\\
r''(c_1,\dots,c_m)\equiv
r'(w_1(c_1,\dots,c_m),\dots,w_m(c_1,\dots,c_m)), \dots
\label{vved4}
\end{eqnarray}
(the sign "$\equiv$" \label{equiv} means the letter-for-letter
equality) and define ${\cal R}'$ as the set of all $r', r'',\dots,
r^{(i)},...$ for every $r\in {\cal R}$. We consider the group
obtained from (\ref{vved3}) by the following formula

\begin{equation}\label{4.1}
\bar{\cal H}={\cal H}(n)/({\cal R'})^{{\cal H}(n)}=\langle
c_1,\dots,c_m\ |{\cal R}\cup{\cal V}\cup{\cal R}'\rangle
\end{equation}

Since $B(m,n)\cong \langle c_1,\dots,c_m\ |{\cal V}\rangle$, we
have from (\ref{4.1})
\begin{equation}
 \bar{\cal H} \cong B_c(m, n)/({\cal R}\cup {\cal
 R}')^{B_c(m,n)}.\label{vved4.5}
 \end{equation}
where $B_c(m,n)$ is the free Burnside group with $m$ free
generators $c_1,...,c_m$. Let ${\cal R}_b$ and ${\cal R}'_b$
consist of copies of words from ${\cal R}$ and ${\cal R}'$ written
in letters $b_1,...,b_m$. Then from (\ref{vved4.5}), we
immediately deduce

\begin{equation}
 \bar{\cal H} \cong B_c(m, n)/({\cal R}\cup {\cal
 R}')^{B_c(m,n)}\cong B_b(m,n)/({\cal R}_b\cup {\cal
 R}'_b)^{B_b(m,n)}.\label{vved5}
 \end{equation}

Since the embedding of $B_b(m,n)$ into ${\cal H}(n)$ is an
CEP-embedding (see (\ref{vved2})), the right-hand side of
(\ref{vved5}) is isomorphic to the group $$\bar B=B_b(m,n)/(B_b(m,
n)\cap ({\cal R}_b\cup{\cal R}'_b)^{{\cal H}(n)}).$$ However, by
the definition (\ref{vved4}), $r(b_1,\dots,b_m)=r'(c_1,\dots,c_m),
r'(b_1,\dots,b_m)=r''(c_1,\dots,c_m),\dots$ in ${\cal H}(n),$ and
so $({\cal R}_b\cup{\cal R}'_b)^{{\cal H}(n)}=({\cal R'})^{{\cal
H}(n)}.$ Hence $\bar B=B_b(m,n)/(B_b(m, n)\cap ({\cal R'})^{{\cal
H}(n)})$, i.e. $\bar B$ is the canonical image of the subgroup
$B_b(m,n)\le {\cal H}(n)$ in $\bar{\cal H}$. Thus, the left-hand
side $\bar{\cal H}$ of (\ref{vved5}) is isomorphic to the image
$\bar B$ of $B_b(m,n)$ in $\bar {\cal H}$ under the mapping
$c_i\mapsto b_i$, $i=1,\dots m$. (We preserve notations
$c_1,\dots,c_m$, $b_1,\dots,b_m$ for the images of the generators
in $\bar {\cal H}$.)

The isomorphism of $\bar{\cal H}$ with its subgroup $\bar B$
allows us to introduce the ascending HNN-extension

$${\cal G}=\langle \bar {\cal H}, t\ |tc_it^{-1}=b_i,
i=1,\dots,m\rangle=\langle c_1,\dots,c_m\ |{\cal R}\cup{\cal
V}\cup {\cal R}'\cup {\cal U}\rangle, $$ where ${\cal
U}=\{tc_it^{-1}w_i^{-1}, i=1,...,m\}$.\footnote{E.Rips has been
publicizing a simpler idea of constructing ${\cal G}$ for several
years. He suggested to consider a finitely presented group ${\cal
H}=\langle c_1,\dots,c_m\ |{\cal R}\rangle $, containing
$B(2,n)$, subject to relations from ${\cal R}$ and additional
relations of the form $c_i=u_i$, $i=1,\dots m$, where the words
$u_i(c_1,\dots,c_m)$ represent ``long and complicated" elements
in $B(2,n)$. The resulting group is clearly finitely presented
and torsion since the images of ${\cal H}$ and $B(2,n)$ coincide
in it. So if it were infinite, it would be an example of an
infinite finitely presented torsion group (a solution of an
outstanding group theory problem). If such an example were to
exist, it would most certainly be non-amenable. Unfortunately in
all known cases the resulting group is obviously finite, and
Rips' idea has not been implemented yet. Our group (\ref{vved6})
is not torsion but it is torsion-by-cyclic, and it satisfies the
identity $[x,y]^n=1$. } However, relators from ${\cal U}$ and
formulas (\ref{vved4}) imply that
$tr^{(i)}(c_1,\dots,c_m)t^{-1}=r^{(i+1)}(c_1,\dots,c_m)$ for $i\ge
0$ where $r^{(0)}= r, r^{(1)}=r', \dots $ for $r\in {\cal R}.$
Hence ${\cal G}=\langle c_1,\dots,c_m\ | {\cal R}\cup {\cal
V}\cup {\cal U}\rangle.$

Next we observe that for every word $v$ in letters $b_1,...,b_m$,
the relation $v(b_1,\dots,b_m)^n=1$ follows from ${\cal R}$
because $B_b(m,n)$ is a subgroup of ${\cal H}$. Therefore
$v(c_1,\dots,c_m)^n=t^{-1}v(b_1,\dots,b_m)^nt=1$ follows from
${\cal R}\cup {\cal U}$. Hence
\begin{equation}
{\cal G}=\langle c_1,\dots,c_m\ |{\cal R}\cup {\cal U}\rangle
\label{vved6}, \end{equation} and so ${\cal G}$ is a finitely
presented group.

The group ${\cal G}$ is an ascending HNN-extension of the group
$\bar{\cal H}$, which is of exponent $n$ being a homomorphic
image of the group ${\cal H}(n)$.

Finally, we notice that the canonical image $B$ of the subgroup
$B_b(m+2,n)\le {\cal H}(n)$ in $\bar{\cal H}$ is isomorphic to
$B_b(m+2,n)/({\cal R'}^{{\cal H}(n)}\cap B_b(m+2,n))=
B_b(m+2,n)/({\cal R'})^{B_b(m+2,n)}$ because of the CEP-embedding
$B_b(m+2,n)\le_{CEP}{\cal H}$ (see (\ref{vved2})). The group
$B_b(m+2,n)/({\cal R'})^{B_b(m+2,n)}$ can be homomorphically
mapped onto the free Burnside group $B(2,n)$ because all words
from ${\cal R'}$ are written in the first $m$ of the $m+2$
generators of  $B_b(m+2,n)$. Therefore $B$ has a retract which is
isomorphic to $B(2,n)$. (The retraction $B(2,n)\to  B$ does exist
since the group $B(2,n)$ is free in the variety of periodic group
defined by the identity $x^n=1$.)

Thus, statements 1 and 3 of Theorem \ref{main} (and so statements
2 and 4, as we noticed earlier) follow from Theorems \ref{SQintr}
and \ref{CEP}.

\begin{remark}{\rm
The embedding of $B(s,n)$ into the group ${\cal H}$ in Theorem
\ref{CEP} and the set of defining relation ${\cal R}$ of ${\cal
H}$ will be constructed explicitly. Since the function $f$ in
Lemma \ref{lm01} and the number of generators in Lemma \ref{02}
are obtained explicitly, the CEP-embedding for Theorem
\ref{SQintr} is constructive too. Therefore one can also write
down an explicit presentation (\ref{vved6}) of the group ${\cal
G}$ in Theorem \ref{main}. But an easier, clear and explicit
construction of a {\bf finite} set of the subwords $w_1,\dots w_m$
of relators ${\cal U}$ from (\ref{vved6}) is explained in Remark
\ref{rem1} which goes back to \cite{Sa87}.}

\end{remark}

\subsection{The scheme of the proof of Theorem \ref{CEP}
and the plan of the paper } \label{intro2}

 To prove Theorem \ref{CEP}, we have to construct special embeddings of given
 free Burnside group $B(m,n)$ generated by a set ${\cal B}$ of cardinality $m$
 into a finitely presented group ${\cal H}$ such that ${\cal H}^n\cap
 B(m,n)=1$, and the canonical image of $B(m,n)$ in ${\cal H}(n)={\cal H}/{\cal H}^n$ enjoys CEP.
 The embedding in ${\cal H}$ is obtained
 in a similar way as in our papers \cite{OS},
\cite{talk} but we need a more complicated $S$-machine than in
\cite{OS} ($S$-machines were introduced by Sapir in \cite{SBR}),
and now the proof requires a detailed analysis of the construction
of the group $ {\cal H}(n)$.

This group is a factor-group of the group ${\cal H}$ generated by
${\cal C}=\{c_1,\dots, c_m\}$ subject to the relations  ${\cal R}$
of ${\cal H}$  over the normal subgroup generated by Burnside
relations. In other words, $ {\cal H}(n)$ is the {\em Burnside
factor} of ${\cal H}$. Burnside factors of free groups and free
products have been studied first by Adian and Novikov in
\cite{AN}, \cite{Adian2}. Geometric approach based on the notion
of $A$-map was employed in the study of Burnside factors of these
and more complicated groups in \cite{book}. Papers \cite{Ol95},
\cite{OlIvanov} extends this approach to Burnside factors of
hyperbolic groups. The main problem we face in this paper is that
${\cal H}$ is ``very" non-hyperbolic. In particular, the set of
relations ${\cal R}$ (denoted below by $Z(\sss,\Lambda)$) contains
many commutativity relations, so ${\cal H}$ contains non-cyclic
torsion-free abelian subgroups which cannot happen in a hyperbolic
group.

Nevertheless (and it is one of the main ideas of the paper) one
can make the Cayley graph of ${\cal H}$ look hyperbolic if one
divides the generators from ${\cal C}$ into two sets and consider
letters from one set as zero letters, and the corresponding edges
of the Cayley graph as edges of length 0. Thus the \label{lop}{\em
length of a path} in the Cayley graph or in a \vk diagram over the
presentation of ${\cal H}$ is the number of non-zero edges of the
path.

More precisely, the group ${\cal H}$ is similar to the group
$G(\sss)$ of \cite{SBR}, \cite{talk}, \cite{OS}. As we mentioned
above it corresponds to an $S$-machine $\sss$. The set ${\cal C}$
consists of tape letters (the set $\aaa$), state letters (the set
$\bkk$) and command letters (the set $\rr$). Recall that unlike an
ordinary Turing machine, an $S$-machine works with elements of a
group, not elements of the free semigroup.

It turns out that the most productive point of view is to consider
an $S$-machine as an inverse semigroup of partial transformations
of a set of states which are special ({\em admissible}) words of a
certain group $\hnka$ generated by $\bkk\cup \aaa$. The generators
of the semigroup are the $S$-{\em rules}. The machine $\sss$ is
the set of the $S$-rules. Every computation of the machine
corresponds to a word over $\sss$ which is called the
\label{history}{\em history of computation}, i.e. the string of
commands used in the computation. With every computation $h$
applied to an admissible element $W$, one associates a \vk diagram
$T(W,h)$ (called a trapezium) over the presentation of $\hnka$
(see \setcounter{pdeleven}{\value{ppp}} Figure \thepdeleven; a
more precise picture of a trapezium is Figure 14 in Section
\ref{trapezia}).

The first and the last words of the computation are written on the
bases the trapezium, copies of the history of the computation are
written on the vertical sides. The horizontal strips (bands) of
cells correspond to applications of individual rules in the
computation.

The trapezia $T(W,h)$ play central role in our study of the
Burnside factor ${\cal H}(n)$ of ${\cal H}$. As in \cite{book},
the main idea is to construct a graded presentation of ${\cal
H}(n)$ where longer relations have higher ranks and such that
every \vk diagram over the presentation of ${\cal H}(n)$ has the
so called property A from \cite{book}. In all diagrams over the
graded presentation of ${\cal H}(n)$, cells corresponding to the
relations from ${\cal R}$ are considered as $0$-cells or cells of
rank $1/2$, and cells corresponding to Burnside relations from the
graded presentation are considered as cells of ranks 1, 2,.... So
in these \vk diagrams ``big" Burnside cells are surrounded by
``invisible" 0-cells and ``small" cells.

\unitlength=1.00mm \special{em:linewidth 0.4pt}
\linethickness{0.4pt}
\begin{picture}(108.72,66.94)
\put(53.72,5.94){\line(1,0){36.33}}
\put(90.06,5.94){\line(1,3){15.67}}
\put(105.72,52.94){\line(-1,0){56.00}}
\put(49.72,52.94){\line(-1,0){6.00}}
\put(43.72,52.94){\line(1,-5){9.33}}
\put(59.06,5.94){\line(-1,0){6.33}}
\put(57.39,5.94){\line(-1,5){9.33}}
\put(85.72,5.94){\line(1,3){15.67}}
\put(101.39,52.61){\line(1,0){4.33}}
\put(101.72,52.78){\line(1,0){4.00}}
\put(52.39,9.11){\line(1,0){4.33}}
\put(51.56,12.78){\line(1,0){4.50}}
\put(50.72,17.28){\line(1,0){4.33}}
\put(49.89,22.28){\line(1,0){4.33}}
\put(48.89,27.28){\line(1,0){4.33}}
\put(47.89,31.67){\line(1,0){4.33}}
\put(47.06,36.44){\line(1,0){4.33}}
\put(46.06,40.94){\line(1,0){4.33}}
\put(45.06,46.11){\line(1,0){4.33}}
\put(86.72,9.11){\line(1,0){4.50}}
\put(87.89,12.78){\line(1,0){4.50}}
\put(89.56,17.33){\line(1,0){4.17}}
\put(91.06,22.33){\line(1,0){4.50}}
\put(92.89,27.28){\line(1,0){4.17}}
\put(92.89,27.28){\line(1,0){4.33}}
\put(94.39,31.67){\line(1,0){4.33}}
\put(95.89,36.33){\line(1,0){4.33}}
\put(97.39,41.00){\line(1,0){4.33}}
\put(99.06,46.00){\line(1,0){4.50}}
\put(49.39,46.00){\line(1,0){50.00}}
\put(53.22,46.00){\line(0,1){6.83}}
\put(61.39,46.00){\line(0,1){6.83}}
\put(61.39,53.00){\line(0,0){0.00}}
\put(79.22,46.00){\line(0,1){6.83}}
\put(89.72,46.00){\line(0,1){6.83}}
\put(71.22,46.00){\line(0,1){6.83}}
\put(37.39,29.61){\makebox(0,0)[cc]{History $h$}}
\put(108.72,29.94){\makebox(0,0)[cc]{History $h$}}
\put(70.06,2.94){\makebox(0,0)[cc]{Starting word $W$}}
\put(70.06,56.94){\makebox(0,0)[cc]{Last word $W\cdot h$}}
\put(46.00,41.00){\line(1,0){54.33}}
\put(49.33,36.33){\line(1,0){49.00}}
\put(50.67,31.67){\line(1,0){46.00}}
\put(51.00,27.28){\line(1,0){44.67}}
\put(51.33,22.33){\line(1,0){42.33}}
\put(53.00,17.33){\line(1,0){38.67}}
\put(54.33,12.67){\line(1,0){35.67}}
\put(54.67,9.00){\line(1,0){35.00}}
\put(53.33,45.67){\line(0,-1){4.67}}
\put(61.33,46.00){\line(0,-1){5.00}}
\put(79.33,46.00){\line(0,-1){5.00}}
\put(89.67,46.00){\line(0,-1){5.00}}
\put(57.67,41.00){\line(0,-1){4.67}}
\put(66.00,41.00){\line(0,-1){23.67}}
\put(76.33,53.00){\line(0,-1){47.00}}
\put(79.33,53.00){\line(0,-1){47.00}}
\put(83.33,41.00){\line(0,-1){35.00}}
\put(71.33,40.67){\line(0,-1){34.67}}
\put(61.67,17.33){\line(0,-1){4.67}}
\put(66.00,12.67){\line(0,-1){6.67}}
\put(89.67,41.00){\line(0,-1){13.67}}
\put(61.33,43.33){\line(0,-1){21.00}}
\put(53.33,36.33){\line(0,1){7.33}}
\put(66.00,46.00){\line(0,-1){5.00}}
\end{picture}

\begin{center}
\nopagebreak[4] Fig. \theppp.

\end{center}
\addtocounter{ppp}{1}

The main part of property A is the property that if a diagram over
${\cal H}(n)$ contains two Burnside cells $\Pi_1, \Pi_2$ connected
by a rectangular {\em \ct} subdiagram $\Gamma$ of rank 0 where the
sides contained in the contours of the two Burnside cells are
``long enough"  then these two cells cancel, that is the union of
$\Gamma$, $\Pi, \Pi'$ can be replaced by a smaller subdiagram.
This is a ``graded substitute" to the classic property of small
cancellation diagrams (where \ct subdiagrams contain no cells).

Roughly speaking nontrivial \ct subdiagrams of rank 0 turn out to
be trapezia of the form $T(W,h)$ (after we clean them of Burnside
$0$-cells), so properties of \ct subdiagrams can be translated
into properties of the machine $\sss$, and the inverse semigroup
described above.

As an intermediate step in studying the group ${\cal H}(n)$, we
construct a graded presentation of the Burnside factor of the
subgroup of ${\cal H}$ generated by $\rr\cup\aaa$. To avoid
repeating the same arguments twice, for ${\cal H}$ and for the
subgroup, we formulate certain key properties (Z1), (Z2), (Z3) of
a presentation of a group with a separation of generators into
zero and non-zero generators, so that there exists a graded
presentation of the Burnside factor of the group which satisfies
property $A$ (Section \ref{segregation}).

In order to roughly explain these conditions, consider the
following example. Let $P=F_A\times F_B$ be the direct product of
two free groups of rank $m$. Then the Burnside factor of $P$ is
simply $B(m,n)\times B(m,n)$. Nevertheless the theory of
\cite{book} cannot be formally applied to $P$. Indeed, there are
arbitrarily thick rectangles corresponding to relations $u\iv v\iv
uv=1$ in the Cayley graph of $P$ so diagrams over $P$ are not
A-maps in the terminology of \cite{book} (i.e. they do not look
like hyperbolic spaces). But one can obtain the Burnside factor of
$P$ in two steps. First we factorize $F_A$ to obtain
$Q=B(m,n)\times F_B$. After that we consider all edges labeled by
letters from $A$ in the Cayley graph of $Q$ as edges of length 0.
As a result the Cayley graph of $Q$ becomes a hyperbolic space.
This allows us to apply the theory of A-maps from \cite{book} to
obtain the Burnside factor of $Q$. The real reason for the theory
from \cite{book} to work in $Q$ is that $Q$ satisfies our
conditions (Z1), (Z2), (Z3). But the class of groups satisfying
these conditions is much bigger and includes groups corresponding
to $S$-machines considered in this paper. In particular (Z3) holds
in $Q$ because all 0-letters centralize $F_B$. This does not
happen in more complicated situations. But we associate with every
cyclically reduced non-0-element $w$ a ``personal" subgroup
$\oo(w)$ consisting of $0$-elements which is normalized by $w$.

The plan of the paper is the following. First in Section
\ref{prelim} we discuss some basic facts about \vk diagrams and
the main tools of studying them. In Section \ref{amaps}, we repeat
the main definitions and the main results of the theory of A-maps
from \cite[Chapter 5]{book}. We also give some modifications of
statements from \cite{book} that we need later.

Then in Section \ref{segregation}, we introduce properties (Z1),
(Z2), (Z3) and study Burnside factors. Although we follow the
general scheme of \cite{book}, we encounter new significant
difficulties. One of the main difficulties is that non-zero
elements can be conjugates of zero elements. As a bi-product of
this study, we show that under certain mild conditions, one can
one can construct analogs of HNN extensions in the class of
Burnside groups of sufficiently large odd exponents (see Corollary
\ref{Bernsrem}).

In Section \ref{appendix}, we prove some properties of the free
Burnside group, in particular the existence of a countably
generated free Burnside subgroup with congruence extension
property (Theorem \ref{SQ}). Roughly speaking, we obtain  the
required embedding of the group $B(\infty, n)$ into $B(s,n)$ if we
choose the images of free generators of $B(\infty,n)$ as aperiodic
words $A_1, A_2,\dots$ in generators of $B(s,n)$ satisfying a
strong small cancellation condition. Here we start the inductive
proof with the factor-group of the absolutely free group $F_s$ by
arbitrary relations of form $r(A_1,A_2,\dots)$. (This first step
was described in \cite{Ol95}.)

Then in Section \ref{presentation}, we define the presentation of
our group ${\cal H}$ by listing the rules of an $S$-machine and
showing how to convert the $S$-machine into a group.

In Section \ref{properties} we conduct a detailed study of the
inverse semigroup corresponding to our $S$-machine. Two important
properties reflect the flavor of that section and can be singled
out. Lemma \ref{burns} (iv) says that if there is a computation
starting with an element $W$, with history of the form $h^2$ (that
is if $W$ is in the domain of $h^2$), then there exists a
computation with history $h^s$ for every integer $s$. In other
words, if a sequence of rules can be applied twice in a row, it
can be applied any number of times. That lemma is used in the
geometric part of the paper when we are cleaning \ct subdiagrams
of Burnside 0-cells (see Lemma \ref{starcirc}). Proposition
\ref{propsss} provides an important property of stabilizers of
elements in the inverse semigroup generated by the $S$-rules. It
says, basically, that if $W$ is stabilized by $h$, $ghg\iv$ and
$g\iv hg$ then it is stabilized by $g^shg^{-s}$ for any integer
$s$. This proposition is used when we consider long and narrow \ct
subdiagrams (see Lemma \ref{ra3.21}). Another very important
feature of our $S$-machine is that two words $h_1$ and $h_2$ in
the alphabet $\sss$ define the same partial transformation on the
intersection of their domains provided $h_1$ and $h_2$ are equal
modulo Burnside relations (see Lemma \ref{burns} (ii)). Thus the
work of the $S$-machine is ``compatible" with Burnside relations.

In the remaining sections of the paper, we show that properties
(Z1), (Z2), (Z3) are satisfied by the presentation of the group
$H_{kra}$ corresponding to our $S$-machine.

The proof ends in Section \ref{proof}.

\bigskip

{\bf Acknowledgment.} The authors are grateful to Eliyahu Rips for
inspiring discussions of topics related to this paper. The authors
are also grateful to Goulnara Arjantseva who read the original
version of the paper and sent us numerous comments about it. This
allowed us to significantly improve the exposition.


\section{Maps and diagrams}
\label{amaps}

Planar maps and \vk diagrams are standard tools to study
complicated groups (see, for example, \cite{LS} and \cite{book}).
Here we present the main concepts related to maps and diagrams. A
\vk diagram is a labeled map, so we start by discussing maps.

\subsection{Graphs and maps}
\label{maps}

We are using the standard definition of a \label{graph}{\em
graph}. In particular, every edge has a direction, and every edge
has an inverse edge (having the opposite direction). For every
edge $e$, \label{emin}$e_-$ and $e_+$ are the beginning and the
end vertices of the edge.

A \label{mapp}{\em map} is simply a finite plane graph drawn on a
disc or an annulus which subdivides this surface into polygonal
cells. Edges do not have labels, so no group presentations are
involved in studying maps. On the other hand the theory of maps
helps find the structure of \vk diagrams because every diagram
becomes a map after we remove all the labels.

Let us recall the necessary definitions from \cite{book}.

A map is called \label{graded}{\em graded} if every cell $\Pi$ is
assigned a non-negative number, $r(\Pi)$, its \label{rank}{\em
rank}. A map $\Delta$ is called a \label{rankmap}{\em map of rank}
$i$ if all its cells have ranks $\le i$. Every graded map $\Delta$
has a \label{type}{\it type} $\tau (\Delta)$ which is the vector
whose first coordinate is the maximal rank of a cell in the map,
say, $r$, the second coordinate is the number of cells of rank $r$
in the map, the third coordinate is the number of cells of rank
$r-1$ and so on. We compare types lexicographically.

Let $\Delta$ be a graded map on a surface $X$. The cells of rank
$0$ in $\Delta$ are called \label{ocell}0-{\em cells}. Some of the
edges in $\Delta$ are called \label{oedge}$0$-{\em edges}. Other
edges and cells will be called \label{pedge}{\em positive}. We
define the {\em length of a path} $p$ in a map as the number of
positive edges in it. The length of a path $p$ is denoted by
$|p|$, in particular, $|\partial(\Pi)|$ is the \label{perimeter}
perimeter of a cell $\Pi$.

We assume that \label{m123}

\begin{itemize}
\item[(M1)] the inverse edge of a $0$-edge is also a $0$-edge,
\item[(M2)] either all edges in a 0-cell are 0-edges
or exactly two of these edges are positive,
\item[(M3)] for every cell $\Pi$ of rank $r(\Pi)>0$ the length
$|\partial(\Pi)|$ of its contour is positive.
\end{itemize}

\subsection{Bands}

\label{prelim}

A \label{vk}\vk diagram over a presentation $\la {\cal X}\ |\
{\cal R}\ra$ is a planar map with edges labeled by elements from
${\cal X}^{\pm 1}$, such that the label of the contour of each
cell belongs to ${\cal R}$ up to taking cyclic shifts and inverses
\cite{book}. In $\S 11$ of \cite{book}, 0-edges of \vk diagrams
have label 1, but in this paper 0-edges will be labeled also by
the so called \label{0let} {\em 0-letters}. All diagrams which we
shall consider in this paper, when considered as planar maps, will
obviously satisfy conditions (M1), (M2), (M3) from Section
\ref{maps}. By \label{vkl}\vk Lemma, a word $w$ over ${\cal X}$ is
equal to 1 modulo ${\cal R}$ if and only if there exists a \vk
disc diagram over ${\cal R}$ with boundary label $w$. Similarly
(see \cite{LS}) two words $w_1$ and $w_2$ are conjugate modulo $R$
if and only if there exists an annular diagram with the internal
counterclockwise contour labeled by $w_1$ and external
counterclockwise contour labeled by $w_2$. There may be many
diagrams with the same boundary label. Sometimes we can reduce the
number of cells in a diagram or ranks of cells by replacing a
2-cell subdiagram by another subdiagram with the same boundary
label. The simplest such situation is when there are two cells in
a diagram which have a common edge and are mirror images of each
other. In this case one can {\em cancel} these two cells (see
\cite{LS}). In order strictly define this cancellation (to
preserve the topological type of the diagram), one needs to use
the so called $0$-refinement. The $0$-refinement consists of
adding $0$-cells with (freely trivial) boundary labels of the form
$1aa\iv$ or $111$ (see \cite{book}, $\S 11.5$).

As in our previous papers \cite{SBR}, \cite{BORS} and
\cite{OlDist}, one of the main tools to study \vk diagrams are
bands and annuli.

Let $S$ be a subset of ${\cal X}$. An \label{sband}$S$-band $\bb$
is a sequence of cells $\pi_1,...,\pi_n$, for some $n\ge 0$ in a
\vk diagram such that

\begin{itemize}
\item Each two consecutive cells in this sequence have a common edge
labeled by a letter from $S$.
\item Each cell $\pi_i$, $i=1,...,n$ has exactly two $S$-edges
(i.e. edges labeled by a letter from $S$).
\item If $n=0$, then the boundary of $S$ has form $ee^{-1}$ for
an $S$-edge $e$.
\end{itemize}

\setcounter{pdten}{\value{ppp}} Figure \thepdten\ illustrates this
concept. In this Figure edges $e, e_1,...,e_{n-1},f$ are
$S$-edges, the lines $l(\pi_i,e_i), l(\pi_i,e_{i-1})$ connect
fixed points in the cells with fixed points of the corresponding
edges.

\bigskip
\unitlength=0.90mm \special{em:linewidth 0.4pt}
\linethickness{0.4pt}
\begin{picture}(149.67,30.11)
\put(19.33,30.11){\line(1,0){67.00}}
\put(106.33,30.11){\line(1,0){36.00}}
\put(142.33,13.11){\line(-1,0){35.67}}
\put(86.33,13.11){\line(-1,0){67.00}}
\put(33.00,21.11){\line(1,0){50.00}}
\put(110.00,20.78){\line(1,0){19.33}}
\put(30.00,8.78){\vector(1,1){10.33}}
\put(52.33,8.44){\vector(0,1){10.00}}
\put(76.33,8.11){\vector(-1,1){10.33}}
\put(105.66,7.78){\vector(1,2){5.33}}
\put(132.66,8.11){\vector(-1,1){10.00}}
\put(16.00,21.11){\makebox(0,0)[cc]{$e$}}
\put(29.66,25.44){\makebox(0,0)[cc]{$\pi_1$}}
\put(60.33,25.44){\makebox(0,0)[cc]{$\pi_2$}}
\put(133.00,25.78){\makebox(0,0)[cc]{$\pi_n$}}
\put(145.33,21.44){\makebox(0,0)[cc]{$f$}}
\put(77.00,25.11){\makebox(0,0)[cc]{$e_2$}}
\put(122.00,25.44){\makebox(0,0)[cc]{$e_{n-1}$}}
\put(100.33,25.44){\makebox(0,0)[cc]{$S$}}
\put(96.33,30.11){\makebox(0,0)[cc]{$\dots$}}
\put(96.33,13.11){\makebox(0,0)[cc]{$\dots$}}
\put(26.66,4.78){\makebox(0,0)[cc]{$l(\pi_1,e_1)$}}
\put(52.66,4.11){\makebox(0,0)[cc]{$l(\pi_2,e_1)$}}
\put(78.33,4.11){\makebox(0,0)[cc]{$l(\pi_2,e_2)$}}
\put(104.66,3.78){\makebox(0,0)[cc]{$l(\pi_{n-1},e_{n-1})$}}
\put(134.00,3.11){\makebox(0,0)[cc]{$l(\pi_n,e_{n-1})$}}
\put(88.66,30.11){\makebox(0,0)[cc]{$q_2$}}
\put(88.66,13.11){\makebox(0,0)[cc]{$q_1$}}
\put(49.33,25.11){\makebox(0,0)[cc]{$e_1$}}
\put(33.33,21.11){\circle*{0.94}}
\put(46.66,21.11){\circle*{1.33}}
\put(60.33,21.11){\circle*{0.94}}
\put(73.66,21.11){\circle*{1.33}}
\put(95.00,20.56){\makebox(0,0)[cc]{...}}
\put(19.33,21.00){\line(1,0){14.00}}
\put(129.67,20.67){\line(1,0){12.67}}
\put(19.00,13.00){\vector(0,1){17.00}}
\put(46.67,13.00){\vector(0,1){17.00}}
\put(73.67,13.00){\vector(0,1){17.00}}
\put(117.00,13.00){\vector(0,1){17.00}}
\put(117.00,20.67){\circle*{1.33}}
\put(129.67,20.67){\circle*{0.67}}
\put(142.33,13.00){\vector(0,1){17.00}}
\put(142.33,20.67){\circle*{1.33}}
\put(19.00,21.00){\circle*{1.33}}
\put(14.00,8.33){\vector(1,1){11.00}}
\put(10.33,5.00){\makebox(0,0)[cc]{$l(\pi_1,e)$}}
\put(148.00,8.00){\vector(-1,1){10.33}}
\put(154.67,3.11){\makebox(0,0)[cc]{$l(\pi_n,f)$}}
\end{picture}

\begin{center}
\nopagebreak[4] Fig. \theppp.

\end{center}
\addtocounter{ppp}{1}


The broken line formed by the lines $l(\pi_i,e_i)$,
$l(\pi_i,e_{i-1})$ connecting points inside neighboring cells is
called the \label{cl}{\em connecting line\footnote{In \cite{SBR},
\cite{BORS} this line was called the {\em median} of the band, but
in \cite{book}, it is called the connecting line.}} of the band
$\bb$. The $S$-edges $e$ and $f$ are called the \label{start}{\em
start} and {\em end} edges of the band.

The counterclockwise boundary of the subdiagram formed by the
cells $\pi_1,...,\pi_n$ of $\bb$ has the form $e\iv q_1f q_2\iv$.
We call $q_1$ the \label{bottom}{\em bottom} of $\bb$ and $q_2$
the {\em top} of $\bb$, denoted \label{bottom1}$\bott(\bb)$ and
$\topp(\bb)$.

A band $\pi_1,...,\pi_t$ is called \label{redd}{\em reduced} if
$\pi_{i+1}$ is not a mirror image of $\pi_i$, $i=1,...,t-1$
(otherwise cells $\pi_i$ and $\pi_{i+1}$ {\em cancel} and there
exists a diagram with the same boundary label as $\cup_i\pi_i$ ad
containing fewer cells). Similarly, we define reduce annuli.

We say that two bands \label{cross}{\em cross} if their connecting
lines cross. We say that a band is an \label{annulus}{\em annulus}
if its connecting line is a closed curve. In this case the first
and the last cells of the band coincide (see
\setcounter{pdeleven}{\value{ppp}} Figure \thepdeleven a).

\bigskip
\begin{center}
\unitlength=1.5mm
\linethickness{0.4pt}
\begin{picture}(101.44,22.89)
\put(30.78,13.78){\oval(25.33,8.44)[]}
\put(30.78,13.89){\oval(34.67,18.00)[]}
\put(39.67,9.56){\line(0,-1){4.67}}
\put(25.67,9.56){\line(0,-1){4.67}}
\put(18.56,14.89){\line(-1,0){5.11}}
\put(32.56,7.11){\makebox(0,0)[cc]{$\pi_1=\pi_n$}}
\put(21.00,7.33){\makebox(0,0)[cc]{$\pi_2$}}
\put(15.89,11.78){\makebox(0,0)[cc]{$\pi_3$}}
\put(43.44,12.89){\line(1,0){4.67}}
\put(43.00,8.44){\makebox(0,0)[cc]{$\pi_{n-1}$}}
\put(19.22,10.89){\line(-1,-1){4.00}}
\put(25.89,20.44){\circle*{0.00}}
\put(30.78,20.44){\circle*{0.00}}
\put(30.33,1.56){\makebox(0,0)[cc]{a}}
\put(35.44,20.44){\circle*{0.00}}
\put(62.78,4.89){\line(0,1){4.67}}
\put(62.78,9.56){\line(1,0){38.67}}
\put(101.44,9.56){\line(0,-1){4.67}}
\put(101.44,4.89){\line(-1,0){38.67}}
\put(82.11,16.22){\oval(38.67,13.33)[t]}
\put(62.78,15.78){\line(0,-1){7.33}}
\put(101.44,16.44){\line(0,-1){8.22}}
\put(82.00,12.22){\oval(27.78,14.67)[t]}
\put(67.89,12.89){\line(0,-1){8.00}}
\put(95.89,13.11){\line(0,-1){8.22}}
\put(67.89,13.78){\line(-1,1){4.67}}
\put(95.67,14.00){\line(1,1){5.11}}
\put(73.67,9.56){\line(0,-1){4.67}}
\put(89.67,9.56){\line(0,-1){4.67}}
\put(76.56,21.33){\circle*{0.00}}
\put(81.44,21.33){\circle*{0.00}}
\put(86.11,21.33){\circle*{0.00}}
\put(65.22,12.22){\makebox(0,0)[cc]{$\pi_1$}}
\put(65.22,6.89){\makebox(0,0)[cc]{$\pi$}}
\put(71.00,6.89){\makebox(0,0)[cc]{$\gamma_1$}}
\put(92.78,7.11){\makebox(0,0)[cc]{$\gamma_m$}}
\put(98.56,7.11){\makebox(0,0)[cc]{$\pi'$}}
\put(98.33,12.67){\makebox(0,0)[cc]{$\pi_n$}}
\put(94.11,20.22){\makebox(0,0)[cc]{$S$}}
\put(85.89,6.89){\makebox(0,0)[cc]{$T$}}
\put(76.56,6.89){\circle*{0.00}} \put(79.44,6.89){\circle*{0.00}}
\put(82.33,6.89){\circle*{0.00}}
\put(83.44,1.56){\makebox(0,0)[cc]{b}}
\end{picture}
\end{center}
\begin{center}
\nopagebreak[4] Fig. \theppp.

\end{center}
\addtocounter{ppp}{1}

If a band $\bb$ connects two edges on the boundary of a \vk disc
diagram $\Delta$ then the connecting line of $\bb$ divides the
region bounded by $\partial(\Delta)$ into two parts. Let
$\Delta_1$ and $\Delta_2$ be the maximal subdiagrams of $\Delta$
contained in these parts. Then we say that the (connecting line of
the) band $\bb$ \label{divides}{\em divides $\Delta$ into two
subdiagrams $\Delta_1$ and $\Delta_2$}. Notice that $\bb$ does not
belong to either $\Delta_1$ or $\Delta_2$ and
$\Delta=\Delta_1\cup\Delta_2\cup\bb$. Similarly if $\Delta$ is an
annular diagram and $\bb$ is an annulus surrounding the hole of
$\Delta$ then (the connecting line of) $\bb$ divides $\Delta$ into
two annular subdiagrams.

The maximal subdiagram bounded by the connecting line of an
annulus in a disc diagram is called the \label{inside}{\em inside
diagram of this annulus}.

Let $S$ and $T$ be two disjoint subsets of ${\cal X}$, let ($\pi$,
$\pi_1$, \ldots, $\pi_n$, $\pi'$) be an $S$-band and let ($\pi$,
$\gamma_1$, \ldots, $\gamma_m$, $\pi'$) is a $T$-band. Suppose
that:
\begin{itemize}
\item the connecting lines of these bands form a simple closed
curve,
\item on the boundary  of $\pi$ and on the boundary of $\pi'$ the pairs of $S$-edges separate the pairs of $T$-edges,
\item the start and end edges of these bands are not contained in the
region bounded by the connecting lines of the bands.
\end{itemize}
Then we say that these bands form an \label{stanl}{\em
$(S,T)$-annulus} and the closed curve formed by the parts of
connecting lines of these bands is the {\em connecting line} of
this annulus (see Figure 3b). For every $(S,T)$-annulus we define
the {\em inside} subdiagram diagram of the annulus as the maximal
subdiagram bounded by the connecting line of the annulus.

\subsection{Bonds and contiguity submaps}

 \label{bcs}

Now let us introduce the main concepts from Chapter 5 of
\cite{book}.  Consider a map $\Delta$. A path in $\Delta$ is
called \label{redpath} {\em reduced} if it does not contain
subpaths of length 2 which are homotopic to paths of length 0 by
homotopies involving only 0-cells. We shall usually suppose that
the contour of $\Delta$ is subdivided into several subpaths called
\label{sect}{\em sections}. We usually assume that sections are
reduced paths.

By property M2 if a 0-cell has a positive edge, it has exactly two
positive edges. Thus we can consider bands of 0-cells having
positive edges. The start and end edges of such bands are called
\label{adj}{\em adjacent}.

Let $e$ and $f$ be adjacent edges of $\Delta$ belonging to two
positive cells $\Pi_1$ and $\Pi_2$ or to sections of the contour
of $\Delta$. Then there exists a band of $0$-cells connecting
these edges. The union of the cells of this band is called a
\label{obond}{\em $0$-bond} between the cells $\Pi_1$ and $\Pi_2$
(or between a cell and a section of the contour of $\Delta$, or
between two sections of the contour). The  connecting line of the
band is called the \label{clb}{\em connecting line} of the bond.
The contour of the $0$-bond has the form $p\iv es\iv f\iv$ where
$p$ and $s$ are paths of length 0 because every $0$-cell has at
most two positive edges.

Now suppose we have chosen two pairs $\{e_1, f_1\}$ and
$\{e_2,f_2\}$ of adjacent edges such that $e_1, e_2$ belong to the
contour of a cell $\Pi_1$ and $f_1\iv, f_2\iv$ to the contour of a
cell $\Pi_2$. Then we have two bonds $E_1$ and $E_2$. If $E_1=E_2$
(that is $e_1=e_2, f_1=f_2$) then let $\Gamma=E_1$. Now let
$E_1\ne E_2$ and $z_1e_1w_1f_1\iv$ and $z_2e_2w_2f_2\iv$ be the
contours of these bonds. Further let $y_1$ and $y_2$ be subpaths
of the contours of $\Pi_1$ and $\Pi_2$ (of $\Pi_1$ and $q$, or
$q_1$ and $q_2$, where $q$, $q_1$ and $q_2$ are some segments of
the boundary cotours) where $y_1$ (or $y_2$) has the form
$e_1pe_2$ or $e_2pe_1$ (or $(f_1uf_2)^{-1}$ or $(f_2uf_1)^{-1}$).
If $z_1y_1w_2y_2$ (or $z_2y_1w_1y_2$) is a contour of a disc
submap $\Gamma$ which does not contain $\Pi_1$ or $\Pi_2$, then
$\Gamma$ is called a $0$-{\em contiguity submap} of $\Pi_1$ to
$\Pi_2$ (or $\Pi_1$ to $q$ or of $q_1$ to $q_2$). The contour of
$\Gamma$ is naturally subdivided into four parts. The paths $y_1$
and $y_2$ are called the {\em contiguity arcs}. We write
$y_1=\Gamma \bigwedge \Pi_1$, $y_2=\Gamma \bigwedge \Pi_2$ or
$y_2=\Gamma \bigwedge q$. The other two paths are called the {\em
side arcs} of the \ct submap.

The ratio $|y_1|/|\partial(\Pi_1)|$ (or $|y_2|/|\partial(\Pi_2)|$
is called the {\em degree of contiguity} of the cell $\Pi_1$ to
the cell $\Pi_2$ or to $q$ (or of $\Pi_2$ to $\Pi_1$). We denote
the degree of contiguity of $\Pi_1$ to $\Pi_2$ (or $\Pi_1$ to $q$)
by $(\Pi_1,\Gamma,\Pi_2)$ (or $(\Pi_1,\Gamma,q)$). Notice that
this definition is not symmetric and when $\Pi_1=\Pi_2=\Pi$, for
example, then $(\Pi,\Gamma,\Pi)$ is a pair of numbers.

A {\em connecting line} of $\Gamma$ is chosen as the connecting
line of one of the bonds $E_1$ or $E_2$.

We say that two contiguity submaps are {\em disjoint} if none of
them has a common point with the interior of the other one, and
their contiguity arcs do not have common edges.

Now we are going to define $k$-bonds and $k$-contiguity submaps
for $k>0$. In these definitions we need a fixed number
$\varepsilon$, $0<\varepsilon< 1$.

Let $k>0$ and suppose that we have defined the concepts of
$j$-bond and $j$-contiguity submap for all $j<k$. Consider three
cells $\pi, \Pi_1, \Pi_2$ (possibly with $\Pi_1=\Pi_2$) satisfying
the following conditions:

\begin{enumerate}
\item $r(\pi)=k, r(\Pi_1)>k, r(\Pi_2)>k$,
\item there are disjoint submaps $\Gamma_1$ and $\Gamma_2$ of
$j_1$-contiguity of $\pi$ to $\Pi_1$ and of $j_2$-contiguity of
$\pi$ to $\Pi_2$, respectively, with $j_1<k, j_2<k$, such that
$\Pi_1$ is not contained in $\Gamma_2$ and $\Pi_2$ is not
contained in $\Gamma_1$,
\item $(\pi,\Gamma_1,\Pi_1)\ge\varepsilon$,
$(\pi,\Gamma_2,\Pi_2)\ge\varepsilon$.
\end{enumerate}

Then there is a minimal submap $E$ in $\Delta$ containing $\pi,
\Gamma_1, \Gamma_2$. This submap is called a \label{princ}{\em
$k$-bond between $\Pi_1$ and $\Pi_2$ defined by the contiguity
submaps $\Gamma_1$ and $\Gamma_2$ with principal cell $\pi$} (see
\setcounter{pdeleven}{\value{ppp}} Figure \thepdeleven).

The {\em contiguity arc} $q_1$ of the bond $E$ to $\Pi_1$ is
$\Gamma_1\bigwedge \Pi_1$. It will be denoted by $E\bigwedge
\Pi_1$. Similarly $E\bigwedge\Pi_2$ is by definition the arc
$q_2=\Gamma_2\bigwedge\Pi_2$. The contour $\partial(E)$ can be
written in the form $p_1q_1p_2q_2$ where $p_1$, $p_2$ are called
the {\em side arcs} of the bond $E$.

\unitlength=1mm \special{em:linewidth 0.4pt} \linethickness{0.4pt}
\begin{picture}(90.33,55.67)
\put(67.33,40.33){\oval(46.00,9.33)[b]}
\put(44.33,42.00){\line(0,1){3.00}}
\put(44.33,46.67){\line(0,1){3.33}}
\put(90.33,42.00){\line(0,1){3.00}}
\put(90.33,46.67){\line(0,1){3.33}}
\put(47.17,52.33){\oval(5.00,5.33)[lt]}
\put(86.17,51.67){\oval(8.33,8.00)[rt]}
\put(67.50,23.17){\oval(29.67,11.00)[]}
\put(58.67,28.67){\line(0,1){7.00}}
\put(74.33,28.67){\line(0,1){7.00}}
\put(67.33,5.33){\oval(46.00,10.67)[t]}
\put(60.67,17.67){\line(0,-1){7.00}}
\put(71.67,10.67){\line(0,1){7.00}}
\put(58.00,47.00){\makebox(0,0)[cc]{$\Pi_1$}}
\put(68.67,38.33){\makebox(0,0)[cc]{$q_1$}}
\put(62.33,7.33){\makebox(0,0)[cc]{$q_2$}}
\put(80.00,4.00){\makebox(0,0)[cc]{$\Pi_2$}}
\put(71.83,10.67){\line(0,1){7.00}}
\put(71.67,17.50){\line(1,0){5.50}}
\put(75.42,23.08){\oval(14.17,11.17)[r]}
\put(74.50,28.83){\line(1,0){2.67}}
\put(77.17,28.83){\line(0,0){0.00}}
\put(74.50,28.83){\line(0,1){6.67}}
\put(58.50,35.67){\line(0,-1){7.00}}
\put(58.50,28.67){\line(-1,0){0.83}}
\put(58.17,23.17){\oval(11.33,11.33)[l]}
\put(58.33,17.50){\line(1,0){2.17}}
\put(60.50,17.50){\line(0,-1){6.83}}
\put(60.50,15.33){\vector(0,-1){2.00}}
\put(74.50,30.33){\vector(0,1){2.67}}
\put(66.00,10.67){\vector(-1,0){2.67}}
\put(65.00,35.67){\vector(1,0){2.83}}
\put(70.00,27.17){\makebox(0,0)[cc]{$v_1$}}
\put(68.83,19.83){\makebox(0,0)[cc]{$v_2$}}
\bezier{108}(63.67,35.67)(60.33,20.50)(67.00,10.67)
\put(61.00,23.00){\makebox(0,0)[cc]{$t$}}
\put(76.00,23.00){\makebox(0,0)[cc]{$\pi$}}
\end{picture}

\begin{center}
\nopagebreak[4] Fig. \theppp.

\end{center}
\addtocounter{ppp}{1}

For each $k$-bond we fix a connecting line obtained by joining the
terminal points of connecting lines of $\Gamma_1$ and $\Gamma_2$
inside $\pi$. The connecting line divides $\pi$ (and the whole
$E$) into two connected components.

Bonds between a cell and a section of the contour of $\Delta$ or
between two sections of the contour are defined in a similar way.

Now let $E_1$ be a $k$-bond and $E_2$ be a $j$-bond between two
cells $\Pi_1$ and $\Pi_2$, $j\le k$ and either $E_1=E_2$ or these
bonds are disjoint. If $E_1=E_2$ then $\Gamma=E_1=E_2$ is called
the $k$-contiguity submap between $\Pi_1$ and $\Pi_2$ determined
by the bond $E_1=E_2$. If $E_1$ and $E_2$ are disjoint then the
corresponding \label{csm}{\em contiguity submap} $\Gamma$ is
defined as the smallest submap containing $E_1$ and $E_2$ and not
containing $\Pi_1$ and $\Pi_2$ (see
\setcounter{pdeleven}{\value{ppp}} Figure \thepdeleven).

\unitlength=1mm \special{em:linewidth 0.4pt} \linethickness{0.4pt}
\begin{picture}(135.00,76.00)
\put(73.50,68.00){\oval(123.00,16.00)[b]}
\put(73.50,13.00){\oval(123.00,20.00)[t]}
\put(35.00,41.00){\oval(26.00,16.00)[]}
\put(107.00,41.00){\oval(28.00,16.00)[]}
\put(30.00,54.00){\oval(6.00,4.00)[]}
\put(43.00,54.00){\oval(6.00,4.00)[]}
\put(102.50,54.00){\oval(7.00,4.00)[]}
\put(114.50,54.00){\oval(7.00,4.00)[]}
\put(29.00,27.00){\oval(6.00,4.00)[]}
\put(41.00,27.00){\oval(6.00,4.00)[]}
\put(102.00,27.00){\oval(6.00,4.00)[]}
\put(114.00,27.00){\oval(6.00,4.00)[]}
\put(28.00,29.00){\line(0,1){4.00}}
\put(31.00,33.00){\line(0,-1){4.00}}
\put(39.00,29.00){\line(0,1){4.00}}
\put(42.00,33.00){\line(0,-1){4.00}}
\put(28.00,25.00){\line(0,-1){2.00}}
\put(31.00,23.00){\line(0,1){2.00}}
\put(39.00,25.00){\line(0,-1){2.00}}
\put(43.00,23.00){\line(0,1){2.00}}
\put(28.00,49.00){\line(0,1){3.00}}
\put(32.00,52.00){\line(0,-1){3.00}}
\put(44.00,52.00){\line(1,-2){3.00}}
\put(41.00,56.00){\line(-1,4){1.00}}
\put(45.00,60.00){\line(-1,-4){1.00}}
\put(29.00,56.00){\line(-1,4){1.00}}
\put(33.00,60.00){\line(-1,-2){2.00}}
\put(101.00,56.00){\line(-1,4){1.00}}
\put(107.00,60.00){\line(-3,-4){3.00}}
\put(113.00,56.00){\line(0,1){4.00}}
\put(116.00,60.00){\line(0,-1){4.00}}
\put(99.00,33.00){\line(1,-2){2.00}}
\put(104.00,29.00){\line(1,4){1.00}}
\put(100.00,25.00){\line(0,-1){2.00}}
\put(100.00,23.00){\line(0,0){0.00}}
\put(104.00,23.00){\line(0,1){2.00}}
\put(113.00,23.00){\line(0,1){2.00}}
\put(115.00,25.00){\line(0,-1){2.00}}
\put(113.00,29.00){\line(-1,2){2.00}}
\put(114.00,33.00){\line(1,-4){1.00}}
\put(101.00,52.00){\line(0,-1){3.00}}
\put(104.00,49.00){\line(0,1){3.00}}
\put(41.00,52.00){\line(-5,-3){5.00}}
\put(113.00,52.00){\line(-4,-3){4.00}}
\put(116.00,52.00){\line(2,-3){4.00}}
\end{picture}
\begin{center}
\nopagebreak[4] Fig. \theppp.

\end{center}
\addtocounter{ppp}{1}
 The
\label{car}{\em contiguity arcs} $q_1$ and $q_2$ of $\Gamma$ are
intersections of $\partial(\Gamma)$ with $\partial(\Pi_1)$ and
$\partial(\Pi_2)$. The contour of $\Gamma$ has the form
$p_1q_1p_2q_2$ where $p_1$ and $p_2$ are called the
\label{sar}{\em side arcs} of $\Gamma$. The \label{cls}{\em
connecting line} of $\Gamma$ is one of the connecting lines of
$E_1$ or $E_2$. The ratio $|q_1|/|\partial(\Pi_1)|$ is called the
\label{cde}{\em contiguity degree} of $\Pi_1$ to $\Pi_2$ with
respect to $\Gamma$ and is denoted by
\label{cde1}$(\Pi_1,\Gamma,\Pi_2)$. If $\Pi_1=\Pi_2=\Pi$ then
$(\Pi,\Gamma,\Pi)$ is a pair of numbers.

Contiguity submaps of a cell to a section of the contour and
between sections of the contour are defined in a similar way. (In
this paper we do not use notion  "degree of contiguity of a
section of a contour to" anything.)

We are going to write ``bonds" and ``\ct submaps" in place of
``$k$-bonds" and ``$k$-contiguity submaps". Instead of writing
``the \ct submap $\Gamma$ of a cell $\Pi$ to ...", we sometimes
write ``the $\Gamma$-\ct of $\Pi$ to ...".

The above definition involved the standard decomposition of the
contour of a \ct submap $\Gamma$ into four sections $p_1q_1p_2q_2$
where \label{wedge}$q_1=\Gamma\bigwedge \Pi_1$,
$q_2=\Gamma\bigwedge \Pi_2$ (or $q_2=\Gamma\bigwedge q$ if $q$ is
a section of $\partial(\Delta)$), we shall write
\label{par}$p_1q_1p_2q_2=\partial(\Pi_1,\Gamma,\Pi_2)$, and so on.

As in \cite{book} (see $\S 15.1$) we fix certain real numbers
\label{param}$\iota << \zeta << \varepsilon <<\delta << \gamma <<
\beta << \alpha$ between 0 and 1 where "$<<$" means ``much
smaller". Here ``much" means enough to satisfy all the
inequalities in Chapters 5, 6 of \cite{book}. We also set
$\bar\alpha=\frac{1}{2}+\alpha$, $\bar\beta=1-\beta$,
$\bar\gamma=1-\gamma, h=\delta\iv, n=\iota\iv$.

\subsection{Condition A}

\label{conditionA}

The set of positive cells of a diagram $\Delta$ is denoted by
$\Delta(2)$ and as before the length of a path in $\Delta$ is the
number of positive edges in the path. A path $p$ in $\Delta$ is
called {\em geodesic} if $|p|\le |p'|$ for any path $p'$
combinatorially homotopic to $p$.

The condition A has three parts:\label{A123}

\begin{itemize}
\item[A1] If $\Pi$ is a cell of rank $j>0,$ then
$|\partial(\Pi)|\ge nj$.

\item[A2] Any subpath of length $\le \max(j,2)$ of the contour of
an arbitrary cell of rank $j$ in $\Delta(2)$ is geodesic in
$\Delta$.

\item[A3] If $\pi, \Pi\in \Delta(2)$ and $\Gamma$ is a \ct submap
of $\pi$ to $\Pi$ with $(\pi,\Gamma,\Pi)\ge \varepsilon$, then
$|\Gamma\bigwedge\Pi|<(1+\gamma)k$ where $k=r(\Pi)$.
\end{itemize}

A map satisfying conditions A1, A2, A3 will be called an
\label{amp}$A$-{\em map}. As in \cite{book}, Section 15.2, a
(cyclic) section $q$ of a contour of a map $\Delta$ is called a
\label{smooth}{\em smooth section of rank} $k > 0$ (we write
\label{rq}$r(q)=k$) if:

1) every subpath of $q$ of length $\le \max(k,2)$ of $q$ is
geodesic in $\Delta$;

2) for each \ct submap $\Gamma$ of a cell $\pi$ to $q$ satisfying
$(\pi, \Gamma,q)\ge \varepsilon$, we have $|\Gamma\bigwedge
q|<(1+\gamma)k$.

It is shown in \cite{book}, $\S\S 16,  17$, that $A$-maps have
several ``hyperbolic" properties. We summarize these properties
here. (They show that $A$-maps of rank $i$ are hyperbolic spaces
with hyperbolic constant depending on $i$ only.)

The next three lemmas are Lemmas 15.3, 15.4 and 15.8 in
\cite{book}.

Lemma \ref{anal15.3} shows that in a \ct diagram between two cells
side arcs are ``much shorter" than perimeters of these cells (so
the cells are close to each other).

\begin{lm}\label{anal15.3} Let $\Delta$ be an A-map,
$\Gamma$ a contiguity submap of either:

(i)  a cell $\Pi_1$ to a cell $\Pi_2$ or

(ii) a cell $\Pi_1$ to a boundary section $q^1$ of $\Delta$, or

(iii) a boundary section $q^1$ to a boundary section $q^2$.

Let $p_1$ and $p_2$ be the side arcs of $\Gamma$, and
$c=\max(|p_1|,|p_2|)$. Then $c<2\varepsilon^{-1}r(\Pi_1)<\zeta n
r(\Pi_1)$ in cases (i) and (ii), and $c<\zeta n r(\Pi_2)$ in case
(i). Furthermore $c<\zeta n r(q^1)$  in cases (ii) and (iii)
provided $q^1$ is smooth.
\end{lm}

Lemma \ref{anal15.4} says that \ct arcs in a \ct subdiagram have
``almost the same lengths".

\begin{lm}\label{anal15.4} Let $\psi$ be the degree of
$\Gamma$-contiguity of a cell $\Pi_1$ to a cell $\Pi_2$ (or to a
section $q$ of a contour) in an A-map $\Delta$, and $q_1=\Gamma
\bigwedge\Pi_1$, $q_2$ are the \ct arcs of $\Gamma$. Then
$$(\bar\beta-2\zeta\psi^{-1})|q_1|<|q_2|.$$ In particular,
$|q_2|>(\psi-2\beta)|\partial(\Pi_1)|,$ and if also
$\psi\ge\varepsilon,$ then $|q_1|<(1+2\beta)|q_2|.$

Furthermore if $q_2=\Gamma\bigwedge\Pi_2$ or $q_2=\Gamma\bigwedge
q$ and $q$ is a smooth section, then
$$|q_1|>\bar\beta(1+2\zeta\psi^{-1})^{-1}|q_2|.$$ In particular,
if $\psi\ge\varepsilon$, then $$|q_1|>(1-2\beta)|q_2|.$$
\end{lm}

Lemma \ref{anal15.8} shows that the degree of \ct of one cell to
another cell cannot be significantly bigger than $1/2$ and the
degree of \ct of a cell of a bigger rank to a cell of a smaller
rank is very small.

\begin{lm}\label{anal15.8} In an arbitrary A-map $\Delta$
the degree of contiguity of an arbitrary cell $\pi$ to an
arbitrary cell $\Pi$ or to an arbitrary smooth section $q$ of the
contour via an arbitrary contiguity submap $\Gamma$ is less than $
\bar\alpha$, and the degree of contiguity of $\pi$ to $\Pi$ (or to
$q$) with $r(\Pi)\le r(\pi)$ (or $r(q)\le r(\pi)$) is less than
$\varepsilon.$
\end{lm}

The next lemma is a generalization of the Greenlinger lemma from
the theory of small cancellation groups \cite{LS}. It shows that
any $A$-map contains a cell whose boundary is almost entirely
contiguous to the boundary of the diagram.

\begin{lm}\label{cor16.1,16.2}  a) Let
$\Delta$ be a disc A-map whose contour is subdivided into at most
$4$ sections $q^1, q^2,q^3,q^4$, and $r(\Delta)>0$. Then there
exists a positive cell $\pi$ and disjoint \ct submaps
$\Gamma_1,...,\Gamma_4$ of $\pi$ to $q^1,...,q^4$ (some of these
submaps may be absent) such that the sum of corresponding \ct
degrees $\sum_{i=1}^4 (\pi,\Gamma_i,q^i)$ is greater than
$\bar\gamma=1-\gamma$.

b) Let $\Delta$ be an annular A-map whose contours are subdivided
into at most 4 sections $q^1, q^2,q^3,q^4$ regarded as cyclic or
ordinary paths, and $r(\Delta)>0$. Then $\Delta$ has a positive
cell $\pi$ and disjoint \ct submaps $\Gamma_1$,..., $\Gamma_4$ of
$\pi$ to $q^1$ $q^2$, $q^3$, $q^4$ (up to three of these submaps
may be absent) such that the sum of \ct degrees of these \ct
submaps is greater than $\bar\gamma$.
\end{lm}

\proof Both statements are immediate corollaries of Theorem 16.1
from \cite{book}.
\endproof

The next lemma follows from Lemma \ref{cor16.1,16.2} in the same
way as Theorem 16.2  follows from Corollary 16.1 in \cite{book}.

\begin{lm} \label{th16.2} If $\Delta$
is a disc or annular A-map with contour subdivided into at most
$4$ sections, and $r(\Delta)>0$, then there exists a positive cell
$\Pi$ and a \ct submap $\Gamma$ of $\Pi$ to one of the sections
$q$ of the contour of $\Delta$ such that $r(\Gamma)=0$ and
$(\Pi,\Gamma,q)\ge \varepsilon$.
\end{lm}

Here we add an additional property whose proof is contained in the
proof of Theorem 16.2 in \cite{book}.

\begin{lm}\label{Anal16.2} Let $\Delta$ be a disc or annular A-map with contour
subdivided into at most $4$ sections. Assume there is a cell $\Pi$
in $\Delta$ and a contiguity submap $\Gamma$ of $\Pi$ to one of
the sections $q$ of the boundary with
$(\Pi,\Gamma,q)\ge\varepsilon$. Then there is a cell $\pi$ of
$\Delta$ and a contiguity submap $\Gamma'$ of $\pi$ to $q$ such
that $r(\Gamma')=0$ and $(\pi,\Gamma',q)\ge\varepsilon.$
\end{lm}

The next four lemmas are Theorem 17.1, Corollary 17.1, Lemma 17.1
and Lemma 17.2 in \cite{book}.

The following lemma shows that smooth sections of contours of
$A$-maps are almost geodesic.

\begin{lm} \label{analteom17.1} Let $\Delta$ be a
disc map with contour $qt$ or an annular $A$-map with contours $q$
and $t.$ If $q$ is a smooth section then $\bar\beta |q|\le |t|.$
\end{lm}

The next nine lemmas demonstrate different variants of the fact
that in $A$-maps, quasigeodesic quadrangles are thin. But we need
more precise statements than in the usual theory of hyperbolic
spaces. The constants appearing in these statements do not depend
on the rank of diagrams while the hyperbolic constants of $A$-maps
tends to infinity with the rank.

\begin{lm} \label{cor17.1} If a disc $A$-map $\Delta$
contains a cell $\Pi$ of positive rank, then
$|\partial(\Delta)|>\bar\beta|\partial(\Pi)|.$
\end{lm}

\begin{lm} \label{anallemma17.2} Let $\Delta$ be a
disc A-map with contour $p_1q_1p_2q_2$ where $q_1$ and $q_2$ are
smooth sections such that $|q_1|+|q_2|>0$,
$|p_1|+|p_2|\le\gamma(|q_1|+|q_2|).$ Then there is either a 0-bond
between $q_1$ and $q_2$ or a cell $\Pi\in\Delta(2)$ and disjoint
contiguity submaps $\Gamma_1$ and $\Gamma_2$ of $\Pi$ to $q_1$ and
$q_2$, respectively, such that
$(\Pi,\Gamma_1,q_1)+(\Pi,\Gamma_2,q_2)>\bar\beta=1-\beta.$
\end{lm}

We also need the following version of Lemma \ref{anallemma17.2}.

We say that a bond between two boundary sections $q_1$ and $q_2$
of a map is \label{nb}{\it narrow} if it is either a 0-bond or, in
the inductive definition of the bond, we add that
$(\pi,\Gamma_1,q_1)+(\pi, \Gamma_2,q_2)>1-\alpha/2$ for the
principal cell $\pi$ of the bond. A \ct submap between $q_1$ and
$q_2$ is \label{ncs}{\it narrow} if it is defined by narrow bonds
(according to the inductive definition).

\begin{lm}\label{narrow}

Let $\Delta$ be a disc A-map with contour $p_1q_1p_2q_2$ where
$q_1$ and $q_2$ are smooth sections such that $|q_1|+|q_2|>0$,
$|p_1|+|p_2|\le \alpha^2(|q_1|+|q_2|)$. Then there is either a
$0$-bond between $q_1$ and $q_2$ or a cell $\Pi\in \Delta(2)$ and
disjoint \ct submaps $\Gamma_1$ and $\Gamma_2$ of  $\Pi$ to $q_1,
q_2$ (respectively) such that
$(\Pi,\Gamma_1,q_1)+(\Pi,\Gamma_2,q_2)\ge 1-\alpha/2.$ Hence
$\Pi$, $\Gamma_1$ and $\Gamma_2$ form a narrow bond between $q_1$
and $q_2$.
\end{lm}

\proof The difference between the formulations of Lemma
\ref{narrow} and Lemma 17.2 of \cite{book} is that the parameters
$\gamma$ and $\bar\beta$ are now substituted by $\alpha^2$ and
$1-\alpha/2$, respectively. The proof differs from the proof of
Lemma 17.2 \cite{book} only in some insignificant details.

There is no change in Case 1) of the proof of Lemma 17.2
\cite{book}. The only alteration in Case 2) is the substitution of
the first factor $\gamma$ for $\alpha^2$ in the left-hand side of
the long inequality at the end. It is easy to see that the
resulting inequality is also true by the choice of the parameters.

In Case 3), we should assume now that
$(\Pi,\Gamma,p_1)>\bar\gamma-(1-\alpha/2)=\alpha/2-\gamma$
(instead of $(\Pi,\Gamma,p_1)>\bar\gamma-\bar\beta=\beta-\gamma$
in \cite{book}). This leads to the similar substitution of the
term $\beta-\gamma$ for $\alpha/2-\gamma$ in the right-hand side
of inequality (7) of the proof. It remains to notice that the last
long inequality of Case 3) is still true after the substitution of
$\beta-\gamma$ by $\alpha/2-\gamma$ in the right-hand side, and
the substitution of the factor $\gamma$ by $\alpha^2$ in the
left-hand side.

The difference $\beta-\gamma$ should be replaced by
$\alpha/2-\gamma$ in the hypothesis of Case 4). This leads to the
same substitution in inequality (9) of the proof of Lemma 17.2
\cite{book}. Then again, it remains to replace the term
$\beta-\gamma$ by $\alpha/2-\gamma$ in the right-hand side of the
concluding inequality, and to replace the factor $\gamma$ by
$\alpha^2$ in the left-hand side of it.
\endproof

\begin{lm}\label{Anal17.3.2}
Let $\Delta$ be a disc A-map with contour $p_1q_1p_2q_2,$ where
$q_1$ and $q_2$ are smooth sections of some ranks $k$ and $l,$
$k\le l,$ and $|p_1|+|p_2|\le \gamma (|q_1|+|q_2|)$. Then there is
a narrow \ct submap between $q_1$ and $q_2$ with contour
$p'_1q'_1p'_2q'_2,$ where $q_1=u'q'_1u''$, $q_2 =v''q'_2 v'$ and
$|u'|, |v'|< \gamma^{-1}(|p_1|+\alpha k),$ $|u''|, |v''|<
\gamma^{-1}(|p_2|+\alpha k).$
\end{lm}

\proof Up to the term ``narrow", the proof is the same as the
proof of items 1) and 2) of Lemma 17.3 \cite{book} (Lemma
\ref{Anal17.3} below). In fact Lemma \ref{Anal17.3.2} follows from
the proof of Lemma 17.3 \cite{book} because $\alpha/2
> \beta$.
\endproof

\begin{lm}\label{Anal17.3} (Lemma 17.3 \cite{book}) Let $\Delta$
be a disc A-map with contour $p_1q_1p_2q_2$ where $q_1$ and $q_2$
are smooth sections of ranks $k$ and $l$, $k\le l$, such that

\begin{equation}
|p_1|+|p_2|\le\gamma(|q_1|+|q_2|).\label{(?.1)}
\end{equation}

Then:

(1) there exist vertices $o_1$ and $o_2$ on $q_1$ and $q_2$ and a
path $x$ joining them in $\Delta$ such that $|x|<\alpha k.$

(2) one can choose $o_1$ and $o_2$ in such a way that the lengths
of the initial segments of $q_1$ and $q_2^{-1}$ (or of $q_1^{-1}$
and $q_2$) ending at $o_1$ and $o_2,$ are less than
$\gamma^{-1}(|p_1|+\alpha k)$ (or than $\gamma^{-1}(|p_2|+\alpha
k)).$

(3)if (\ref{(?.1)}) is replaced by $$|p_1|, |p_2|<\alpha k,$$ then
any vertex $o$ of $q_1$ (or of $q_2$) can be joined in $\Delta$ by
a path $y$ to a vertex $o'$ on $q_2$ (or on $q_1$) in such a way
that $|y|<\gamma^{-1}k.$
\end{lm}

Consider an A-map $\Delta$ with contour $p_1q_1p_2q_2,$ where
$q_1,q_2$ are smooth sections and $\Delta$ is a narrow contiguity
submap between $q_1$ and $q_2.$ We define an
\label{alpha}$\alpha$-{\em series} as a maximal series of disjoint
narrow bonds between $q_1$ and $q_2.$ In particular, the maximal
series contains the bonds from the definition of narrow contiguity
submap. By $\Pi_0$ we denote a cell of maximal perimeter among
principal cells of the bonds over all $\alpha$-series.

\begin{lm} \label{Nearly}
With the preceding notation, any vertex of $q_1$ is at the
distance at most $ \alpha^{-1}|\partial(\Pi_0)|$ from the path
$q_2$.
\end{lm}

\proof Assume that $E_1,\dots,E_s$ is the $\alpha$-series. Let
$x_1^ly_1^lx_2^ly_2^l$ be the standard contour of $E_l,$ where
$y_j^l$ is a subpath of $q_j$ for $j=1,2.$ By Lemma \ref{anal15.3}
and the definition of narrow bonds,
$|x_j^l|<(\alpha+2\zeta)|\partial(\Pi_0)|.$ By Lemmas
\ref{anal15.4} and \ref{anal15.8},
$|y_j^l|<(1-2\beta)^{-1}\bar\alpha|\partial(\Pi_0)|.$ Therefore
the distance between any vertex of some $y_1^l$ and $q_2$ is less
than
$(\alpha+2\zeta+(1-2\beta)^{-1}\bar\alpha/2)|\partial(\Pi_0)|<|\partial(\Pi_0)|/3
<\alpha^{-1}|\partial(\Pi_0)|.$

Then we consider the submaps $\Delta_1,\dots,\Delta_{s-1}$ with
contours $(x_2^l)^{-1}z_1^l(x_1^{l+1})^{-1}z_2^l$, where $z_j^l$
is a subpath of $q_j$ for $j=1,2.$ One may suppose that
$|z_1^l|+|z_2^l|<\alpha^{-2}(|x_2^l|+|x_1^{l+1}|),$ because
otherwise by Lemma \ref{narrow} there would be an additional
narrow bond in $\Delta_l,$ which contradicts the definition of the
$\alpha$-series. (Indeed $(\Pi,\Gamma_j,
q_j)\ge1-\alpha/2-\bar\alpha >\varepsilon$ for $j=1,2$ in the
Lemma \ref{narrow} by Lemma \ref{anal15.8}, and therefore that
cell $\Pi$ is a principle cell of a narrow bond between $q_1$ and
$q_2.$) Since
$|x_2^l|+|x_1^{l+1}|<(\alpha+4\zeta)|\partial(\Pi_0)|,$ we have
$|z_1^l|<\alpha^{-1}(1+4\zeta\alpha^{-1})|\partial(\Pi_0)|,$ and
the distance between any vertex of $z_1^l$ and $q_2$ is less than
$|z_1^l|/2+(\alpha+2\zeta)|\partial(\Pi_0)|<\alpha^{-1}|\partial(\Pi_0)|,$
and the lemma is proved. \endproof

\begin{lm}\label{anal17.4} (Lemma 17.4 \cite{book}.) Let $\Delta$
be a disc A-map with contour $p_1q_1p_2q_2$, where $q_1$ and $q_2$
are smooth sections of ranks $k$ and $l,$ $k\le l$, and $|p_1|,
|p_2|<\zeta n k$. Then the perimeter of every cell of $\Delta$ is
less than $3\gamma^{-1}\zeta n k<nk$ and
$r(\Delta)<3\gamma^{-1}\zeta k < k.$
\end{lm}

\begin{lm}\label{anal17.5} (Lemma 17.5 \cite{book}). Let $\Delta$ be a
disc A-map with contour $p_1q_1p_2q_2$, where $q_1$ and $q_2$ are
smooth sections of ranks $k$ and $l,$ $k\le l$, and $|q_j|-M\ge 0$
for $j=1,2$ and $M=\gamma^{-1}(|p_1|+|p_2|+k).$ Then it is
possible to excise a submap from $\Delta$ with contour
$p'_1q'_1p'_2q'_2$ where $q'_1$ ($q'_2$) is a subpath of $q_1$ (of
$q_2$), $|p'_1|, |p'_2| <\alpha k,$ and $|q'_j|>|q_j|-M$ for
$j=1,2.$
\end{lm}

The next  lemma shows that annular $A$-maps can be cut by short
paths to obtain a disc map.

\begin{lm}\label{anallemma17.1}  (Lemma 17.1, \cite{book}) Let $\Delta$ be an
annular A-map with contours $p$ and $q$ such that any loop
consisting of 0-edges of $\Delta$ is 0-homotopic. Then there is a
path $t$ connecting some vertices $o_1$ and $o_2$ of the paths $p$
and $q$ respectively, such that $|t|<\gamma(|p|+|q|).$
\end{lm}


\section{Burnside quotients}

\label{segregation}

\subsection{Axioms}
\label{axioms}

Recall that in Chapter 6 of \cite{book}, it is proved that for any
sufficiently large odd $n$, there exists a presentation of the
free Burnside group $B(m,n)$ with $m\ge 2$ generators, such that
every reduced diagram over this presentation is an $A$-map.

Our goal here is to generalize the theory from \cite{book} to
Burnside factors of non-free groups in the following way.

In \cite{book}, 0-edges are labeled by 1 only, and (trivial)
0-relators have the form $1\cdot\dots\cdot 1$ or $a\cdot 1\cdot
... \cdot 1\cdot a^{-1} \cdot 1\cdot...\cdot 1$. We need to
increase the collection of labels of 0-edges and the collection of
0-relations. In other words, we divide the generating set of our
group into the set of 0-letters and the set of positive letters.
Then in the diagrams, 0-letters label 0-edges. We shall also add
defining relators of our group which have the form $uavb^{-1}$
where $u,v$ are words consisting of 0-letters only (we call such
words 0-{\it words}), and $a,b$ are a positive letters, as well as
defining relations without positive letters, to the list of zero
relations. The corresponding cells in diagrams will be called
$0$-cells.

Let ${\bf G}$ be a group given by a presentation $\la {\cal X}\ |
\ {\cal R}\rangle$. Let ${\cal X}$ be a symmetric set (closed
under taking inverses). Let ${\cal Y}\subseteq {\cal X}$ be a
symmetric subset which we shall call the {\it set} of {\em
non-zero letters} or $Y$-{\em letters}. All other letters in
${\cal X}$ will be called $0$-{\em letters}.

By the \label{yl}$Y$-{\em length} of a word $A$, written
\label{ay}$|A|_Y$ or just $|A|$, we mean the number of occurrences
of \label{yt}$Y$-letters in $A$. The $Y$-length $|g|_Y$ of an
element $g\in G$ is the length of the $Y$-shortest word
representing element $g$.

Let ${\bf G}(\infty)$ be the group given by the presentation $\la
{\cal X}\ |\ {\cal R}\ra$ inside the variety of Burnside groups of
exponent $n$. We will choose a presentation of the group ${\bf
G}(\infty)$ in the class of all groups, that is we add enough
relations of the form $A^n=1$ to ensure that ${\bf G}(\infty)$ has
exponent $n$. This infinite presentation will have the form
$\langle {\cal X} \ | \ {\cal R}(\infty)= \cup {\cal S}_i\rangle$
where sets ${\cal S}_i$ will be disjoint, and relations from
${\cal S}_i$ will be called relations of {\it rank } $i$. Unlike
\cite[Chapter 6]{book}, $i$ runs over $0, 1/2, 1, 2,3,\dots.$
($i=0,1,2,\dots$ in \cite[Chapter 6] {book}.) The set ${\cal R}$
will be equal to ${\cal S}_0\cup {\cal S}_{1/2}$. The detailed
description of sets ${\cal S}_i$ will be given below. We shall
call this presentation of ${\bf G}(\infty)$ a \label{gp}{\it
graded presentation}. Our goal is the same as in \cite{book}: we
want every diagram over ${\cal R}(\infty)$ to be an A-map.

We are going to prove that this goal is possible to achieve if
${\cal R}$ satisfies conditions \label{z123}(Z1), (Z2), (Z3)
presented below. While listing these conditions, we also fix some
notation and definitions.

\begin{itemize}
\item[(Z1)] The set ${\cal R}$ is the union of two disjoint subsets
${\cal S}_0={\cal R}_0$ and ${\cal S}_{1/2}$. The group $\langle
{\cal X}\ |\ {\cal S}_0\rangle$ is denoted by ${\bf G}(0)$ and is
called the group of rank $0$. We call relations from ${\cal S}_0$
relations of rank 0. The relations from ${\cal S}_{1/2}$ have rank
$1/2.$ The group $\langle {\cal X}\ |{\cal R}\rangle$ is denoted
by ${\bf G}(1/2).$
\begin{itemize}
\item[(Z1.1)]The set ${\cal S}_0$ consists of all relations from ${\cal R}$ which have
$Y$-length $0$ and all relations of ${\cal R}$ which have the form
(up to a cyclic shift) $ay_1by_2^{-1}=1$ where $y_1,y_2\in Y^{+}$,
$a$ and $b$ are of $Y$-length $0$ (i.e. \label{0word}$0$-words).

\item[(Z1.2)] The set ${\cal S}_0$ implies all Burnside relations $u^n=1$ of $Y$-length 0.
\end{itemize}
\end{itemize}
The subgroup of ${\bf G}(0)$ generated by all $0$-letters will be
called the \label{0sub}{\em $0$-subgroup} of ${\bf G}(0)$.
Elements from this subgroup are called \label{0el}$0$-{\em
elements}. Elements which are not conjugates of elements from the
$0$-subgroup will be called \label{ess}{\em essential}. An
essential element $g$ from ${\bf G}(0)$ is called \label{cyre}{\em
cyclically $Y$-reduced} if $|g^2|_Y=2|g|_Y$. We shall show that
$g$ is cyclically $Y$-reduced if and only if it is cyclically
minimal (in rank 0), i.e. it is not a conjugate of a $Y$-shorter
element in ${\bf G}(0)$ (Lemma \ref{minimal}).

Notice that condition (Z1.1) allows us to consider $Y$-bands in
\vk diagrams over ${\cal S}_0$. Maximal $Y$-bands do not
intersect.

\begin{itemize}
\item[(Z2)] The relators of the set ${\cal S}_{1/2},$ will be called
\label{hub1}{\it hubs}. The corresponding cells in \vk diagrams
are also called hubs. They satisfy the following properties
\begin{itemize}
\item[(Z2.1)] The $Y$-length of every hub is at
least $n$.
\item[(Z2.2)] Every hub is linear in $Y,$ i.e. contains
at most one occurrence $y^{\pm 1}$ for every letter $y\in Y.$
\item[(Z2.3)]
Assume that each of words $v_1w_1$ and $v_2w_2$ is a cyclic
permutation of a hub or of its inverse, and
$|v_1|\ge\varepsilon|v_1w_1|.$ Then an equality $u_1v_1=v_2u_2$
for some 0-words $u_1,$ $u_2,$ implies in ${\bf G}(0)$ equality
$u_2w_1=w_2u_1$ (see \setcounter{pdeleven}{\value{ppp}} Figure
\thepdeleven).

\unitlength=.85mm \special{em:linewidth 0.4pt}
\linethickness{0.4pt}
\begin{picture}(113.00,74.00)
\put(40.00,25.00){\framebox(60.00,20.00)[cc]{}}
\put(71.50,57.00){\oval(83.00,24.00)[]}
\put(71.50,14.00){\oval(83.00,22.00)[]}
\put(40.00,38.00){\vector(0,-1){7.00}}
\put(100.00,40.00){\vector(0,-1){8.00}}
\put(70.00,25.00){\vector(1,0){10.00}}
\put(62.00,45.00){\vector(1,0){12.00}}
\put(85.00,69.00){\vector(-1,0){10.00}}
\put(86.00,3.00){\vector(-1,0){12.00}}
\put(37.00,34.00){\makebox(0,0)[cc]{$u_1$}}
\put(71.00,22.00){\makebox(0,0)[cc]{$v_1$}}
\put(69.00,42.00){\makebox(0,0)[cc]{$v_2$}}
\put(103.00,36.00){\makebox(0,0)[cc]{$u_2$}}
\put(70.00,65.00){\makebox(0,0)[cc]{$w_2$}}
\put(102.00,6.00){\makebox(0,0)[cc]{$w_1$}}
\end{picture}

\begin{center}
\nopagebreak[4] Fig. \theppp.

\end{center}
\addtocounter{ppp}{1}

\end{itemize}
\end{itemize}

A word $B$ over ${\cal X}$ is said to be \label{cmwe}{\it
cyclically minimal in rank} $1/2$ if it is not a conjugate in
${\bf G}(1/2)$  of a word of smaller $Y$-length.  An element $g$
of ${\bf G}(0)$ is called {\em cyclically minimal in rank} $1/2$
if it is represented by a word which is cyclically minimal in rank
$1/2$.

\begin{itemize}
\item[(Z3)]  With every essential element $g\in {\bf G}(0)$ we
associate a subgroup $\oo(g)\le {\bf G}(0)$, normalized by $g$. If
$g$ is cyclically $Y$-reduced, let \label{0g}$\oo(g)$ be the
maximal subgroup consisting of $0$-elements which contains $g$ in
its normalizer. By Lemmas  \ref{reduced} (2) and \ref{cyclically}
(3) below arbitrary essential element $g$ is equal to a product
$vuv^{-1}$ where $u$ is cyclically $Y$-reduced. In this case we
define $\oo(g)=v\oo(u)v^{-1}.$ We shall prove in Lemma
\ref{welldef} that $\oo(g)$ is well defined.

\begin{itemize}
\item[(Z3.1)]
Assume that $g$ is an essential cyclically minimal in rank $1/2$
element of ${\bf G}(0)$ and for some $x\in {\bf G}(0)$, both $x$
and $g^{-4}xg^4$ are 0-elements. Then $x\in \oo(g).$

\item[(Z3.2)] For every essential cyclically minimal in rank $1/2$
element $g\in {\bf G}(0)$ there exists a $0$-element $r$ such that
$gr\iv$ commutes with every element of $\oo(g)$.\footnote{At the
beginning of May 2002, influenced by the talk by Sergei Ivanov at
the Geometric Group Theory conference in Montreal, we realized
that Condition (Z3.2) can be replaced by the following weaker
condition:

(Z3.2') For every essential cyclicly minimal in rank $1/2$ element
$g\in {\bf G}(0)$, the extension of $\oo(g)$ by the automorphism
induced by $g$ (acting by conjugation) satisfies the identity
$x^n=1$.

It is clear, that in the proof of Lemma \ref{star} below  (the
only lemma in this paper, where (Z3.2) is used), property (Z3.2)
can be replaced by (Z3.2').}
\end{itemize}
\end{itemize}

In this section, we show how the theory from \cite[Chapter
6]{book} can be adapted to diagrams over a presentation satisfying
(Z1), (Z2), (Z3).


\subsection{Corollaries of (Z1), (Z2), (Z3)}
\label{corollaries}

Here we are going to deduce some useful facts about groups ${\bf
G}(0)$ whose presentations satisfy (Z1), (Z2), (Z3). (Some of them
can be deduced from a presentation of ${\bf G}(0)$ as an
HNN-extension.)

The following lemma is obvious.

\begin{lm} \label{ZZ}
If the set ${\cal R}$ satisfies property (Z1.1) then there exists
a homomorphism from ${\bf G}(0)$ to the group of integers
$\mathbf{Z}$ which sends every $Y$-letter to $1$ and every
$0$-letter to $0$.
\end{lm}

\begin{lm} \label{uchastok}
Suppose that ${\cal R}$ satisfies (Z1.1). Let $\Delta$ be a \vk
diagram over ${\cal S}_0$ and $p$ be a subpath of the boundary
$\partial(\Delta)$. Suppose that a $Y$-band starting on $p$ ends
on $p\iv$. Then $\Lab(p)$ contains a subword of $Y$-length 2 which
is equal to a $0$-word in ${\bf G}(0)$.
\end{lm}

\proof Since maximal $Y$-bands cannot intersect, there exists a
$Y$-band $\ttt$ starting on edge $e$ and ending on edge $f\iv$,
where $e,f\in p$, such that the shortest subpath $p_1$ of $p$
containing $e$ and $f$ has length $2$. Let $\Delta_1$ be the
subdiagram of $\Delta$ bounded by $p_1$ and one of the sides $q_1$
of $\ttt$. The word $\Lab(q_1)$ is a $0$-word, so $\Lab(p_1)$ is a
subword of $A$ of $Y$-length $2$ which is equal to a $0$-word in
${\bf G}(0)$.
\endproof

A word $A$ which does not contain a subword of $Y$-length 2 that
is equal to a 0-word in ${\bf G}(0)$, will be called \label{yr}
$Y$-{\em reduced}. A word $A$ is called \label{crw} {\em
cyclically $Y$-reduced} if every cyclic shift of $A$ is
$Y$-reduced.

\begin{lm} \label{reduced}
Suppose that ${\cal R}$ satisfies (Z1.1). Then the following
conditions hold.

(1)  Let $A$ be a word representing an element $g\in {\bf G}(0)$,
but $|A|_Y
> |g|_Y$. Then $A$ is not $Y$-reduced.

(2) If $A$ is a $Y$-reduced word that represents a conjugate of an
essential element $g$, but it is not a shortest word representing
a conjugate of $g$, then the word $A$ has the form as a product
$UVU'$ where $|U|_Y>0$ and $U'U$ is equal to a $0$-word in ${\bf
G}(0)$, $V$ is a cyclically $Y$-reduced word.
\end{lm}

\proof (1) Let $B$ be a $Y$-shortest word representing $g$. Then
there exists a diagram over ${\cal S}_0$ with boundary $pq$ where
$\Lab(p)\equiv A$, $\Lab(q)\equiv B\iv$.

Condition (Z1.1) imply that every $Y$-edge on $A$ is the start
edge of a maximal $Y$-band. Each of these $Y$-bands ends on the
boundary of $\Delta$. Since $|A|_Y>|B|_Y$, one of the $Y$-bands
starting on $p$ ends on $p$. By Lemma \ref{uchastok} there exists
a subword of $A$ of length $2$ which is equal to a $0$-word in
${\bf G}(0)$.

To prove part (2), one can take an annular diagram instead of a
disc diagram as in part (1), and prove that a cyclic shift of $A$
has a subword of $Y$-length 2 which is equal to a $0$-word in
${\bf G}(0)$. Hence $A\equiv UVU'$ where $|U|_Y>0$, $U'U$ is equal
to a $0$-word $W$ in ${\bf G}(0)$. Since the word $VW$ represents
a conjugate of $g$ and is shorter than $A$, we can conclude the
proof by induction on $|A|_Y$.
\endproof

\begin{lm} \label{cyclically}
Suppose that ${\cal R}$ satisfies (Z1.1). Then the following
statements are true.

\begin{enumerate}
\item[(1)] If $g\in {\bf G}(0)$ and $A$ is a $Y$-reduced word representing $g$
then there exists a decomposition $A\equiv  A_1A_2$ such that the
element represented by $A_2A_1$ is cyclically $Y$-reduced in ${\bf
G}(0)$.

\item[(2)] If $g\in {\bf G}(0)$ is cyclically $Y$-reduced then for every integer $m$,
$|g^m|_Y= |m||g|_Y$ and $g^m$ is cyclically $Y$-reduced.

\item[(3)] If $g$ is cyclically $Y$-reduced and $A$ is a $Y$-shortest word
representing $g$ in ${\bf G}(0)$ then $A$ is cyclically
$Y$-reduced.
\end{enumerate}
\end{lm}

\proof (1) By Lemma \ref{reduced} (1) we can suppose that $A$ is a
$Y$-shortest word representing $g$. Let $B$ be a shortest word
representing $g^2$. If $|B|_Y=2|A|_Y$ then $g$ is cyclically
$Y$-reduced. Suppose that $2|A|_Y>|B|_Y$. By Lemma \ref{reduced}
there exists a subword $C$ of $A^2$ of $Y$-length $2$ which is
equal to a $0$-word $D$ in ${\bf G}(0)$. By Lemma \ref{ZZ} we can
assume that $C\equiv y_1uy_2\iv$ (or $C\equiv y_1^{-1}uy_2$) where
$u$ is a $0$-word and $y_1,y_2\in Y^+$.

Since $A$ is $Y$-reduced, one can represent $A$ as $u_2y_2\iv
A_1y_1u_1$ where $u_1u_2\equiv  u$. Then the cyclic shift
$A_1y_1uy_2$ of $A$ can be represented by $A_1D$ which is  a
reduced word. Since this word is $Y$-shorter than $A$ we can apply
induction and conclude that there exists a cyclic shift of $A_1D$
which represents a cyclically reduced element.

(2) Let $A$ be a shortest word representing $g$. Suppose that for
some $m>0$, $A^m$ is not a shortest word representing $g^m$. Then
by Lemma \ref{reduced} $A^m$ contains a subword of $Y$-length 2
which is equal to a $0$-word in ${\bf G}(0)$. But then the same
subword occurs in $A^2$, so $A^2$ is not a shortest word
representing $g^2$, a contradiction.

(3) Notice that for every cyclic shift $A_1$ of $A$, $A_1^2$ is a
subword of $A^3$. It remains to apply part (2) of this lemma.
\endproof

\begin{lm} \label{minimal}
Suppose that ${\cal R}$ satisfies (Z1.1). Let $A$ be a
$Y$-shortest word representing an essential element $g$ in ${\bf
G}(0)$.

Then the following conditions are equivalent:

(i)  $A$ is cyclically $Y$-reduced;

(ii) $g$ is cyclically minimal in ${\bf G}(0)$ (that is $g$ is not
conjugate in ${\bf G}(0)$ of an element of a smaller $Y$-length);

(iii) $g$ is cyclically $Y$-reduced in ${\bf G}(0)$ (that is
$|g^2|_Y=2|g|_Y$).
\end{lm}

\proof (i)$\to $(ii). Suppose that every cyclic shift of $A$ is
$Y$-reduced. Then $A$ cannot be represented in the form $UVU'$
where $|U|_Y>0$ and $U'U$ is equal to a $0$-word in ${\bf G}(0)$.
Therefore by Lemma \ref{reduced} (2), $A$ is a shortest word among
all words representing conjugates of $g$. Hence $g$ is not a
conjugate of a shorter word, so $g$ is cyclically minimal.

(ii)$\to $(iii). Suppose that $g$ is cyclically minimal but
$|g^2|_Y<2|g|_Y$. Then $A^2$ is not a shortest word representing
$g^2$. By Lemma \ref{reduced} (1), $A^2$ is not $Y$-reduced, which
means that $A^2$ contains a subword $W$ of $Y$-length $2$ which is
equal to a 0-word in ${\bf G}(0)$. Since $A$ is reduced (by
assumption it is a shortest word representing $g$), $W\equiv U'U$
where $U'$ is a suffix of $A$, $U$ is a prefix of $A$. By Lemma
\ref{ZZ}, $|A|_Y\ge 2$. Hence $A\equiv UVU'$ for some word $V$.
The cyclic shift $VU'U$ of $A$ is not $Y$-reduced, so it
represents an element of ${\bf G}(0)$ which is shorter than $g$.
Hence $g$ is a conjugate of a shorter element, a contradiction.

(iii)$\to$ (i). This follows directly from Lemma \ref{cyclically}
(3).
\endproof

\begin {lm} \label {sdvig}
Suppose that ${\cal R}$ satisfies (Z1.1). Let $x$ and $g^{-1}xg$
be $0$-elements and $g=hf$ where $|g|=|f|+|h|.$ Then $h^{-1}xh$ is
also an $0$-element.
\end{lm}
\proof Let $X, G, H$ be $Y$-reduced words representing the
elements $x,g,h,$ respectively. By Lemma \ref{reduced}, the only
possible reduction of the word $G^{-1}XG$ can be fulfilled in the
middle of it. Then, after $|h|$ such reductions, we will see that
the subword $HXH^{-1}$ is equal to a $0$-word.
\endproof

The next several lemmas deal with the subgroups $\oo(g)$ in
property (Z3).

\begin {lm}\label{welldef} Suppose that ${\cal R}$ satisfies
(Z1.1). Then:

(1) the subgroup $\oo(g)$ is well-defined, i.e. $v\oo(u)v^{-1}$
does not depend on the presentation $vuv^{-1}$ of an element $g$
where $u$ is cyclically $y$-reduced.

(2) For every $h\in {\bf G}(0)$ and every essential $g,$ we have
$\oo(hgh^{-1})=h\oo(g)h^{-1}.$
\end{lm}
\proof Assume that $g=v'u'(v')^{-1}$ and $v'\oo(u')(v')^{-1}\ne
v\oo(u)v^{-1}.$ Then $$\oo(u')\ne (v')^{-1}v\oo(u)v^{-1}v'.$$
Thus, to prove the statement, one may assume that $g$ is
cyclically $Y$-reduced. Therefore one has to obtain equality
$\oo(g)=v\oo(u)v^{-1}$ for the maximal 0-subgroups $\oo(g)$ and
$\oo(u)$ normalized by $g$ and $u,$ respectively. By the symmetry,
it suffices to prove that $v\oo(u)v^{-1}$ is a 0-subgroup, because
it is normalized by $g=vuv^{-1}.$

Let $G, U, V$ be shortest words representing elements $g,u,v.$
Then there is an annular diagram diagram over ${\bf G}(0)$ whose
boundary labels are $G\equiv \Lab(q_1)$ and $U\equiv
\Lab(q_2^{-1})$ and a path $t$, connecting the initial vertices of
$q_1$ and $q_2$ is labeled by $V.$ Arbitrary $Y$-band starting
with $q_1$ ends on $q_2,$ and vise versa, since $g$ and $u$ are
cyclically $Y$-reduced. Since these words are essential, $(q_1)_-$
can be connected with a vertex of $q_2$ by a path $p$ with label
$P$ where $|P|=0$. It defines a decomposition $U\equiv  U_1U_2.$
Hence $V=PU_1^{-1}U^{k}$ in ${\bf G}(0)$ for some integer $k$ by
Lemma 11.4 from \cite{book}.

Now, $U^k\oo(U)U^{-k}=\oo(U),$ by the definition of
$\oo(U)=\oo(u).$ Then $U_1^{-1}\oo(U)U_1$ is also a 0-subgroup by
Lemma \ref{sdvig} since the word $U\equiv  U_1U_2$ normalizes
$\oo(U).$ Finally, $V\oo(U)V^{-1}= P(U_1^{-1}\oo(U)U_1)P^{-1}$ is
a 0-subgroup because $P$ is a $0$-word.

The second statement of the lemma immediately follows from the
first one.
\endproof

\begin{lm} \label{Z3.5a}
Suppose that ${\cal R}$ satisfies (Z1.1). If $g$ is an essential
cyclically $Y$-reduced element of ${\bf G}(0)$ and $x\in \oo(g)$
then $g^sx$ is an essential and cyclically $Y$-reduced element for
$s\ne 0$ .
\end{lm}

\proof Suppose that $g^sx=h^{-1}uh$ for $0$-element $u$ and $h\in
{\bf G}(0)$. Then the $Y$-length of $(g^sx)^m$ is at most $2|h|_Y$
for every integer $m\ge 0$. On the other, we have
$(g^sx)^m=g^{sm}x_1$ for some $0$-element $x_1\in \oo(g)$ since
the 0-subgroup $\oo(g)$ is normalized by $g$ . Therefore
$g^{sm}=h\iv u^mhx_1\iv$ and the $Y$-length of $g^{sm}$ is at most
$2|h|_Y$ for every integer $m$. This contradicts the fact that
$g^{sm}$ is cyclically $Y$-reduced by the second part of Lemma
\ref{cyclically}.

Similarly if $g^sx$ is not cyclically $Y$-reduced then $|(g^s
x)^2|_Y<2|g^sx|_Y=2|g^s|_Y$. Therefore for some $0$-element $x_1$
we have $|g^{2s}|_Y=|g^{2s}x_1|_Y=|(g^s x)^2|<2|g^s|_Y$ which
contradicts Lemma \ref{cyclically} (2).
\endproof

\begin{lm} \label{4.8.1}
Suppose  that ${\cal R}$ satisfies (Z1.1). If $g$ is an essential
element and $x\in \oo(g)$ then $\oo(gx)=\oo(g)$.
\end{lm}

\proof By Lemma \ref{welldef} we can assume that $g$ is cyclically
$Y$-reduced. Since $x\in \oo(g)$, the subgroup $\oo(g)$ which is
normalized by $g$, is also normalized by $gx$. Since by Lemma
\ref{Z3.5a}, $gx$ is also cyclically $Y$-reduced, we have
$\oo(g)\subseteq \oo(gx)$. In particular $x\in \oo(gx)$, so
$x\iv\in\oo(gx)$. Therefore, as above,
$\oo(gx)\subseteq\oo(gxx\iv)=\oo(g)$. Hence $\oo(g)=\oo(gx)$.
\endproof

\begin{lm} \label{Z3.5c}
Suppose  that ${\cal R}$ satisfies (Z1.1). If $g$ is an essential
element then for every integer $m\ne 0,$ element $g^m$ is
essential and $\oo(g^m)=\oo(g)$
\end{lm}

\proof By Lemma \ref{welldef} we can assume that $g$ is cyclically
$Y$-reduced. By Lemma \ref{cyclically} (2) one may assume that
$g^m$ is also cyclically $Y$-reduced. Then  $g^m$ is essential by
lemmas \ref{cyclically} (2) and \ref{reduced}, and subgroups
$\oo(g)$ and $\oo(g^m)$ are 0-subgroups. It is clear that
$\oo(g)\subseteq \oo(g^m)$. Let $x\in \oo(g^m)$. Then $x\in
\oo(g^{ms})$ for any $s\ne 0$. Take any integer $t\ne0$. We need
to show that $g^{-t}xg^{t}$ is a $0$-word. But we may choose $s$
so that $g^{ms}=g^tg^{t'}$ for some $t'$ and
$|g^{ms}|=|g^t|+|g^{t'}|$ by Lemma \ref{reduced}. Then
$g^{-t}xg^t$ is a 0-word by Lemma \ref{sdvig}.
\endproof

\begin{lm} \label{Z3.5d}
Suppose that ${\cal R}$ satisfies (Z1.1). If $g$ is an essential
element of ${\bf G}(0)$ and $x\in \oo(g)$, $l\ne 0$, then $g^l x$
is an essential element and $\oo(g^l x)=\oo(g)$.
\end{lm}

\proof The claim follows immediately from lemmas \ref{Z3.5c},
\ref{4.8.1}, \ref{Z3.5a} and \ref{welldef}.
\endproof

\begin{lm} \label{star}
Suppose that ${\cal R}$ satisfies conditions (Z1.1), (Z1.2) and
(Z3.2). Then for every essential element $g\in {\bf G}(0)$ which
is cyclically minimal in rank $1/2$, and every $x\in \oo(g)$, we
have $(gx)^n=g^n$.
\end{lm}

\proof By the definition of $\oo(g)$ and Lemma \ref{welldef}, we
can assume that $g$ is cyclically reduced and $x$ is a 0-element.
We have $(gx)^n=g^n (g^{-n+1}xg^{n-1})...(g\iv xg)x$. By by
property (Z3.2), there exists a $0$-element $r$ such that $gr\iv$
commutes with every element of $\oo(g)$. Therefore $r$ also
normalizes $\oo(g)$ and the actions of $r$ and $g$ on $\oo(g)$ (by
conjugation) coincide. Hence for every integer $s$ the actions of
$r^s$ and $g^s$ on $\oo(g)$ also coincide, so $g^{-s}xg^s =
r^{-s}xr^s$. Therefore $$ (g^{-n+1}xg^{n-1})...(g\iv
xg)x=(r^{-n+1}xr^{n-1})...(r\iv xr)x=r^{-n}(rx)^{n}.$$ By property
(Z1.2) $(rx)^n=r^n=1$ in ${\bf G}(0)$. Hence $(gx)^n=g^n$.
\endproof

\begin{lm} \label{star1}
Suppose that ${\cal R}$ satisfies conditions (Z1.1), (Z1.2),
(Z3.2). Then for every essential element $g\in {\bf G}(0)$ which
is cyclically minimal in rank $1/2$, and every $x\in \oo(g)$ we
have $g^nx=xg^n$.
\end{lm}

\proof It follows from Lemma \ref{star} that $g^n$ commutes with
both $g$ and $gx$. Therefore it commutes with $x$.
\endproof

\begin{lm}\label{C}
Suppose that ${\cal R}$ satisfies (Z1.1). Let $A$ be a $Y$-reduced
essential word over ${\cal X}$ (i.e. it represents an essential
element in ${\bf G}(0)$) and let $C$ be a 0-word. Suppose that the
subgroup generated by $A^{-s}CA^s$ for all integers $s$ consists
of $0$-elements. Then $C\in \oo(A)$.
\end{lm}

\proof If $A$ is cyclically $Y$-reduced then $C\in \oo(A)$ by
definition. Otherwise by Lemma \ref{reduced} (2), $A$ can be
represented as $UVU'$ where $U'U$ is equal to a 0-word $W$ in
${\bf G}(0)$ and $V$ is cyclically $Y$-reduced, so every cyclic
shift of $V$ is $Y$-reduced. Since we assume that $A$ is
essential, $VW$ is an essential word. By Lemma \ref{reduced} (1)
then every cyclic shift of $VW$ is $Y$-reduced. Hence $VW$ is
cyclically $Y$-reduced. By Lemma \ref{sdvig}, $U\iv C U$ is a
$0$-element. Since the subgroup generated by $A^{-s}CA^{s}$, $s\in
\mathbb{Z}$ consists of $0$-elements and is normalized by $A$,
$U(VW)^{-s}(U\iv CU) (VW)^s)U\iv$ is a $0$-element for every $s$.

Notice that $U(VW)^sU\iv$ is a $Y$-reduced form of $A^s$. Indeed,
$V$ is cyclically $Y$-reduced subword of $A$ by definition, $VW$
is a conjugate of $A$ of the same $Y$-length as $V$. Hence by
Lemma \ref{reduced} (2), $VW$ is cyclically $Y$-reduced. By Lemma
\ref{cyclically} (2) then $(VW)^s$ is cyclically $Y$-reduced. By
Lemma \ref{reduced}(1), $U(VW)^sU\iv$ is reduced because every
subword of $Y$-length $2$ of $U(VW)^sU\iv$ is either a subword of
$A$ or a subword of $(VW)^s$.

Since  $A^{-s}CA^s = U(VW)^{-s}(U\iv CU) (VW)^s)U\iv$ is a
$0$-element for every $s$, and $U'\equiv WU\iv$ is a suffix of
$A^s$, by Lemma \ref{sdvig}, the subgroup generated by $(VW)^{-s}
(U\iv CU) (VW)^s$ consists of $0$-elements and is normalized by
$U\iv AU=VW$. Hence $U\iv CU$ is in $\oo(U\iv AU)$. Therefore
$C\in \oo(A)$ by Lemma \ref{welldef}.
\endproof

\begin{lm} \label{Z}
Suppose that ${\cal R}$ satisfies (Z1.1), (Z1.2), (Z3.1) and
(Z3.2). Let $g$ be an essential element of ${\bf G}(0)$
represented by a cyclically $Y$-reduced word $A$ in ${\bf G}(0)$.
Suppose that $A$ contains a subword $B^6$, where $B$ is a
cyclically minimal in rank $1/2$ word of a positive $Y$-length.
Let $A'$ be the word obtained by replacing that occurrence of
$B^6$ by $B^{6-n}$. Then $\oo(A)\le \oo(A'),$ provided the word
$A'$ represents an essential element of ${\bf G}(0)$.
\end{lm}

\proof Let $x\in \oo(A)$. Then for every integer $s$, $g^{-s}xg^s$
is a $0$-element. Fix $s>0$.  Let $C$ be a $0$-word representing
$x$, $D$ be the word representing $g^{-s}xg^s$. Notice that $B^6$
contains the fifth power of every cyclic shift of $B$. Let $B_1$
be a cyclic shift of $B$ starting with a letter from $Y$. Let
$A\equiv A_1B_1^5A_2$. Clearly then the word $A'$ obtained from
$A$ by replacing $B^6$ by $B^{6-n}$ is equal to $A_1B_1^{5-n}A_2$
in the free group.

Consider a \vk diagram over ${\bf G}(0)$ with boundary label
$p_1q_1p_2\iv q_2\iv$ where $\Lab(p_1)=C$, $\Lab(p_2)=D$,
$\Lab(q_1)=\Lab(q_2)=A^s$. The path $q_1$ can be represented as
$$q_1=t_{1,1}t_{1,2}t_{1,3}t_{1,4}t_{2,1}t_{2,2}t_{2,3}t_{2,4}...t_{s,1}t_{s,2}t_{s,3}t_{s,4}$$
where $\Lab(t_{i,2})\equiv B_1^4$, $\Lab(t_{i,3})\equiv B_1$,
$i=1,...,s$ ($\Lab(t_{i,2}t_{i,3})$ is the fixed occurrence of
$B_1^5$ in $A$).

Since $A^s$ is $Y$-reduced and $p_1,p_2$ do not contain $Y$-edges,
the maximal $Y$-bands starting on $q_1$ provide a one-to-one
correspondence between $Y$-edges of $q_1$ and $Y$-edges of $q_2$.
Since $B_1$ starts with a $Y$-letter, $t_{i,2}$ and $t_{i,3}$
start with a $Y$-edges, $e_i$ and $f_i$ respectively, $i=1,...,s$.

For every $i=1,...,s$, consider the maximal $Y$-band $\ttt_i$
starting on $e_i$ and the maximal $Y$-band $\bb_i$ starting on
$f_i$. Consider the smallest subdiagram $\Delta_i$ bounded by a
side of $\ttt_i$, a side of $\bb_i$, $q_1$, $q_2$, containing
$\ttt_i$ but not containing $\bb_i$. The boundary of each
$\Delta_i$ has the form $p_i't_{i,2}(p_i'')\iv (t'_i)\iv$ where
$\Lab(t_{i,2})\equiv \Lab(t_i')\equiv B_1^4$, $v_i\equiv
\Lab(p_i')$ and $v_i'\equiv \Lab(p''_i)$ are $0$-words. Therefore

$$B_1^{-4}v_iB_1^4=v_i'$$ for every $i=1,...,s$.

By Property (Z3.1) (all premises of (Z3.1) hold), $v_i\in
\oo(B_1)$. By Lemma \ref{star1}, $v_i$ commutes with $B_1^n$. For
every $i=1,...,s$ consider a \vk diagram $\Gamma_i$ over ${\bf
G}(0)$ with boundary $\bar p_i \bar q_i (\bar p_i')\iv (\bar
q_i')\iv$ where $\Lab(\bar p_i)\equiv  \Lab(\bar p_i')\equiv v_i$,
$\Lab(\bar q_i)\equiv \Lab(\bar q_i')\equiv B_1^{-n}$. Now cut the
diagram $\Delta$ along each path $p_i$. Thus each path $p'_i$
turns into two copies of $p'_i$, which we shall call the left copy
and the right copy: the right copy is contained in the diagram
containing the band $\ttt_i$. For every $i=1,...,s$, identify
paths $\bar p_i$ of $\Gamma_i$ with the left copy of $p_i$, $\bar
p_i'$ with the right copy of $p_i$, $i=1,...,s$. As a result, we
get a \vk diagram with boundary label $(A_1B_1^{-n}B_1^5A_2)^{-s}
C(A_1B_1^{-n}B_1^5A_2)^sD\iv$. Since as we have noticed,
$A_1B_1^{5-n}A_2=A'$, we have that $(A')^{-s}C(A')^s=D$ in ${\bf
G}(0)$ for every positive $s$. The same result for negative $s$
can be proved similarly.

Hence the subgroup generated by $(A')^{-s}C(A')^{s}$, $s\in
\mathbb{Z}$, consists of $0$-elements and is normalized by $A'$.
By Lemma \ref{C}, $C\in \oo(A')$ as required. \endproof

\begin{lm} \label{Z3.6}
Suppose that ${\cal R}$ satisfies (Z1.1), (Z2.1), (Z2.2), (Z2.3).
Let $g$ be an essential element of ${\bf G}(0)$ represented by a
cyclically $Y$-reduced word $A$ in ${\bf G}(0)$. Suppose that $A$
contains a subword $B$ which is a prefix of a cyclic shift $C$ of
a hub, $C\equiv BB'$, and $|B|\ge \varepsilon |C|$. Let $A'$ be
the word obtained by replacing that occurrence of $B$ by
$(B')\iv$. Then, provided $A'$ represents an essential element of
${\bf G}(0),$ we have inclusion $\oo(A)\le \oo(A')$.
\end{lm}

\proof The proof is similar to the proof of Lemma \ref{Z}. Let
$x\in \oo(A)$. Then for every integer $s$, $g^{-s}xg^s$ is a
$0$-element. Fix $s>0$. Let $C$ be a $0$-word representing $x$,
$D$ be the word representing $g^{-s}xg^s$. Since letters not from
$Y$ have $Y$-length $0$, we can assume that $B$ starts and ends
with $Y$-letters. Let $A\equiv A_1BA_2$. Clearly then the word
$A'$ obtained from $A$ by replacing $B$ by $(B')\iv$ is equal to
$A_1(B')\iv A_2$ in the free group.

Consider a \vk diagram with boundary label $p_1q_1p_2\iv q_2\iv$
where $\Lab(p_1)\equiv C$, $\Lab(p_2)\equiv D$, $\Lab(q_1)\equiv
\Lab(q_2)\equiv A^s$. The path $q_1$ can be represented as
$$q_1=t_{1,1}t_{1,2}t_{1,3}t_{2,1}t_{2,2}t_{2,3}...t_{s,1}t_{s,2}t_{s,3}$$
where $\Lab(t_{i,2})\equiv B$, $i=1,...,s$ ($\Lab(t_{i,2})$ is the
fixed occurrence of $B$ in $A$).

Since $A^s$ is $Y$-reduced by Lemma \ref{cyclically} (2), and
$p_1,p_2$ do not contain $Y$-edges, the maximal $Y$-bands starting
on $q_1$ provide a one-to-one correspondence between $Y$-edges of
$q_1$ and $Y$-edges of $q_2$. Since $B$ starts and ends with a
$Y$-letter, $t_{i,2}$ starts and ends with a $Y$-edges, $e_i$ and
$f_i$ respectively, $i=1,...,s$.

For every $i=1,...,s$, consider the maximal $Y$-band $\ttt_i$
starting on $e_i$ and the maximal $Y$-band $\bb_i$ starting on
$f_i$. Consider the smallest subdiagram $\Delta_i$ bounded by a
side of $\ttt_i$, a side of $\bb_i$, $q_1$, and $q_2$, containing
$\ttt_i$ and $\bb_i$. The boundary of each $\Delta_i$ has the form
$p_i't_{i,2}(p_i'')\iv (t'_i)\iv$ where $\Lab(t_{i,2})\equiv
\Lab(t_i')\equiv B$, $v_i\equiv \Lab(p_i')$ and $v_i'\equiv
\Lab(p''_i)$ are $0$-words. Therefore

$$v_iB=Bv_i'$$ for every $i=1,...,s$.

By Property (Z2.3) (all premises of (Z2.3) hold), we have

$$v_i(B')\iv=(B')\iv v_i', i=1,...,s.$$

Thus for every $i=1,...,s$ there exists a \vk diagram $\Gamma_i$
over ${\cal R}$ with boundary $p_i' \bar t_{i,2} (p_i'')\iv (\bar
t_i')\iv$ where $\Lab(\bar t_{i,2})\equiv (B')\iv$,  $\Lab(\bar
t_i')\equiv (B')\iv$. Now we can substitute the subdiagrams
$\Delta_i$ in $\Delta$ by $\Gamma_i$. The boundary label of the
resulting diagram is equal $(A')^{-s}C(A')^sD\iv$. Hence
$(A')^{-s}C(A')^s=D$ in ${\bf G}(0)$ for every positive $s$. The
same result for negative $s$ can be proved similarly.

Hence the subgroup generated by $(A')^{-s}C(A')^{s}$, $s\in
\mathbb{Z}$, consists of $0$-elements and is normalized by $A'$.
By Lemma \ref{C} $C\in \oo(A')$ as required. \endproof

\subsection{The construction of a graded presentation}

\label{construction}

In this section, we show how the theory from \cite[Chapter
6]{book} can be adapted to diagrams over a presentation satisfying
(Z1), (Z2), (Z3).

We are defining several concepts analogous to concepts in
\cite{book}, chapter 6, by induction on the rank $i$.

We say that a certain statement {\it holds in rank $i$} if it
holds in the group \label{gi}${\bf G}(i)$ given by ${\cal
R}_i=\cup_{j=0}^i
 {\cal S}_j$.

Recall that  a word $w$  is \label{per} $A$-{\em periodic} if it
is a subword of some power of the word $A.$ Any $A$-periodic word
can be read (say, clockwise) on a cyclic graph with $|A|$ edges,
labeled, starting with a vertex $o,$ by the word $A.$ Then every
$A$-periodic word corresponds to a path on this graph. Let
$W\equiv  W_1W_2$ be a decomposition of an $A$-periodic word $W,$
such that the word $W_1$ ends at $o,$ when reading the word $W$ on
the cyclic graph, $|W|\ge|A|$. Then we say that the decomposition
is a \label{phase}{\it phase} decomposition of $W.$ If the label
of a simple path $p$ of an arbitrary labeled graph has an
$A$-periodic label $W\equiv \phi(p),$ then a vertex of $p,$
defining a phase decomposition of $p$ is said to be a {\it phase}
vertex of $p.$

{\em Up to the end of this section, the $Y$-length of a word $W$
will be called simply {\em length} and will be denoted by $|W|$.}

  Let us define simple words of rank $i=0,1/2,1,2,...$ and periods of rank
$i=1,2,3,\dots$. Let ${\cal S}_0={\cal R}_0$ and ${\cal S}_{1/2}$
be the set of hubs. Thus the groups ${\bf G}(0)$ and ${\bf
G}(1/2)={\bf G}$ are defined by relations of the sets ${\cal R}_0$
and ${\cal R}_{1/2}={\cal S}_0\cup{\cal S}_{1/2},$ respectively. A
word $A$ of positive length is said to be {\it simple} in rank 0
(resp. 1/2) if it is not a conjugate in rank 0 (resp. 1/2) either
of a word of length 0 or of a product of the form $B^l P$, where
$|B|<|A|$ and $P\in \oo(B)$ in rank 0.

Suppose that $i\ge 1$, and we have defined the sets of relators
\label{ri}${\cal R}_j$, $j<i$, and the corresponding groups ${\bf
G}(j)$ of ranks $j<i$. Suppose also that we have defined simple
words of ranks $j$, $j<i$, and periods of rank $j$, $1\le j< i$.

For every $i = 1/2,1,2,3,..$ let $i_-$ denote $i-1$ if $i>1$,
$1/2$ if $i=1$ or $0$ if $i=1/2$. Similarly $i_+$ is $1/2$ if
$i=0$ or $1$ if $i=1/2$ or $i+1$ if $i\ge 1$.

For every $i\ge 1$ consider the set ${\cal X}'_i$ of all words of
length $i$ which are simple in rank $i_-$, and the
\label{sim}equivalence $\sim_{i_-}$ given in Lemma \ref{Equival}
(we apply it for the smaller rank here).  Now choose a set of
representatives \label{xi}${\cal X}_i$ of the equivalence classes
of ${\cal X}'_i$. The words of ${\cal X}_i$ are said to be
\label{period1}{\it periods} of rank $i.$ The set of words ${\cal
R}_i$ (defining the group ${\bf G}(i)$ or rank $i$) is the union
of ${\cal R}_{i_-}$ and \label{si}${\cal S}_i= \{A^n, \ A\in {\cal
X}_i\}$.

Let $A$ be a word of positive length. We say that $A$ is
\label{simple1}{\it simple} in rank $i\ge 1$ if it is not a
conjugate in rank $i$ either of a word of length 0 or of a word of
the form $B^lP,$ where $B$ is a period of rank $j\le i$ or an
essential word with $|B|<|A|$, and $P$ represents a word of the
subgroup $\oo(B).$

Let \label{rr}${\cal R}(\infty)=\cup_0^\infty {\cal R}_i$. If
$\Pi$ is a cell of a diagram over ${\cal R}(\infty)$ with the
boundary labeled by a word of the set ${\cal S}_{j},$ then, by
definition, $r(\Pi)=j.$ The boundary label of every cell of a rank
$\ge 1$ (read counterclockwise), has a period $A$ defined up to
cyclic shifts.

\medskip

For $j\ge 1$, a pair of distinct cells $\Pi_1$ and $\Pi_2$ of rank
$j$ of a diagram $\Delta$ is said to be a \label{jp1}$j$-{\it
pair}, if their counterclockwise contours $p_1$ and $p_2$ are
labeled by $A^n$ and $A^{-n}$ for a period $A$ of rank $j$ and
there is a path $t=(p_1)_- -(p_2)_-$ without self-intersections
such that $\phi(t)$ is equal, in rank $j_-$, to an element of
$\oo(A).$ Then $\phi(t)$ and $A^n$ commute in rank $j_-$ by Lemma
\ref{star1}, and so the subdiagram with contour $p_1tp_2t^{-1}$
can be replaced in $\Delta$ by a diagram of rank $j_-$. As a
result, we obtain a diagram with the same boundary label as
$\Delta$ but of a smaller type.

Similarly, a pair of hubs $\Pi_1,$ $\Pi_2$ with counterclockwise
contours $p_1$, $p_2$ forms a \label{12p}$\frac 12$-{\it pair}, if
the vertices $(p_1)_-$ and $(p_2)_-$ are connected by a path $t$
without self-intersections such that the label of the path
$p_1tp_2t^{-1}$ is equal to 1 in ${\bf G}(0).$ Consequently, any
diagram can be replaced by a diagram having no $j$-pairs, $j=1/2,
1, 2,...$, i.e. by a \label{gred}{\it g-reduced} diagram with the
same boundary label(s)\footnote{Recall that in \cite{book}, $\S
13.2,$ such diagrams are called {\em reduced}. But in this paper
we shall consider two kinds of reduced diagrams, so we call
diagrams without $j$-pairs g-reduced; ``g" stands for ``graded".}.

Assume that $q_1$ and $q_2$ are disjoint sections of the boundary
$\partial(\Delta)$ of a diagram $\Delta$ of rank $i,$ and
$\phi(q_1), \phi(q_2)$ are $A$- and $A\iv$-periodic words
respectively, where $|q_1|,|q_2|\ge |A|$, $A$ is simple in rank
$i$ or a period of rank $j\le i$ . We call $q_1$ and $q_2$
\label{comp1}{\it compatible} if there exist phase vertices $o_1$
and $o_2$ on $q_1$ and $q_2,$ respectively, and a path $t=o_1-o_2$
without self-intersections, such that $\phi(t)$ is equal in rank
$i$ to an element of $\oo(A)$; moreover the equality must be true
in rank $j_-$ if $A$ is a period of rank $j.$

Similarly one defines the compatibility of an $A$-periodic cell
$\Pi$ of rank $j\ge 1$ with a $A^{\pm 1}$-periodic section $q$ of
$\partial(\Delta)$.

\medskip

\begin{lm}\label{anal13.2}  Let $\Pi_1$ be a cell of rank $j\ge1$ in a
g-reduced diagram $\Delta$ of rank $i$. Assume it has a boundary
label $A^n$. Let $\Gamma$ be a subdiagram of $\Delta$ with contour
$p_1q_1p_2q_2$, where $q_1$ is a subpath of $\partial(\Pi_1)$ and
$q_2$ is a subpath of the contour of a cell $\Pi_2$ with the
boundary label $A^{-n}$, ($\Pi_1,\Pi_2$ do not belong to
$\Gamma$). Then

1) $q_1$ and $q_2$ cannot be compatible in $\Gamma$;

2) the subdiagram $\Gamma$ contains no cells compatible with
$q_1.$
\end{lm}

\proof If one of the statements were false, we could include the
cell $\Pi_1$ into a $j$-pair and come to an obvious contradiction
with the assumption that $\Delta$ is g-reduced diagram, as in
Lemma 13.2 from \cite{book}. \endproof

\begin{lm}\label{anal13.3} Let $p_1,p_2\dots$ be a partition of the
boundary of a diagram $\Delta$ of rank $i$ into several sections,
and $\phi(p_1)\equiv  A^l$ for some period of rank $j$. Assume
that a cell $\Pi$ of $\Delta$ is compatible with $p_1$. Then there
exists a g-reduced diagram $\Delta'$ with topologically identical
partition of the boundary into $p'_1,p'_2,\dots,$
 such that $\tau(\Delta')<\tau(\Delta),$ $\phi(p'_k)\equiv  \phi(p_k)$
for $k\ge 2,$ $\phi(p'_1)\equiv  A^{l+sn}$ for some integer $s,$
and $\Delta'$ has no cells compatible with $p_1.$
\end{lm}

\proof Let $t$ be the path connecting a phase vertex $o_1$ of
$p_1$ and a phase vertex $o_2$ of the $A^{\pm 1}$-periodic cell
$\Pi,$ according to the definition of the compatibility, and
$T\equiv \phi(t)$ is equal in rank $i$ to an element of $\oo(A).$
The vertex $o_1$ defines a decomposition $A^{l_1}A^{l_2},$ and
there is a path $p$ homotopic to $p_1$ with label $W\equiv
A^{l_1}TA^{\pm n}T^{-1}A^{l_2}.$ By Lemma \ref{star1}, we have
$W=A^{l\pm n}$ in rank $j_-.$ Thus one can cut $\Pi$ out of
$\Delta,$ by cutting along $p$. Then one can replace the deleted
subdiagram by the diagram of the equality $W=A^{l\pm n}.$ By
continuing this process, we prove the lemma.
\endproof

The following proposition is the main statement of this section.

\begin{prop}\label{mainprop} If the presentation ${\cal R}$ of a group ${\bf G}$ satisfies
properties (Z1), (Z2), (Z3) then the graded presentation ${\cal
S}_0\cup {\cal S}_{1/2}\cup {\cal S}_1\cup {\cal S}_2\cup...$
satisfies statements of Lemmas \ref{$V_2UV_1$} through
\ref{anal19.5} below and is such that ${\cal R} = {\cal S}_0\cup
{\cal S}_{1/2}$. In particular, all g-reduced diagrams over this
presentation are $A$-maps. Every relator from ${\cal S}_i$, $i\ge
1$, is of the form $A^n$ for a cyclically reduced word $A$ of
$Y$-length $i$. The group defined by this presentation is ${\bf
G}/\la g^n, g\in {\bf G}\ra$.
\end{prop}

{\bf Lemmas from \ref{$V_2UV_1$} to \ref{anal19.5} will be proved
by simultaneous induction on rank $i\ge \frac{1}{2}.$}

\medskip

\begin{lm}\label{$V_2UV_1$} Let $V$ be a periodic word with a
simple in rank $i$ period. Assume that  for some word $U$ with
$2|U|\le |V|$, the word $VU$ is a conjugate in rank $i$ of a
0-word. Then there is a decomposition $V\equiv  V_1V_2,$ such that
$||V_1|-|V_2||<\beta|V|+|U|$ and for some word $S$ of length at
most $3\gamma |V|,$ $V_2UV_1=S$ in rank $i.$
\end{lm}

\proof Consider a g-reduced annular diagram $\Delta$ for the
conjugacy of $VU$ to a 0-word $W.$ We may assume that $\Delta$ has
the smallest possible type. Therefore there are no $Y$-bands
surrounding the hole of $\Delta,$ since otherwise their boundary
labels were 0-words. Hence, by lemmas \ref{anal19.4} and
\ref{anallemma17.1}, there are vertices $o_1$ on the contour $q_1$
labeled by $VU$ and $o_2$ on the contour $q_2$ labeled by $W,$
which can be connected by a simple path $t$ of length at most
$\gamma(|V|+|U|)\le 3\gamma|V|/2$ (see
\setcounter{pdeleven}{\value{ppp}} Figure \thepdeleven).

One may suppose that the vertex $o_2$ does not cut the word $W,$
because $|W|=0.$ The vertex $o_1$ cannot cut the word $|U|.$
Indeed, in case it defines a decomposition $U\equiv  U_1U_2,$ we
obtain the following equality in rank $i$:
$U_1TW^{-1}T^{-1}U_2=V^{-1}, $ where $T\equiv \phi(t).$ It follows
by lemmas \ref{anal19.5} and \ref{analteom17.1} that
$\bar\beta|V|<|U|+2\gamma(|V|+|U|),$ whence
$|V|<(\bar\beta-2\gamma)^{-1}(1+2\gamma)|U|\le 2|U|,$ contrary to
the hypothesis of the lemma.

\unitlength=1mm \special{em:linewidth 0.4pt} \linethickness{0.4pt}
\begin{picture}(98.00,37.00)
\put(42.00,7.00){\line(1,0){53.00}}
\put(95.00,7.00){\line(0,1){23.00}}
\put(95.00,30.00){\line(-1,0){53.00}}
\bezier{120}(42.00,30.00)(24.00,30.00)(23.00,18.00)
\bezier{120}(23.00,18.00)(23.00,7.00)(42.00,7.00)
\put(26.00,25.00){\circle*{1.00}} \put(34.00,18.00){\circle{8.94}}
\put(34.00,22.00){\circle*{1.00}}
\bezier{76}(26.00,25.00)(38.00,27.00)(34.00,22.00)
\put(40.00,15.00){\makebox(0,0)[cc]{$W$}}
\put(30.00,23.00){\makebox(0,0)[cc]{$t$}}
\put(23.00,27.00){\makebox(0,0)[cc]{$o_1$}}
\put(34.00,20.00){\makebox(0,0)[cc]{$o_2$}}
\put(64.00,33.00){\makebox(0,0)[cc]{$V_1$}}
\put(98.00,18.00){\makebox(0,0)[cc]{$U$}}
\put(75.00,10.00){\makebox(0,0)[cc]{$V_2$}}
\put(95.00,16.00){\vector(0,1){8.00}}
\put(83.00,30.00){\vector(-1,0){18.00}}
\put(58.00,7.00){\vector(1,0){12.00}}
\end{picture}

\begin{center}
\nopagebreak[4] Fig. \theppp.

\end{center}
\addtocounter{ppp}{1}

Thus the vertex $o_1$ defines a decomposition $V\equiv  V_1V_2,$
and $TW^{-1}T^{-1}V_2U=V_1^{-1}$ in rank $i.$ Then, as above,
$\bar\beta|V_1|< |V_2|+|U|+2\gamma(|V|+|U|)$. Therefore
$|V_1|-|V_2|<\beta|V_1|+|U|+3\gamma|V|$. Similarly,
$|V_2|-|V_1|<\beta|V_1|+|U|+3\gamma|V|,$ and so
$||V_1|-|V_2||<\beta |V|+|U|$, as $2|U|\le |V|.$

To complete the proof we observe that the product $V_2UV_1$ is
equal in rank $i$ to the word $S\equiv  TWT^{-1}$ of length at
most $3\gamma |V|$. \endproof

Let $q_1,$ $q_2$ be two paths labeled by $A^{\pm 1}$-periodic
words for some period $A.$ We say that vertices $o_1$, $o_2$ of
these paths are in the \label{sp}{\it same phase} (with respect to
$A$) if they are phase vertices for the paths $q_1,$ $q_2$ with
respect to a cyclic permutation $A'$ of the word $A.$

\begin{lm}\label{1/3|Pi|} Let $\Pi$ be a cell of positive rank
 in a g-reduced disc
diagram $\Delta$ of rank $i$ with contour $p_1q_1p_2q_2$, where
$q_1$ and $q_2^{-1}$ are labeled by $A$-periodic words where $A$
is a simple in rank $i$ word. Suppose there are disjoint
contiguity subdiagrams $\Gamma_j$ of $\Pi$ to $q_j$ of degrees
$\psi_j$ for $j=1,2,$ with $\psi_1+\psi_2>1-\alpha/2.$ Denote
$s_1^jt_1^js_2^jt_2^j=\partial(\Pi, \Gamma_j,q_j)$ (see
\setcounter{pdeleven}{\value{ppp}} Figure \thepdeleven). Let $p$
be a subpath of $q_2^{-1},$ such that vertex $o=(s_2^2)_+$ belongs
to $p,$ where the lengths of subpaths $p'=p_--o$ and $p''=o-p_+$
are less than $\frac13 |\partial(\Pi)|.$ Then any vertex $o_2$ of
$p$ is not in the same phase with  vertex $o_1=(s_1^1)_-$.

\end{lm}
\proof By lemmas \ref{anal19.4}, \ref{anal19.5} and
\ref{anal15.8},

\begin{equation}
\psi_1,\psi_2\in [1/2-2\alpha,\bar\alpha], \label{(??.1)}
\end{equation}
and by Lemma \ref{anal15.4},

\begin{equation}
|t_2^1|, |t_2^2|\in [(1+2\beta)^{-1}(1/2-2\alpha)|\partial(\Pi)|,
(1-2\beta)^{-1}\bar\alpha|\partial(\Pi)|] \label{(??.2)}
\end{equation}

Assume that $o_2$ and $o_1$ are in the same phase and $o_2\in p'.$
First, let $\Pi$ be a hub. Then $r(\Gamma_1)=r(\Gamma_2)=0$ by
Lemma \ref{HubCont} in the smaller rank, hence $|s_l^j|=0$,
$l,j=1,2$. Recall that by (Z2.2), every hub is a linear word. By
definition, every simple word in rank $i$ is cyclically
$Y$-reduced.

Therefore by Lemmas \ref{uchastok} and \ref{cyclically} every
maximal $Y$-band starting on $t_l^j$ ends on $t^j_{3-l}$ and forms
a $0$-bond between these paths. Hence any subpath of $t_l^j$
starting and ending with $Y$-letters, is a contiguity arc for some
contiguity subdiagram between $t_l^j$ and $t_{3-l}^j$ for
$l,j=1,2.$ By our assumption and inequality (\ref{(??.2)}), there
are subpaths $u^1$ and $u^2$ of $t_2^1$ and $(t_2^2)^{-1}$ with
the same label $U$, such that $$|U|>((1+2\beta)^{-1}(1/2-2\alpha)
-1/3)|\partial(\Pi)|>\frac 17|\partial(\Pi)|$$ Gluing along $u_1$
and $u_2$ one can construct a contiguity diagram of rank 0 between
disjoint arcs of the cell $\Pi$ of length greater than $\frac 17
|\partial(\Pi)|.$ But this contradicts condition (Z2.2).

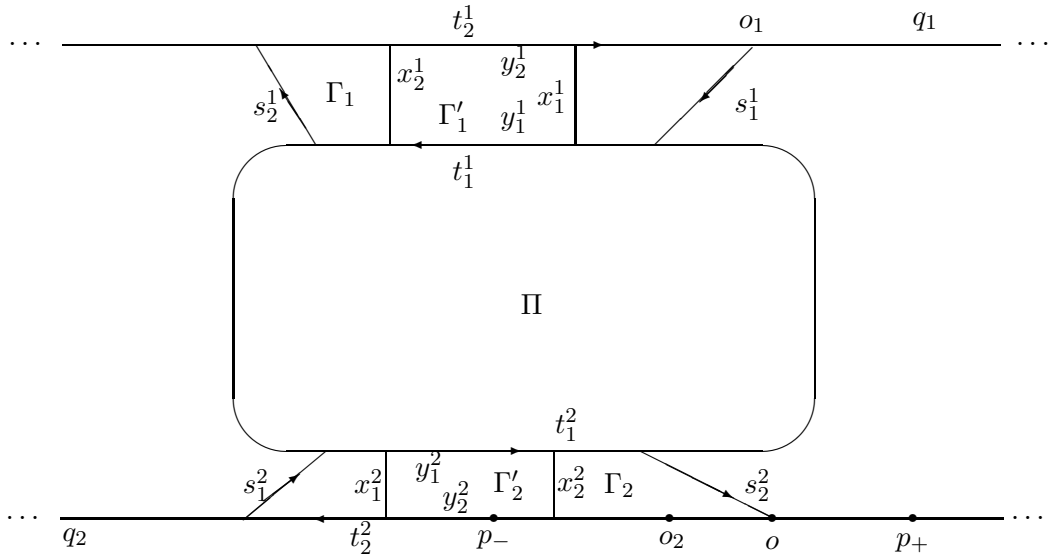
\unitlength=1.00mm \special{em:linewidth 0.4pt}
\linethickness{0.4pt}
\begin{picture}(136.95,80.11)
\put(7.61,7.94){\line(1,0){125.33}}
\put(7.95,70.94){\line(1,0){124.67}}
\put(69.28,37.28){\oval(77.33,40.67)[]}
\put(41.61,57.61){\line(-3,5){8.00}}
\put(86.61,57.61){\line(1,1){13.00}}
\put(42.95,16.94){\line(-6,-5){11.00}}
\put(84.61,16.94){\line(2,-1){18.00}}
\put(40.28,59.94){\vector(-2,3){3.50}}
\put(34.61,10.28){\vector(4,3){4.83}}
\put(97.11,67.94){\vector(-1,-1){4.50}}
\put(88.28,15.11){\vector(2,-1){8.83}}
\put(52.11,16.94){\vector(1,0){16.83}}
\put(71.61,57.61){\vector(-1,0){17.33}}
\put(58.11,70.94){\vector(1,0){21.83}}
\put(51.61,7.94){\vector(-1,0){10.33}}
\put(51.45,57.61){\line(0,1){13.33}}
\put(76.11,57.61){\line(0,1){13.33}}
\put(50.95,16.94){\line(0,-1){9.00}}
\put(73.28,7.94){\line(0,1){9.00}}
\put(70.28,36.44){\makebox(0,0)[cc]{$\Pi$}}
\put(61.28,54.11){\makebox(0,0)[cc]{$t_1^1$}}
\put(54.28,67.11){\makebox(0,0)[cc]{$x_2^1$}}
\put(72.95,64.11){\makebox(0,0)[cc]{$x_1^1$}}
\put(67.95,68.77){\makebox(0,0)[cc]{$y_2^1$}}
\put(67.95,61.44){\makebox(0,0)[cc]{$y_1^1$}}
\put(59.94,61.45){\makebox(0,0)[cc]{$\Gamma_1'$}}
\put(44.95,64.44){\makebox(0,0)[cc]{$\Gamma_1$}}
\put(34.95,62.78){\makebox(0,0)[cc]{$s_2^1$}}
\put(98.95,63.11){\makebox(0,0)[cc]{$s_1^1$}}
\put(61.28,74.11){\makebox(0,0)[cc]{$t_2^1$}}
\put(74.95,20.44){\makebox(0,0)[cc]{$t_1^2$}}
\put(33.61,12.44){\makebox(0,0)[cc]{$s_1^2$}}
\put(100.28,12.44){\makebox(0,0)[cc]{$s_2^2$}}
\put(47.61,5.44){\makebox(0,0)[cc]{$t_2^2$}}
\put(9.61,
5.44){\makebox(0,0)[cc]{$q_2$}}
 \put(65.28,8.11){\circle*{1.00}}
\put(88.61,8.11){\circle*{1.00}} \put(102.28,8.11){\circle*{1.00}}
\put(65.28,5.11){\makebox(0,0)[cc]{$p_-$}}
\put(88.95,5.11){\makebox(0,0)[cc]{$o_2$}}
\put(102.28,4.78){\makebox(0,0)[cc]{$o$}}
\put(120.95,8.11){\circle*{1.00}}
\put(120.95,4.78){\makebox(0,0)[cc]{$p_+$}}
\put(48.61,12.44){\makebox(0,0)[cc]{$x_1^2$}}
\put(75.61,12.78){\makebox(0,0)[cc]{$x_2^2$}}
\put(56.61,14.78){\makebox(0,0)[cc]{$y_1^2$}}
\put(60.28,10.78){\makebox(0,0)[cc]{$y_2^2$}}
\put(67.28,12.44){\makebox(0,0)[cc]{$\Gamma_2'$}}
\put(81.95,12.44){\makebox(0,0)[cc]{$\Gamma_2$}}
\put(99.61,73.78){\makebox(0,0)[cc]{$o_1$}}
\put(122.61,74.11){\makebox(0,0)[cc]{$q_1$}}
\put(2.95,71.11){\makebox(0,0)[cc]{$\dots$}}
\put(2.61,8.11){\makebox(0,0)[cc]{$\dots$}}
\put(136.95,71.11){\makebox(0,0)[cc]{$\dots$}}
\put(136.28,8.11){\makebox(0,0)[cc]{$\dots$}}
\end{picture}

\begin{center}
\nopagebreak[4] Fig. \theppp.

\end{center}
\addtocounter{ppp}{1}

Now assume that $\Pi$ is a cell with a period $B$ of rank $k\le
i$. In this case, $r(\Gamma_1), r(\Gamma_2)<k$ by lemmas
\ref{anal15.3}, \ref{anal17.4} and \ref{anal19.5}. As in the
previous paragraph, there are subpaths $u_1$ and $u_2$ of $t_2^1$
and $(t_2^2)^{-1}$ with the same label $U,$ where $|U|>\frac 17
|\partial(\Pi)|.$

  Lemma \ref{anal15.3} and inequality (\ref{(??.1)}) show that Lemma
  \ref{Anal17.3} is applicable to each of the
diagrams $\Gamma_1, \Gamma_2.$ Since $|s_1^j|<\zeta
|\partial(\Pi)|$, it follows from Lemma \ref{Anal17.3} that an
arbitrary vertex of $t_l^j$ ($l,j \in \{1,2\}$) can be connected
in $\Gamma_j$ with a vertex of $t_{3-l}^j$ by a path of length at
most $3\gamma^{-1}\zeta|\partial(\Pi)|.$ Therefore, for $j=1,2,$
there exist subdiagrams $\Gamma'_j$ of $\Gamma_j$ with contours
$x_1^jy_1^jx_2^jy_2^j,$ where
$|x_l^j|<3\gamma^{-1}\zeta|\partial(\Pi)|,$ $\phi(y_2^1)\equiv
\phi(y_2^2)^{-1}\equiv  U$ and $y_1^j$ are subpaths of $t_1^j.$

  By gluing diagrams $\Gamma'_1$ and $\Gamma'_2$ along $y_2^1$ and
$(y_2^2)^{-1},$ one can obtain (after deleting $j$-pairs) a
g-reduced diagram $\Gamma_0$ with contour $z_1y_1^1z_2y_1^2$ where
$|z_1|, |z_2|<6\gamma^{-1}\zeta|\partial(\Pi)|.$ Notice, that by
property A1 and Lemma \ref{analteom17.1},
$$|y_1^j|\ge\bar\beta|y_2^j|-|x_1^j|-|x_2^j|>\frac 17
\bar\beta|\partial(\Pi)| -6\gamma^{-1}\zeta|\partial(\Pi)|>\frac
18 |\partial(\Pi)|\ge \frac n8 |B|$$ Therefore one can apply Lemma
\ref{anal17.5} to $\Gamma_0$ and  obtain a subdiagram $\Gamma'_0$
of $\Gamma_0$ with contour $v_1w_1v_2w_2,$ where $w_j$ is a
subpath of $y_1^j$ (for $j=1,2$), $|v_1|, |v_2|<\alpha |B|,$ and
$$|w_1|, |w_2|> \frac 18 |\partial(\Pi)|
-\gamma^{-1}(12\gamma^{-1}\zeta|\partial(\Pi)|+|B|)
>\frac n9 |B|$$
since $|\partial(\Pi)|\ge n|B|$ and the rank of $\Pi$ is $k=|B|$.
Since the labels of both $w_1$ and $w_2$ are $B$-periodic words,
and period $B$ is simple in rank $r(\Gamma'_0)<k,$ we obtain a
contradiction with Lemma \ref{anal18.9} in a smaller rank.

 Assume now that $o_2\in p''.$ Then we consider $o'_1= o=(s_2^2)_+$. There
is a vertex $o'_2$ on the path $q_1,$ such that $o'_2-o_1$ is a
subpath of $q_1$ of length $|o_2-o|<\frac 13 |\partial(\Pi)|$.
Then $o'_1$ and $o'_2$ are in the same phase. Thus, this case is
completely analogous to the case when $o_2\in p'$ up to the
interchange of the paths $q_1$ and $q_2^{-1}.$ The lemma is
proved. \endproof

\begin{lm}\label{$o_1-o_2$}  Let $\Delta$ be a narrow
contiguity diagram of rank $i$ between $q_1$ and $q_2$ with
contour $p_1q_1p_2q_2,$ where $q_1$ and $q_2^{-1}$ are paths with
$A$-periodic labels for some simple in rank $i$ period $A.$ Assume
that there are vertices $o_1$ and $o_2$ on $q_1$ and $q_2$
respectively, which are in the same phase, the label
$\phi(o_1-o_2)$ of a path $o_1 - o_2$ is a conjugate of a 0-word
in rank $i$ and $||(q_1)_- - o_1|-|(q_2)_+-o_2||\le \frac 13
\min(|q_1|,|q_2|),$ where $(q_1)_--o_1$ and $o_2-(q_2)_+$ are
subpaths of $q_1$ and $q_2$. Then $r(\Delta)=0$.
\end{lm}

\proof As in the proof of Lemma \ref{Nearly}, $\Delta$ can be
partitioned into several narrow bonds $E_l$ of an $\alpha$-series
and several subdiagrams $\Delta_1, \Delta_2, \dots$ between
successive bonds. If $r(E_l)=r(E_{l+1})=0,$ then $r(\Delta_l)=0$
by Lemma \ref{Anal17.3.2} and the maximality of the
$\alpha$-series. Thus it suffices to prove that $r(E_l)=0$ for
every $E_l.$  Proving by contradiction, we choose the bond $E_l$
whose principal cell $\Pi_0$ has maximal perimeter (see
\setcounter{pdeleven}{\value{ppp}} Figure \thepdeleven).

\unitlength=1.00mm \special{em:linewidth 0.4pt}
\linethickness{0.4pt}
\begin{picture}(119.33,52.33)
\put(0.00,2.67){\line(0,1){41.67}}
\put(0.00,44.33){\line(1,0){116.33}}
\put(116.33,44.33){\line(0,-1){41.33}}
\put(116.33,3.00){\line(-1,0){116.33}}
\put(85.17,23.83){\oval(35.00,23.00)[]}
\put(76.67,35.33){\line(-3,5){5.33}}
\put(99.00,44.33){\vector(-3,-4){6.67}}
\put(76.67,12.33){\line(-3,-5){5.67}}
\put(88.67,12.33){\line(1,-1){9.00}}
\put(71.33,3.00){\circle*{1.00}} \put(97.67,3.00){\circle*{1.00}}
\put(99.00,44.33){\circle*{1.00}}
\put(71.67,44.33){\circle*{1.00}}
\put(30.67,44.33){\circle*{1.00}} \put(52.33,3.00){\circle*{1.00}}
\bezier{104}(30.67,44.33)(38.33,36.33)(39.67,21.00)
\bezier{112}(39.67,20.67)(38.00,8.33)(52.33,3.00)
\put(0.00,10.67){\vector(0,1){19.33}}
\put(15.33,44.33){\vector(1,0){9.67}}
\put(81.00,44.33){\vector(1,0){12.33}}
\put(92.33,3.00){\vector(-1,0){13.00}}
\put(116.33,34.67){\vector(0,-1){15.33}}
\put(37.67,31.33){\vector(-1,3){1.00}}
\put(102.67,22.67){\vector(0,1){4.67}}
\put(77.00,35.33){\vector(-2,3){3.00}}
\put(72.00,4.67){\vector(2,3){3.67}}
\put(91.33,9.67){\vector(1,-1){5.00}}
\put(37.00,3.00){\vector(-1,0){20.33}}
\put(113.67,3.00){\vector(-1,0){12.67}}
\put(67.67,28.00){\vector(0,-1){8.67}}
\put(13.67,20.33){\makebox(0,0)[cc]{$\Delta$}}
\put(-2.67,23.67){\makebox(0,0)[cc]{$p_1$}}
\put(30.67,47.33){\makebox(0,0)[cc]{$o'$}}
\put(52.33,47.33){\makebox(0,0)[cc]{$q_1$}}
\put(84.67,41.67){\makebox(0,0)[cc]{$t_2^1$}}
\put(99.33,47.00){\makebox(0,0)[cc]{$o_1$}}
\put(99.33,39.00){\makebox(0,0)[cc]{$s_1^1$}}
\put(71.00,39.00){\makebox(0,0)[cc]{$s_2^1$}}
\put(86.00,32.33){\makebox(0,0)[cc]{$t_1^1$}}
\put(70.00,24.33){\makebox(0,0)[cc]{$u_1$}}
\put(100.00,24.33){\makebox(0,0)[cc]{$u_2$}}
 \put(85.33,22.67){\makebox(0,0)[cc]{$\Pi_0$}}
\put(81.00,15.00){\makebox(0,0)[cc]{$t_1^2$}}
\put(84.67,5.33){\makebox(0,0)[cc]{$t_2^2$}}
\put(80.67,12.33){\vector(1,0){5.67}}
\put(76.67,0.33){\makebox(0,0)[cc]{$v_2$}}
\put(43.33,0.33){\makebox(0,0)[cc]{$v_1$}}
\put(53.33,5.67){\makebox(0,0)[cc]{$o$}}
\put(21.67,41.00){\makebox(0,0)[cc]{$q_1'$}}
\put(41.00,29.67){\makebox(0,0)[cc]{$t$}}
\put(119.33,25.33){\makebox(0,0)[cc]{$p_2$}}
\put(108.67,5.33){\makebox(0,0)[cc]{$q_2''$}}
\put(109.00,41.67){\makebox(0,0)[cc]{$q_1''$}}
\put(71.00,7.67){\makebox(0,0)[cc]{$s_1^2$}}
\put(95.67,8.67){\makebox(0,0)[cc]{$s_2^2$}}
\put(59.33,3.00){\vector(1,0){8.00}}
\put(33.67,3.00){\circle*{1.00}}
\put(33.33,-0.33){\makebox(0,0)[cc]{$o_2$}}
\put(14.67,5.67){\makebox(0,0)[cc]{$q_2'$}}
\put(12.33,46.33){\makebox(0,0)[cc]{$q_1'$}}
\put(40.00,3.00){\vector(1,0){7.33}}
\end{picture}
\begin{center}
\nopagebreak[4] Fig. \theppp.

\end{center}
\addtocounter{ppp}{1}

By Lemma \ref{analteom17.1} (and lemmas \ref{anal19.4},
\ref{anal19.5}) $|q_2|\le\bar\beta^{-1}(|q_1|+ |p_1|+|p_2|).$ We
have $|p_1|+|p_2|<(\alpha+4\zeta)|\partial(\Pi_0)|$ as in Lemma
\ref{Nearly} for the side arcs of narrow bonds. Denote by
$\Gamma_1$ and $\Gamma_2$ the contiguity subdiagrams of $\Pi_0$ to
$q_1$ and $q_2$ with contours  $s_1^jt_1^js_2^jt_2^j$ ($j=1,2$) as
in Lemma \ref{1/3|Pi|}. Applying lemmas \ref{anal15.4},
\ref{anal15.8}, to $\Gamma_1$, we have
$|q_1|>(1-\alpha-\bar\alpha)|\partial(\Pi_0)|,$ because $q_1$
contains $t_2^1.$ Hence $|p_1|+|p_2|<\frac52\alpha|q_1|$ and
$|q_2|<(1+3\alpha) |q_1|.$

Then we decompose $q_1=q'_1t_2^1q''_1$,  $q_2=q''_2t_2^2q'_2$ and
$\partial(\Pi)=t_1^1u_1t_1^2u_2$.

  Observe that since $E_l$ is also a bond of a narrow contiguity
subdiagram between the paths $q'_1t_2^1$ and $t_2^2q'_2,$ we have
the similar inequality as for $|q_1|$ and $|q_2|$ above:
$|q'_1t_2^1|<(1+3\alpha)|t_2^2q'_2|.$ This implies that
$|q'_1|+|q''_2|<(1+3\alpha)(|q'_2|+|t_2^2|+|q''_2|)=(1+3\alpha)|q_2|.$
Similarly, $|q'_1|+|q''_2|<(1+3\alpha)|q_1|.$

Assume first that $o_1\in q'_1$ and $o_2\in q''_2$. We may also
suppose by symmetry that $|q'_1|<\frac{1+3\alpha}{2}\min (|q_1|,
|q_2|).$ Notice that equal shifts of the vertices $o_1, o_2$ along
$q_1$ and $q_2$ preserve the assumption of the lemma. Now assume
that, after shifts, $o_2$ becomes equal to $(q_2)_-$ and $o_1$
still belongs to $q'_1$. The last assumption leads to a
contradiction, because we have for the new $o_1$ and $o_2$ that
$|o_2-(q_2)_+| - |o_1-(q_1)_-| >
|q_2|-\frac{1+3\alpha}{2}|q_2|>|q_2|/3$ contrary to the assumption
of the lemma. Therefore one can move $o_1$ to $(t_2^1)_-$ and
$o_2$ to a vertex laying on $q_2$, in this case.

If $o_1\in q'_1$ and $o_2\in t_2^2$ ($o_2\in q'_2$), then
obviously, one can make equal shifts of these vertices so that
either $o_1$ becomes equal to $(t_2^1)_-$ or $o_2$ becomes equal
to $(t_2^2)_-$ . Similarly, one of the vertices $o_1, o_2$ can be
made equal to one of the $(t_2^1)_{\pm}, (t_2^2)_{\pm}$ after
parallel shifts, if $o_1\in t_2^1$ and $o_2\in t_2^2.$  The
remaining cases are symmetric to the examined ones.

Thus, to prove the lemma, one may assume that either  $o_1=
(t_2^1)_-$ or $o_1=(t_2^1)_+,$ or $o_2=(t_2^2)_-,$  or
$o_2=(t_2^2)_+.$ We will assume that $o_1=(t_2^1)_+=(s_1^1)_-.$

In addition, one may assume that the vertex $o_2$ belongs to the
subpath $t_2^2q'_2,$ because otherwise it belongs to $q''_2$ and
$|o_2-(t_2^2)_-|\ge \frac13 |\partial(\Pi_0)|$ by Lemma
\ref{1/3|Pi|}. In this case, by Lemma \ref{anal15.4}, one might
shift $o_1$ to $(t_2^1)_-$ so that corresponding vertex $o_2$
shifts to a vertex of the path $q''_2t_2^2$ because
$$\begin{array}{l}|o_2-(t_2^2)_+|>\frac13 |\partial(\Pi_0)|+|t_2^2|>
(1/3+(1-2\beta)(1-\alpha-\bar\alpha)) |\partial(\Pi_0)|>\\ \frac23
|\partial(\Pi_0)|>(1+2\beta)\bar\alpha|\partial(\Pi_0)|>|t_2^1|\end{array}$$
Then one gets a situation which is symmetric to the one considered
above with respect to changing of notations $p_1 \to p_2^{-1}$,
$p_2\to p_1^{-1}$, $q_1\to q_1^{-1}$, $q_2\to q_2^{-1}$ (i.e. we
replace $\Delta$ by its mirror copy).

Thus, up to symmetries, one has to consider only one case when
$o_1=(t_2^1)_+$ and $o_2\in t_2^2q'_2.$ The vertices $o_2$ and
$o_1$ can be connected by a path $vu,$ where $v= o_2-(t_2^2)_-$
and $u=(s_2^2)^{-1}u_2(s_1^1)^{-1}.$ Set $V\equiv  \phi(v)$ and
$U\equiv  \phi(u).$ As has been mentioned, $|U|\le
(\alpha+2\zeta)|\partial(\Pi_0)|.$ On the other hand, $|V|\ge
\frac13 |\partial(\Pi_0)|$ by Lemma \ref{1/3|Pi|}. Therefore one
gets the decomposition $V\equiv  V_1V_2$ guaranteed by Lemma
\ref{$V_2UV_1$}. Denote by $o$ the vertex of $v$ which gives the
corresponding decomposition $v=v_1v_2.$ There are two subcases.

(1)
  $|V|<3\alpha^{-1}|\partial(\Pi_0)|.$ Recall, that by Lemma \ref{$V_2UV_1$}
    the word
$V_2UV_1$ is equal in rank $i$ to a word $S$ of length at most
$3\gamma |V|$ and

\begin{equation}
||V_1|-|V_2||<\beta|V|+|U|<(3\alpha^{-1}\beta+\alpha+2\zeta)|\partial(\Pi_0)|
<2\alpha|\partial(\Pi_0)| \label{(?.3)}.
\end{equation}
 By
Lemma \ref{1/3|Pi|}, $|V|\ge |\partial(\Pi_0)|/3$, whence $|V_1|,
|V_2|>(1/6-\alpha) |\partial(\Pi_0)|
>1/7|\partial(\Pi_0)|$ by (\ref{(?.3)}). Since
$$0.45|\partial(\Pi_0)|<(1+2\beta)^{-1}(1-\bar\alpha-\alpha)|\partial(\Pi_0)|<|t_2^j|<
(1-2\beta)^{-1}\bar\alpha|\partial(\Pi_0)|$$ by lemmas
\ref{anal15.8} and \ref{anal15.4}, we have in case $o\in t_2^2$
that $$|o-(s_1^2)_-|<|t_2^2|-|V_2|
<((1+2\beta)\bar\alpha-1/6+\alpha)|\partial(\Pi_0)|<0.35|\partial(\Pi_0)|$$

If $o\in q'_2$, then
$$|o-(s_1^2)_-|=|V_2|-|t_2^2|<|V_1|+2\alpha|\partial(\Pi_0)|-0.45|\partial(\Pi_0)|.$$
 Hence in any case the path between $o$ and
$(s_1^2)_-$ has length at most $\max (0.35|\partial(\Pi_0)|,
|V_1|-2\alpha|\partial(\Pi_0)|),$ because of (\ref{(?.3)}) and the
inequality $4\alpha<0.45.$ Since
$|s_1^2u_1^{-1}s_2^1|<(\alpha+2\zeta)|\partial(\Pi_0)|,$ the
vertex $o$ can be connected with $(t_2^1)_-$ by a path $p$ of
length at most $\max (0.4|\partial(\Pi_0)|, |V_1|).$

The path $pt_2^1$ is homotopic to the path $v_2u.$ Thus the word
$S$ is equal in rank $i$ to the word $P\phi(t_2^1)V_1$, where
$P\equiv  \phi(p).$ The word $W\equiv \phi(t_2^1)V_1$ is
A-periodic because the vertices $o_2$ and $o_1$ are in the same
phase. It has length at least
$|V_1|+|t_2^1|>|V_1|+0.45|\partial(\Pi_0)|.$  Then $W=P^{-1}S$ in
rank $i,$ which implies, by lemmas \ref{analteom17.1} and
\ref{anal19.4}, \ref{anal19.5} applied to the diagram for this
equality, that
\begin{equation}
\max (0.4|\partial(\Pi_0)|, |V_1|)+3\gamma|V|>|P|+|S|>
\bar\beta(|V_1|+|t_2^1|)>\bar\beta(|V_1|+0.45|\partial(\Pi_0)|)\label{dopoln}
\end{equation}
On the other hand, since $|V_1|\le |V|\le
3\alpha^{-1}|\partial(\Pi_0)|,$  we have
$(\beta+3\gamma)|V|<(0.45\bar\beta-0.4)|\partial(\Pi_0)|$ which
contradicts inequality (\ref{dopoln}).

(2) Let $|V|\ge 3\alpha^{-1}|\partial(\Pi_0)|.$ Recall that
$|V_1|-|V_2|<\beta|V|+ (\alpha+2\zeta)|\partial(\Pi_0)|$, and
since $|V|=|V_1|+ |V_2|\ge 3\alpha^{-1}|\partial(\Pi_0)|$, the
length of $v_2$ is greater than $\frac{7}{15}|V|\ge 1.4\alpha
^{-1}|\partial(\Pi_0)|$. By Lemma \ref{Nearly}, the vertex $o$ can
be connected with a vertex $o'\in q'_1t_2^1$ by a path $t$ of
length at most $\alpha^{-1}|\partial(\Pi_0)|<\frac 13 |V|.$ By
Lemma \ref{analteom17.1}, the length of the path $q=o'-o_1$ is at
least
$$\bar\beta|v_2|-|t|-|u|>\frac {7}{15}|V|-\frac 13 |V|-
\alpha^2|V|>|V|/8$$ because
$|u|\le(\alpha+2\zeta)|\partial(\Pi_0)|<\alpha^2|V|$ by the
assumption of case (2).
 Therefore the length of the $A$-periodic word
$W\equiv  QV_1$ (where $Q\equiv \phi(q)$) is greater than
$(\frac{7}{15}+\frac{1}{8})|V|\ge\frac {7}{12}|V|.$ By Lemma
\ref{analteom17.1}, $W$ cannot be equal in rank $i$ to a word of
length $\le |V|/2$.  But it is equal in rank $i$ to the label of
$t$ times $S$ which is equal in rank $i$ to $V_2UV_1$. Therefore
$W$ is equal in rank $i$ to a word of length $<(\frac 13 +3\gamma)
|V|,$ a contradiction. \endproof

\begin{lm}\label{Hvost} Assume that $A$ is a word which is not
a conjugate in rank $i$ of a 0-word. Let $A=XB^lRX^{-1}$ in rank
$i$, where $B$ is a simple in rank $i$ word or a period of rank
$j\le i$ and $R\in \oo(B)$. Then $\oo(A)\le X\oo(B)X^{-1}$ in rank
$i$ (i.e. the canonical image of $\oo(A)$ in ${\bf G}(i)$ is a
subgroup of the image of $X\oo(B)X^{-1}).$
\end{lm}

\proof One may suppose that $A$ is a $Y$-reduced word, moreover
that $A$ is cyclically $Y$-reduced, by Lemma \ref{cyclically}(1)
and Lemma \ref{welldef}.

Consider a diagram $\Delta_1$ of rank $i$ for the equality
$A=XB^lRX^{-1}$. Let $\Delta_2$ be an annular diagram which arises
after identifying the boundary sections of $\Delta_1$ labeled by
$X.$ It has contours $q_1$ and $q_2$ labeled by $A$ and $B^lR,$
respectively. Initial vertices $o_1$ of $q_1$ and $o_2$ of $q_2$
are connected in $\Delta_2$ by a path $x$ labeled by the word $X.$
If $\Delta_2$ is not g-reduced, one can reduce it and obtain a
diagram $\Delta_3$ of a smaller type in which the vertices $o_1$
and $o_2$ are connected by a path $x'$ with $X'\equiv \phi(x')=X$
in rank $i.$ (Lemma 11.3 \cite{book}.) We will induct on the type
of $\Delta_3.$

If $r(\Delta_3)=0,$ then $\oo(A)= X'\oo(B^lR)(X')^{-1}$ in rank 0,
by Lemma \ref{welldef}. Since $X'=X$ in rank $i,$ the statement of
the lemma follows from Lemma \ref{Z3.5d}.

Assume now that $r(\Delta_3)>0.$ One may suppose that there are no
cells compatible with $q_2$ in $\Delta_3$ (if $B$ is a period of
rank $j\le i),$ since otherwise one can decrease the type of
$\Delta_3$ by Lemma \ref{anal13.3} and replace the label of $q_2$
for $B^{\ell \pm sn}R.$ Therefore the cyclic section $q_2$ is
smooth by Lemma \ref{anal19.5}, and so the degree of contiguity of
any cell of $\Delta_3$ to $q_2$ is less $\bar\alpha$ by Lemma
\ref{anal15.8}. By Lemma \ref{cor16.1,16.2} there is a positive
cell $\Pi$ and its contiguity subdiagram $\Gamma'$ to $q_1$ such
that
$(\Pi,\Gamma',q_1)>\frac{1}{2}(\bar\gamma-\bar\alpha)>\varepsilon.$

It follows, from Lemma \ref{Anal16.2} that there is a positive
cell $\pi$ of $\Delta_3$ and a contiguity subdiagram $\Gamma$ of
rank 0 of $\pi$ to $q_1,$ such that $(\pi,\Gamma,q_1)\ge
\varepsilon.$

Let $\partial(\pi,\Gamma,q_1) = s_1t_1s_2t_2.$ Then
$|s_1|=|s_2|=0$ since $r(\Gamma)=0$. There are no $Y$-bands
connecting different edges of subpath $t_1$ of $\partial(\pi)$ by
axiom A2. Therefore $|t_1|\le |t_2|,$ and the word $A'$ obtained
by the exchange in $A$ of the subword $\phi (t_2)$ by
$\phi(s_2t_1s_1)^{-1},$ is $Y$-reduced (and equal to $A$ in rank
0.)

The subword $\phi(t_1)^{-1}$ of $A'$ is a subword of the boundary
label of $\pi$ and $|\phi(t_1)^{-1}|\ge \varepsilon
|\partial(\pi)|.$ By Lemmas \ref{Z}, \ref{Z3.6}, if one
substitutes $\phi(t_1)^{-1}$ in $A'$ by $\phi(t)$, where $t_1t$
is the contour of $\pi,$ and replaces the resulting word by an
equal in rank $0$ $Y$-reduced word $A'',$ then $\oo(A'')\ge
\oo(A')=\oo(A).$

The word $A''$ is the label of the contour of a g-reduced diagram
$\Delta_4$ which is obtained from $\Delta_3$ by cutting the cell
$\pi$ out and (possibly) adding several 0-cells. Of course, such
surgeries (cutting off a \ct diagram or a cell) can touch the path
$x',$ but there is a homotopic path $x''$ with label $X''=X'=X$ in
rank $i$ which is not touched by the surgeries (i.e. this path
does not pass through the \ct subdiagram and does not have common
edges with the cell). By the inductive hypothesis for the diagram
$\Delta_4$ of a smaller type, $\oo(A'')\le X''\oo(B)(X'')^{-1}$ in
rank $i$. Hence $\oo(A)\le X\oo(B)X^{-1},$ as desired. \endproof

We need the statement of the following lemma to define periods of
rank $i_+.$

\begin{lm}\label{Equival} The following relation $\sim_i$ is an
equivalence on the set of all simple in rank $i$ words:  $A\sim_i
B$ by definition, if there are words $X, P, R$, where $P\in
\oo(A)$, $R\in\oo(B)$ (in rank 0), such that $AP=XB^{\pm
1}RX^{-1}$ in rank $i$.
\end{lm}

\proof Obviously, $\sim_i$ is a reflexive relation.

The equality $AP=XB^{\pm 1}RX^{-1}$ implies $B^{\pm 1}R=X^{-1}APX$
in rank $i$. Taking the inverses, if necessary, we have equality
$BR_{\pm}=X^{-1}A^{\pm 1}P_{\pm}X$ for a suitable $P_{\pm}\in
\oo(A)$, $R_{\pm}\in \oo(B)$, because $\oo(A)$ and $\oo(B)$ are
normalized by $A$ and $B$ respectively. By Lemma \ref{4.8.1}
$\oo(B)=\oo(BR_{\pm})$, and so $R_{\pm}\in X^{-1}\oo(A)X$ by Lemma
\ref{Hvost}. Then $Q\equiv  XR_{\pm}X^{-1}\in\oo(A)$ by Lemma
\ref{welldef}. Hence $B=X^{-1}A^{\pm 1}VX$ in rank $i$, where
$V\equiv  P_{\pm}Q\in \oo(A)$, and so the relation $\sim_i$ is
symmetric.

We have also proved that $B\sim_i A$ implies that $B=X^{-1}A^{\pm
1}VX$ in rank $i$ for some $V\in\oo(A)$ and a word $X.$

Assume now that $B\sim_i C$, that is $B=YC^{\pm 1}QY^{-1}$ in rank
$i$ for a word $Y$ and some $Q\in\oo(C).$ Then, in rank $i$, $AP=
X(YC^{\pm 1}Q_{\pm}Y^{-1})RX^{-1}$ for suitable $Q_{\pm}\in
\oo(C)$, and $V\equiv  Y^{-1}RY\in\oo(C)$ by lemmas \ref{Hvost}
and \ref{welldef}. Hence, in rank $i$, $AP=(XY)C^{\pm
1}(Q_{\pm}V)(XY)^{-1}$, where $Q_{\pm }V\in\oo(C).$ Therefore the
relation $\sim_i$ is transitive.
\endproof

\begin{lm}\label{anal18.1} Every word $X$ is a conjugate
in rank $i\ge 0$ either of a 0-word or of a word $A^lP,$ where
$|A|\le |X|,$ $A$ is either period of rank $j\le i$ or a simple in
rank $i$ word and $P$ represents in rank $0$ an element of the
subgroup $\oo(A).$
\end{lm}

\proof If the statement were false for $X$, then, in rank $i$,
$X=YA^sRY^{-1}$ for some ($Y$-reduced) essential word $A$ with
$|A|<|X|$ and word $R\in \oo(A)$. By the inductive hypothesis, we
may assume that $A=ZB^tPZ^{-1}$ in rank $i$ where $B$ is either a
period of rank $j\le i$ or a simple in rank $i$ word, and $P\in
\oo(B).$

We have that $A^sR=Z(B^tP)^sZ^{-1}R$ in rank $i$. The word
$Z^{-1}RZ$ is equal in rank $i$ to a word $Q\in \oo(B^tP)=\oo(B)$
by Lemma \ref{Hvost} and Lemma \ref{Z3.5d}. Therefore
$A^sR=Z((B^tP)^sQ)Z\iv$. Since $\oo(B)$ is normalized by $B$ (by
definition), $X$ is a conjugate in rank $i$ of a word of the form
$B^{st}S$ for some $S\in \oo(B),$ as desired.
\endproof

\begin{lm} \label{anal18.3} If a word $X$ has a finite order in rank $i,$
then it is a conjugate in rank $i$ either of a 0-word or of a word
$A^lR,$ where $A$ is a period of rank $j\le i$ and $R\in \oo(A).$
\end{lm}

\proof If the claim were false, then, by Lemma \ref{anal18.1},
since $\oo(A)$ is normalized by $A$, $A^lR=1$ in rank $i$ for
some simple in rank $i$ word $A,$ $R\in \oo(A)$ and $|l|>0.$ Then
$A^l=R^{-1}$ in rank $i$ and $\bar\beta|l||A|<|R|=0$ by lemmas
\ref{anal19.4}, \ref{anal19.5} and Lemma \ref{analteom17.1}. The
statement is proved by contradiction.
\endproof

\begin{lm}\label{anal18.4} If $A$ and $B$ are simple in rank $i$ words
and $A$ is equal in rank $i$ to $XB^lRX^{-1}$ for some word $X$,
and some word $R\in \oo(B),$ then $l=\pm 1$ and $|A|=|B|.$
\end{lm}

\proof When proving the first assertion by contradiction, one can
assume that $l\ge 2.$ By Lemma \ref{anal19.5} and Lemma
\ref{analteom17.1} for a g-reduced annular diagram for the
conjugacy of $A$ and $B^lR$, we have $2\bar\beta|B|<|A|,$ whence
$|B|<|A|$ contrary to the simplicity of the word $A$ in rank $i.$
Then we have $P\equiv  XRX^{-1}\in \oo(A)$ in rank $i$ by Lemma
\ref{welldef}, Lemma \ref{Z3.5d}, and $B= X^{-1}(AP^{-1})^{\pm
1}X.$ Hence $|A|=|AP^{-1}|\ge |B|$ since $|P|=0$ by definition of
$\oo(A)$ (see property (Z3)) and the simplicity of $B$ in rank
$i.$
\endproof

\begin{lm}\label{anal18.5} If words $X$ and $Y$ are conjugate in rank
$i,$ and $X$ is not a conjugate in rank $i$ of a $0$-word, then
there is a word $Z,$ such that $|Z|\le\bar\alpha (|X|+|Y|)$ and
$ZYZ^{-1}$ is equal to $X$ in rank $i.$
\end{lm}

\proof The proof is analogous to the proof of Lemma 18.5
\cite{book}. One may apply Lemma \ref{anallemma17.1} (which is an
analog of Lemma 17.1 \cite{book}) since $X$ is not a conjugate of
a 0-word.
\endproof

\begin{lm}\label{Phase}
 Let $\Delta$ be a narrow contiguity diagram between sections $q_1$
 and $q_2$ of its contour. Assume that it is of rank 0 and has
contour $p_1q_1p_2q_2,$ where $q_1$ and $q_2^{-1}$ are paths
labeled by $A$-periodic words with a simple in rank $i$ period
$A.$ Assume that $|q_1|, |q_2|\ge 2|A|,$ and $o_1$ and $o_2$ are
vertices of $q_1$ and $q_2$ which are in the same phase,
$|(q_1)_-- (o_1)|\le |q_1|/2$, $|(q_2)_+- (o_2)|\le |q_2|/2,$ and
$\phi(o_1-o_2)$ is a conjugate in rank $i$ of a 0-word. Then the
vertices $(q_1)_-$ and $(q_2)_+$ are in the same phase too.
\end{lm}

\proof Without loss of generality, one may assume that
$(q_2)_+=o_2$ since the vertices $o_1$ and $o_2$ can be
simultaneously shifted along the paths $q_1$ and $q_2$.

Every $Y$-band of $\Delta$ is a 0-bond between $q_1$ and $q_2$
because the narrow diagram $\Delta$ of rank 0 is defined by
0-bonds and by Lemma \ref{uchastok}, there are no bands connecting
different $Y$-edges of $q_1$ (of $q_2$) as the period $A$ is
simple in rank 0. Hence $|p_1|=0,$ and the vertex $o_1$ can be
connected with a vertex $o_0$ of $q_2$ by a path of length 0.
Proving by contradiction, one may suppose that $|V'|>0$ for the
word $V'\equiv \phi((q_1)_--o_1).$ Then $|V|=|V'|>0$ for the word
$V\equiv  \phi(o_2-o_0).$

Denote $U\equiv  \phi(o_0-o_1)$. Then $|U|=0,$ and the word $VU$
is a conjugate of a 0-word in rank $i$ by the hypothesis of the
lemma. Lemma \ref{$V_2UV_1$} provides us with a decomposition
$V\equiv V_1V_2,$ where $V_2UV_1$ is equal in rank 0 to a word $S$
of length at most $3\gamma|V|.$

Consider the vertex $o$ of the subpath $o_2-o_0$ that defines the
decomposition $V_1V_2$ of the label $V$ of $o_2-o$. Then the word
$V_2UV_1$ is readable in $\Delta$ starting with $o$: the suffix
$V_1$ can be read in the beginning of $o_1-(q_1)_+$ since $o_1$
and $o_2$ are in the same phase. The vertex $o$ can be connected
with a vertex $o'$ of $q_1$ by a path of length $0$ (see
\setcounter{pdeleven}{\value{ppp}} Figure \thepdeleven).

\unitlength=1mm \special{em:linewidth 0.4pt} \linethickness{0.4pt}
\begin{picture}(130.44,28.00)
\put(126.00,7.00){\line(-1,0){83.00}}
\put(43.00,7.00){\vector(0,1){15.00}}
\put(43.00,22.00){\line(1,0){83.00}}
\put(67.00,7.00){\circle*{1.00}} \put(67.00,22.00){\circle*{1.00}}
\put(93.00,7.00){\circle*{1.00}} \put(93.00,22.00){\circle*{1.00}}
\put(43.00,7.00){\circle*{1.00}}
\put(68.00,7.17){\line(1,0){24.00}}
\put(93.17,7.00){\line(0,1){15.00}}
\put(93.33,22.00){\line(0,-1){14.83}}
\put(93.33,21.83){\line(1,0){17.50}}
\bezier{64}(67.00,7.00)(64.83,13.00)(67.00,22.00)
\put(93.33,21.67){\line(1,0){20.00}}
\put(110.33,21.78){\line(1,0){2.89}}
\put(81.22,7.11){\vector(1,0){6.00}}
\put(93.44,10.67){\vector(0,1){6.44}}
\put(97.67,21.78){\vector(1,0){7.56}}
\put(41.11,13.44){\makebox(0,0)[cc]{$p_1$}}
\put(42.78,4.44){\makebox(0,0)[cc]{$o_2$}}
\put(54.44,3.78){\makebox(0,0)[cc]{$V_1$}}
\put(55.78,25.11){\makebox(0,0)[cc]{$V'$}}
\put(63.78,9.78){\makebox(0,0)[cc]{$o$}}
\put(69.44,19.44){\makebox(0,0)[cc]{$o'$}}
\put(93.11,25.78){\makebox(0,0)[cc]{$o_1$}}
\put(94.78,3.78){\makebox(0,0)[cc]{$o_0$}}
\put(95.44,14.44){\makebox(0,0)[cc]{$U$}}
\put(108.11,25.44){\makebox(0,0)[cc]{$V_1$}}
\put(78.78,3.78){\makebox(0,0)[cc]{$V_2$}}
\put(130.44,22.11){\makebox(0,0)[cc]{$q_1$}}
\put(130.11,7.11){\makebox(0,0)[cc]{$q_2$}}
\put(126.11,7.11){\vector(-1,0){6.67}}
\put(119.11,22.11){\vector(1,0){7.00}}
\end{picture}

\begin{center}
\nopagebreak[4] Fig. \theppp.

\end{center}
\addtocounter{ppp}{1}

Thus, up to a factor of length 0, the word $V_2UV_1$ is equal to a
subword of $q_1$ of the length $|V_2|+|V_1|=|V|,$ which is not
equal in rank $i$ to a word of length $\le\bar\beta|V|$, by Lemma
\ref{anal19.5} and Lemma \ref{analteom17.1}. The inequality
$\bar\beta|V|\le 3\gamma|V|$ contradicts the Lowest Parameter
Principle. This completes the proof of the lemma. \endproof

\medskip

The following lemmas \ref{anal18.6}--\ref{anal18.8} are analogous
to lemmas 18.6 - 18.8 \cite{book}. {\bf For a given rank $i\ge
1/2,$ they are proved by simultaneous induction on the sum $L$ of
lengths of periods.}

\begin{lm}\label{anal18.6} Let $\Delta$ be a g-reduced diagram of rank $i$
with contour $p_1q_1p_2q_2$, where $\phi(q_1)$ and
$\phi(q_2)^{-1}$ are periodic words with a simple in rank $i$
period $A.$ Then the paths $q_1$ and $q_2$ are $A$-compatible in
$\Delta,$ provided $|p_1|, |p_2| <\alpha |A|$ and $|q_1|, |q_2|>
(\frac56 h+1)|A|.$ (The inductive parameter $L$ is equal to
$|A|+|A|.$)
\end{lm}

\proof The same argument, as in the proof of Lemma 18.6
\cite{book} allows us to assume that the vertices $(p_1)_+$ and
$(p_1)_-$ are in the same phase, $|q_1|,|q_2|
>\frac56 h|A|$ and $|p_1|<(\alpha+1/2)|A|=\bar\alpha|A|.$ In view
of Lemma \ref{welldef}, to prove the lemma, one may replace period
$A$ by a cyclic permutation. Therefore we suppose that both $q_1$
and $q_2^{-1}$ start with the prefix $A.$

Suppose that the word $P_1\equiv \phi(p_1)$ is not a conjugate in
rank $i$ of a 0-word. In this case the annular diagram
$\bar\Delta$ obtained after identification of the the long
subpaths of $q_1$ and $q_2$ with equal labels, has no paths of
length 0 surrounding the hole. Therefore we can apply Lemma
\ref{anallemma17.1} (which is Lemma 17.1 from \cite{book}). The
rest of the proof of Lemma 18.6 from \cite{book} carries out with
one minor change: by Lemma \ref{anal18.1}, the word $D\equiv
\phi(p_1)$ is equal in rank $i$ to the product $YB^mRY^{-1},$
$R\in \oo(B)$ instead of simply $YB^mY\iv$. The word $D^l$ is
equal in rank $i$ to $YB^{ml}QY^{-1},$ for some 0-word $Q\in
\oo(B)$ since $\oo(B)$ is normalized by $B$. Then we have to
define $Z\equiv QY^{-1}\phi(s')$, instead of $Z\equiv
Y^{-1}\phi(s')$ in \cite{book}. This does not effect on the
further proof because $|Q^l|=0$ by the definition of $\oo(B)$.

Assume now that the word $P$ is a conjugate in rank $i$ of a
0-word. Since $q_1$ and $q_2$ are long enough ($h=\delta^{-1})$,
we can excise a narrow contiguity subdiagram $\Delta^0$ of
$\Delta$ with contour $p'_1q'_1p'_2q'_2$ according to Lemma
\ref{Anal17.3.2} (which can be applied because $q_1$ and $q_2$ are
smooth sections of rank $|A|$ by Lemma \ref{anal19.5} in the
smaller rank). Moreover by Lemma \ref{Anal17.3.2}, we may suppose
that $|q'_1|,|q'_2|> \frac12 h|A|$ and $|p_1'|,|p_2'|<\alpha |A|$.
Again, by Lemma \ref{Anal17.3.2}, there are vertices $o_1$ and
$o_2$ of $q'_1$ and $q'_2$ respectively, such that the paths
$(p_1)_+-o_1$ and $(p_1)_--o_2$ have the same label (say, $X$) and
$$|(q'_1)_--o_1|, |(q'_2)_+-o_2|<
\gamma^{-1}(\bar\alpha+\alpha)|A|< \frac h6 |A|.$$ Then the label
of any path $o_1-o_2$ is equal in rank $i$ to $X^{-1}P^{-1}X,$ and
therefore it is a conjugate in rank $i$ of a 0-word. By Lemma
\ref{$o_1-o_2$}, $r(\Delta^0)=0,$ and the vertices $(q'_1)_-$ and
$(q'_2)_+$ are in the same phase by Lemma \ref{Phase} (see
\setcounter{pdeleven}{\value{ppp}} Figure \thepdeleven).

\unitlength=1mm \special{em:linewidth 0.4pt} \linethickness{0.4pt}
\begin{picture}(130.78,42.44)
\put(10.44,2.78){\line(0,1){39.67}}
\put(10.44,42.44){\line(5,-3){29.00}}
\put(39.44,25.11){\line(1,0){55.67}}
\put(95.11,25.11){\line(2,1){33.00}}
\put(128.11,41.44){\line(0,-1){38.67}}
\put(128.11,2.78){\line(-5,3){31.67}}
\put(96.44,21.78){\line(-1,0){56.67}}
\put(39.78,21.78){\line(-3,-2){29.33}}
\put(8.11,21.78){\makebox(0,0)[cc]{$p_1$}}
\put(130.78,22.11){\makebox(0,0)[cc]{$p_2$}}
\put(77.78,28.11){\makebox(0,0)[cc]{$q_1'$}}
\put(78.11,18.44){\makebox(0,0)[cc]{$q_2'$}}
\put(39.44,21.44){\line(0,1){3.67}}
\put(95.44,25.11){\line(0,-1){3.33}}
\put(36.44,23.11){\makebox(0,0)[cc]{$p_1'$}}
\put(99.44,23.44){\makebox(0,0)[cc]{$p_2'$}}
\put(60.78,25.11){\circle*{1.00}}
\put(55.11,21.78){\circle*{1.00}}
\put(61.78,27.78){\makebox(0,0)[cc]{$o_1$}}
\put(55.11,18.44){\makebox(0,0)[cc]{$o_2$}}
\put(45.11,23.44){\makebox(0,0)[cc]{$\Delta^0$}}
\put(40.11,14.11){\makebox(0,0)[cc]{$\Delta$}}
\put(10.44,25.11){\vector(0,1){11.00}}
\put(39.44,21.78){\vector(0,1){3.33}}
\put(95.44,25.11){\vector(0,-1){3.33}}
\put(128.11,28.78){\vector(0,-1){13.00}}
\put(81.78,25.11){\vector(1,0){8.67}}
\put(72.44,21.78){\vector(-1,0){10.00}}
\put(92.94,25.11){\line(1,0){2.50}}
\end{picture}

\begin{center}
\nopagebreak[4] Fig. \theppp.

\end{center}
\addtocounter{ppp}{1}

Every $Y$-band starting on $q'_1$ must end on $q'_2$ by Lemma
\ref{uchastok}. Hence one may suppose that $\phi(q'_1)\equiv
\phi(q'_2)^{-1}\equiv  U$ for some word $U.$ Then
$U^{-1}\phi(p'_1)U=\phi(p'_2)^{-1}$ in ${\bf G}(0)$, and since
$|U|>\frac h2|A|>4|A|$, we have $A_1^{-4}\phi(p_1)A_1^4\in
\oo(A_1)$ by Lemma \ref{sdvig}, where $A_1$ is the cyclic shift of
$A$ which is the prefix of $\phi(q'_1).$
 Thus the paths $q'_1$ and $q'_2$ are compatible by
(Z3.1) since $A$ is simple in rank $1/2$ (recall that $i\ge 1/2$).
Therefore the paths $q_1, q_2$ are also compatible by Lemma
\ref{welldef}, and the lemma is proved.
\endproof

\begin{lm}\label{anal18.7} Let $Z_1A^{m_1}Z_2=A^{m_2}$ in rank $i$,
$m=\min(m_1,m_2),$ and $A$ be a simple word in rank $i.$ Then each
of the words $Z_1,Z_2$ is equal in rank $i$ to a product of a
power of $A$ and a word representing an element of the subgroup
$\oo(A)$ provided $|Z_1|+|Z_2|<(\gamma(m-\frac 56 h -1)-1)|A|.$
(The inductive parameter $L$ is equal to $|A|+|A|$.)
\end{lm}
\proof As in Lemma 18.7 \cite{book} we can consider a g-reduced
disc diagram $\Delta$ of rank $i$ with the contour $p_1q_1p_2q_2,$
where $\phi(q_1)\equiv  A^{m_1}$, $\phi(q_2^{-1})\equiv  A^{m_2}$,
$\phi(p_1)=Z_1$, $\phi(p_2)=Z_2$. By lemmas \ref{anal19.4},
\ref{anal19.5}, $q_1$ and $q_2$ are smooth of rank $|A|$. By
lemmas \ref{anal19.4} and \ref{anal17.5}, we can excise from
$\Delta$ a subdiagram $\Delta'$ with contour $p_1'q_1'p_2'q_2'$
where $$
\begin{array}{l}
|p_1'|,|p_2'| < \alpha|A|, |q_j'|>|q_j|-\gamma\iv
(|Z_1|+|Z_2|+|A|)>\\ m|A|-(m-\frac 56 h -1)|A| = (\frac 56
h+1)|A|, j=1,2.
\end{array}
$$ By Lemma \ref{anal18.6}, $q_1'$ and $q_2'$ are compatible in
$\Delta'$. So there exist phase vertices $o_1\in q_1$ and $o_2\in
q_2$ such that the label $T=\phi(t)$ of a path $t$ connecting
these vertices represents an element of $\oo(A)$ in rank $i$. Then
$\phi((q_1)_- - o_1)$ and $\phi((q_2)_+-o_2)$ are powers of $A$,
we have an equality of the following form in rank $i$:
$Z_1=A^{s_1}TA^{s_2}$. Since $\oo(A)$ is normalized by $A$ in rank
0, the lemma is proved.
\endproof

\begin{lm}\label{anal18.8} Let $\Delta$ be a g-reduced disc diagram of
rank $i$ with the contour $p_1q_1p_2q_2,$ where $\phi(q_1)$ and
$\phi(q_2)$ are periodic words with simple in rank $i$ periods $A$
and $B,$ respectively, and $|A|\ge |B|.$ Assume that
$|p_1|,|p_2|<\alpha |B|,$ $|q_1|>\frac34 h|A|$ and $|q_2|>h|B|.$
Then the word $A$ is a conjugate in rank $i$ of a product $B^{\pm
1}R,$ where $R$ represents an element of the subgroup $\oo(B),$
and $|A|=|B|.$ Moreover, if the words $\phi(q_1)$ and
$\phi(q_2)^{-1}$ start with $A$ and $B^{-1},$ respectively, then
$A$ is equal in rank $i$ to $\phi(p_1)^{-1}B^{\pm 1}R\phi(p_1).$
(The inductive parameter $L$ is equal to $|A|+|B|.$)
\end{lm}

\proof The proof is similar to the proof of Lemma 18.8
\cite{book}, but, at the very end, by Lemma \ref{anal18.7}, we
obtain the equality $Z_1=B^sR$ (instead of $Z_1=B^s$) with $R\in
\oo(B)$. Then $s=\pm 1$ and $|A|=|B|$ by Lemma
\ref{anal18.4}.\endproof

\begin{lm} \label{anal18.9} Let $\Delta$ be a disc diagram of rank $i$
with the contour $p_1q_1p_2q_2,$ where both $\phi(q_1)$ and
$\phi(q_2)$ are periodic words with a simple in rank $i$ period
$A$ and $|p_1|, |p_2|<\alpha |A|.$ Then $|q_1|,|q_2|\le h |A|.$
\end{lm}
\proof We prove the statement by contradiction. The proof is
parallel to the proof of Lemma 18.9 \cite{book} up to the
application of Lemma \ref{anal18.8} (instead of Lemma 18.8). This
lemma gives the equality $AR=Z^{-1}(A^{-1})^lZ$ in rank $i$, where
$l=\pm 1, R\in \oo(A), |Z|<(1+\alpha)|A|.$ Then, in case $l=-1,$
$Z^{-1}A^kZ=A^kQ$ in rank $i,$ where $Q\in \oo(A)$, because of
(Z3) where $\phi(q_1)\equiv A^k\bar A$, $|\bar A|<|A|$. The rest
of the argument of the proof of Lemma 18.9 from \cite{book} in
this case remains valid because $\oo(A)$ is normalized by $A$ in
rank $0$ and all elements in $\oo(A)$ have length $0$ (so the
corresponding factors do not affect length estimates).

In case $l=1,$ the equality $AR= Z^{-1}A^{-1}Z$, in rank $i,$
implies $A^n=  Z^{-1}A^{-n}Z$ by Lemma \ref{star}. Hence $Z^2$
commutes with $A^n$ in rank $i$. Therefore $Z^2$ commutes with any
$A^{sn}$, $s=1,2,...$. Taking $s$ sufficiently large, we can apply
Lemma \ref{anal18.7} and conclude that $Z^2=A^dP$ in rank $i,$
with $P\in \oo(A).$ Then again, $Z^{2n}=A^{dn}$ in rank $i$ by
Lemma \ref{star}. Therefore the word $Z$ commutes with $A^{dn}$ in
rank $i.$ But the equality $A^n= Z^{-1}A^{-n}Z$ implies that
$Z^{-1}A^{dn}Z=A^{-dn}$ in rank $i.$ Hence $A^{2dn}=1$ in rank
$i.$ This would contradict the definition of a simple word and
Lemma \ref{anal18.3}, if $d\ne 0.$

 Consequently $d=0$ and
$Z^2=P$ in rank $i.$ Then $Z$ has  finite order in rank $i$ by
property (Z1.2). This order must divide $n$ by property (Z1.2),
Lemma \ref{star} and Lemma \ref{anal18.3} since $B^n=1$ in rank
$i$  for every period $B$ of rank $\le i$. Therefore $Z^2=P\in
\oo(A)$ implies that $Z\in \oo(A)$ in rank $i$ as $n$ is an odd
number.

By Lemma \ref{star1} $Z$ commutes with $A^n$ in rank $i$. It
follows that $A^{2n}=1$ in rank $i.$ This is the final
contradiction, as simple in rank $i$ word $A$ has infinite order
in rank $i$ by Lemma \ref{anal18.3}. \endproof

\begin{lm}\label{anal19.1} Let $\Delta$ be a g-reduced disc diagram of
rank $i$ with the contour $p_1q_1p_2q_2,$ where $\phi(q_1)$ and
$\phi(q_2)$ are periodic words with simple in rank $i$ periods $A$
and $B$. Assume that $|p_1|,|p_2|<\alpha|B|,$
$|q_1|>(1+\frac{\gamma}{2})|A|,$ $|q_2|> \frac12 \varepsilon
n|B|.$ Then the word $A$ is a conjugate in rank $i$ of a word
$B^{\pm 1}R,$ where $R\in \oo(B).$ Moreover it is equal in rank
$i$ to $\phi(p_1)^{-1}B^{\pm 1}R\phi(p_1)$ if the word $\phi(q_1)$
starts with $A$ and $\phi(q_2)^{-1}$ starts with $B^{-1}.$
\end{lm}

\proof To prove the lemma, one should repeat the argument from
Lemma 19.1 \cite{book} up to the reference to Lemma
\ref{anal18.7}, instead of Lemma 18.7 \cite{book} which provides
us with the additional factor $R$ as compared with \cite{book}.
\endproof

\begin{lm}\label{anal19.2} Let $\Delta$ be a g-reduced disc diagram of
rank $i$ with the contour $p_1q_1p_2q_2$, where $\phi(q_1)$ and
$\phi(q_2)$ are periodic words with simple in rank $i$ periods $A$
and $B;$ $|q_2|\ge \varepsilon n|B|,$ $|p_1|, |p_2|<\zeta nc$ for
$c=\min(|A|,|B|).$ Then either $|q_1|<(1+\gamma)|A|,$ or $A$ is a
conjugate in rank $i$ of a word $B^{\pm 1}R$ where $R\in \oo(B)$.
Furthermore, equality $A\equiv B^{\pm 1}$ implies $A\equiv B^{-1}$
and the paths $q_1$, $q_2$ are $A$-compatible in $\Delta$.
\end{lm}

\proof The statement can be reduced to the claim of Lemma
\ref{anal19.1} in the similar way as in Lemma 19.2 \cite{book}.
\endproof

Recall that $i_+$ means 1 if $i=1/2$, and $i+1$ if $i\ge 1$.

\begin{lm}\label{anal19.3} 1) If $|X|\le i_+$ then $X^n=1$ in rank $i_+.$

2)  period $A$ of rank $i_+$ is simple in any rank $j\le i.$

3) Let $A$ and $B$ be periods of ranks $k,l\le i_+$, and $A\sim_j
B$ for some $j<k,l.$ Then $A\equiv  B.$
\end{lm}

\proof 1) If $X$ is a conjugate of a word of length 0 in rank
$i_+$, the statement follows from (Z1.2). Otherwise by  Lemma
\ref{anal18.1} $X$ is a conjugate in rank $i$ of a word of the
form $A^lP$ with $|A|\le |X|, P\in \oo(A)$ and $A$ is either a
period of rank $\le i$ or a simple word in rank $i$. Then
$(A^lP)^n=A^{ln}$ by Lemma \ref{star}. If $|A|<i_+,$ then
$A^{ln}=1$ in rank $i$ by the definition of relators of rank $\le
i.$ Hence $X^n=1$ in rank $i$, and consequently, in rank $i_+.$ If
$|A|=i_+$ and so, $A$ is simple in rank $i$, then $AV=YB^{\pm 1}
RY^{-1}$ in rank $i$ for a period $B$ of rank $i_+$, some $V\in
\oo(A)$, $R\in \oo(B)$ and a word $Y$, by the definition of the
set ${\it X}_{i_+}.$ Hence equality $A^n=1$ in rank $i_+$ follows
from equality $B^n=1$ in rank $i_+$ and Lemma \ref{star}.

2) Follows from the inclusion of the set of relators of rank $j$
into the set of relators of rank $i$ for $j\le i.$

3) If $k$ were less than $l$, for example, then $B=XA^{\pm 1
}PX^{-1}$ in rank $l_-$ for some $X$ and $P\in \oo(A)$ as it was
shown at Lemma \ref{Equival}, but this contradicts the choice of
periods of rank $l$. If $k=l$, then $A\equiv  B$ by the definition
of the set ${\it X}_l$.
\endproof

{\bf Lemmas \ref{HubCont} -- \ref{anal19.5} will be proved by
simultaneous induction on the type of diagrams.}

\begin{lm}\label{HubCont} Let $\pi$ be a hub in a g-reduced diagram
$\Delta$ of rank $i_+.$ Then (1) the degree of contiguity of
arbitrary cell $\Pi$ to $\pi$ is less than $\varepsilon;$
(2)arbitrary contiguity subdiagram $\Gamma$ of $\pi$ to a cell
$\Pi$ or to a section $p$ of the boundary $\partial(\Delta),$
labeled by a $A$-periodic word where $A$ is a simple in rank $i_+$
word or a period of rank $\le i_+,$ is of rank 0, provided in the
last case that there are no cells compatible with $p,$ in
$\Delta.$
\end{lm}

\proof Consider a contiguity subdiagram $\Gamma$ between a hub
$\pi$ and another cell $\Pi.$ Let $E$ be one of the bonds defining
$\Gamma.$ If $r(E)>0,$ then the principal cell of $E$ has the
contiguity degree to $\pi$ at least $\varepsilon$, and this
contradicts the statement of the lemma for a subdiagram of a
smaller type (the subdiagram does not contain $\Pi$). Thus $E$ is
a 0-bond.

Furthermore the same argument shows that a cell $\pi'$ of positive
rank from $\Gamma$ possesses no contiguity subdiagram to $\pi$ (to
$\Pi$ if $\Pi$ is also a hub), with degree $\ge\varepsilon$. If
$\Pi$ is a cell of rank $\ge 1$, then by lemmas \ref{anal19.4},
\ref{anal19.5} for $\Gamma$ and Lemma \ref{anal13.2}, one may
apply Lemma \ref{anal15.8} and conclude that a contiguity degree
of $\pi'$ to $\Pi$ is less than $\bar\alpha$. Since side arcs of
$\Gamma$ are of length 0, this contradicts Lemma
\ref{cor16.1,16.2} because $\varepsilon+\bar\alpha<\bar\gamma.$
Thus $r(\Gamma)=0$ if $\pi\ne \Pi$ .

If $\Pi$ is a hub, then $(\Pi,\Gamma,\pi)<\varepsilon$ by property
(Z2.3) because otherwise a subdiagram containing hubs $\pi$ and
$\Pi$ could be replaced by a subdiagram of rank 0, but $\Delta$ is
a g-reduced diagram. If $\Pi$ has rank $\ge 1$, then every maximal
$Y$-band of $\Gamma,$ starting on $\Pi,$ must end on $\pi$ by
Lemma \ref{uchastok}, and so $\Gamma\bigwedge\Pi$ cannot contain a
subword of length $>|A|$ since every letters of
$\phi(\Gamma\bigwedge\pi)$ are distinct by $(Z2.2)$. Therefore
$(\Pi,\Gamma,\pi)\le \iota <\varepsilon $ by A1, and the lemma
claim is proved for a contiguity subdiagram between a hub and
another cell.

If $\pi=\Pi$, then again the bond $E$ has rank 0, because
otherwise the principal cell $\pi'$ of this bond differs from
$\pi$, and, by previous case, $\pi'$ possesses no contiguity
subdiagram of degree $\ge\varepsilon$ to $\pi$ contrary to the
definition of a bond. Then $r(\Gamma)=0$ as above, and the
assumption $(\pi,\Gamma,\pi)>\varepsilon$ would contradict
property (Z2.2).

 Finally, if $\Gamma$ is the contiguity subdiagram of the hub
 $\pi$ to the boundary section $p$ given in the lemma hypothesis,
 then again the rank of the bond $E$ is 0 as above, and
 then the same argument shows that $r(\Gamma)=0.$ \endproof

\begin{lm}\label{anal19.4} A g-reduced diagram $\Delta$ of rank
$i_+$ is an A-map.
\end{lm}

\proof 1) Hubs have rank 1/2 and perimeter at least $n$ by
property (Z2.1). The perimeter of a cell with a period of a
positive integer rank $j\le i_+$ is equal to $nj$. Thus $\Delta$
satisfies condition A1.

2) Let $q$ be a subpath of $\partial(\Pi)$, where $r(\Pi)=j\le
i_+,$ and let $|q|\le\max(j,2).$ Let $p$ be a path homotopic to
$q$ in $\Delta,$ and $\Gamma$ denote the subdiagram with contour
$p^{-1}q.$ If $\Pi$ does not occur in $\Gamma$, then the type of
$\Gamma$ is smaller than that of $\Delta$, and $\Gamma$ is an
$A$-map by the inductive hypothesis. If also $|q|\le 2$ and
$|p|<2$, then $|\partial(\Gamma)|\le 3$, and $r(\Gamma)=0$ by
Lemma \ref{cor17.1}. Then $\phi(p)=\phi(q)$ in ${\bf G}(0)$,
which is impossible by the cyclic irreducibility of periods and,
if $\Pi$ is a hub, by (Z2.2).

We suppose next that $|q|\le j$ and $j>1.$ If $|p|<|q|,$ then
$|\partial(\Gamma)|<2j,$ and by Lemma \ref{cor17.1}, the
perimeters of all cells of positive rank in $\Gamma$ are less than
$2j\bar\beta^{-1},$ and so $r(\Gamma)<2j\bar\beta^{-1}n^{-1}<j$ by
A1. In this case, one can find a subword $\phi(q)$ of the cyclic
permutation of the period $A$ of $\Pi$, equal in rank $j-1$ to the
shorter word $\phi(p)$. Thus, one has a  simple in rank $j-1$ word
$A$ that is a conjugate in rank $j-1$ of a word of length $<|A|$;
a contradiction.

We now assume that $\Pi$ occurs in $\Gamma$. Then $\Pi$ does not
occur in the subdiagram $\Gamma'$ with contour $pq'$ where $q'$ is
the complement of $q$ in the contour of $\Pi.$ Thus $\Gamma'$  is
an $A$-map by the inductive hypothesis. If $r(\Pi)\ge 1,$  then,
by Lemma \ref{anal13.2} 2), no cell of $\Gamma'$ is  compatible
with $q'$. Hence, applying Lemma \ref{anal19.5} to  $\Gamma',$
$q'$ is a smooth section in $\partial(\Gamma')$, and by  Lemma
\ref{analteom17.1}, $\bar\beta|q'|\le |p|.$ Therefore  $$|q|\le
2(n-2)^{-1}|q'|\le 2\bar\beta^{-1}(n-2)^{-1}|p|<|p|.$$ It remains
to suppose that $\Pi$ is a hub. Then, by Lemma   \ref{HubCont},
there is no cell in $\Gamma'$ with a contiguity subdiagram to $q'$
of degree $\ge\varepsilon.$ Hence $q'$ is a  smooth section in
$\partial(\Gamma')$. Then again  $|q|\le
2(n-2)^{-1}\bar\beta^{-1}|p|<|p|,$ which completes the
verification of condition A2.

3) Let $\Gamma$ be a contiguity subdiagram of a cell of rank $j$
to a cell of rank $k$ in $\Delta$, with $p_1q_1p_2q_2=\partial(\pi,\Gamma,\Pi)$
and $(\pi,\Gamma,\Pi)\ge\varepsilon.$ Then $k\ge 1$ by Lemma
\ref{HubCont}, and $\Gamma$ is an $A$-map by the inductive
hypothesis.

Assume first that $j=1/2.$ Then $r(\Gamma)=0$ by Lemma
\ref{HubCont}. By Lemma \ref{uchastok} every maximal $Y$-band of
$\Gamma$ starting on $q_2$ must end on $q_1,$ and therefore,
by property (Z2.2), all the edges of $q_2$ have distinct labels.
Hence $|q_2|\le|A|=k<(1+\gamma)k$ where $A$ is the period of $\Pi.$

Then suppose that $j\ge 1.$ In this case, $q_1$ and $q_2$ are
smooth sections of ranks $j$ and $k$, respectively, in
$\partial(\Gamma)$ by Lemma \ref{anal13.2} 2) and Lemma
\ref{anal19.5} applied to $\Gamma.$ Thus, applying Lemma
\ref{anal15.3} to $\Gamma,$  $|p_1|, |p_2|<\zeta n \min(j,k).$
Hence, by Lemma \ref{anal17.4}, $r=r(\Gamma)<\min (j,k),$ and, by
Lemma \ref{anal19.3}, the periods $B,A$ corresponding to $\pi$ and
$\Pi$ are simple in rank $r.$ Thus, applying Lemma \ref{anal19.2}
to $\Gamma$, either $|q_2|<(1+\gamma)|A|=(1+\gamma)k$ (that is,
condition A3 holds), or $A=XB^{\pm 1}RX^{-1}$ in rank $r$ for
$R\in \oo(B)$, and therefore $A\equiv  B$ by Lemma \ref{anal19.3}.
Hence, by Lemma \ref{anal19.2}, $q_1$ and $q_2$ are compatible,
contradicting Lemma \ref{anal13.2} 1). This completes the proof.
 \endproof

\begin{lm}\label{anal19.5}
Let $p$ be a section of a contour of a g-reduced diagram $\Delta$
of rank $i_+$ and one of the following two conditions holds.

(1) The path $p$ is not a cyclic section and $\phi(p)$ is
 an $A$-periodic word
 where $A$ is a period of a rank
$k\le i_+,$ and there are no cells in $\Delta$, compatible with
$p$, or a simple word in rank $i_+$.

(2)The path $p$ is a cyclic section and $\phi(p)\equiv  A^mP$
where $A$ is a simple in rank $i_+$ word and $P\in\oo(A)$.

Then $p$ is a smooth section of rank $|A|$ in $\partial(\Delta)$.
 \end{lm}

\proof {\bf Case 1.} Suppose first that assumption (1) from the
lemma holds. Let us check two conditions from the definition of a
smooth section.

 1)  Let $q$ be a subpath in $p$, $|q|\le\max (|A|,2)$, and $t$
a path homotopic to $q$ in $\Delta.$ The subdiagram $\Gamma$ with
contour $qt^{-1}$ is an A-map by Lemma \ref{anal19.4}. If
$|t|<|q|\le 2,$ then $|\partial(\Gamma)|\le 3$, which gives a
contradiction as in the proof of Lemma \ref{anal19.4}. We can thus
assume that $|q|\le|A|.$ If $A$ is a simple in rank $i_+,$ word
then it is not a conjugate in rank $i_+$ of  a shorter word, and
so we have, for the subword $\phi(q)$ of a cyclic permutation of
$A$, $|\phi(q)|\le|\phi(t)|,$ since $\phi(q)=\phi(t)$ in rank
$i_+.$ This enables us to assume that $r(A)=k\le i_+.$

If we now suppose that $|t|<|q|,$ then $|\partial(\Gamma)|<2k,$
and by Lemma \ref{cor17.1}, $r(\Gamma)=j<k.$ Then by Lemma
\ref{anal19.3}, $A$ is simple in rank $j$, and we arrive at the
previous case. Thus, the first condition in the definition of
smoothness of section of rank $|A|$ is verified.

2) Let $\Gamma$ be a contiguity subdiagram of a cell $\Pi$ of rank
$j$ to $q,$ where $\partial(\Pi,\Gamma,q)=p_1q_1p_2q_2$ and
$(\Pi,\Gamma,q)\ge\varepsilon.$ If $j=1/2,$ then $|q_2|\le
|A|<(1+\gamma)|A|$ as in the proof of Lemma \ref{anal19.4}. Hence
we suppose that $j\ge 1$ and $B$ is the period for $\Pi.$

Subdiagram $\Gamma$ is an A-map by Lemma \ref{anal19.4}, and by
the inductive hypothesis $q_2$ is a smooth section of rank $|A|$
in $\partial(\Gamma).$ The section $q_1$ is smooth of rank $j$ by
lemmas \ref{anal13.2} 2) and \ref{anal19.5} applied to $\Gamma.$
Applying Lemma \ref{anal15.3} to $\Gamma$, $|p_1|,|p_2|<\zeta n
\min (j,|A|).$ Hence, by Lemma \ref{anal17.4},
$r=r(\Gamma)<\min(j,|A|),$ and by lemma \ref{anal19.3}, the
periods $A$ and $B$ are simple in rank $r.$ So, applying Lemma
\ref{anal19.2} to $\Gamma,$ either $|q_2|<(1+\gamma)|A|$ (that is,
the second condition in the definition of smooth section holds),
or $A=XB^{\pm 1}RX^{-1}$ in rank $r$, where $R\in\oo(B).$ Then, by
Lemma \ref{anal19.3}, $A\equiv  B^{\pm 1}$, and by Lemma
\ref{anal19.2}, $q_1$ and $q_2$ are $A$-compatible in $\Gamma$,
that is, $\Pi$ and $q$ are $A$-compatible in $\Delta$, contrary to
the hypothesis of the lemma

{\bf Case 2.} Suppose now the assumption (2) of the lemma holds.
The label of the contiguity arc $q_2$ in item 2) (and $\phi(q)$ in
item 1)) can be of the form $U_1PU_2$ where $U_1U_2$ is a phase
decomposition of an $A$-periodic word. By the definition of
$\oo(A)$, by Lemma \ref{sdvig} and by Lemma \ref{welldef}
$PU_2=U_2P'$ in ${\bf G}(0)$, where $P'\in\oo (A')$,$|P'|=0$, $A'$
is a cyclic permutation of $A$ and $U\equiv  U_1U_2$ is an
$A'$-periodic word with phase decomposition $U\cdot 1.$ Hence one
can transform $\Gamma$ by adding a subdiagram of rank $0$ whose
contour labelled by $UP'(U_1PU_2)^{-1}.$ The new diagram $\Gamma'$
can be considered as in Case 1, because $|P'|=0$ and $P'$ can be
can be considered a part of the label of the side arc of
$\Gamma'.$

The lemma is proved.
\endproof

This together with Lemma \ref{anal19.3} 1) complete the proof of
Proposition \ref{mainprop}.

Proposition \ref{mainprop} immediately implies the following
statement that shows that one can construct analogs of HNN
extensions in the class of groups of odd exponent $n>>1$.

\begin{cy}\label{Bernsrem} Let $G=\la A\ |\ R\ra$ be a group satisfying
the identity $x^n=1$ for some odd $n>>1$. Let $P_i, Q_i$ be
subgroups of $G$, $i=1,2,...,k$, $k\ge 1$, and for each
$i=1,..,k$, let $\phi_i: P_i\too Q_i$ be an isomorphism. Let $H$
be the (multiple) HNN extension of $G$ corresponding to these
pairs of isomorphic subgroups, $$H=\la A, t_1,...,t_k \ | \ R,
t_iut_i\iv =\phi_i(u), \hbox{ for all } u\in P_i, i=1,...,k\ra.$$
Consider $A$ as a set of $0$-generators and $Y=\{t_1,...,t_k\}$ as
the set of non-zero generators of $H$. Suppose that $H$ satisfies
the conditions (Z3.1) and (Z3.2) (or (Z3.1) and (Z3.2')). Then the
natural homomorphism $G\too H/H^n$ is an embedding.
\end{cy}

\proof Indeed, notice that $H$ obviously satisfies conditions
(Z1). Conditions (Z2) also are satisfied because there are no
hubs. Conditions (Z3) are satisfied by our assumption. Hence by
Proposition \ref{mainprop} every $g$-reduced diagram over some
graded presentation $R\cup S_1 \cup S_2\cup...$ of $H/H^n$ is an
A-map, and every relator in $S_i$, $i=1,2,...$ is of the form
$A^n$ where $A$ is a cyclically reduced word of $Y$-length $i$.

Let $W$ be a word in $A$ which is equal to $1$ in $H/H^n$. Let
$\Delta$ be a $g$-reduced \vk diagram with boundary label $W$. If
the rank of $\Delta$ is at least 1 then by Lemma \ref{th16.2}
$\Delta$ contains a cell $\Pi$ of rank $\ge 1$ and a contiguity
subdiagram of rank $0$ of  $\Pi$ to $\partial(\Delta)$ with
non-zero contiguity degree. This implies that the label of
$\partial(\Delta)$ must contain a letter from $Y$, a contradiction
(we assumed that $W$ contains only letters from $A$). Hence
$\Delta$ has rank 0, so $\Delta$ is a diagram over $R$. Hence
$W=1$ in $G$.
\endproof

\begin{remark}{\rm Sergei Ivanov also proved Corollary \ref{Bernsrem}\cite{Ivanov}.
We added this corollary in our paper after the second author
attended Ivanov's talk at the Geometric Group Theory conference in
Montreal (May, 2002). A particular case of that corollary has been
obtained by K.V.Mikhajlovskii \cite{Mikh}. He considered the case
when $k=1$, each of the subgroups $P_1$ and $Q_1$ is malnormal
(i.e. $P_1\cap xP_1x\iv=\{1\}$ if $x\in G\backslash P_1$, $Q_1\cap
xQ_1x\iv=\{1\}$ if $x\in G\backslash Q_1$), and $xP_1x\iv \cap
Q_1=\{1\}$ for every $x\in G$. In that case $\oo(h)=\{1\}$ for
every non-zero cyclically minimal element $h\in H$. Hence
condition (Z3) trivially holds.}\end{remark}

\section{Subgroups of the free
Burnside groups.} \label{appendix}

In this section we prove certain properties of free Burnside
groups. Although Theorem \ref{SQ} below is of independent
interest, this section plays an auxiliary role in the paper and
can be skipped in the first reading.

\subsection{Sets of words with small cancellation}

 Let ${\cal A}=\{a_1,\dots,a_m\}$ be an alphabet. For a constant $c>1$,  we say that a
word $v$ is \label{caper} $c$-{\it aperiodic} if for every word
$A$, there are no non-empty $A$-periodic subwords $w$ of $v$ with
$|w|\ge c|A|.$

\begin{lm} \label{lm01}  For  every  $\varepsilon>0$,  one  can  find
$m=m(\varepsilon)$ such that the set of (positive)
$(1+\varepsilon)$-aperiodic words in the alphabet ${\cal A}_m =
\{a_1,\dots, a_m\}$ is infinite. Moreover, the number $f(i)$ of
such words of length $i\ge1$ is greater than
$m^{(1-\varepsilon/2)i}.$
\end{lm}

\proof First notice that one can assume that $\varepsilon$ is
sufficiently small and fix a number $m$ such that

\begin{equation}
m-m^{\frac{\varepsilon}{2}+\frac{1}{1+\varepsilon}}
-m^{-1+\frac{\varepsilon}{2}+\frac{2}{1+\varepsilon}}>m^{1-\frac
{\varepsilon}{2}}-m^{\frac{1}{1+\varepsilon}} \label{00}
\end{equation}

Obviously, $f(1)=m>m^{1-\varepsilon/2}.$ We will prove the
inequality $f(i+1)>m^{1-\varepsilon/2}f(i)$, for $i\ge 1,$ by
induction on $i.$

     By adding a letter to a $(1+\varepsilon)$-aperiodic word $w'$ of length
$i$ on the right, one obtains $mf(i)$ words of length $i+1.$ But
some of these words $w$ may contain non-empty $A$-periodic
subwords $v$ with $|v|\ge (1+\varepsilon)|A|$ for some $|A|.$ We
need an upper estimate for the number of such ``bad"
possibilities.

     Clearly, the subword $v$ must be a suffix of  $w$,
$w\equiv  uv,$ since otherwise $v$ would be a subword of the
prefix $w'$ of length $i$ in $w.$ Further we may assume that the
word $v$ starts with the period $A$, and the previous argument
shows that $A$ is uniquely determined by $v$, because otherwise
word $w'$ were not $(1+\varepsilon)$-aperiodic. Since $|v|\le
(1+\varepsilon)|A|,$ the total number of the words $v$ of length
$k$, which can occur, is at most $m^{\frac{k}{1+\varepsilon}}$.
Therefore the number of ``bad" products $w\equiv uv$ with $|v|=k$
and $|u|=n+1-k$, is at most $f(n+1-k)m^{\frac{k}{1+\varepsilon}}.$
By using the inductive hypothesis, one can substitute the first
factor by $f(n)m^{(1-k)(1-\varepsilon/2)}.$ Then, by summing over
$k\ge 2,$ we obtain

$$f(n+1)>mf(n)-f(n)\sum_{k=2}^{\infty}
m^{(1-\varepsilon/2)(1-k)+\frac{k}{1+\varepsilon}}=$$
$$f(n)\left(m-\frac{m^{-1+\frac{\varepsilon}{2}+\frac{2}{1+\varepsilon}}}
{1-m^{-1+\frac{\varepsilon}{2} +\frac{1}{1+\varepsilon}}}\right)>
m^{1-\varepsilon/2}f(n)$$ by (\ref{00}), as required. \endproof

\begin{lm} \label{02} For any $\varepsilon>0,$ there exists
$m=m(\varepsilon)$ and an infinite set $S$ of (positive) words
$A_1, A_2,\dots$ in the alphabet
 ${\cal A}_m
=\{a_1,\dots,a_m\}$ such that:

(1) every reduced product, whose factors belong to $S,$ has no
non-empty $A$-periodic subwords of length $\ge(1+\varepsilon)|A|$
unless the word $A$ is freely conjugate to a product of some words
of $S$;

(2) suppose that $A'\equiv  UV'$ and $A''\equiv  UV''$, or
$A'\equiv  V'U$ and $A''\equiv  V''U,$ are distinct cyclic
permutations of the words of $S$ (that is, $U$ is a common prefix
or a common suffix). Then
$|U|<\frac{\varepsilon}{10}\min(|A'|,|A''|).$

(3) $|A_i|\ge n$, $i=1,2,...$.
\end{lm}

\proof We can assume that $\varepsilon<1/3$ and choose a number
$m_1=m_1(\varepsilon_1)$ for $\varepsilon_1=\varepsilon/10$ and an
infinite set of positive words $B_1, B_2,\dots$ in alphabet ${\cal
A}_{m_1}$, according to Lemma \ref{lm01}. Then define $m=sm_1$ for
a number $s\ge 2\varepsilon_1^{-1}.$

  Set $A_i=B_{i,1}\dots B_{i,s}$ for $i=1,2,\dots$, where $B_{i,j}$ are the copies
of the word $B_i$ in the alphabet ${\cal
A}_j=\{a_{1j},\dots,a_{m_1j}\}$. Throwing out, if necessary, a
finite subset of $\{A_1, A_2,...\}$ we can assume that $|A_i|\ge
n$. Hence property (3) holds. Property (2) follows from the
inequality $s\ge 2\varepsilon_1^{-1}$ because the alphabets ${\cal
A}_j$ are disjoint.

  In particular, in any product $A_iA_j^{-1},$ we can cancel at most
$\varepsilon_1$ (in length) of $A_i$ or $A_j$ if $i\ne j,$ and no
letter of $B_{i,s-1}$ or $B_{j,s-1}^{-1}$ can be cancelled in the
product. Similar remark concerns the product $A_i^{-1}A_j.$

  To prove (1), suppose that $v$ is a non-empty $A$-periodic subword of a reduced form
$w\equiv  \bar A_{i_1}^{\pm 1}\dots \bar A_{i_l}^{\pm 1}$ of a
product $A_{i_1}^{\pm 1}\dots A_{i_l}^{\pm 1},$ and $|v|\ge
(1+\varepsilon)|A|.$ We may assume that $v\equiv  A'A''A',$ where
$A\equiv  A'A''$ and $|A'|\ge\varepsilon |A|,$ because otherwise
we can just shorten the subword $v.$

Assume first that both occurrences $A'$ are contained in the same
factor $\bar A_{i_r}^{\pm 1}$ of $w.$ Then $v^{\pm 1}$ is a
subword of the product $B_{i_r1}\dots B_{i_rs}$, and the both
occurrences of $A'$ are present in some $B_{i_rj}.$ Hence $v$ is a
subword of $B_{i_rj},$ contrary to the choice of the words $B_1,
B_2\dots.$

  Now assume, the first $A'$ of the decomposition of $v$ starts in some
$\bar A_{i_l}^{\pm 1}$, and the last $A'$ ends in some $\bar
A_{i_r}^{\pm 1}$
 with $r>l.$ Then the definition of the words $A_i$ implies that $v$
 must contain at least $\frac{s-3}{2}>\frac 25 s$ factors $B_{i_lj}$ of some $A_{i_k}$
 where $l\le k\le r.$  Then

$$|A'|\ge
\varepsilon(1+\varepsilon)^{-1}|v|>\varepsilon(1+\varepsilon)^{-1}
\frac 25 |A_{i_k}|>3\varepsilon_1|A_{i_k}| $$

By the choice of $s$, $A'$ contains a subword
$(B_{i_k,t}B_{i_k,t+1}B_{i_k,t+2}B_{i_k,t+3})^{\pm 1}$ for some
$t.$ We will suppose that the exponent is equal to $+1$. Hence a
cyclic permutation of $A$ is freely equal to a word
$(B_{i_k,t+2}B_{i_k,t+3})\dots(B_{i_k,t}B_{i_k,t+1}).$ This
subword of $w$ must be freely equal to a cyclic permutation of a
product of a number of factors $A_i^{\pm 1}$, because there is
only one occurrence of the form $B_{j,q}B_{j,q+1}$ in all the
words $A_1, A_2,\dots$. The lemma is proved.
\endproof

\begin{remark} \label{rem1}
{\rm In the proof of Theorem \ref{main}, we do not need {\bf
infinite} collections of words satisfying the conditions of Lemma
\ref{02}. We only need a collection of words satisfying that
conditions and consisting of $Cm$ words where $C$ is a constant
depending on $n$ (but not on $m$), $C>n>\frac{20}{\varepsilon}$,
and $m$ is a sufficiently large integer.

An explicit and easy construction of such collections of words can
be found in \cite{Sa87}. Let us repeat this construction here. Let
us fix integers $r>C^2$ and $m=rn$. Consider an $r^2\times
n$-matrix $M$ where every odd-numbered column is equal to
$$(1,2,...,r,1,2,...,r,...,1,2,...,r)$$ and every even-numbered
column is equal to $$(1,1,...,1,2,2,...,2,...,r,r,...,r).$$ Now
let us replace every number $i$ in column $j$ of $M$ by letter
$a_{i,j}$. Rows of the resulting matrix can be viewed as words of
length $n$  over the alphabet with $m=rn$ letters $\{a_{i,j},
i=1,...,r, j=1,...,n\}$. For example the first word is equal to
$a_{1,1}a_{1,2}...a_{1,n}$, etc. Let $S=\{w_1,...,w_{r^2}\}$ be
the collection of these words. Let us call words from $S$ and
their inverses {\em blocks}.

It is mentioned in \cite{Sa87} (and obvious) that no two different
words in $S$ have a common subword of length $2$. Condition (3) of
Lemma \ref{02} holds by definition of the blocks. Since
$n>\frac{20}{\varepsilon}$, condition (2) of Lemma \ref{02} holds.

To prove (1), consider a reduced product
$W=w_{i_1}^{s_1}...w_{i_k}^{s_k}$ of words from $S$, $s_i=\pm 1$,
where $w_{i_j}^{s_j}w_{i_{j+1}}^{s_{j+1}}\ne 1$, $j=1,...,k-1$.
Suppose that $W$ contains an $A$-periodic subword $V$ where
$|V|>(1+\varepsilon)|A|$. Since $n>4/\varepsilon$, we can assume
that $V$ contains a subword $AB$ where $B$ is the beginning of $A$
of length $4$. Then either the first two letters of $B$ are in one
block of $W$ or the last two letters of $B$ are in one block of
$W$. Since every block is completely determined by any of its
2-letter subwords, we can shift $A$ to the left or to the right in
$W$, replacing $A$ by a cyclic shift $A'$ of $A$, such that $A'$
starts with the beginning of a block and ends with the end of a
block. Then $A'$ is a product of blocks as required. }\end{remark}

\subsection{Subgroups of $B(m,n)$ satisfying the congruence\\
extension property} \label{cep}

Recall that a subgroup $H\le G$ satisfies the {\em congruence
extension property} if for every normal in $H$ subgroup $N$ there
is a normal subgroup $L$ of $G$, such that $N=L\cap H.$

In this section, we prove

\begin{theo} \label{SQ} For any sufficiently large odd $n,$ there
exists $m=m(n)$ such that the free Burnside group $B(m,n)$ has a
subgroup $H$ with three properties:

     (1) $H$ is isomorphic to the free Burnside group $B(\infty,n)$
of the infinite countable rank;

     (2) $H$ satisfies the congruence extension property in $B(m,n)$.

     (3) the words $A_1, A_2,\dots$ representing the free generators of
$H$ can be effectively chosen in accordance with Lemma \ref{02},
and for every $l\ge 1$, the free Burnside subgroup generated by
$A_1,\dots,A_l$, also possesses the congruence extension property
in $B(m,n).$
\end{theo}

The statements (1) and (2) constitute Theorem \ref{SQintr}. The
proof of Theorem \ref{SQ} follows very closely Section 6 of
\cite{book} although we repeat a few notation and definitions from
\cite{book}.

In particular, we use  the positive parameters $\alpha,
\beta,\gamma,\delta,\varepsilon,\zeta,\iota, n$ from Section
\ref{amaps}. Let $m=m(\varepsilon)$ be taken according to Lemma
\ref{02}. Then we fix the infinite set of words $A_1, A_2,\dots$
which exists by Lemma \ref{02}. By Lemma \ref{02} (3) the words
$A_k$ are sufficiently long, say, $|A_k|\ge n.$ Let $H$ be the
subgroup generated by all words $A_1,A_2,\dots$ in the free
Burnside group $B(m,n)$ with the free basis ${\cal
A}=\{a_1,..,a_m\}.$

Suppose that $v$ is a (cyclically) reduced word of length $t$ in
certain variables $x_1,\dots,x_s$. The (cyclic) word
$v(A_1,\dots,A_s)\equiv  A_{i_1}^{\eta_1}\dots A_{i_t}^{\eta_t}$
obtained by replacing $A_i$ for $x_i$ will be called a (cyclic)
\label{hword}$H$-{\it word} and will be considered together with
its decomposition into factors $A_{i_1}^{\eta_1},\dots,
A_{i_t}^{\eta_t}.$

 \begin{lm} \label{03} The words $A_1,A_2\dots,$ satisfying conditions (1)
 and (2) of Lemma \ref{02},
freely generate a free Burnside subgroup (of infinite countable
rank) of exponent $n$ in the group $B(m,n).$
\end{lm}

\proof To prove the lemma, we consider the graded presentation of
$B(m,n)$ from Section 18 in \cite{book}. Assume there is a
non-trivial relation $v(A_1,\dots,A_l)=1$ . Then there is a
g-reduced diagram $\Delta$ of some rank $i>0$ whose boundary label
$\Lab(q)$ is freely equal to an $H$-word. Proving by
contradiction, we may suppose that the number of cells of $\Delta$
is minimal.

By Lemma \ref{th16.2}, there is a cell $\Pi$ of $\Delta$ and a
contiguity subdiagram $\Gamma$ with
$(\Pi,\Gamma,q)\ge\varepsilon.$ Since $n\varepsilon>2,$ we obtain,
by Lemma \ref{02}(1), (3) that $\Pi$ corresponds to a period $A$
where $A$ is freely conjugate to an $H$-word $\bar A$. By Lemma
\ref{02}(2) the boundary label $\Lab(q')$ of the subdiagram
$\Delta',$ obtained from $\Delta$ by cutting off the cell $\Pi$,
is freely equal to an $H$-word also, and $\Lab(q)$ is freely equal
to the product of $\Lab(q')$ and a word $X\bar A^nX^{-1}$ for some
$H$-word $X.$ By the minimality of $\Delta,$ diagram $\Delta'$
corresponds to a relation $v'(A_1,\dots,A_l)=1$ such that
$v'(x_1,...,x_l)$ follows from the Burnside identity $x^n=1$. Then
the relation $v(A_1,\dots,A_l)=1$ follows from Burnside relations
too.
\endproof

Let $\Delta$ be a map. A smooth section $q$ of rank 1 in the
contour $\partial(\Delta)$ is said to be
\label{epss}$\varepsilon$-{\it section} if for any cell $\pi$ of
$\Delta$ and any contiguity subdiagram $\Gamma$, $(\pi,\Gamma,
q)<\varepsilon.$

Let us now fix an arbitrary normal subgroup $N$ in $H$. We would
like to define a group presentation of the factor group of
$B(m,n)$ by the normal closure of $N$ in $B(m,n)$.

We denote by $T$ the set of all $H$-words (in the alphabet
$a_1,\dots, a_m$) representing elements of $N$, that are
cyclically reduced as elements of the free subgroup $H$ (of the
absolutely free group) with generators $A_1,A_2,\dots.$ In
particular, $T$ contains all powers of the form
$v(A_1,\dots,A_s)^n.$ The following alteration of the definitions
from Section 18 of \cite{book} depends on $N$.

We define ${\bf G}(0)=F({\cal A})$ to be the absolutely free group
with the empty set ${\cal R}_0$ of defining words. For every $i\ge
0$ a non-trivial word $A$ is called \label{simple2}{\it simple} in
rank $i$ if it is not conjugate in rank $i$ (i.e. in group ${\bf
G}(i)$) either to a power of a period $B$ of rank $k$, where $1\le
k\le i,$ or to a power of a word $C$, where $|C|<|A|$, or to
 an $H$-word.

Every word of $T$ is included in the system ${\cal R}_1$ of
relators of rank 1. (Let us call them \label{trel}$T$-{\it
relators}.) We also include words $a_1^n,\dots, a_m^n$ in the
system ${\cal R}_1.$ Every word $a_k$ of length 1 is, by
definition, a period of rank 1, and ${\bf G}(1)=\langle {\cal A}\
|\ {\cal R}_1\rangle .$

For ranks $i\ge 2$, the definitions of \label{period2} periods of
rank $i$ and the group \label{ginf} ${\bf G}(i$)=$\langle{\cal
A}|{\cal R}_i\rangle, {\bf G}(\infty),$ given in Section 18 of
\cite{book} remain valid.

Instead of Lemma 18.1 \cite{book} we claim now that the following
statement holds.

\begin{lm} \label{18.1}
Every word $X$ over ${\cal A}$ is a conjugate in rank $i\ge 0$
either of a power of some period of rank $j\le i$ or of a power of
a simple in rank $i$ word, or of an $H$-word.
\end{lm}

The proof can be easily derived from the definition as in
\cite{book}.

Analogs of Lemma 18.2, Corollaries 18.1 and 18.2 are not used in
the proof, so we do not need them now.

The following statement replaces Lemma 18.3 \cite{book}.

\begin{lm} \label{18.3}
If a word $X$ has a finite order in rank $i$, then it is a
conjugate in rank $i$ of a power of a period of rank $k\le i$ or
of an $H$-word.
\end{lm}

There are no changes in the proof as compared with \cite{book}.

\medskip

Assume that $\Pi$ is a $T$-cell (i.e. it corresponds to a
$T$-relator) in a g-reduced annular diagram $\Delta$ of rank $i,$
and $\Gamma$ is a contiguity diagram of the cell $\Pi$ to itself,
such that the hole of $\Delta$ is ``surrounded" by $\Pi$ and
$\Gamma$, i.e. the union of $\Pi$ and $\Gamma$ does not belong to
any disc subdiagram of $\Delta$ (see
\setcounter{pdeleven}{\value{ppp}} Figure \thepdeleven).

\unitlength=1mm \special{em:linewidth 0.4pt} \linethickness{0.4pt}
\begin{picture}(102.00,53.00)
\put(72.00,29.00){\oval(60.00,48.00)[]}
\put(71.50,28.00){\oval(11.00,10.00)[]}
\put(63.50,28.00){\oval(15.00,30.00)[l]}
\put(84.50,28.00){\oval(11.00,30.00)[r]}
\put(63.00,13.00){\line(1,0){21.00}}
\put(63.00,43.00){\line(1,0){8.00}}
\put(84.00,43.00){\line(-1,0){10.00}}
\put(58.50,27.50){\oval(15.00,37.00)[l]}
\put(88.50,27.50){\oval(13.00,37.00)[r]}
\put(58.00,9.00){\line(1,0){31.00}}
\put(58.00,46.00){\line(1,0){13.00}}
\put(74.00,46.00){\line(1,0){14.00}}
\put(71.00,43.00){\line(0,1){3.00}}
\put(74.00,46.00){\line(0,-1){3.00}}
\put(71.00,43.00){\line(1,0){3.00}}
\put(74.00,46.00){\line(-1,0){3.00}}
\put(72.00,46.00){\line(0,-1){3.00}}
\put(73.00,43.00){\line(0,1){3.00}}
\put(71.00,43.00){\rule{3.00\unitlength}{3.00\unitlength}}
\put(69.00,44.00){\makebox(0,0)[cc]{$\bar p$}}
\put(76.00,44.00){\makebox(0,0)[cc]{$p$}}
\put(49.00,37.00){\makebox(0,0)[cc]{$q_1$}}
\put(58.00,34.00){\makebox(0,0)[cc]{$q_2$}}
\put(51.00,32.00){\vector(0,-1){9.00}}
\put(56.00,23.00){\vector(0,1){7.00}}
\put(55.00,14.00){\makebox(0,0)[cc]{$\Pi$}}
\put(48.00,47.00){\makebox(0,0)[cc]{$\Delta$}}
\end{picture}

\begin{center}
\nopagebreak[4] Fig. \theppp.

\end{center}
\addtocounter{ppp}{1}

The subdiagram $\Gamma$ is of rank $0$ by Lemma \ref{07} (1)
applied for a smaller rank. If the degree of $\Gamma$-contiguity
of $\Pi$ to $\Pi$ is at least $\varepsilon, $ we call $\Pi$ a
\label{hoop}{\it hoop}. The reduced boundary $q$ of a hoop can be
decomposed as $q_1pq_2{\bar p}^{-1}$, where $p, \bar p$ are the
\ct arcs of $\Gamma,$, $\Lab(\bar p)$ is freely equal to
$\Lab(p)$, and so $|p|, |\bar p|\ge \varepsilon |q|$. In view of
property (2) of Lemma \ref{02}, there is a subpath $p'$ of $p$,
such that $\Lab(p')$ is a unique subword of a unique word
$A_k^{\pm 1}.$ (More details of such an argument can be found in
\S 2 of \cite{Ol95}.) This implies that the words $\Lab(q_1)$ and
$\Lab(q_2)$ are freely conjugate to some $H$-words. Besides, it
follows from $T$-relations that $X\Lab(q_1)X^{-1}=\Lab(q_2)^{-1}$
for some $H$-word $X.$

\medskip
Lemma 18.4 in \cite{book} remains unchanged:

\begin{lm} \label{18.4} If $A$ and $B$ are simple in rank $i$ and $A$ is equal
in rank $i$ to $XB^lX\iv$ for some $X$ then $l = \pm 1$.
\end{lm}

The proof does not change much also. Notice only that by
definition, a simple word of rank $i\ge 1$ cannot be conjugate in
rank $i$ of an $H$-word. Therefore the diagram $\Delta$ from the
proof of Lemma 18.4 \cite{book}, cannot have hoops, and so the
analog of Lemma 19.4 \cite{book} (see Lemma \ref{19.4} below) is
applicable.

\medskip

The formulation of Lemma 18.5 changes as follows.

\begin{lm}\label{18.5}  If words $X$ and
$Y$ are conjugate in rank $i$ and $X$ is not a conjugate of any
$H$-word in rank $i$ then there exists a word $Z$ such that $X$ is
equal to $ZYZ\iv$ in rank $i$ and $|Z|\le \bar\alpha(|X|+|Y|)$.
\end{lm}

The additional assumption that $X$ is not a conjugate of an
$H$-word again allows us to apply Lemma \ref{19.4}, the analog of
Lemma 19.4 \cite{book}.

\medskip

Consider a section $q$ of a contour of some diagram $\Delta$,
where $\Lab(q)$ is freely equal to an $H$-word. We call $q$ an
\label{hsec}$H$-{\it section}. Assume that there is a $T$-cell
$\Pi$ in $\Delta$ with boundary path $p$ starting with a vertex
$o$ of $\partial(\Pi),$  and there is a path $t$ having no
self-intersections, which connects $o$ and a vertex $o'$ cutting
$q$ into a product $q_1q_2.$ We say that $\Pi$ is
\label{compq2}{\it compatible} with $q$, if the boundary label of
the path $q_1t^{-1}ptq_2,$ without self-intersections, is freely
equal to an $H$-word. Obviously, one can cut such a cell out of
$\Delta$ replacing $q$ by another $H$-section. The definition of
compatibility can be generalized in the case when $\Lab(q)$ is a
subword of a reduced form $W$ of an $H$-word with $|\Lab(q)|\ge
\varepsilon |W|$ in view of Lemma \ref{02}(2). (Such a path $q$ is
also an $H$-section by definition.)

Similarly, one defines the compatibility of two $T$-cells in a
diagram. A pair of two compatible $T$-cells can be substituted by
(at most) one $T$-cell. Therefore a g-reduced diagram of some rank
$i$ cannot contain pairs of compatible $T$-cells or $j$-{\em
pairs} in the terminology of \cite{book} for $j=1$. For
\label{bercel}{\em Burnside cells} (i.e. cells corresponding to
relations of the form $A^n=1$ where $A$ is a period of rank $j>0$)
we are using the same definition of $j$-pairs as in \cite{book}.

\begin{lm} \label{05} Assume that a word $W$ is a conjugate of
an $H$-word in rank $i.$ Then $W=UVU^{-1}$, in rank $i,$ for some
$H$-word $V$ and some word $U$, such that $|V|\le\bar\beta^{-1}
|W|$ and $|U|\le 2|W|.$
\end{lm}

\proof Consider a g-reduced annular diagram $\Delta$ for the
conjugacy of $W$ to some $H$-word $V$. We may assume that there
are no hoops in $\Delta,$ because otherwise one can replace
$\Delta$ by a diagram of a smaller type, that also demonstrates
the conjugacy of $W$ to a subword of the boundary label of a hoop.
For the same reason, we may assume that $\Delta$ has no $T$-cells
compatible with the contour, say $p,$ of $\Delta$ labeled by $V$.
Then $\Delta$ is an $A$-map by Lemma \ref{19.4}, and $p$ is a
smooth section of $\partial(\Delta)$ by Lemma \ref{19.5}. Then
$|V|\le \bar\beta^{-1}|W|$ by Lemma \ref{analteom17.1}, and the
word $U$ can be taken in such a way that
$$|U|\le\bar\alpha(1+\bar\beta^{-1})|W|\le 2|W|$$ by Lemma
\ref{18.5}. \endproof

The following Lemma implies malnormality of the subgroup $H$ in
$B(m,n).$

\begin{lm} \label{06} Let $V$ and $W$ be two non-trivial in
rank $i$ $H$-words and $W=UVU^{-1}$ in rank $i$. Then the word $U$
is itself equal in rank $i$ to an $H$-word.
\end{lm}

\proof Let $\Delta$ be a g-reduced annular diagram of rank $i$
with contours $p$ and $q$ labeled by $V$ and $W$ respectively. One
can assume that there is a path $t=q_- - p_-$ in $\Delta$, such
that $\phi(t)= U$ in rank $i$. Then it suffices to show that
$\phi(t)$ is freely equal to an $H$-word.

One can suppose that $\Delta$ has no $T$-cells compatible with the
contours, because otherwise our problem reduces to a similar
problem for an annular diagram of smaller type. If $\Delta$ has a
hoop, then $t$ can be constructed as a product $t_1t_2t_3$ where
the path $t_2$ cuts the hoop and $t_1,  t_3$ are the cutting paths
for annular subdiagrams of smaller types whose labels are
$H$-words by induction. Thus, we may assume that there are no
hoops in $\Delta.$ Then both $p$ and $q$ are
$\varepsilon$-sections by Lemma \ref{08}. It follows that $\Delta$
has rank 0, because otherwise the inequality
$\bar\gamma>2\varepsilon,$ Lemma \ref{cor16.1,16.2} and Lemma
\ref{19.4} lead to a contradiction. Then the desired property of
$\Delta$ follows from the malnormality of the subgroup generated
by the words $A_1,A_2,\dots$ in the absolutely free group $F({\cal
A}),$ which, in turn, can be easily derived from the statement of
Lemma \ref{02}(2). \endproof

Next as in \cite[Section 18]{book} we prove analogs of lemmas
18.6,18.7, 18.8 by a simultaneous induction on the sum $L$ of the
two periods involved in the formulations of these lemmas. The
formulations of these lemmas do not change.

\begin{lm} \label{18.6} (Analog of Lemma 18.6 \cite{book}.) Let $\Delta$ be a
g-reduced disc diagram of rank $i$ with contour $p_1q_1p_2q_2$
where $\Lab(q_1)$ and $\Lab(q_2\iv)$ are periodic words with
period $A$ simple in rank $i$. If $|p_1|, |p_2|<\alpha|A|$ and
$|q_1|, |q_2|> (\frac 56 h+1)|A|$, then $q_1$ and $q_2$ are
$A$-compatible in $\Delta$. (The inductive parameter is
$L=|A|+|A|$.)
\end{lm}

\proof As in the proof of Lemma 18.6 from \cite{book}, assume that
$A$ is a prefix of both $\phi(q_1)$ and $\phi(q_2\iv)$ by making
if necessary $p_1$ a little longer (by at most $1/2|A|$). The new
path is also denoted by $p_1$.

Notice first that if $\Lab(p_1)$ is not a conjugate of an $H$-word
in rank $i$ then in any g-reduced annular diagrams of rank $\le i$
where one of the contours is labeled by $\phi(p_1)$, there are no
hoops. So the proof proceeds as in \cite{book}.

Otherwise, as in \cite{book}, consider the annular diagram
$\bar\Delta$ of rank $i$ with boundary labels $W\equiv  \Lab(p_1)$
and $W'\equiv  \Lab(q')\Lab(p_2)^{-1}$.  Then $W, W'$ are
conjugate in rank $i$ of some $H$-word. If $W=1$ in rank $i$ then
there is nothing to prove. Thus we assume that the word $W$ is
non-trivial in rank $i$.

Inequalities (2), (3) from Section 18 of \cite{book} are still
valid: $$|W|<\bar\alpha |A|,\;\;
|W'|<(\bar\beta^{-1}-1)(\frac{5}{6} \delta^{-1}|A|+|A|)
+2|A|<(\bar\beta^{-1}-1)(|q_2|+2|A|)$$

By Lemma \ref{narrow}, one can find $H$-words $V$ and $V'$ such
that $$|V|\le \bar\beta^{-1}|W|<
\bar\beta^{-1}\bar\alpha|A|<|A|,\;\;
|V'|<\bar\beta^{-1}(\bar\beta^{-1}-1)
(5/6\delta^{-1}|A|+|A|)<\beta\delta^{-1}|A|,$$ and, in rank $i,$
$W=UVU^{-1}, W'=U'V'U'^{-1},$ where $|U|<2|A|,
|U'|<2\beta\delta^{-1}|A|$ by Lemma \ref{05}.

On the other hand, $QWQ^{-1}=W'$, in rank $i,$ for $Q\equiv
\Lab(q_2),$ $|Q|>\frac 56 \delta^{-1}|A|$ (recall that
$h=\delta\iv$). Therefore $H$-words $V$ and $V'$ are conjugate in
rank $i$ by the product $U'^{-1}QU.$ By Lemma \ref{06}, this
product must be equal, in rank $i$ to an $H$-word $X$.

Let $\Delta_0$ be a g-reduced diagram of rank $i$ for the equality
$(U')^{-1}QU=X.$ Its contour has a natural decomposition
$(u')^{-1}qux^{-1}.$ By changing $H$-word $X,$ one can suppose
that there are no $T$-cells compatible with $x$ in $\Delta.$ Hence
this section is an $\varepsilon$-section of $\Delta_0$ by Lemma
\ref{08}. Therefore the sum of the degrees of contiguity of an
arbitrary cell $\Pi$ of $\Delta$ to $q$ and to $x$, is less than
$\bar\alpha+\varepsilon<2/3$ by Lemmas \ref{anal15.8} and lemma
\ref{19.5} below. Recall that $Q$ is a periodic word with the
simple in rank $i$ period $A.$

Notice that $$|u'|+|u|<2(\beta\delta^{-1}+1)|A| <\frac 56
\alpha^2\delta^{-1}|A|<\alpha^2|q|.$$ Thus, there is an 0-bond
between $q$ and $x$ by Lemma \ref{narrow} (we have mentioned above
that the other possibility from that lemma is not possible). Then
$q=q(1)eq(2),$ for the contiguity edge $e$ of the 0-bond. If the
length of $q(1)$ or $q(2)$ is greater than
$2\alpha^{-2}(\beta\delta^{-1}+1)|A|,$ then one is able to get  a
new  $0$-bond  between $q$ and $x$ by Lemma \ref{narrow}.  By
keeping obtaining new and new such $0$-bonds, we conclude, that
there is a subpath $q(0)$ of $q$ and a contiguity subdiagram
$\Gamma(0)$ of rank 0 between $q(0)$ and $x$, where
$|q(0)|>|q|/2>2|A|>(1+\varepsilon)|A|.$

Consequently, a reduced form of the $H$-word $X^{-1}$ has an
$A$-periodic subword of length $>(1+\varepsilon)|A|.$ By Lemma
\ref{02}(1), the word $A$ is itself freely conjugate to an
$H$-word. This contradicts the definition of simple words and
completes the proof of our lemma.
\endproof

\medskip

The proofs of Lemmas 18.7, 18.8 remain unchanged.

The formulations of lemmas 18.9 - 19.3 \cite{book} are left
unchanged. The proofs are also unchanged. Only when we apply our
analog of Lemma 18.3 (Lemma \ref{18.3} above) in the proof of
Lemma 18.9, we take into account that by $T$-relations every
$H$-word has order dividing the odd number $n$ in rank $i\ge 1$.
In the proof of the first statement of Lemma 19.3, one should
remember that the reduced form of $W^n$ is included in the list of
$T$-relators for an arbitrary $H$-word $W.$

In order to prove analogs of Lemmas 19.4 and 19.5 \cite{book}, we
need two more Lemmas \ref{07} and \ref{08}. {\bf These two lemmas
and our analogs of Lemmas 19.4, 19.5 \cite{book} will be proved by
a simultaneous induction on the number of positive cells.}

\begin{lm} \label{07} Let $\Gamma$ be a contiguity subdiagram
of a cell $\pi$ to a $T$-cell $\Pi$ in a g-reduced diagram
$\Delta$ of rank $i+1.$ Then:

(1) the rank $r(\Gamma)$ is $0$ and

(2) if $\pi\ne\Pi,$ then $(\pi,\Gamma,\Pi)< \varepsilon.$
\end{lm}

\proof We can assume that the statement of the lemma is not true
and the triple $(\pi, \Gamma,\Pi)$ is a minimal counterexample in
the sense that $|\Gamma(2)|$ is minimal possible.

Since $|\Gamma(2)|<|\Delta(2)|$, the subdiagram $\Gamma$ is an
$A$-map by Lemma \ref{19.4}. Let $p_1q_1p_2q_2=
\partial(\pi, \Gamma,\Pi).$  Notice that $|p_1|=|p_2|=0,$ because
otherwise the principal cell of one of the bonds defining $\Gamma$
would be a positive cell with degree of \ct to $\Pi$ at least
$\varepsilon$, which would contradict the minimality of our
counterexample.

There are no cells in $\Gamma$ which are compatible with either
$q_1$ or $q_2$ because otherwise these cells would be compatible
with cells in $\Delta,$ and $\Delta$ would be a non-g-reduced
diagram. Then $q_2$ is an $\varepsilon$-section by Lemma
\ref{08}(2) for $\Gamma.$ Also $q_1$ is a smooth section in
$\partial(\Gamma)$ by Lemma \ref{19.5} or \ref{08} (2). Therefore
the sum of contiguity degrees of any cell of $\Gamma$ to $q_1$ and
$q_2$ is less than $\bar\alpha+\varepsilon<\bar\gamma,$ and so
$r(\Gamma)=0$ by Lemma \ref{cor16.1,16.2}.

Assume that $(\pi,\Gamma,\Pi)\ge\varepsilon.$ It follows from
Lemma \ref{02}(2) that the cells $\pi$ an $\Pi$ are compatible if
$\pi$ is a $T$-cell. But this is impossible in the g-reduced
diagram $\Delta.$ If $\pi$ is not a $T$-cell, we obtain a
contradiction to the definition of periods by Lemma \ref{02}(1).
Thus the lemma is proved by contradiction. \endproof

\begin{lm} \label{08} Let $\Gamma$ be a contiguity diagram of a
cell $\pi$ to an $H$-section $q$ of the boundary of a g-reduced
diagram $\Delta$ of rank $i+1.$ Assume there are no $T$-cells
compatible with $q$ in $\Delta.$ Then:

(1) $r(\Gamma)=0$ and

(2) $q$ is an $\varepsilon$-section of $\partial(\Delta).$
\end{lm}

\proof The proof is analogous to the proof of Lemma \ref{07}.
\endproof

The formulation of Lemma 19.4 of \cite{book} is modified:

\begin{lm}\label{19.4}
A g-reduced diagram of rank $i+1,$ having no hoops, is an $A$-map.
\end{lm}

\proof The perimeter of any $T$-cell is at least $n$ by Lemma
\ref{02} (3), and so condition A1 is true. The proof of part A2
remains the same except at the end the argument related to
$A$-compatibility should be replaced by a similar argument related
to the compatibility of a $T$-cell with an $H$-section, if the
cell $\Pi$ is a $T$-cell. The assumption $(\pi,\Gamma,\Pi)\ge
\varepsilon,$ in part A3 means that $\Pi$ is not a $T$-cell by
Lemma \ref{07}(2), since $\pi\ne\Pi$ (the diagram contains no
hoops). Therefore there are no changes in the proof of this part
if $\pi$ is not a $T$-cell as well. Otherwise the inequality
$|q_2|<(1+\gamma)|A|$ follows from Lemma \ref{07}(1) and Lemma
\ref{02}(1) (because $\gamma>\varepsilon$).
\endproof

We add one more possibility to the hypothesis of the analog of
Lemma 19.5 in \cite{book}. This lemma is now formulated as
follows.

\begin{lm} \label{19.5}  Let $p$ be a section of the
contour of a g-reduced diagram $\Delta$ of rank $i+1$ and one of
the following possibilities holds:
\begin{itemize}
\item The label of $p$ is an $A$-periodic word where $A$ is simple in rank $i+1$,
\item The label of $p$ is an $A$-periodic word where $A$ is a period of rank $k\le i+1$,
and $\Delta$ has no cells of rank $k$ $A$-compatible with $p$.
\item The label of $p$ is an $H$-word and there are no $T$-cells
in $\Delta$ which are compatible with $p$.
\end{itemize}
(If $p$ is a cyclic section then in the first two cases we further
require that $\Lab(p)\equiv A^m$ for some integer $m$.)  Then $p$
is a smooth section of rank $|A|$ in the first and the second
cases and of rank 1 in the third case, in $\partial(\Delta)$.
\end{lm}

\proof We need to prove the two properties from the definition of
a smooth section.

1. The proof that every subpath of $p$ of length $\le \max(k,2)$
is geodesic in $\Delta$ coincides with part 1) in the proof of
Lemma 19.5 \cite{book} if $\phi(p)$ is not an $H$-word. In the
case when $\Lab(p)$ is an $H$-word, this property follows from
Lemmas \ref{07}, \ref{cor17.1} and \ref{02}.

2. The proof that for every \ct submap $\Gamma$ of a cell $\pi$ to
$p$ satisfying $(\pi, \Gamma,p)\ge \varepsilon$ we have $|\Gamma
\bigwedge p| < (1+\gamma)(i+1)$ proceeds along the lines of part
2) of Lemma 19.5 \cite{book}. But we need to consider a new case
when $\pi$ is a $T$-cell. In this case, $r(\Gamma)=0$ as in the
proof of Lemma \ref{07} (one should just replace $\Pi$ by $p$).
Then the inequality $|q_2|\le (1+\gamma)(i+1)$ follows from Lemma
\ref{02}(1) and from the definition of $T$-relations.

If $\Lab(p)$ is an $H$-word,  then the second property of smooth
sections follows from Lemma \ref{08} (2). \endproof

\begin{lm} \label{09} The constructed group ${\bf G}(\infty)$ satisfies identity
$x^n=1.$
\end{lm}

\proof This immediately follows from Lemma 19.3 (1) \cite{book}
(its formulation remains the same in our situation). That lemma
states that the order of every word of length at most $i$, in rank
$i$, is a divisor of $n$.
\endproof

\medskip

{\bf Proof of Theorem \ref{SQ}.} By Lemma \ref{03}, it remains to
show that for arbitrary normal in $H$ subgroup $N,$ the
intersection of the normal closure $L$ of $N$ in $B(m,n)$ and $H$
is equal to $N,$ i.e. any $H$-word $w$ representing an element of
$L,$ represents an element of $N.$

Consider the group ${\bf G}(\infty)$ constructed above for this
particular $N$. Since $w$ represents an element from $L$ in
$B(m,n)$, it belongs to the normal closure, in the absolutely free
group, of set $T$ plus the set of all words of the form $A^n=1$.
By Lemma \ref{09} and the definition of ${\bf G}(1)$ then $w$
represents $1$ in ${\bf G}(\infty)$.
 Therefore there
is a g-reduced diagram $\Delta$ of some rank $i$ whose boundary
label $\Lab(q)$ is freely equal to $w.$

By contradiction assume that the $H$-word $w$ does not represent
an element of $N$, and the diagram $\Delta$ has minimal possible
number of positive cells.

By Lemma \ref{th16.2}, there are a cell $\Pi$ of $\Delta$ and a
contiguity subdiagram $\Gamma$ of rank 0 with $(\Pi,\Gamma,q)\ge
\varepsilon.$ If $\Pi$ is a $T$-cell, then it must be compatible
with $q$ by Lemma \ref{02}(2). Therefore when we cut the cell
$\Pi$ from $\Delta$, we multiply (in the free group) the boundary
label of $\Delta$ by a word $U$ of the form $VWV\iv$ where $W$ is
the boundary label of $\Pi$, and $V$ is an $H$-word. Since $W$
represents an element of $N$,  $V$ represents an element of $H$,
and $N$ is normal in $H$, the word $U$ represents an element of
$N$. Hence the boundary label of the resulting diagram $\Delta_1$
is an $H$-word which does not represent an element from $N$. So
$\Delta_1$ is a counterexample with a smaller number of positive
cells. If $\Pi$ is not a $T$-cell, i.e. it corresponds to a
relation $A^n=1,$ then we get a contradiction with Lemma
\ref{08}(2).

It is clear that if $H$ is a subgroup of $B(m,n)$ satisfying
properties (1) and (2) from Theorem \ref{SQ} then any subgroup
generated by finitely many of the free generators of $H$ is also
free in the variety of Burnside groups of exponent $n$, and has
the congruence extension property because retracts of subgroups
with the congruence extension property also have this property and
because the congruence extension property is transitive (if $H_1$
has the congruence extension property in $H$ and $H$ has it in
${\bf G}$ then $H_1$ has the congruence extension property in
${\bf G}$).

Theorem \ref{SQ} (and Theorem \ref{SQintr}) are proved.

\subsection{Burnside groups and free products}
\label{bgfp}

The following lemma is Theorem VI.3.2 from \cite{Adian1} (or
Theorem 19.5 from \cite{book}).

\begin{lm} \label{VI} The centralizer of every non-trivial element
of $B(m,n)$ is cyclic of order $n$.
\end{lm}

We are going to use the following Lemma repeatedly in Section
\ref{properties}.

\begin{lm} \label{additional} If $a,b,c\in B(m,n)$, $n>>1$ odd,
$a\ne 1$, $c \ne 1$, $[a,c]=1$, $[bab\iv, c]=1$ then
$[b,c]=[b,a]=1$.
\end{lm}

\proof Indeed, suppose that $[a,c]=1$, $[bab\iv, c]=1$. By Lemma
\ref{VI}, the centralizer of every non-trivial element of $B(m,n)$
is cyclic of order $n$. Hence $a, bab\iv,c$ belong to the same
cyclic subgroup $C$ of $B(m,n)$. Therefore $bCb\iv$ and $C$ are
non-trivially intersecting cyclic subgroups of order $n$. By Lemma
\ref{VI}, $bCb\iv =C$. Hence $\la b, C\ra$ is a finite subgroup of
$B(m,n)$ (or order at most $n^2$). By Theorem 19.6 from
\cite{book}, $\la b, C\ra$ is cyclic of order $n$, so $b\in C$.
Hence $[b,c]=1$ as required.
\endproof

In order to prove Proposition \ref{propsss} below, we will need
some properties of free products. Let $P$ and $Q$ be two free
groups. Elements of the free product $P*Q$ are represented by
words of the form $p_1q_1p_2...q_sp_{s+1}$ where $p_i\in P$,
$i=1,...,s+1$, is called a $P$-syllable, $q_i\in Q$, $i=1,...,s$,
is called a $Q$-syllable. By definition, a $P$-syllable (resp.
$Q$-syllable) is a maximal subword over the set of generators of
$P$ (resp. $Q$). Elements of $P$ and $Q$ are called 1-syllable
elements.

A word $p_1q_1...p_sq_sp_{s+1}$ from $P *Q$, where $p_i, q_j$ are
syllables, is called {\em \label{str}*-reduced in $P*Q$} if none
of its syllables except possibly for $p_1, p_{s+1}$ is equal to 1.
If a product is *-reduced, it is not equal in $P*Q$ to such a
product with smaller number of syllables. The $P$-syllable
$p_{s+1}$ is called the {\em last} syllable of
$p_1q_1...p_sq_sp_{s+1}$ (even if $p_{s+1}$ is empty). If in
$p_1q_1...p_sq_sp_{s+1}$, a syllable, say, $p_i$, is freely equal
to 1, we can remove it and {\em glue} $q_{i-1}$ with $q_i$ to form
a new $Q$-syllable. Applying this operation several times, we get
a *-reduced in $P*Q$ form of the product which is unique. Notice
that a word which is *-reduced in $P*Q$ may not be freely reduced
(as a word written in the free generators of $P$ and $Q$).

Let $\psi:P\to \bar P$ be a homomorphism of $P$ onto a group $\bar
P$ satisfying the identity $x^n=1$.  We say that a word $h$ of $P
* Q$ is \label{ptriv}{\em $\psi$-trivial} if the $\psi$-image of every syllable
of $h$ is 1 in $\bar P$.

\begin{lm} \label{bernsyl}
Let $h_1, h_2$ be words in $P*Q$ which are equal modulo Burnside
relations. Then the following properties hold.

(i) If $h_1$ is $\psi$-trivial, and $h_2$ is of minimal length
among all words which are equal to $h_1$ modulo Burnside relations
then $h_2$ is $\psi$-trivial.

(ii) Suppose that $h_1$ has exactly two $P$-syllables $p,p'$ with
non-trivial $\psi$-images, $h_1\equiv xpyp'z$, and $y$ is not
equal to a word in $P$ modulo Burnside relations. Then $h_2$ is
not $\psi$-trivial.

(iii) If $h_1$ is $\psi$-trivial, $h_1\ne 1$ modulo Burnside
relations, $h_2\in P$, and modulo Burnside relations
$h_1=xh_2x\iv$ for some $x\in P*Q$, then there exists a word
$x_1\in P*Q$ such that $x_1h_2x_1\iv=h_1$ modulo Burnside
relations and the word $x_1h_2x_1\iv$ is $\psi$-trivial.
\end{lm}

\proof (i) Obviously we can assume that $h_1$ is *-reduced.
Consider a g-reduced diagram $\Delta$ over the graded presentation
of the free Burnside quotient of $P*Q$ from Section 18 of
\cite{book} corresponding to the equality $h_1=h_2$. If the rank
of $\Delta$ is 0 then $h_1=h_2$ in $P*Q$, and each syllable of the
word $h_2$ (of minimal length) must be a product of several
syllables of $h_1$. Hence $h_2$ is $\psi$-trivial.

Suppose that the rank of $\Delta$ is $> 0$, and
$\partial(\Delta)=s_1s_2$ where $\Lab(s_1)=h_1$,
$\Lab(s_2\iv)=h_2$. Since $s_2$ is geodesic, $s_2$ is a smooth
section of rank $|s_2|$. By Lemma \ref{cor16.1,16.2}, part a),
and Lemma \ref{anal15.8} there exists a cell $\Pi$ of rank $>0$
and its \ct subdiagram $\Gamma$ to $s_1$ of \ct degree
$>1/2-\gamma-\alpha>3\varepsilon$. By Lemma \ref{Anal16.2},
$\Delta$ contains a cell $\pi$ and a rank 0 \ct subdiagram
$\Gamma'$ of $\pi$ to $s_1$, such that the \ct degree of $\Gamma'$
is $>\varepsilon$. Let $t$ be the \ct arc of $\Gamma'$ on
$\partial(\pi)$, $|t|>\varepsilon|\partial(\pi)|$. There are two
possibilities.

\medskip

{\bf Case 1.} The period of $\Lab(\partial(\pi))$ is a 1-syllable
word. Then $\Lab(\partial(\pi))$ is a 1-syllable word, and one can
cut the subdiagram $\pi\cup\Gamma'$ from $\Delta$ preserving the
$\psi$-triviality of $\Lab(s_1)$. Here we use the fact that $\bar
P$ satisfies the identity $x^n=1$.

\medskip
{\bf Case 2.} The period $A$ of $\Lab(\partial(\pi))$ contains at
least 2 syllables. Since $\Gamma'$ is a diagram of rank 0, and
$\varepsilon n>2$, every syllable of $A$ is a syllable of
$\Lab(s_1)^{\pm 1}=h_1^{\pm 1}$. Hence $A$ is $\psi$-trivial
(since $h_1$ is $\psi$-trivial). Therefore, again, one can cut the
subdiagram $\Gamma'\cup\pi$ from $\Delta$ preserving the
$\psi$-triviality of $\Lab(s_1)$: we replace subword $A$ of
$\Lab(s_1)^{\pm 1}$ by $A^{1-n}$ and reduce the resulting word.

\medskip

Thus in each of these cases we replace $\Delta$ by a diagram of
smaller type, and complete the proof by induction on the type of
$\Delta$.

\medskip

(ii) By contradiction assume that $h_2$ is $\psi$-trivial. We have
that $pyp'=x\iv h_2z\iv$. By part (i) of this lemma we can replace
$x\iv h_2z\iv$ by a $\psi$-trivial word $h$ of minimal possible
length. Consider a g-reduced diagram $\Delta$ with
$\partial(\Delta)=s_1t_1s_2t_2$ where $\Lab(s_1)=p,$
$\Lab(t_1)=y$, $\Lab(s_2)=p'$, $\Lab(t_2\iv)=h$. If the rank of
$\Delta$ is 0 then $pyp'=h$ in $P*Q$. Since $y$ is not equal to a
word from $P$ and is $\psi$-trivial, and $\psi(p)\ne 1$, we have
that the *-reduced in $P*Q$ form of $pyp'$ is not $\psi$-trivial,
a contradiction with the $\psi$-triviality of $h$.

If the rank of $\Delta$ is $>0$ then as in part (i), we get a cell
$\pi$ and a \ct subdiagram $\Gamma'$ of $\pi$ to one of the
subpaths $s_1,t_1,s_2$ with \ct degree $\ge \varepsilon$, such
that $\Gamma'$ has rank 0. Since $\Lab(s_1)$ and $\Lab(s_2)$ are
1-syllable words, the period $A$ of $\pi$ is a 1-syllable word if
$\Gamma'$ is a \ct subdiagram to $s_1$ or $s_2$. If $\Gamma'$ is a
\ct subdiagram to $t_1$, $A$ is a $\psi$-trivial word. As in part
(i) we can remove $\Gamma'\cup\pi$ from $\Delta$ preserving the
structure of the boundary of $\Delta$ and the condition that
$\Lab(y)$ is a $\psi$-trivial, and does not belong to $P$ modulo
Burnside relations, $\Lab(s_1)$, $\Lab(s_2)$ are 1-syllable words
from $P$, $\psi(\Lab(s_1)), \psi(\Lab(s_2))\ne 1$. We can complete
the proof by induction of the type of $\Delta$.

\medskip

(iii) Consider a g-reduced annular diagram $\Delta$ with pointed
contours $s_1$, $s_2$, where $\Lab(s_1)\equiv h_1, \Lab(s_2)\equiv
h_2$. If the rank of $\Delta$ is 0 then $h_1=xh_2x\iv$ in $P*Q$
and $xh_2x\iv$ is $\psi$-trivial.

Suppose that the rank of $\Delta$ is $>0$. As in parts (i), (ii),
we can remove a cell of rank $>0$ from $\Delta$ and change $s_1$
or $s_2$ preserving $\Lab(s_1), \Lab(s_2)$ modulo Burnside
relations. The labels of the contours $s_1'$, $s_2'$ of the new
annular diagram satisfy all the conditions of part (iii) of the
lemma, and again we can complete the proof by induction on the
type of $\Delta$.
\endproof

\begin{lm} \label{phitriv}
Suppose that $h, g\in P*Q$, and there is no $x\in P*Q$ such that
$g=xpx\iv, h=xp'x\iv$ modulo Burnside relations for some $p,p'\in
P$. Suppose also that $h$, the *-reduced forms of $ghg\iv$, and a
word $b$ which is equal to $g\iv hg$ modulo Burnside relations are
$\psi$-trivial. Then $g$ is $\psi$-trivial.
\end{lm}

\proof Suppose that $g$ is not $\psi$-trivial. Let $g\equiv
g_1pg_2$ where $g_2$ is $\psi$-trivial, $p$ is a $P$-syllable,
$\psi(p)\ne 1$. Consider the *-reduced form of $ghg\iv\equiv
g_1pg_2hg_2\iv p\iv g_1\iv$. Since this word is $\psi$-trivial,
and $g_2$ is also $\psi$-trivial, the syllable $p$ is glued with a
$P$-syllable $p_h$ of $h$ and with the $P$-syllable $p\iv$ of
$g\iv$. This means that $h=g_2\iv p_hg_2$ in the free group $P*Q$,
and the *-reduced form of $ghg\iv$ is $g_1pp_hp\iv g_1\iv$. This
implies that $g_1$ is $\psi$-trivial, so $g$ has only one
$P$-syllable with non-trivial $\psi$-image. Now consider the
*-reduced form $b'$ of the word $g\iv hg=g_2\iv p\iv g_1\iv g_2\iv
p_hg_2g_1pg_2$.

Since $b$ is equal to $b'$ modulo Burnside relations and $b$ is
$\psi$-trivial, we can conclude by Lemma \ref{bernsyl}(ii) that
$g_1\iv g_2\iv p_hg_2g_1$ is equal to a word $p''\in P$ modulo
Burnside relations. Taking the projection of the equality
$(g_2g_1)\iv p_h g_2g_1=p''$ onto the subgroup $P$, we get
$P(g_2g_1)\iv p_hP(g_2g_1)=p''$ modulo Burnside relations, where
$P(g_2g_1)$ is the $P$-projection of $g_2g_1$. Therefore
$P(g_2g_1)\iv (g_2g_1)$ centralizes $p_h$ modulo Burnside
relations. Since the centralizers of $p_h$ in the free Burnside
quotient of $P$ and $P*Q$ have order $n$ (by Lemma \ref{VI}),
these centralizers coincide. Hence $P(g_2g_1)\iv g_2g_1$ belongs
to $P$ modulo Burnside relations, so $g_2g_1\in P$ modulo Burnside
relations. Then $g_2gg_2\iv\in P $, $g_2hg_2\iv =g_2g_2\iv
p_hg_2g_2\iv=p_h\in P$ modulo Burnside relations, which is
impossible by our assumption.
\endproof

\section{The presentation of the group}
\label{presentation}

In this section, we define an $S$-machine consisting of a set of
(positive) $S$-rules \label{s+1}$\sss^+$. Using this $S$-machine,
we define our groups ${\cal H}$ (Theorem \ref{CEP}).

 Fix some numbers first. We use the
same odd number $n>>1$ as in the previous sections. Let $N>n$, and
let $m=m(n)$ be the number which exists by Theorem \ref{SQ}.

\subsection{$S$-machines and their interpretation}
\label{smachines}

We start with definition of the generating set ${\bf X}$ of the
group $H$ and the alphabet of the $S$-machine.

Let \label{aaaa}${\cal A}$ be the set $\{a_1,...,a_m\}$. Let
\label{omeg}$\Omega$ be the set $\{0,...,2n+1\}$.

The set $\sss^+$ that we shall define in the next section, will be
the union of disjoint sets \label{s+o1}$\sss^+_\omega$, $\omega\in
\Omega$, and the set $\{\bc_0,...,\bc_{2n}\}$ called the
\label{cor1}{\em set of connecting rules}. For each $\omega\in
\{0,...,2n\}$, the set $\sss^+_\omega$ will consist of $m$ rules,
one for each letter in ${\cal A}$. The rules from
$\cup_{\omega=0}^{2n}\sss_\omega\cup\{\bc_0,...,\bc_{2n}\}$ will
be called \label{wor1}{\em working rules}.

The rules from $\sss_{2n+1}^+$ will be called {\em cleaning
rules}. The set $\sss_{2n+1}^+$ will be in one-to-one
correspondence with the set of all working rules (thus the set of
cleaning rules and the set of working rules are disjoint but have
the same cardinality).

Consider the alphabet \label{bkk}$$\bkk=\{\tk(i,\omega),
\kappa(i), \bk(i),\rk(i,j,\omega), \lk(i,j,\omega),  i=1,...,N,
j=1,...,n, \omega\in \Omega\}.$$ Elements of $\bkk$ will be called
$\bkk$-letters. Let also \label{kk}$\kk$ be the set of letters
obtained from $\bkk$ by removing the $\Omega$-coordinate $\omega$,
and let us call the corresponding map from $\bkk$ to $\kk$ the
$\kk$-projection. The preimage of every letter $z\in\kk$ under the
$\kk$-projection will be denoted by \label{kz}$\kk(z)$. For every
$z\in (\kk)\iv$ we denote $\kk(z\iv)\iv$ by $\kk(z)$. Sometimes it
will be convenient to assume that letters with $\kk$-projection
$\kappa(i)$ or $\bk(i)$ also have $\Omega$-coordinates but we
identify all letters $\kappa(i,\omega)$ with $\kappa(i)$ and all
letters $\bk(i,\omega)$ with $\bk(i)$.

If $U$ is any word containing letters from $\kk$, and $\omega\in
\Omega$, then we define \label{uo}$U(\omega)$ as the word,
obtained from $U$ by replacing every letter $z$ from
$\kk\backslash\{\kappa(i), \bk(i), i=1,...,N\}$ by $z(\omega)$ (in
other words, we add the coordinate $\omega$ to each letter from
$\kk$ in $U$ except $\kappa(i), \bk(i)$). In this case we shall
call $\omega$ the $\Omega$-coordinate of $U(\omega)$.

For every $z\in \kk\backslash \{\tk(i), \bk(i), \kappa(j), 1\le
i\le N, j=2,...,N \}$ consider a copy \label{az}$A(z)$ of the
alphabet ${\cal A}$. For each $z\in \{\tk(i),\bk(i), 1\le i\le
N\}$, consider the set $A(z)$ which is a copy of the set of
working rules: with every working rule $\tau$ we associate a
letter \label{btz}$\bt{z}$ in $A(z)$ (all sets $A(z)$ are
disjoint). Notice that we have not defined the sets
$A(\kappa(j))$, $j=2,...,N$. For convenience we assume that all
these sets are empty.

Let \label{aa}$\aaa$ be the union of all the sets $A(z)$, $z\in
\kk$. Elements of this set will be called $\aaa$-letters.

Consider also the set of $\rr$-letters \label{RR}$\rr=\{r(\tau,z),
\tau\in \sss^+, z\in \kk\}$.

Now let \label{xxx}$$\xxx=\bkk\cup \aaa\cup\rr.$$ This is the
generating set of our group ${\cal H}$. We rank some of the
letters in $\xxx$ in the following way: $\bkk>\rr>\aaa$, that is
letters from $\bkk$ have rank higher than all other letters, etc.
For every word over $\xxx$ we can consider its $\bkk$-length, etc.

Consider the following word \label{lambda}$\Lambda$ over $\kk$:

\begin{equation} \label{hub}
\begin{array}{l}
\tk(1)\kappa(1)\bk(1)\rk(1,1)\lk(1,1)...\rk(1,n)\lk(1,n)...\\
\tk(N)\kappa(N)\bk(N)\rk(N,1)\lk(N,1)...\rk(N,n)\lk(N,n)
\end{array}
\end{equation}

The word $\Lambda(0)$ will be called the \label{hub2}{\em hub}.

We define a circular order on $\kk$ according to the appearance in
$\Lambda$ if we consider the word $\Lambda$ as written on a
circle. For every $z\in \kk$ let \label{zmp}$z_-$ be the letter
immediately preceding $z$ in our order, and let $z_+$ be the
letter next to $z$ in that order. We extend that order also on
letters from $\kk\iv$. If $z\in \kk$ then $(z\iv)_-\equiv
(z_+)\iv$, $(z\iv)_+\equiv (z_-)\iv$. We define a corresponding
circular order on $\bkk\cup\bkk\iv$: $$z(\omega)_+\equiv
z_+(\omega).$$

For every $z\in \kk\iv$ we set $A(z)=A(z\iv_-)$. Also by
definition for every $\bar z\in \bkk$ with $\kk$-projection $z$ we
let $A(\bar z^{\pm 1})=A(z^{\pm 1})$.

The language of \label{adm}{\em admissible words} consists of all
words of the form $y_1u_1y_2u_2...y_t$ where $y_i\in \bkk$, $u_i$
are words in $A(y_i)$, $i=1,2,...,t-1$, and for every $i=1,2...$,
either $y_{i+1}\equiv (y_i)_+$ or $y_{i+1}\equiv y_i\iv$. Thus all
$\bkk$-letters in an admissible word have the same \label{oca}{\em
$\Omega$-coordinates}. (Notice that $\Lambda(\omega)$ is an
admissible word for every $\omega$.) The subword $y_iu_iy_{i+1}$
is called the \label{sector}$y_i$-{\em sector} of the admissible
word, $i=1,2....$. Two admissible words $y_1u_1y_2u_2...y_t$ and
$y_1'u_1'y_2'u_2'...y_{t'}'$ are considered {\em equal} if $t=t'$,
$y_i\equiv y_i'$, $i=1,...,t$, and $u_i=u_i'$, $i=1,...,t-1$, in
the free Burnside group (later on we will identify these
admissible words if $y_i\equiv y_i'$ and $u_i=u_i'$ in the group
$\hnka$ defined below).

The set of defining relations of the group ${\cal H}$ consists of
the hub $\Lambda(0)$ and relations corresponding to the $S$-rules.
We are going to present these relations now.

We start with defining the set of \label{srule}$S$-rules (an
$S$-machine in the terminology of \cite{SBR}) and the
corresponding defining relations.

Each rule means simultaneous replacing certain subwords in
admissible words by other subwords (the rule specifies the
subwords which need be replaced and the replacements).

We shall write a rule $\tau$ in the form $[z_1\too
v((z_1)_-)z_1u(z_1),..., z_s\too
v((z_s)_-)z_su(z_s);\omega\too\omega']$ where $z_i\in \kk$, each
of $u(z),v(z)$ is either a letter from $A(z)^{\pm 1}$ or empty,
$z\in\kk$.

Some of the arrows in a rule can have the form
\label{tool}$\tool$. This means that the rule can be applied if
the corresponding sector does not have $\aaa$-letters in it. We
shall say that the corresponding sectors are \label{lock}{\em
locked} by the rule and the rule {\em locks} these sectors. If a
$z$-sector is locked by the rule, $u(z)$ and $v(z)$ must be empty.

If $z\in\kk$ does not appear in this rule then we set $u(z)$,
$v(z_-)$ to be empty. Thus with every rule $\tau$ and every
$z\in\kk$ we associate two words over $\aaa$: $u(z)$ and $v(z)$.
We can extend this definition to $z\in \kk\iv$ by saying that
$v(z\iv)\equiv u(z)\iv$, $u(z\iv)\equiv v(z_-)\iv$.

If for some $z\in \kk$, $u(z)$ is not empty then we call $\tau$
\label{lra}{\em left active} for $z$-sectors, if $v(z)$ is not
empty then $\tau$ is {\em right active} for $z$-sectors.

If $\tau$ is a rule as above then for every $z\in\kk\cup\kk\iv$ we
define two words \label{lrtz}$L(\tau,z)\equiv r(\tau,z_-)v(z_-)$
and $R(\tau,z)\equiv r(\tau,z)u(z)\iv$. If $h\equiv
\tau_1\tau_2...\tau_t$ is a word in the alphabet $\sss$, we set
$L(h,z)\equiv L(\tau_1,z)L(\tau_2,z)...L(\tau_t,z)$, $R(h,z)\equiv
R(\tau_1,z)R(\tau_2,z)...R(\tau_t,z)$.

The rule $\tau$ is \label{app}{\em applicable} to a one-letter
admissible word $y$ if and only if the $\Omega$-coordinate of $y$
is $\omega$.

Let $W\equiv y_1uy_2$ be an admissible word with two
$\bkk$-letters. The following two conditions must hold for $\tau$
to be applicable to $W$:

(i) $y_1$ and $y_2$ belong to $\kk(\omega)$;

(ii) if $y_1$ is one of the letters $z_i$, $i=1,...,s$, $\tau$
locks $y_1$-sectors then $u$ must be empty.

If conditions (i) and (ii) hold then the result of the application
of $\tau$ to $W$ is the word $$y_1'u(y_1)uv((y_2)_-)y_2'$$ where
$y_i'$ are obtained from $y_i$ by replacing the
$\Omega$-coordinate $\omega$ by $\omega'$ (recall that $\tau$
contains the instruction  $\omega\to\omega'$).

If $W\equiv y_1u_1y_2...u_ty_t$ is an admissible word with $t\ge
2$ $\kk(\omega)$-letters then the rule $\tau$ is applicable to $W$
if and only if it is applicable to every sector of $W$. In order
to apply an $S$-rule to $W$ we apply it to every sector of $W$.

Let $W'$ be the resulting word. Notice that $W'$ is clearly again
an admissible word.

With every rule $\tau=[z_1\too v((z_1)_-)z_1u(z_1),..., z_s\too
v((z_s)_-)z_su(z_s);\omega\too\omega']$ we associate the
\label{invr} {\em inverse} rule $\tau\iv=[z_1\too v((z_1)_-)\iv
z_1u_(z_1)\iv,...,z_s\too v((z_s)_-)\iv z_su(z_s)\iv; \omega'\to
\omega]$. Clearly if $\tau$ is applicable to an admissible word
$W$ and $W'$ is the result of application of $\tau$ to $W$, then
$\tau\iv$ is applicable to $W'$ and $W$ is the result of
application of $\tau\iv$ to $W'$.

The rules from $\sss^+$, described below, will be called
\label{pnr}{\em positive}, the inverses of these rules will be
called {\em negative}.

Now let us describe the procedure of \label{crr}{\em converting}
an $S$-rule into a set of relations. After that we shall present
the list of rules $\sss$ corresponding to our group ${\cal H}$.

Let $\tau=[z_1\too v((z_1)_-)z_1u(z_1),...,z_s\to
v((z_s)_-)z_su(z_s); \omega\too\omega']$ be an $S$-rule. Then for
every $z\in \kk$, we introduce a relation

$$r(\tau,z_-)\iv z(\omega)r(\tau,z) = v(z_-)z(\omega')u(z).$$

If (and only if) $\tau$ does not lock $z$-sectors, we also include
all relations of the form $$r(\tau,z)b = br(\tau,z), b\in A(z).$$

For every $z\in \kk$ let \label{rrz}$\rr(z)$ denote the set
$\{r(\tau, z)\ | \ \tau\in \sss\}$.

To make formulas look less frightening let us agree on the
following:

{\bf
\begin{quotation}
In every word of the form $u_{ar}zv_{ar}$ appearing in an $S$-rule
or in a defining relation of our group, or in an admissible word,
where $z\in\kk(\omega)$, $u_{ar}$, $v_{ar}$ are words in the
alphabet $\aaa\cup\rr$, we shall always assume that the letters in
$u_{ar}$ are from $A(z_-)\cup\rr(z_-)$ only, and the letters in
$v_{ar}$ are from $A(z)\cup\rr(z)$ only. So we shall omit the
$\kk$-coordinates in these words. We shall also omit the
$\Omega$-coordinate in $\bkk$-letters if the value of it is clear.
\end{quotation}
}

{\bf Example.} Suppose, for example, that the rule $\tau$ has the
form $$[\kappa(1)\too \kappa(1) a, \rk(1,1)\too \rk(1,1)a,
\lk(1,2)\tool a\lk(1,2), 1\too 2 ]$$ (this rule is not in $\sss$,
we use it just as an example). Then the set of relations
corresponding to this rule consists of the following relations.

\begin{itemize}
\item $r(\tau)\iv \kappa(1) r(\tau) = \kappa(1)a;$
\item $r(\tau)\iv \rk(1,1,1) r(\tau)= \rk(1,1,2)a;$
\item $r(\tau)\iv\lk(1,2,1) r(\tau) = a\lk(1,2,2);$
\item $r(\tau)\iv z(1)r(\tau)=z(2)$ for every $z\in\kk, z\nee
\kappa(1), \rk(1,1)$.
\item $r(\tau)b=br(\tau)$ for every $b\in A(z)$ and $z\nee
\lk(1,2)$.
\end{itemize}

Notice that, for example in the relation $r(\tau)\iv \rk(1,1,1)
r(\tau)= \rk(1,1,2)a$, the two letters $r(\tau)$ are different
according to our agreement: the first letter is
$r(\tau,\kappa(1))$ and the second letter is $r(\tau,\bk(1))$.
This rule $\tau$ is left active for $\kappa(1)$-sectors and for
$\rk(1,1)$-sectors, right active for $\lk(1,2)$-sectors, and
locking for $\lk(1,1)$-sectors.

Here is one of the main  properties of the systems of relations
corresponding to $S$-rules. The lemma immediately follows from
definitions.

\begin{lm} \label{lm1}
Let $W$ be an admissible word, starting with $z\in \bkk$ and
ending with $z'\in \bkk$, let $\tau=[z_1\too
v((z_1)_-)z_1u(z_1),...,z_s\to v((z_s)_-)z_su(z_s);
\omega\too\omega']$ be an $S$-rule applicable to $W$ (some of the
arrows $\too$ may be of the form $\tool$). Then:

\begin{enumerate}
\item  $L(\tau,z)\iv WR(\tau,z')$ is equal modulo the relations corresponding to $\tau$
to the admissible word $W'$ obtained by applying $\tau$ to $W$.

\item If $\tau$ locks $z$-sectors then the
$z$-sectors of $W$ and $W'$ have trivial retraction to the
subgroup $\la\aaa\ra$ in the free group generated by
$\aaa\cup\bkk$.

\item For every $z\in \kk$,
$L(\tau,z)z(\omega')=z(\omega)R(\tau,z)$ is one of the defining
relations corresponding to $\tau$.
\end{enumerate}
\end{lm}

\subsection{The $S$-machine $\sss$}
\label{thesmachine}

Now let us define our collection of $S$-rules \label{s+2}$\sss^+$.
As we have announced in the Section \ref{smachines}, the set
$\sss^+$ will be divided into $2n+3$ subsets: $$\sss_0^+,
...,\sss_{2n+1}^+, \{ {\bf c}_i; i=0,...,2n\}.$$ The sets $\sss_i$
will be called {\em steps}, $\bc_i$ will be called
\label{cor2}{\em connecting rules}. Rules from $\sss_{2n+1}$ are
called \label{cler2}{\em cleaning rules}, other rules are called
\label{wor2}{\em working}.

In each set $\sss_i^+$, $i\le 2n$, each rule corresponds to a
letter $a\in {\cal A}$, it will be denoted by
\label{tia}$\tau(i,a)$.

The rule $\tau(0,a)$ in the set $\sss_0^+$ is the following

$$\left[\begin{array}{l}\tk(i)\tool\tk(i), \kappa(i)\tool
\kappa(i), \bk(i)\tool \bk(i), \rk(i,j)\too \rk(i,j)a\iv,
\lk(i,j)\tool \lk(i,j),\\ i=1,...,N, j=1,...,n;0\to
0\end{array}\right].$$ This rule inserts copies of the letter
$a\iv$ to the right of $\rk(i,j)$. A sequence
$$\tau(0,a_{i_1})^{\pm 1},...,\tau(0,a_{i_s})^{\pm 1}$$ of these
rules inserts copies of the word $(a_{i_1}^{\pm 1}...a_{i_s}^{\pm
1})\iv$ in $\rk(i,j)$-sectors. These rules lock all other sectors.
The connecting rule $\bc_0$ has the form

$$[\tk(i)\tool\tk(i), \kappa(i)\tool \kappa(i), \bk(i)\tool
\bk(i), \lk(i,j)\tool\lk(i,j), i=1,...,N, j=1,...,n; 0\to 1].$$ It
locks all sectors except for $\rk(i,j)$-sectors.

The rule $\tau(1,a)$ has the form

$$\left[\begin{array}{l}\kappa(i)\tool \bt{1,a}\kappa(i),
\rk(i,1)\too \bt{1,a}\,\rk(i,1)a,\\ \rk(i,j)\too \rk(i,j)a,
\lk(i,j)\too \lk(i,j)a,\\ 1\le i\le N, 2\le j\le n; 1\too
1\end{array}\right].$$ This rule moves letters from
$\rk(i,j)$-sectors to $\lk(i,j)$-sectors (i.e. the rule removes a
letter from $\rk(i,j)$-sectors, and adds a copy of that letter to
each of the $\lk(i,j)$-sectors). It writes a copy of itself in
the $\tk(i)$-sectors and $\bk(i)$-sectors. It locks all
$\kappa(i)$-sectors.

The connecting rule $\bc_1$ has the form

$$\left[\begin{array}{l}\kappa(i)\tool\bbc{1}\kappa(i),
\rk(i,1)\tool\bbc{1}\rk(i,1), \rk(i,j)\tool \rk(i,j),
\\ 1\le i\le N, 2\le j\le n; 1\too 2
\end{array}\right].$$
It writes a copy of itself in the $\tk(i)$-sectors and
$\bk(i)$-sectors. It locks $\rk(i,j)$-sectors and
$\kappa(i)$-sectors.

For every $t=2,...,n$ and every $a\in {\cal A}$ the rule
$\tau(2t-1,a)$ has the form:

$$\left[\begin{array}{l}\kappa(1)\too \bt{2t-1,a}\kappa(1),
\kappa(i')\tool \bt{2t-1,a}\kappa(i'),\\
\rk(i,1)\tool\bt{2t-1,a}\,\rk(i,1), \lk(i,1)\tool \lk(i,1),
\\ \rk(i,s)\tool
\rk(i,s), \lk(i,s)\tool \lk(i,s), \rk(i,j)\too \rk(i,j)a,
\lk(i,j)\too \lk(i,j)a,\\ 1\le i\le N, 2\le i'\le N, 2\le s<t,
t\le j\le n;\quad 2t-1\to 2t-1\end{array}\right].$$ This rule
moves letters from $\rk(i,j)$-sectors to $\lk(i,j)$-sectors,
$j=t,...,n$, and writes a copy of itself in the $\tk(i)$-sectors
and $\bk(i)$-sectors. It locks $\kappa(i)$-sectors, $i\ge 2$, all
$\rk(i,j)$-sectors, $j<t$, and all $\lk(i,j)$-sectors, $j<t$.

The connecting rule $\bc_{2t-1}$ has the form

$$\left[\begin{array}{l}\kappa(1)\too\bbc{2t-1}\kappa(1),
\kappa(i')\tool \bbc{2t-1}\kappa(i'),\\
\rk(i,1)\tool\bbc{2t-1}\rk(i,1), \lk(i,1)\tool\lk(i,1),
\\ \rk(i,j)\tool \rk(i,j), \lk(i,s)\tool \lk(i,s),
\\ 1\le i\le N, 2\le i'\le N, 2\le s<t, 2\le j\le n; 2t-1\to
2t\end{array}\right].$$ It locks $\kappa(i)$-sectors, $i\ge 2$,
all $\rk(i,j)$-sectors, and all $\lk(i,j)$-sectors, $j<t$. It also
writes a copy of itself in the $\tk(i)$-sectors and
$\bk(i)$-sectors.

For every $t=1,...,n$ the rule $\tau(2t,a)$ has the form:

$$\left[\begin{array}{l}  \kappa(1)\too \bt{2t,a}\kappa(1)a\iv,
\rk(i,1)\tool \bt{2t,a}\rk(i,1), \lk(i,1)\tool\lk(i,1),
\\ \kappa(i')\tool\bt{2t,a}\kappa(i'), \rk(i,s)\tool \rk(i,s),
\lk(i,s')\tool \lk(i,s'), \rk(i,j)\too \rk(i,j)a\iv,\\
\lk(i,j')\too \lk(i,j')a\iv,\\ 1\le i\le N, 2\le i'\le N, 2\le
s\le t, t< j\le n, 1\le s'<t, t\le j' \le n; 2t\to
2t\end{array}\right].$$ This rule moves letters from
$\lk(i,t)$-sector to $k(1)$-sector, letters from
$\lk(i,j)$-sectors, $j>t$, to $\rk(i,j)$-sectors, and writes a
copy of itself in the $\tk(i)$-sectors and $\bk(i)$-sectors. It
locks all other sectors.

The connecting rule $\bc_{2t}$, $t=1,...,n-1$, has the form

$$\left[\begin{array}{l}
\kappa(i)\tool \bbc{2t}\kappa(i),
\lk(i,1)\tool\lk(i,1)\\
\rk(i,s)\tool \rk(i,s), \lk(i,j)\tool \lk(i,j),\\ 1\le i\le N,
2\le s\le t, 1\le j\le n;\quad 2t\to 2t+1\end{array}\right].$$ It
writes a copy of itself in the $\tk(i)$-sectors and
$\bk(i)$-sectors. It locks all $\kappa(i)$-sectors, $i>1$, all
$\lk(i,j)$-sectors and all $\rk(i,j)$-sectors, $i\le t$.

The connecting rule $\bc_{2n}$, $t=1,...,n-1$, has the form

$$\left[\begin{array}{l}\kappa(i)\tool\bbc{2n}\kappa(i),
\rk(i,1)\tool\bbc{2n}\rk(i,1), \lk(i,j)\tool\lk(i,j)\\
\rk(i,j')\tool \rk(i,j'),\\ 1\le i\le N, 1\le j\le n, 2\le j'\le
n;\quad  2n\to 2n+1\end{array}\right].$$ It writes a copy of
itself in the $\tk(i)$-sectors and $\bk(i)$-sectors. It locks all
other sectors.

The cleaning rules are in one-to-one correspondence with working
rules. For every working rule $r$, the corresponding cleaning rule
$\tau(2n+1,r)$ has the form

$$\left[\begin{array}{l}\kappa(i)\tool \overline{r}\kappa(i),
\rk(i,1)\tool \overline{r}\rk(i,1), \lk(i,1)\tool\lk(i,1),
\\ \rk(i,j)\tool \rk(i,j), \lk(i,j)\tool\lk(i,j), 1\le i\le N, 2\le
j\le n; 2n+1\to 2n+1\end{array}\right].$$ These rules remove the
content of $\tk(i,j)$-sectors and $\bk(i,j)$-sectors (which
contain copies of the history of the computation) and locks all
other sectors.

Notice that only connecting rules change the $\Omega$-coordinates
of $\bkk$-letters. Notice also that for every $a\in {\cal A}$ and
every $i=1,...,2n$, one can define the rule $\tau(i,a\iv)$ as
before. Similarly for every working rule $r$ one can define
$\tau(2n+1,r\iv)$. It is easy to see that
$$\tau(i,a\iv)=\tau(i,a)\iv.$$

The following picture (\setcounter{pdeleven}{\value{ppp}}Figure
\thepdeleven) shows which rules of the $S$-machine $\sss$ are
locking, active or neither for which sectors of admissible words.

\unitlength=1.00mm \special{em:linewidth 0.4pt}
\linethickness{0.4pt}
\begin{picture}(155.67,165.99)
\put(4.33,15.00){\makebox(0,0)[cc]{$\sss_0$}}
\put(4.33,25.00){\makebox(0,0)[cc]{$\bc_0$}}
\put(4.33,34.67){\makebox(0,0)[cc]{$\sss_1$}}
\put(4.33,44.67){\makebox(0,0)[cc]{$\bc_1$}}
\put(4.33,134.67){\makebox(0,0)[cc]{$\sss_{2n}$}}
\put(4.33,144.67){\makebox(0,0)[cc]{$\bc_{2n}$}}
\put(4.33,154.67){\makebox(0,0)[cc]{$\sss_{2n+1}$}}
\put(4.33,55.00){\makebox(0,0)[cc]{$\sss_2$}}
\put(4.33,65.00){\makebox(0,0)[cc]{$\bc_2$}}
\put(10.00,7.00){\makebox(0,0)[cc]{${\tk(1)}$}}
\put(20.00,7.00){\makebox(0,0)[cc]{${\kappa(1)}$}}
\put(30.00,7.00){\makebox(0,0)[cc]{${\bk(1)}$}}
\put(40.00,7.00){\makebox(0,0)[cc]{$_{\rk(1,1)}$}}
\put(50.00,7.00){\makebox(0,0)[cc]{$_{\lk(1,1)}$}}
\put(60.00,7.00){\makebox(0,0)[cc]{$_{\rk(1,2)}$}}
\put(70.00,7.00){\makebox(0,0)[cc]{$_{\lk(1,2)}$}}
\put(95.00,7.00){\makebox(0,0)[cc]{$_{\rk(1,n)}$}}
\put(105.00,7.00){\makebox(0,0)[cc]{$_{\lk(1,n)}$}}
\put(115.00,7.00){\makebox(0,0)[cc]{${\tk(2)}$}}
\put(125.00,7.00){\makebox(0,0)[cc]{${\kappa(2)}$}}
\put(135.00,7.00){\makebox(0,0)[cc]{${\bk(2)}$}}
\put(10.00,10.33){\circle{1.33}} \put(10.00,10.33){\circle*{1.33}}
\put(20.00,10.33){\circle*{1.33}} \put(30.00,10.33){\circle{1.33}}
\put(30.00,10.33){\circle*{1.33}}
\put(40.00,10.33){\circle*{1.33}}
\put(50.00,10.33){\circle*{1.33}} \put(60.00,10.33){\circle{1.33}}
\put(60.00,10.33){\circle*{1.33}}
\put(70.00,10.33){\circle*{1.33}} \put(95.00,10.33){\circle{1.33}}
\put(95.00,10.33){\circle*{1.33}}
\put(105.00,10.33){\circle*{1.33}}
\put(115.00,10.33){\circle*{1.33}}
\put(125.00,10.33){\circle{1.33}}
\put(125.00,10.33){\circle*{1.33}}
\put(135.00,10.33){\circle*{1.33}}
\put(10.00,20.33){\circle{1.33}} \put(10.00,20.33){\circle*{1.33}}
\put(20.00,20.33){\circle*{1.33}} \put(30.00,20.33){\circle{1.33}}
\put(30.00,20.33){\circle*{1.33}}
\put(40.00,20.33){\circle*{1.33}}
\put(50.00,20.33){\circle*{1.33}} \put(60.00,20.33){\circle{1.33}}
\put(60.00,20.33){\circle*{1.33}}
\put(70.00,20.33){\circle*{1.33}} \put(95.00,20.33){\circle{1.33}}
\put(95.00,20.33){\circle*{1.33}}
\put(105.00,20.33){\circle*{1.33}}
\put(115.00,20.33){\circle*{1.33}}
\put(125.00,20.33){\circle{1.33}}
\put(125.00,20.33){\circle*{1.33}}
\put(135.00,20.33){\circle*{1.33}}
\put(10.00,30.33){\circle{1.33}} \put(10.00,30.33){\circle*{1.33}}
\put(20.00,30.33){\circle*{1.33}} \put(30.00,30.33){\circle{1.33}}
\put(30.00,30.33){\circle*{1.33}}
\put(40.00,30.33){\circle*{1.33}}
\put(50.00,30.33){\circle*{1.33}} \put(60.00,30.33){\circle{1.33}}
\put(60.00,30.33){\circle*{1.33}}
\put(70.00,30.33){\circle*{1.33}} \put(95.00,30.33){\circle{1.33}}
\put(95.00,30.33){\circle*{1.33}}
\put(105.00,30.33){\circle*{1.33}}
\put(115.00,30.33){\circle*{1.33}}
\put(125.00,30.33){\circle{1.33}}
\put(125.00,30.33){\circle*{1.33}}
\put(135.00,30.33){\circle*{1.33}}
\put(10.00,40.33){\circle{1.33}} \put(10.00,40.33){\circle*{1.33}}
\put(20.00,40.33){\circle*{1.33}} \put(30.00,40.33){\circle{1.33}}
\put(30.00,40.33){\circle*{1.33}}
\put(40.00,40.33){\circle*{1.33}}
\put(50.00,40.33){\circle*{1.33}} \put(60.00,40.33){\circle{1.33}}
\put(60.00,40.33){\circle*{1.33}}
\put(70.00,40.33){\circle*{1.33}} \put(95.00,40.33){\circle{1.33}}
\put(95.00,40.33){\circle*{1.33}}
\put(105.00,40.33){\circle*{1.33}}
\put(115.00,40.33){\circle*{1.33}}
\put(125.00,40.33){\circle{1.33}}
\put(125.00,40.33){\circle*{1.33}}
\put(135.00,40.33){\circle*{1.33}}
\put(10.00,50.33){\circle{1.33}} \put(10.00,50.33){\circle*{1.33}}
\put(20.00,50.33){\circle*{1.33}} \put(30.00,50.33){\circle{1.33}}
\put(30.00,50.33){\circle*{1.33}}
\put(40.00,50.33){\circle*{1.33}}
\put(50.00,50.33){\circle*{1.33}} \put(60.00,50.33){\circle{1.33}}
\put(60.00,50.33){\circle*{1.33}}
\put(70.00,50.33){\circle*{1.33}} \put(95.00,50.33){\circle{1.33}}
\put(95.00,50.33){\circle*{1.33}}
\put(105.00,50.33){\circle*{1.33}}
\put(115.00,50.33){\circle*{1.33}}
\put(125.00,50.33){\circle{1.33}}
\put(125.00,50.33){\circle*{1.33}}
\put(135.00,50.33){\circle*{1.33}}
\put(10.00,60.33){\circle{1.33}} \put(10.00,60.33){\circle*{1.33}}
\put(20.00,60.33){\circle*{1.33}} \put(30.00,60.33){\circle{1.33}}
\put(30.00,60.33){\circle*{1.33}}
\put(40.00,60.33){\circle*{1.33}}
\put(50.00,60.33){\circle*{1.33}} \put(60.00,60.33){\circle{1.33}}
\put(60.00,60.33){\circle*{1.33}}
\put(70.00,60.33){\circle*{1.33}} \put(95.00,60.33){\circle{1.33}}
\put(95.00,60.33){\circle*{1.33}}
\put(105.00,60.33){\circle*{1.33}}
\put(115.00,60.33){\circle*{1.33}}
\put(125.00,60.33){\circle{1.33}}
\put(125.00,60.33){\circle*{1.33}}
\put(135.00,60.33){\circle*{1.33}}
\put(10.00,70.33){\circle{1.33}} \put(10.00,70.33){\circle*{1.33}}
\put(20.00,70.33){\circle*{1.33}} \put(30.00,70.33){\circle{1.33}}
\put(30.00,70.33){\circle*{1.33}}
\put(40.00,70.33){\circle*{1.33}}
\put(50.00,70.33){\circle*{1.33}} \put(60.00,70.33){\circle{1.33}}
\put(60.00,70.33){\circle*{1.33}}
\put(70.00,70.33){\circle*{1.33}} \put(95.00,70.33){\circle{1.33}}
\put(95.00,70.33){\circle*{1.33}}
\put(105.00,70.33){\circle*{1.33}}
\put(115.00,70.33){\circle*{1.33}}
\put(125.00,70.33){\circle{1.33}}
\put(125.00,70.33){\circle*{1.33}}
\put(135.00,70.33){\circle*{1.33}}
\put(10.00,130.33){\circle{1.33}}
\put(10.00,130.33){\circle*{1.33}}
\put(20.00,130.33){\circle*{1.33}}
\put(30.00,130.33){\circle{1.33}}
\put(30.00,130.33){\circle*{1.33}}
\put(40.00,130.33){\circle*{1.33}}
\put(50.00,130.33){\circle*{1.33}}
\put(60.00,130.33){\circle{1.33}}
\put(60.00,130.33){\circle*{1.33}}
\put(70.00,130.33){\circle*{1.33}}
\put(95.00,130.33){\circle{1.33}}
\put(95.00,130.33){\circle*{1.33}}
\put(105.00,130.33){\circle*{1.33}}
\put(115.00,130.33){\circle*{1.33}}
\put(125.00,130.33){\circle{1.33}}
\put(125.00,130.33){\circle*{1.33}}
\put(135.00,130.33){\circle*{1.33}}
\put(10.00,140.33){\circle{1.33}}
\put(10.00,140.33){\circle*{1.33}}
\put(20.00,140.33){\circle*{1.33}}
\put(30.00,140.33){\circle{1.33}}
\put(30.00,140.33){\circle*{1.33}}
\put(40.00,140.33){\circle*{1.33}}
\put(50.00,140.33){\circle*{1.33}}
\put(60.00,140.33){\circle{1.33}}
\put(60.00,140.33){\circle*{1.33}}
\put(70.00,140.33){\circle*{1.33}}
\put(95.00,140.33){\circle{1.33}}
\put(95.00,140.33){\circle*{1.33}}
\put(105.00,140.33){\circle*{1.33}}
\put(115.00,140.33){\circle*{1.33}}
\put(125.00,140.33){\circle{1.33}}
\put(125.00,140.33){\circle*{1.33}}
\put(135.00,140.33){\circle*{1.33}}
\put(10.00,150.33){\circle{1.33}}
\put(10.00,150.33){\circle*{1.33}}
\put(20.00,150.33){\circle*{1.33}}
\put(30.00,150.33){\circle{1.33}}
\put(30.00,150.33){\circle*{1.33}}
\put(40.00,150.33){\circle*{1.33}}
\put(50.00,150.33){\circle*{1.33}}
\put(60.00,150.33){\circle{1.33}}
\put(60.00,150.33){\circle*{1.33}}
\put(70.00,150.33){\circle*{1.33}}
\put(95.00,150.33){\circle{1.33}}
\put(95.00,150.33){\circle*{1.33}}
\put(105.00,150.33){\circle*{1.33}}
\put(115.00,150.33){\circle*{1.33}}
\put(125.00,150.33){\circle{1.33}}
\put(125.00,150.33){\circle*{1.33}}
\put(135.00,150.33){\circle*{1.33}}
\put(10.00,160.33){\circle{1.33}}
\put(10.00,160.33){\circle*{1.33}}
\put(20.00,160.33){\circle*{1.33}}
\put(30.00,160.33){\circle{1.33}}
\put(30.00,160.33){\circle*{1.33}}
\put(40.00,160.33){\circle*{1.33}}
\put(50.00,160.33){\circle*{1.33}}
\put(60.00,160.33){\circle{1.33}}
\put(60.00,160.33){\circle*{1.33}}
\put(70.00,160.33){\circle*{1.33}}
\put(95.00,160.33){\circle{1.33}}
\put(95.00,160.33){\circle*{1.33}}
\put(105.00,160.33){\circle*{1.33}}
\put(115.00,160.33){\circle*{1.33}}
\put(125.00,160.33){\circle{1.33}}
\put(125.00,160.33){\circle*{1.33}}
\put(135.00,160.33){\circle*{1.33}}
\put(10.67,10.33){\rule{8.67\unitlength}{20.00\unitlength}}
\put(10.67,30.33){\framebox(8.67,40.33)[cc]{}}
\put(15.00,44.67){\makebox(0,0)[cc]{$<$}}
\put(14.67,35.67){\makebox(0,0)[cc]{$<$}}
\put(15.33,64.00){\makebox(0,0)[cc]{$<$}}
\put(15.00,55.00){\makebox(0,0)[cc]{$<$}}
\put(10.67,130.33){\framebox(8.67,30.00)[cc]{}}
\put(15.00,135.33){\makebox(0,0)[cc]{$<$}}
\put(15.00,145.33){\makebox(0,0)[cc]{$<$}}
\put(15.00,155.33){\makebox(0,0)[cc]{$<$}}
\put(20.67,10.33){\rule{8.67\unitlength}{30.00\unitlength}}
\put(20.67,40.33){\rule{8.67\unitlength}{10.00\unitlength}}
\put(20.67,50.33){\framebox(8.67,20.33)[cc]{}}
\put(24.67,55.33){\makebox(0,0)[cc]{$>$}}
\put(20.67,130.33){\framebox(8.67,10.00)[cc]{}}
\put(25.00,135.67){\makebox(0,0)[cc]{$>$}}
\put(20.67,140.33){\rule{8.67\unitlength}{20.00\unitlength}}
\put(30.67,10.33){\rule{8.67\unitlength}{20.00\unitlength}}
\put(30.67,30.33){\framebox(8.67,40.33)[cc]{}}
\put(34.33,36.00){\makebox(0,0)[cc]{$<$}}
\put(34.33,45.00){\makebox(0,0)[cc]{$<$}}
\put(34.33,55.33){\makebox(0,0)[cc]{$<$}}
\put(34.33,64.33){\makebox(0,0)[cc]{$<$}}
\put(30.67,130.33){\framebox(8.67,30.00)[cc]{}}
\put(35.00,135.33){\makebox(0,0)[cc]{$<$}}
\put(35.00,145.33){\makebox(0,0)[cc]{$<$}}
\put(35.00,155.33){\makebox(0,0)[cc]{$<$}}
\put(41.00,10.33){\framebox(8.33,30.00)[cc]{}}
\put(45.00,16.00){\makebox(0,0)[cc]{$>$}}
\put(45.00,36.00){\makebox(0,0)[cc]{$>$}}
\put(41.00,40.33){\rule{8.33\unitlength}{10.00\unitlength}}
\put(40.67,130.33){\rule{8.67\unitlength}{30.00\unitlength}}
\put(41.00,50.33){\rule{8.33\unitlength}{20.33\unitlength}}
\put(50.67,10.33){\rule{8.67\unitlength}{20.00\unitlength}}
\put(50.67,30.33){\framebox(8.67,30.00)[cc]{}}
\put(55.00,36.00){\makebox(0,0)[cc]{$>$}}
\put(55.00,55.67){\makebox(0,0)[cc]{$>$}}
\put(50.67,60.33){\rule{9.00\unitlength}{10.33\unitlength}}
\put(50.67,130.33){\rule{8.67\unitlength}{30.00\unitlength}}
\put(60.67,10.33){\framebox(8.67,30.00)[cc]{}}
\put(65.00,16.33){\makebox(0,0)[cc]{$>$}}
\put(65.00,36.00){\makebox(0,0)[cc]{$>$}}
\put(60.67,40.33){\rule{8.67\unitlength}{10.00\unitlength}}
\put(60.67,50.33){\framebox(8.67,20.33)[cc]{}}
\put(65.00,56.00){\makebox(0,0)[cc]{$>$}}
\put(60.67,130.33){\rule{8.67\unitlength}{30.00\unitlength}}
\put(95.67,10.33){\framebox(8.67,30.00)[cc]{}}
\put(100.00,16.33){\makebox(0,0)[cc]{$>$}}
\put(100.00,36.00){\makebox(0,0)[cc]{$>$}}
\put(95.67,40.33){\rule{8.67\unitlength}{10.00\unitlength}}
\put(95.67,50.33){\framebox(8.67,20.33)[cc]{}}
\put(100.00,56.00){\makebox(0,0)[cc]{$>$}}
\put(95.67,130.33){\rule{8.67\unitlength}{30.00\unitlength}}
\put(145.00,7.00){\makebox(0,0)[cc]{$_{\rk(2,1)}$}}
\put(145.00,10.33){\circle*{1.33}}
\put(145.00,20.33){\circle*{1.33}}
\put(145.00,30.33){\circle*{1.33}}
\put(145.00,40.33){\circle*{1.33}}
\put(145.00,50.33){\circle*{1.33}}
\put(145.00,60.33){\circle*{1.33}}
\put(145.00,70.33){\circle*{1.33}}
\put(145.00,130.33){\circle*{1.33}}
\put(145.00,140.33){\circle*{1.33}}
\put(145.00,150.33){\circle*{1.33}}
\put(145.00,160.33){\circle*{1.33}}
\put(135.67,10.33){\rule{8.67\unitlength}{20.00\unitlength}}
\put(135.67,30.33){\framebox(8.67,40.33)[cc]{}}
\put(139.33,36.00){\makebox(0,0)[cc]{$<$}}
\put(139.33,45.00){\makebox(0,0)[cc]{$<$}}
\put(139.33,55.33){\makebox(0,0)[cc]{$<$}}
\put(139.33,64.33){\makebox(0,0)[cc]{$<$}}
\put(135.67,130.33){\framebox(8.67,30.00)[cc]{}}
\put(140.00,135.33){\makebox(0,0)[cc]{$<$}}
\put(140.00,145.33){\makebox(0,0)[cc]{$<$}}
\put(140.00,155.33){\makebox(0,0)[cc]{$<$}}
\put(10.67,40.33){\line(1,0){58.67}}
\put(10.67,50.33){\line(1,0){58.67}}
\put(10.67,60.33){\line(1,0){58.67}}
\put(10.67,140.33){\line(1,0){28.67}}
\put(39.33,140.33){\line(0,0){0.00}}
\put(10.67,150.33){\line(1,0){28.67}}
\put(95.67,20.33){\line(1,0){18.67}}
\put(114.33,20.33){\line(0,0){0.00}}
\put(95.67,30.33){\line(1,0){18.67}}
\put(95.67,60.33){\line(1,0){48.67}}
\put(95.67,40.33){\line(1,0){48.67}}
\put(95.67,50.33){\line(1,0){48.67}}
\put(115.67,140.33){\line(1,0){28.67}}
\put(115.67,150.33){\line(1,0){28.67}}
\put(115.67,10.33){\rule{8.67\unitlength}{20.00\unitlength}}
\put(115.67,30.33){\framebox(8.67,40.33)[cc]{}}
\put(120.00,44.67){\makebox(0,0)[cc]{$<$}}
\put(119.67,35.67){\makebox(0,0)[cc]{$<$}}
\put(120.33,64.00){\makebox(0,0)[cc]{$<$}}
\put(120.00,55.00){\makebox(0,0)[cc]{$<$}}
\put(115.67,130.33){\framebox(8.67,30.00)[cc]{}}
\put(120.00,135.33){\makebox(0,0)[cc]{$<$}}
\put(120.00,145.33){\makebox(0,0)[cc]{$<$}}
\put(120.00,155.33){\makebox(0,0)[cc]{$<$}}
\put(105.67,10.33){\rule{8.67\unitlength}{20.00\unitlength}}
\put(105.67,30.33){\framebox(8.67,30.00)[cc]{}}
\put(110.00,36.00){\makebox(0,0)[cc]{$>$}}
\put(110.00,55.67){\makebox(0,0)[cc]{$>$}}
\put(105.67,60.33){\rule{9.00\unitlength}{10.33\unitlength}}
\put(105.67,140.33){\rule{8.67\unitlength}{20.00\unitlength}}
\put(105.67,130.33){\framebox(8.67,10.00)[cc]{}}
\put(110.00,135.67){\makebox(0,0)[cc]{$>$}}
\put(125.67,10.33){\rule{8.67\unitlength}{60.00\unitlength}}
\put(125.67,130.33){\rule{8.67\unitlength}{30.00\unitlength}}
\put(82.00,40.00){\makebox(0,0)[cc]{$\dots$}}
\put(82.00,79.67){\makebox(0,0)[cc]{$\dots$}}
\put(82.00,145.00){\makebox(0,0)[cc]{$\dots$}}
\put(151.00,145.33){\makebox(0,0)[cc]{$\dots$}}
\put(150.67,40.33){\makebox(0,0)[cc]{$\dots$}}
\put(155.67,7.33){\makebox(0,0)[cc]{$\dots$}}
\put(10.67,20.33){\line(1,0){58.67}}
\put(10.67,30.33){\line(1,0){58.67}}
\put(10.00,80.33){\circle{1.33}} \put(10.00,80.33){\circle*{1.33}}
\put(20.00,80.33){\circle*{1.33}} \put(30.00,80.33){\circle{1.33}}
\put(30.00,80.33){\circle*{1.33}}
\put(40.00,80.33){\circle*{1.33}}
\put(50.00,80.33){\circle*{1.33}} \put(60.00,80.33){\circle{1.33}}
\put(60.00,80.33){\circle*{1.33}}
\put(70.00,80.33){\circle*{1.33}} \put(95.00,80.33){\circle{1.33}}
\put(95.00,80.33){\circle*{1.33}}
\put(105.00,80.33){\circle*{1.33}}
\put(115.00,80.33){\circle*{1.33}}
\put(125.00,80.33){\circle{1.33}}
\put(125.00,80.33){\circle*{1.33}}
\put(135.00,80.33){\circle*{1.33}}
\put(145.00,80.33){\circle*{1.33}}
\put(82.00,93.34){\makebox(0,0)[cc]{$\dots$}}
\put(10.00,90.33){\circle{1.33}} \put(10.00,90.33){\circle*{1.33}}
\put(20.00,90.33){\circle*{1.33}} \put(30.00,90.33){\circle{1.33}}
\put(30.00,90.33){\circle*{1.33}}
\put(40.00,90.33){\circle*{1.33}}
\put(50.00,90.33){\circle*{1.33}} \put(60.00,90.33){\circle{1.33}}
\put(60.00,90.33){\circle*{1.33}}
\put(70.00,90.33){\circle*{1.33}} \put(95.00,90.33){\circle{1.33}}
\put(95.00,90.33){\circle*{1.33}}
\put(105.00,90.33){\circle*{1.33}}
\put(115.00,90.33){\circle*{1.33}}
\put(125.00,90.33){\circle{1.33}}
\put(125.00,90.33){\circle*{1.33}}
\put(135.00,90.33){\circle*{1.33}}
\put(145.00,90.33){\circle*{1.33}}
\put(82.00,104.34){\makebox(0,0)[cc]{$\dots$}}
\put(10.00,110.33){\circle{1.33}}
\put(10.00,110.33){\circle*{1.33}}
\put(20.00,110.33){\circle*{1.33}}
\put(30.00,110.33){\circle{1.33}}
\put(30.00,110.33){\circle*{1.33}}
\put(40.00,110.33){\circle*{1.33}}
\put(50.00,110.33){\circle*{1.33}}
\put(60.00,110.33){\circle{1.33}}
\put(60.00,110.33){\circle*{1.33}}
\put(70.00,110.33){\circle*{1.33}}
\put(95.00,110.33){\circle{1.33}}
\put(95.00,110.33){\circle*{1.33}}
\put(105.00,110.33){\circle*{1.33}}
\put(115.00,110.33){\circle*{1.33}}
\put(125.00,110.33){\circle{1.33}}
\put(125.00,110.33){\circle*{1.33}}
\put(135.00,110.33){\circle*{1.33}}
\put(145.00,110.33){\circle*{1.33}}
\put(82.00,119.67){\makebox(0,0)[cc]{$\dots$}}
\put(10.00,120.33){\circle{1.33}}
\put(10.00,120.33){\circle*{1.33}}
\put(20.00,120.33){\circle*{1.33}}
\put(30.00,120.33){\circle{1.33}}
\put(30.00,120.33){\circle*{1.33}}
\put(40.00,120.33){\circle*{1.33}}
\put(50.00,120.33){\circle*{1.33}}
\put(60.00,120.33){\circle{1.33}}
\put(60.00,120.33){\circle*{1.33}}
\put(70.00,120.33){\circle*{1.33}}
\put(95.00,120.33){\circle{1.33}}
\put(95.00,120.33){\circle*{1.33}}
\put(105.00,120.33){\circle*{1.33}}
\put(115.00,120.33){\circle*{1.33}}
\put(125.00,120.33){\circle{1.33}}
\put(125.00,120.33){\circle*{1.33}}
\put(135.00,120.33){\circle*{1.33}}
\put(145.00,120.33){\circle*{1.33}}
\put(10.67,70.67){\framebox(8.67,19.67)[cc]{}}
\put(10.67,110.33){\framebox(8.67,20.00)[cc]{}}
\put(20.67,70.67){\framebox(8.67,19.67)[cc]{}}
\put(30.67,70.67){\framebox(8.67,19.67)[cc]{}}
\put(20.67,110.33){\framebox(8.67,20.00)[cc]{}}
\put(30.67,110.33){\framebox(8.67,20.00)[cc]{}}
\put(125.67,110.00){\rule{8.67\unitlength}{20.33\unitlength}}
\put(125.33,70.67){\rule{9.00\unitlength}{19.67\unitlength}}
\put(115.67,70.67){\framebox(9.00,19.67)[cc]{}}
\put(135.67,70.67){\framebox(8.67,19.67)[cc]{}}
\put(115.67,110.33){\framebox(8.67,20.00)[cc]{}}
\put(135.67,110.33){\framebox(8.67,20.00)[cc]{}}
\put(4.33,75.00){\makebox(0,0)[cc]{$\sss_3$}}
\put(4.33,85.33){\makebox(0,0)[cc]{$\bc_3$}}
\put(4.33,115.33){\makebox(0,0)[cc]{$\sss_{2n-1}$}}
\put(4.33,125.00){\makebox(0,0)[cc]{$\bc_{2n-1}$}}
\put(10.67,70.67){\line(1,0){58.67}}
\put(10.67,80.33){\line(1,0){58.67}}
\put(10.67,90.33){\line(1,0){58.67}}
\put(95.67,80.33){\line(1,0){48.67}}
\put(95.33,90.33){\line(1,0){49.00}}
\put(10.67,110.33){\line(1,0){58.67}}
\put(10.67,120.33){\line(1,0){58.67}}
\put(95.67,110.33){\line(1,0){48.67}}
\put(95.67,120.33){\line(1,0){48.67}}
\put(40.67,70.67){\rule{8.67\unitlength}{19.67\unitlength}}
\put(40.67,110.33){\rule{8.67\unitlength}{20.00\unitlength}}
\put(50.67,70.67){\rule{8.67\unitlength}{19.67\unitlength}}
\put(50.67,110.33){\rule{8.67\unitlength}{20.00\unitlength}}
\put(60.67,110.33){\rule{8.67\unitlength}{20.00\unitlength}}
\put(60.67,70.67){\framebox(8.67,9.67)[cc]{}}
\put(60.67,80.33){\rule{8.67\unitlength}{10.00\unitlength}}
\put(65.00,75.33){\makebox(0,0)[cc]{$>$}}
\put(34.33,75.00){\makebox(0,0)[cc]{$<$}}
\put(34.33,84.00){\makebox(0,0)[cc]{$<$}}
\put(34.33,115.00){\makebox(0,0)[cc]{$<$}}
\put(34.33,124.00){\makebox(0,0)[cc]{$<$}}
\put(24.67,64.33){\makebox(0,0)[cc]{$>$}}
\put(24.67,75.00){\makebox(0,0)[cc]{$>$}}
\put(24.67,84.00){\makebox(0,0)[cc]{$>$}}
\put(15.33,84.00){\makebox(0,0)[cc]{$<$}}
\put(15.00,75.00){\makebox(0,0)[cc]{$<$}}
\put(24.67,115.00){\makebox(0,0)[cc]{$>$}}
\put(24.67,124.00){\makebox(0,0)[cc]{$>$}}
\put(15.33,124.00){\makebox(0,0)[cc]{$<$}}
\put(15.00,115.00){\makebox(0,0)[cc]{$<$}}
\put(139.66,75.00){\makebox(0,0)[cc]{$<$}}
\put(139.66,84.00){\makebox(0,0)[cc]{$<$}}
\put(120.66,83.67){\makebox(0,0)[cc]{$<$}}
\put(120.33,74.67){\makebox(0,0)[cc]{$<$}}
\put(139.66,115.33){\makebox(0,0)[cc]{$<$}}
\put(139.66,124.33){\makebox(0,0)[cc]{$<$}}
\put(120.66,124.00){\makebox(0,0)[cc]{$<$}}
\put(120.33,115.00){\makebox(0,0)[cc]{$<$}}
\put(95.67,70.67){\framebox(8.67,19.67)[cc]{}}
\put(105.67,71.00){\framebox(8.67,19.33)[cc]{}}
\put(95.67,80.33){\rule{8.67\unitlength}{10.00\unitlength}}
\put(100.00,75.67){\makebox(0,0)[cc]{$>$}}
\put(110.00,76.00){\makebox(0,0)[cc]{$>$}}
\put(95.67,110.33){\framebox(8.67,10.00)[cc]{}}
\put(105.67,110.33){\framebox(8.67,10.00)[cc]{}}
\put(100.00,115.00){\makebox(0,0)[cc]{$>$}}
\put(110.00,115.33){\makebox(0,0)[cc]{$>$}}
\put(95.67,120.33){\rule{8.67\unitlength}{10.00\unitlength}}
\put(105.67,120.33){\framebox(8.67,10.00)[cc]{}}
\end{picture}

\begin{center}
\nopagebreak[4] Fig. \theppp.

\end{center}
\addtocounter{ppp}{1}

A filled box means that the sector is locked, a box with ``$>$" in
it means that the rules are left active for that sector, a box
with ``$<$" in it means that the rules are right active, and an
empty box means that the rules from the corresponding set of
$S$-rules are neither active for that sector nor locking it. For
example, this picture shows that rules from $\sss_2$ are right
active for $\tk(i)$-sectors and $\bk(i)$-sectors, lock
$\rk(1,1)$-sectors but are left active for $\rk(1,2)$-sectors. It
also shows that $\kappa(i)$-sectors, $i\ge 2$, are locked by all
rules of $\sss$.

\subsection{The presentation}
\label{thepresentation}

The set of relations corresponding to the rules from $\sss$ will
be denoted by \label{zss}$Z(\sss)$. By
\label{zsl}$Z(\sss,\Lambda)$ we shall denote the set $Z(\sss)$
together with the hub $\Lambda(0)$.

Our main result is that the group \label{ash1} ${\cal H}$ defined
by the set of relations $Z(\sss,\Lambda)$ satisfies the two
conditions of Theorem \ref{CEP}.

\subsection{The subgroup $\la {\cal A}\ra$ of ${\cal H}$ is a torsion group of exponent $n$}

For every group word $u$ in the alphabet ${\cal A}$ and every
$\bkk$-letter $z\not\in\{\tk(i),\bk(i,1), i=1,...,N\}$ we denote a
copy of $u$ in the alphabet $A(z)$ by $u(z)$. We may identify
${\cal A}$ with $A(\kappa(1))$.

\begin{lm} \label{lmbk}
For every word $u$ in ${\cal A}\cup {\cal A}\iv$ the relation
$u(\kappa(1))^n=1$ follows from $Z(\sss,\Lambda)$.
\end{lm}

\proof Let $u\equiv a_{i_1}a_{i_2}...a_{i_s}$ be a group word in
the alphabet ${\cal A}$. Start with the admissible word
$\Lambda(0)\tk(1,0)$ and apply $S$-rules $\tau(0,a_{i_1})$,...,
$\tau(0,a_{i_s})$ from $\sss_0$ to $\Lambda(0)\tk(1,0)$. The
resulting word is equal to

\begin{equation}
\begin{array}{l}
\tk(1)\kappa(1)\bk(1)\rk(1,1)u\iv\lk(1,1)\rk(1,2)
u\iv\lk(1,2)...\tk(N)\kappa(N)\bk(N)\\
\rk(N,n)u\iv\lk(N,n) \tk(1)
\end{array}
\end{equation}
the $\Omega$-coordinate of each $\bkk$-letter is $0$.

Applying now the connecting rule $\bc_0$, we change the
$\Omega$-coordinates of all $\bkk$-letters to $1$.

After that, applying rules $\tau(1,a_{i_s})$,
...,$\tau(1,a_{i_1})$, we get the word
\begin{equation}
\begin{array}{l}
\tk(1)v_1\kappa(1)\bk(1)v_1\rk(1,1)\lk(1,1)u\rk(1,2)...\tk(N)v_1
\kappa(N)\bk(N)v_1\rk(N,1)...\\ \rk(N,n)\lk(N,n)u\tk(1)
\end{array}
\end{equation}
(the $\Omega$-coordinates are $1$) where $v_1\equiv
\overline{\tau(1,a_{i_s})}...\overline{\tau(1,a_{i_1})}$.

The connecting  rule $\bc_1$ makes all $\Omega$-coordinates $2$
and adds letter $\bbc{1}$ at the end of $v$.

Now we can apply rules $\tau(2,a_{i_1})$, ..., $\tau(2,a_{i_s})$,
$\bc_2$, $\tau(3,a_{i_s})$,..., $\tau(3,a_{i_1})$, $\bc_3$,...,
$\tau(2n, a_{i_1})$, ..., $\tau(2n,a_{i_s})$, and obtain the
following word:

\begin{equation}\label{eqsecond}
\begin{array}{l}
\tk(1)
v_{2n}\kappa(1)u^{-n}\bk(1)v_{2n}\rk(1,1)\lk(1,1)...\tk(N)v_{2n}\kappa(N)\bk(N)v_{2n}\rk(N,1)\lk(N,1)\\...
\rk(N,n)\lk(N,n)\tk(1)
\end{array},
\end{equation}
where the $\Omega$-coordinate of every $\bkk$-letter is $2n$,
$v_{2n}$ is a copy of the history of the computation.

Let $W$ be the word (\ref{eqsecond}) without the last letter.

Let $h\equiv \tau(0,a_{i_1})...\tau(2n,a_{i_s})$ be the string of
all rules from $\sss$ applied so far to $\Lambda(0)\tk(1)$. Then
by Lemma \ref{lm1} (part 3), we have the following equality:

$$L(h,\tk(1))W\tk(1,2n)=\Lambda(0)\tk(1,0)R(h,\tk(1))=
\Lambda(0)L(h,\tk(1))\tk(1,2n).$$ So
$$W=L(h,\tk(1))\iv\Lambda(0)L(h,\tk(1))=1$$ modulo
$Z(\sss,\Lambda)$.

Notice that if we apply connecting rules $\bc_0,...,\bc_{2n}$ to
$\Lambda(0)\tk(1,0)$, we obtain $$\Lambda(2n+1)\tk(1,2n+1).$$ Now
take the word $\Lambda(2n+1)\tk(1,2n+1)$ and apply the cleaning
rules $\tau(2n+1,t_1)...\tau(2n+1,t_s)$ where
$\overline{t_1}...\overline{t_s}\equiv v_{2n}$, and the rule
$\bc_{2n}\iv$. We get the word

\begin{equation}\label{eqthird}
\begin{array}{l}
\tk(1)
v_{2n}\kappa(1)\bk(1)v_{2n}\rk(1,1)\lk(1,1)...\tk(N)v_{2n}\kappa(N)\bk(N)v_{2n}\rk(N,1)\lk(N,1)\\...
\rk(N,n)\lk(N,n)\tk(1)
\end{array},
\end{equation}
with $\Omega$-coordinate $2n$. Let $W'$ be the word
(\ref{eqthird}) without the last letter. Let $h'$ be the word of
rules applied to $\Lambda(0)\tk(1,0)$ to get (\ref{eqthird}). Then
by Lemma \ref{lm1}, we have
$$L(h',\tk(1))W'\tk(1,2n)=\Lambda(0)\tk(1,0)R(h',\tk(1))=\Lambda(0)
L(h',\tk(1))\tk(1,2n)$$ modulo $Z(\sss,\Lambda)$. This implies
$W'=1$ modulo $Z(\sss,\Lambda)$.
 But the words
$W$ and $W'$ differ only by the factor $u(\kappa(1))^{-n}$. Thus
$u(\kappa(1))^n=1$ modulo $Z(\sss,\Lambda)$. Hence the $n$-th
power of every word over $A(\kappa(1))$ is 1 modulo
$Z(\sss,\Lambda)$.
\endproof

\section{Properties of the $S$-machine $\sss$}
\label{properties}

\subsection{The inverse semigroup $P(\sss$)}

We fix an arbitrary set \label{rkappa1} $R(\kappa(1))$ of words in
the alphabet $A(\kappa(1))={\cal A}.$ (We choose
$R(\kappa(1))=\emptyset$ for the construction of the group ${\cal
H}$, but later we should consider arbitrary $R(\kappa(1))$ to
prove part (2) of Theorem \ref{CEP}.)

 Let  \label{ha}$\maa$ be the group generated by the set $\aaa$
subject all relations of the form $u^n=1$ ($u$ is any word in
$\aaa$) and relations $r=1$ for all $r\in R(\kappa(1))$. The group
$\maa$ is isomorphic to the free product in the variety of groups
of exponent $n$ of subgroups $\la A(z)\ra$, $z\in \kk$, where all
$\la A(z)\ra$, $z\ne\kappa(1)$, are free Burnside groups of
exponent $n$, and $\la A(\kappa(1))\ra $ is a group of exponent
$n$ with additional relators from $R(\kappa(1)).$

Let \label{hka}$\hnka$ be the free product of the group $\maa$ and
the free group freely generated by the set $\bkk$.

>From now on we will not distinguish admissible words which are
equal in $\hnka$. Thus a word of the form $y_1u_1...y_t$ where
$y_i\in \bkk$, $u_i$ are words in $A(z_i)$, is \label{adm1.5}{\em
admissible} if for every $i=1,...,t-1$, either $y_{i+1}\equiv
(y_i)_+$ or $y_{i+1}\equiv y_i\iv$.

If no sector of an admissible word $W$ is equal in $\hnka$ to a
word without $\bkk$-letters then $W$ is called a \label{radm}{\em
reduced admissible word}. By Lemma \ref{reduced} this property of
an admissible word is equivalent to being $\bkk$-reduced in
$\hnka$.

The result of an application of an $S$-rule to an admissible word
is again an admissible word.

Let $W'$ be a word obtained from $W$ by an application of a rule
$\tau$ from $\sss$. Then the $\bkk$-letters in $W$ differ from the
corresponding $\bkk$-letters in $W'$ by their $\Omega$-coordinates
only.

Therefore we have the following simple fact.

\begin{lm}\label{p1}
 If $W'$ is obtained from $W$ by a nonempty sequence of
applications of $S$-rules from $\sss$ then the word $W'$ is also
admissible and the $\bkk$-letters of $W'$ differ from the
corresponding $\bkk$-letters in $W$ only by the
$\Omega$-coordinate. The $\Omega$-coordinate of $W'$ is $\omega$
if and only if the last rule applied was from $\sss_\omega$ or
equal to $\bc_{\omega-1}$, $\bc_{\omega}\iv$.
\end{lm}

Thus every $\tau\in\sss$ defines a partial one-to-one
transformation of the set of admissible words. Let
\label{pss}$P(\sss)$ be the subsemigroup  generated by these
partial transformations in the inverse semigroup of all partial
one-to-one transformations of the set of admissible words. Then
$P(\sss)$ is an inverse semigroup because for every rule $\tau$ in
$\sss$, the inverse rule $\tau\iv$ is also in $\sss$. The zero of
$P(\sss)$ is a transformation with empty domain. Removing zero
from $P(\sss)$ we get what is usually called a {\em pseudogroup}
of partial transformations. Every word in $\sss$ induces a partial
one-to-one transformation from $P(\sss)$.

The following obvious lemma is a general property of inverse
semigroups \cite{Cliff-Prest}.

\begin{lm} \label{pseudo}
Let a word $h$ over $\sss$ be graphically equal to $h_1aa\iv h_2$.
Then, considered as a transformation of the set of admissible
words, the \label{domain} domain of $h$ is contained in the domain
of $h_1h_2$, but the partial transformations induced by $h$ and
$h_1h_2$ coincide on the domain of $h$. Thus for every word $h$
over $\sss$ the domain of the transformation induced by $h$ is
contained in the domain of transformation induced by the reduced
form of $h$.
\end{lm}

Notice that the domain of the reduced form of $h$ may be strictly
bigger than the domain of $h$.

If $h$ is a word in $\sss$ then denote by \label{wh}$W\cdot h$
($W$ is an admissible word), the result of the application of the
freely reduced form of $h$ to $W$ (if it exists). Thus the
equality $W\cdot h=W'$ means that $W$ is in the domain of the
freely reduced form of $h$ and the transformation induced by the
reduced form of $h$ takes $W$ to $W'$.

The following easy but important lemma immediately follows from
the definition of application of an $S$-rule to an admissible
word.

\begin{lm}\label{p2}
Let $W$ be an admissible word. Let $h\equiv \tau_1...\tau_t$ be a
word in $\sss$. Then the following conditions are equivalent:

(i) $W$ belongs to the domain of $h$;

(ii) Every sector of $W$ belongs to the domain of $h$.

Moreover $W'=W\cdot h$ if and only if for every $z$-sector $W_1$
of $W$ and $z$-sector $W_1'$ of $W'$ we have $W_1'=W_1\cdot h$,
$z\in \kk\cup\kk\iv$.
\end{lm}

For every $\tau\in\sss$ and $z\in\kk\cup\kk\iv$ we have defined
two words $L(\tau,z)$ and $R(\tau,z)$. These words belong to the
free group generated by $\aaa\cup\rr$. Let
\label{lara}$L^a(\tau,z)$ and $R^a(\tau,z)$ be the
$\aaa$-projections of $L(\tau,z)$ and $R(\tau,z)$. For every word
$h$ in the alphabet $\sss$, we define the words $L(h,z)$,
$R(h,z)$, $L^a(h,z)$ and $R^a(h,z)$ in the natural way ($L(h,z)$
and $R(h,z)$ have been already defined before). Clearly $L^a(.,z)$
and $R^a(.,z)$ are homomorphisms from the free semigroup
generated by $\sss$ into the subgroups $\la A(z_-)\ra$ and $\la
A(z)\ra$ respectively.

We say that a word $w$ is a \label{copy}{\em copy} of a word $w'$
if $w'$ is obtained from $w$ by replacing equal letters for equal
letters and different letters for different letters. For example,
$xyx$ is a copy of $aba$.

The following lemma immediately follows from the definition of the
words $L(\tau,z)$ and $R(\tau,z)$ and the definition of an
application of an $S$-rule.

\begin{lm}\label{landr}
(i) Let $h$ be a word over $\sss$ applicable (in $P(\sss)$) to a
one-sector admissible word $z_1wz_2$ then $z_1wz_2\cdot h
=z_1'(R^a(h,z_1))\iv wL^a(h,z_2)z_2'$ where $z_i'$ differs from
$z_i$ only by the $\Omega$-coordinate, $i=1,2$.

(ii) For every $\omega\in\Omega$ and every $z\in\kk\cup\kk\iv$
either $L^a(\tau,z)$ (resp. $R^a(\tau,z)$) is empty for every
$\tau\in\sss_\omega$ or $L^a(\tau,z)$ (resp. $R^a(\tau,z)$) is not
empty for every $\tau\in\sss_\omega$. If for some $\tau\in\sss$
and some $z\in \kk\cup\kk\iv$, $L^a(\tau,z)$ is not empty (resp.
$R^a(\tau,z_+)$ is not empty) then for every $\tau\in\sss$,
$L^a(\tau,z_+)$ is empty (resp. $R^a(\tau,z)$ is empty).

(iii) Suppose that either none of the rules appearing in a word
$h$ over $\sss$ are from $\sss_{2n+1}$ or all of them are from
$\sss_{2n+1}$. Then if all rules of a word $h$ over $\sss$ are
left active for $z$-sectors for some $z\in \kk$ then $R^a(h,z)$ is
a copy of $h$ written in the alphabet $A(z)\cup A(z)\iv$. If all
rules of a word $h$ are right active for $z$-sectors for some
$z\in\kk$ then $L^a(h,z_+)$ is a copy of $h$ written in the
alphabet $A(z)$.
\end{lm}

\begin{lm} \label{burns}
(i) If $h\equiv h_1^s$ is a reduced word in $\sss$, $s> 1$, $W$
belongs to the domain of $h$ then $W\cdot h$ and $W$ have the same
$\Omega$-coordinates.

(ii) Let $h_1$ and $h_2$ be two words in $\sss$ which are equal in
the free group modulo Burnside relations. Suppose that $W$ is in
the domain of both $h_1$ and $h_2$. Then $W\cdot h_1=W\cdot h_2$.

(iii) Let $h$ be a word in $\sss$ containing subword $g$ which is
equal to 1 modulo Burnside relations: $h=h_1gh_2$. Suppose $W$ is
in the domain of $h$. Then $W$ is in the domain of $h_1h_2$ and
$W\cdot h_1h_2=W\cdot h$.

(iv) If $W$ belongs to the domain of $h^2$, where $h$ is a
cyclically reduced word, then it belongs to the domain of $h^s$
for every integer $s$ (in particular, $h^s$ has a non-empty
domain).
\end{lm}

\proof (i) Recall that non-connecting rules do not change
$\Omega$-coordinates and each connecting rule $\bc_\omega$
replaces $\omega$ by $\omega+1$, $\omega\in \Omega$. Hence for
every word $g$ over $\sss$, there exists $\omega$ in $\Omega$ such
that the domain of $g$ consists of admissible words with
$\Omega$-coordinate $\omega$.  Since $W\cdot h_1$ is in the domain
of $h_1$, the $\Omega$-coordinate of $W\cdot h_1$ is the same as
the $\Omega$-coordinate of $W$. Hence for every $\omega\in\Omega$,
the product of occurrences of $\bc_\omega^{\pm 1}$ in $h_1$ is 1.
The same is true for any power $h=h_1^s$, $s\in \mathbb{Z}$, of
$h$. Therefore for $W\cdot h$ has the same $\Omega$-coordinate as
$W$.

(ii) Let $W$ be in the domain of $h_1, h_2$. Suppose that $h_1$ is
equal to $h_2$ modulo the Burnside relations in the free group.
Clearly it is enough to assume that $W$ has only one sector:
$W\equiv z_1wz_2$. Then by Lemma \ref{landr}(i) $W\cdot
h_i=z_{1,i}(R^a(h_i,z_1))\iv wL^a(h_i,z_2)z_{2,i}$, $i=1,2$. Since
$R^a(.,z_1)$ and $L^a(.,z_2)$ are homomorphisms and $h_1$ and
$h_2$ are equal modulo Burnside relations,
$R^a(h_1,z_1)=R^a(h_2,z_1)$ and $L^a(h_1,z_2)=L^a(h_2,z_2)$ modulo
Burnside relations. But all Burnside relations hold in
$\la\aaa\ra$, so $(R^a(h_1,z_1))\iv
wL^a(h_1,z_2)=(R^a(h_2,z_1))\iv wL^a(h_2,z_2)$.

It remains to show that the $\Omega$-coordinates of $W\cdot h_1$
and $W\cdot h_2$ are the same. These $\Omega$-coordinates are
determined by the projections of $h_1$ and $h_2$ on the alphabet
$\{\bc_1,...,\bc_{2n+1}\}$ (because other rules do not change the
$\Omega$-coordinate). Each connecting rule $\bc_\omega$ induces a
partial transformation on $\Omega$ of the form $\omega\to
\omega+1$. Therefore every product of $\bc_\omega$,
$\omega\in\Omega$ which has a non-empty domain considered as a
transformation of $\Omega$, must be freely equal to a word of the
form $(\bc_i\bc_{i+1}...\bc_s)^{\pm 1}$ for some $i\le s$.

Since $h_1$ and $h_2$ are equal modulo Burnside relations, their
projections are also equal modulo Burnside relations. But words of
the form $(\bc_i\bc_{i+1}...\bc_s)^{\pm 1}$ can be equal modulo
Burnside relations only if they are graphically equal. Thus
projections of $h_1$ and $h_2$ onto the alphabet
$\{\bc_1,...,\bc_{2n+1}\}$ are freely equal. It remains to observe
that if two words in $\bc_\omega$, $\omega\in \Omega$, are freely
equal and induce non-empty transformations of $\Omega$, then they
induce the same transformation of $\Omega$. Hence $W\cdot h$ and
$W\cdot h'$ have the same $\Omega$-coordinates.

(iii) Indeed by part (ii) $W\cdot h_1$ is equal to $W\cdot h_1g$
and is in the domain of $h_2$. Therefore $W$ is in the domain of
$h_1h_2$ and $W\cdot h_1h_2=W\cdot h$.

(iv) Let $W$ belong to the domain of $h^2$ where $h$ is cyclically
reduced. It is enough to consider the case when $W$ consists of
one sector: $W\equiv z_1wz_2$.

Suppose first that $h$ does not contain rules which lock
$z_1$-sectors. By Lemma \ref{landr}(i) $W\cdot
h^2=z_1'(R^a(h^2,z_1))\iv wL^a(h^2,z_2)z_2'$ where $z_1'$ and
$z_2'$ differ from $z_1$ and $z_2$ only by the
$\Omega$-coordinate. By part (i), $z_1'\equiv z_1$ and $z_2'\equiv
z_2$. So $W\cdot h^2$ is in the domain of $h$, because $W\cdot
h^2$ has the same $\Omega$-coordinate as $W$, and $h$ does not
contain rules which lock $z_1$-sectors. Thus $W$ is in the domain
of $h^s$ for all positive integers $s$. By part (ii), $W\cdot
h^n=W$. Therefore $W$ is in the domain of $h^{-2}$ ($n>2$). This
implies as before that $W$ is in the domain of $h^s$ for all
negative integers $s$.

Now suppose that $h$ contains a rule $\tau$ which locks
$z_1$-sectors, $h\equiv h_1\tau h_2$ for some words $h_1,h_2$ in
$\sss$. Then $W\cdot h_1$ has the form $z_1'z_2'$ where $z_1',
z_2'$ differ from $z_1,z_2$ only by the $\Omega$-coordinate,
$h^2\equiv h_1\tau h_2h_1\tau h_2$. Then $W\cdot h_1$ and $W\cdot
h_1\tau h_2h_1$ have the same $\Omega$-coordinates because both
words are in the domain of $\tau$. Since $\tau$ locks
$z_1$-sectors, $W\cdot h_1$ and $W\cdot h_1\tau h_2h_1$ are both
equal to the word of the form $z_1'z_2'$ in $\hnka$. Since $W\cdot
h_1$ is in the domain of $\tau h_2h_1\tau h_2$, $W\cdot h_1\tau
h_2h_1$ must be in the domain of $\tau h_2h_1\tau h_2$ too. Thus
$W$ is in the domain of $h_1\tau h_2h_1\tau h_2h_1\tau h_2\equiv
h^3$. Considering consequently $W\cdot h, W\cdot h^2,...$ instead
of $W$, we conclude that $W$ is in the domain of $h^s$ for every
positive integer $s$. Using part (ii), as before, we deduce the
same fact for negative $s$.
\endproof

\begin{lm} \label{intermediate}
For every admissible word $U$ with $\Omega$-coordinate $\omega$,
and every two non-empty words $h_1,h_2$ over $\sss_\omega$ if $W$
is in the domains of $h_1$, $h_1=h_2$ modulo Burnside relations
then $U\cdot h_1=U\cdot h_2$.
\end{lm}

\proof By Lemma \ref{p2}, we can assume that $U\equiv zuz'$ is a
one-sector admissible word. Recall that either all rules from
$\sss_\omega$ are left (resp. right) active for $z$-sectors or all
these rules lock $z$-sectors or none of these rules are active or
locking for $z$-sectors. Hence we can conclude that since $U$ is
in the domain for $h_1$, and $h_1$ is not empty, it is also in the
domain for $h_2$. The equality $W\cdot h_1=W\cdot h_2$ follows
from Lemma \ref{burns} (ii).\endproof

Words in the alphabet $\sss$, that do not have subwords which are
equal to 1 modulo Burnside relations, will be called
\label{brw}{\em Burnside-reduced}.

\begin{lm} \label{lm:]tap}
Let $W\equiv z_1wz_2$ be an admissible word consisting of one
sector, $z_2\nee z_1\iv$. Let $h$ be a word in $\sss$. Suppose
that $R^a(\tau,z_1)$ (resp. $L^a(\tau,z_2)$) is non-empty for all
$\tau$ in $h$. Suppose that either none of the rules in $h$ is
from $\sss_{2n+1}$ or all of them are from $\sss_{2n+1}$. Suppose
also that $W\cdot h=W$. Finally suppose that $z_1\ne \kappa(1)$,
$z_1\ne \bk(1)\iv$. Then $h$ is equal to 1 modulo Burnside
relations.
\end{lm}

\proof  We only consider the case when $R^a(\tau,z_1)$ is not
empty because the $L^a$-case is similar. Recall that by Lemma
\ref{landr} (ii) then $L^a(\tau,z_1)$ is empty for every $\tau$
from $h$. By Lemma \ref{landr}(i), $W\cdot h$ is graphically equal
to $z_1R^a(h,z_1)\iv wz_2$ and by Lemma \ref{landr} (iii)
$R^a(h,z_1)$ is a copy of the word $h$ written in the alphabet
$A(z_1)\cup A(z_1)\iv$. Since $W\cdot h=W$, $R^a(h,z_1)=1$ in the
subgroup $\la A(z_1)\ra$ of $\hnka$. This subgroup is a free
Burnside group of exponent $n$ (freely generated by $A(z_1)$).
Therefore $h$ is equal to 1 modulo Burnside relations.
\endproof

The presentation of $\hnka$ satisfies conditions (Z1) if we
consider letters from $\aaa$ as $0$-letters and let $Y=\bkk$, so
we can apply Lemmas from Section \ref{corollaries}. In particular,
by Lemma \ref{reduced}, an admissible word $W$ is $\bkk$-reduced
if and only if it does not have sectors $zuz\iv$ where $u=1$
modulo relations of $\maa$.

\begin{lm} \label{p0}
Let $W$ be a $\bkk$-reduced admissible word which is in the domain
of some word $h$ over $\sss$. Then $W\cdot h$ is also
$\bkk$-reduced.
\end{lm}

\proof Indeed, one only needs to check that no sector in $W\cdot
h$ is equal to 1 in $\hnka$. Indeed, it is trivial for sectors of
the form $zwz_+$ because the application of $\tau$ do not change
the $\kk$-projections of $\bkk$-letters and the free group
$\la\bkk\ra$ is a retract of $\hnka$. If a sector of $W$ has the
form $zwz\iv$ then the corresponding sector in $W\cdot h$ has the
form $zuwu\iv z\iv$ for some $u$ by Lemma \ref{landr}. This sector
is equal to 1 only if $uwu\iv=1$. But then $w=1$ and the sector
$zwz\iv$ is also equal to 1 in $\hnka$.
\endproof

{\em From now on we consider only $\bkk$-reduced
\label{adm2}admissible words.}

\begin{lm} \label{lm:badsector}
Let $W\equiv zwz\iv$ be an admissible word consisting of one
sector. Let $h$ be a reduced word in $\sss$ and $W$ is in the
domain of $h$. Then $h$ does not contain rules which lock
$z$-sectors.
\end{lm}

\proof Indeed, suppose $h\equiv h_1\tau h_2$ for some
$\tau\in\sss$ which locks $z$-sectors. Then by Lemma \ref{landr}
$$W\cdot h_1=z'R^a(h_1,z)w(R^a(h_1,z))\iv (z')\iv=z'(z')\iv$$ because
$W\cdot h_1$ is in the domain of $\tau$. This implies that $w=1$
and $W=1$ which contradicts the assumption that $W$ is reduced in
$\hnka$.
\endproof

\begin{lm} \label{lm:badsector1}
Let $W\equiv zwz\iv$ be an admissible word consisting of one
sector. Let $h$ be a word in $\sss$ and all rules in $h$ are left
active for $z$-sectors. Suppose that $W\cdot h=W$. Then

(i) $R^a(h,z)$ commutes with $w$

(ii) if in addition $z\nee \kappa(1)$ and $h$ either does not
contain rules from $\sss_{2n+1}$ or contains only rules from
$\sss_{2n+1}$ then $h$ is equal modulo Burnside relations to the
power of a certain word $\rt_\omega(w)$ depending only on $w$ and
the $\Omega$-coordinate of $z$. Moreover $W\cdot \rt_\omega(w)=W$.
\end{lm}

\proof The first statement immediately follows from Lemma
\ref{landr}(i). To prove (ii), recall that by Lemma
\ref{landr}(iii) $R^a(h,z)$ is a copy of $h$ written in the
alphabet $A(z)\iv$ and the subgroup $\la A(z)\ra$ of $\hnka$ is
the free Burnside group of exponent $n$. By Lemma \ref{VI} then
$R^a(h,z)$ and $w$ belongs to the same maximal cyclic subgroup
$\la w_0\ra$ of $\la A(z)\ra$ where $w_0$ depends only on $w$.
Therefore $h$ is equal modulo Burnside relations to a power of a
copy of $w_0$ written in the alphabet $\sss$. Denote this copy by
$\rt_\omega(w)$. Then all rules in $\rt_\omega(w)$ are left active
for $z$-sectors, and so $W$ is in the domain of $\rt_\omega(w)$.
Since $R^a(\rt_\omega(w),z)=w_0$ commutes with $w$, we have
$W\cdot \rt_\omega(w)=W$.
\endproof

\begin{lm}\label{lm:iji}
Let $h$ be a reduced word in $\sss$ containing a subword of the
form $\bc_\omega^{\pm 1} h'\bc_\omega^{\mp 1}$ where $h'$ is a
word in $\sss_{\omega}$ or $\sss_{\omega+1}$. Let $W\equiv
z_1wz_2$ be an admissible word consisting of one sector, $z_1\nee
\kappa(1), \bk(1)\iv$. Suppose that $\bc_\omega$ locks
$z_1$-sectors, but rules from $h'$ are (left or right) active for
$z_1$-sectors. Finally suppose that $W$ is in the domain of $h$.
Then $h'$ is equal to 1 modulo Burnside relations.
\end{lm}

\proof Indeed, since $\bc_\omega$ locks $z_1$-sectors, by Lemma
\ref{lm:badsector} $z_2\nee z_1\iv$. Thus $z_2\equiv (z_1)_+$ by
the definition of admissible words. From the definition of $\sss$,
it follows that then either $R^a(\tau,z_1)$ is empty for all
$\tau\in\sss$ or $L^a(\tau,z_2)$ is empty for all $\tau\in \sss$.
Without loss of generality assume that $L^a(\tau,z_2)$ is empty
for all $\tau$ from  $h'$. Then $R^a(\tau,z_1)$ are not empty for
$\tau$ from $h'$ and by Lemma \ref{landr} (iii) $R^a(h',z_1)$ is a
copy of $h'$ written in the alphabet $A(z_1)\cup A(z_1)\iv$. By
Lemma \ref{landr} (i), $W\cdot h=z_1R^a(h',z_1)z_2=z_1z_2$ because
$W\cdot h$ is in the domain of $\bc_\omega\iv$ which locks
$z_1$-sectors. Therefore $R^a(h',z_1)=1$, thus $h'$ is 1 modulo
Burnside relations.
\endproof

The argument in the proof of Lemma \ref{lm:iji} can be easily
generalized to prove the following statement.

\begin{lm}\label{lm:iji1}
Let $h$ be a reduced word in $\sss$ of the form $h_1h_2$. Let
$W\equiv z_1wz_2$ be an admissible word consisting of one sector,
$z_2\equiv (z_1)_+$, $z_1\nee \kappa(1), \bk(1)\iv$. Suppose that
$W$ is in the domain of $h$ and $W\cdot h_1=W\cdot h_1h_2$ and
that rules from $h_2$ are active for $z_1$-sectors. Then $h_2=1$
modulo Burnside relations.
\end{lm}

\begin{lm} \label{no2n+1}
Let $W\equiv zuz'$ be a one-sector admissible word and $h,g$ be
words over $\sss$, $W\cdot h=W$, $W\cdot ghg\iv=W$. Suppose that
at least one of the following three possibilities holds for all
rules $\tau$ of $h,g$ simultaneously:

(i) $R^a(\tau,z)$ is not empty;

(ii) $L^a(\tau,z')$ is not empty;

(iii) $L^a(\tau, z)$ and $R^a(\tau,z)$ are empty.

Suppose also that either all rules from $h, g$ are not from
$\sss_{2n+1}$ or all of them are from $\sss_{2n+1}$. Finally
suppose that $z\nee \kappa(1)$, $z\nee \bk(1)\iv$. Then either $g,
h$ commute modulo Burnside relations or $W\cdot g=W$.
\end{lm}

\proof Consider two possible cases.


{\bf Case 1.} $W\equiv zuz_+$. Since  $z\nee \kappa(1),
\bk(1)\iv$, by Lemma \ref{landr}(ii) possibilities (i) and (ii)
from the formulation of the lemma cannot both hold.  By
assumption, we have the following possibilities:

 {\bf Case 1.1.} For every rule $\tau$ from
$h,g$, $R^a(\tau,z)$ is not empty but $L^a(\tau,z_+)$ is empty.

{\bf Case 1.2.} For every rule $\tau$ from $g,h$, $L^a(\tau,z_+)$
is not empty but $R^a(\tau,z)$ is empty.

{\bf Case 1.3.} For every rule $\tau$ from $g,h$ $R^a(\tau,z)$ and
$L^a(\tau,z_+)$ are empty.

In Cases 1.1, 1.2, by Lemma \ref{lm:]tap}, $h=1$ modulo Burnside
relations, so $g, h$ commute modulo Burnside relations, as
required.

In Case 1.3 by Lemma \ref{landr}(i) for every $\tau$ in $g$, we
have $W\cdot \tau=W$, and so $W\cdot g=W$, as required.

 {\bf Case 2.} $W\equiv zuz\iv$  (again $z\nee \kappa(1), \bk(1)\iv$ by assumption).

\medskip

{\bf Case 2.1.} Suppose first that for every rule from $g,h$,
$R^a(\tau,z)$ is not empty.

In this case by Lemma \ref{lm:badsector1} $h$ commutes with
$ghg\iv$ modulo Burnside relations which implies by Lemma
\ref{additional} that $h$ commutes with $g$ modulo Burnside
relations.

{\bf Case 2.2.} Suppose now that every rule $\tau$ of $g,h$,
$R^a(\tau,z)$ (and hence $L^a(\tau,z\iv)$) is empty.

In this case again by Lemma \ref{landr}(i), for every
$\tau\in\sss_\omega$ we have $W\cdot \tau=W$, and so $W\cdot g=W$,
as required.
\endproof

\subsection{Stabilizers in $P(\sss)$ of subwords of the hub}
For every $i=1,...,N$ let \label{lami}$\Lambda_i$ be the subword
of $\Lambda$ starting with $\tk(i)$ and ending with $\tk(i+1)$
(``$N+1$" here is 1).

Suppose $i>1$, and that $\Lambda_i(0)$ is in the domain of a word
$h$ over $\sss$. Let $U_i(\omega)\equiv \Lambda_i(0)\cdot h$, $U_i
\equiv \tk(i)u_1\kappa(i)\bk(i)u_2...\lk(i,n)u_{2n+2}\tk(i+1)$ for
some words $u_1,...,u_{2n+2}$ (notice that the word between
$\kappa(i)$ and $\bk(i)$ is empty, $i=2,...,N$ because every rule
in $\sss$ locks $\kappa(i)$-sectors).

Let $U_j$ be the admissible word obtained from $U_i$ by replacing
all indexes $i$ in $\kk$-letters by $j$:

$$U_j\equiv
\tk(j)u_1\kappa(j)\bk(j)u_2...\lk(j,n)u_{2n+2}\tk(j+1).$$

It is easy to see by inspection that rules of $\sss$ act in the
same way on $U_j$ for all $j>1$. Therefore, an induction on $|h|$
gives us

\begin{equation} \label{gimelj}
\Lambda_j(0)\cdot h=U_j(\omega).
\end{equation}

Let us also define $U_1$:

$$U_1\equiv
\tk(1)u_1\kappa(1)\Phi(u_1)\bk(1)u_2...\lk(1,n)u_{2n+2}\tk(2)$$
where $\Phi$ is the homomorphism from $\la A(\tk(1))\ra$ to $\la
A(\kappa(1))\ra$ which takes every letter $\bt{\omega,a}$ to
$a(\kappa(1))\iv$ if $\omega$ is even and greater than 0, and
every other letter to 1. Notice that the map $\Phi$ from
$A(\tk(1))$ to $A(\kappa(1))\cup\{1\}$ is indeed extendable to a
homomorphism from $\la A(\tk(1))\ra$ to $\la A(\kappa(1))\ra$
because $\la A(\tk(1))\ra$ is a free Burnside group and $\la
A(\kappa(1))\ra$ satisfies all Burnside relations.

\begin{remark}\label{remark10} {\rm
Notice that $U_1$ would be a copy of $U_j$ obtained by replacing
index $1$ for $j$ if not for the $\kappa$-sector: the
$\kappa(j)$-sector of an admissible word contains no
$\aaa$-letters, but $\kappa(1)$-sectors may contain
$\aaa$-letters. If the $\aaa$-subword of the $\kappa(1)$-sector of
$U_1$ is trivial (say, if $U_1$ is in the domain of $\bc_{2n}^{\pm
1}$), then $U_1$ is a copy of $U_j$.}
\end{remark}

For every word $u\equiv a_{i_1}...a_{i_s}$ over ${\cal A}\cup{\cal
A}\iv $ and every even $\omega$ from $\Omega$ let $T_\omega(u)$ be
the word $$\tau(\omega,a_{i_1})...\tau(\omega,a_{i_s}).$$ For
every odd $\omega\ne 2n+1$ from $\Omega$ let $$T_\omega(u)\equiv
\tau(\omega,a_{i_s})...\tau(\omega,a_{i_1}).$$

Let $z,z'$ be two letters from $\kk$, $z'\equiv z_+$ or $z'\equiv
z\iv$ then every sector of an admissible word of the form
$(z(\omega)uz'(\omega))^{\pm 1}$ will be called a \label{z1z2}{\em
sector of the type (form)} $[zz']$.

\begin{lm} \label{copy2}
Let $j>1$. Suppose $\Lambda_j(0)$ is in the domain of a word of
the form  $h\bc_{2n}$ where the word $h$ is Burnside-reduced and
does not contain $\bc_{2n}^{\pm 1}$. Let $U_j(2n+1)\equiv
\Lambda_j(0)\cdot h\bc_{2n}$. Then $\Lambda_1(0)$ is in the domain
of $h\bc_{2n}$ and $U_1(2n+1)\equiv \Lambda_1(0)\cdot h\bc_{2n}$.
\end{lm}

\proof Suppose that $h$ contains a subword of the form $\bc_i
h'\bc_i\iv$ where $h'$ contains no connecting rules: $h\equiv
h_1\bc_i h'\bc_i\iv h_2$ for some even $i$ from $0$ to $2n-2$.
Then consider the sector $W'$ of the type $[\lk(j,n)\tk(j+1)]$ of
$\Lambda_j(0)\cdot h_1\bc_i$. Since $\bc_i$ locks
$\lk(j,n)$-sectors, $W'$ has the form
$\lk(j,n,\omega')\tk(j+1,\omega')$ (the $\aaa$-part of the sectors
is trivial). Since $W'\cdot h'$ is in the domain of $\bc_i\iv$,
$W'=W'\cdot h'$. Since all rules in $h'$ are left active for
$\lk(j,n)$-sectors (all these rules belong to $\sss_{i+1}$, and
$i$ is odd), by Lemma \ref{lm:iji}, $h'=1$ modulo Burnside
relations, a contradiction. Similarly, $h$ does not contain
subwords of the form $\bc_i\iv h'\bc_i$ with even $i>0$ where $h'$
does not contain connecting rules.

Suppose that $h$ contains a subword of the form $\bc_ih'\bc_i\iv$
where $h'$ contains no connecting rules and $i$ is odd, $i\ne
2n-1$. Then applying the argument similar to the one used in the
previous paragraph to the sector of the form $[\rk(j,n)\lk(j,n)]$,
we get a contradiction. Similarly one can prove that $h$ does not
contain subwords of the form $\bc_i\iv h'\bc_i$ where $h'$
contains no connecting rules, $i$ is any odd number between $1$
and $2n-1$.

Suppose now that $h$ contains one of the connecting rules twice.
Since $\Lambda_j(0)$ has $\Omega$-coordinate $0$, and
$\Lambda_j(0)\cdot h$ has $\Omega$-coordinate $2n$, the facts
proved in the previous two paragraphs imply that $h$ contains a
subword of the form
$$h_0\bc_0h_1...h_{2n-2}\bc_{2n-1}h_{2n-1}\bc_{2n-1}\iv
h_{2n-2}'...\bc_0\iv h_0'$$ where the subwords $h_\omega,
h_\omega'$ are words over $\sss_\omega$. Let then $W'$ be the
sector of the form $[\bk(j)\rk(j,1)]$ of $\Lambda_j(0)\cdot h_0$.
The word $W'$ is in the domain of $$h''\equiv
\bc_0h_1...h_{2n-2}\bc_{2n-1}h_{2n-1}\bc_{2n-1}\iv
h_{2n-2}'...\bc_0\iv,$$ $\bc_0$ locks $W'$, all other rules from
$h''$ are right active for $\bk(j)$-sectors. None of these rules
belong to $\sss_{2n+1}$. Hence by Lemma \ref{lm:iji}, $h'=1$
modulo Burnside relations, a contradiction. Therefore $h$ contains
each of its connecting rules only once.

Hence $h$ can be represented in the form
$h_0\bc_0h_1\bc_1...h_{2n}$ for some words $h_i$ over $\sss_i$,
$i=0,...,2n$.

For some words $u_0,...,u_{2n}$ over ${\cal A}\cup{\cal A}\iv$,
$h_i\equiv T_i(u_i)$, $i=0,...,2n$. Let us prove by induction on
$i$ that $u_i=u_0$ modulo Burnside relations. Indeed, this
statement is trivial for $i=0$. Suppose that it is true for some
$\omega$ from $0$ to $2n-1$. Let us prove it for $\omega+1$. Since
$u_0=u_1=...=u_\omega$ modulo Burnside relations, by Lemma
\ref{intermediate}, $$\Lambda_j(0)\cdot
h_0\bc_0h_1\bc_1...h_\omega\bc_{\omega}=\Lambda_j(0)\cdot
T_0(u)\bc_0T_1(u)\bc_1...T_\omega(u)\bc_{\omega}=V.$$ In
particular the $\Omega$-coordinate of $V$ is $\omega+1$.

Suppose that $\omega$ is even. Then as in the proof of Lemma
\ref{lmbk} the sector of the form $[\rk(j,n)\lk(j,n)]$ of $V$ is
equal to $\rk(j,n,\omega+1)u\lk(j,n,\omega+1)$. Since $V\cdot
T_{\omega+1}(u_{\omega+1})$ is in the domain of $\bc_{\omega+1}$
and $\bc_{\omega+1}$ locks $\rk(j,n)$-sectors, we get that the
sector of the form $[\rk(j,n)\lk(j,n)]$ of $V\cdot
T_{\omega+1}(u_{\omega+1})$ is equal to
$\rk(j,n,\omega+1)\lk(j,n,\omega+1)$. By Lemma \ref{landr}, we get
that $u_{\omega+1}=u$ modulo Burnside relations, as required.

If $\omega$ is odd then we can apply the above argument to the
sector of the form $[\lk(j,n)\tk(j+1)]$ of $V$.

Therefore we can apply Lemma \ref{intermediate} $2n+1$ times and
conclude that $\Lambda_1(0)\cdot h_0=\Lambda_1(0)\cdot T_0(u)$,
$\Lambda_1(0)\cdot h_0$ is in the domain of $\bc_0h_1$ and
$\Lambda_1(0)\cdot h_0\bc_0 h_1=\Lambda_1(0)\cdot
T_0(u)\bc_0T_1(u)$,..., $\Lambda_1(0)\cdot
h_1\bc_0...h_{2n-1}=\Lambda_1(0)\cdot T_0(u)\bc_0...T_{2n-1}(u)$
is in the domain of $\bc_{2n-1}h_{2n}$ and

$$\Lambda_1(0)\cdot
h_1\bc_0...h_{2n-1}\bc_{2n-1}h_{2n}=\Lambda_1(0)\cdot
T_0(u)\bc_0...T_{2n-1}(u)\bc_{2n-1}T_{2n}(u)=V'.$$ All sectors of
the latter word contain no $\aaa$-letters except for the
$\tk(1)$-sector, $\kappa(1)$-sector, and $\bk(1)$-sector. The
first and the third of these sectors contain copies
$\overline{h\bc_{2n}}$ of the word $h\bc_{2n}$, the
$\kappa(1)$-sector contains $u^n$ (as in the proof of Lemma
\ref{lmbk}) which is equal to 1 modulo relations of $\la
A(\kappa(1))\ra$. Clearly $\Phi(\overline{h\bc_{2n}})=u^n$. Thus
$\Lambda_1(0)\cdot h\bc_{2n}$ exists and is equal to $U_1(2n+1)$.
\endproof

\begin{prop} \label{copy1}
Suppose that for some Burnside-reduced word $h$ over $\sss$, and
some $j>1$, $\Lambda_j(0)\cdot h=\Lambda_j(0)$. Then for every
$i=1,...,N$, $\Lambda_i(0)\cdot h=\Lambda_i(0)$.
\end{prop}

\proof If $h$ is empty, the statement is obviously true. By
(\ref{gimelj}), the statement is also true if $i>1$. So let $h$ be
non-empty and $i=1$.

Then as in the proof of Lemma \ref{copy2}, Lemma \ref{lm:iji}
implies that $h$ contains no subwords of the form $\bc_i^{\pm
1}f\bc_i^{\mp 1}$ where $f$ does not contain connecting rules,
$i\ne 0, 2n$. Hence $h$ can be represented in the form
$$h_1\bc_{2n}g_1\bc_{2n}\iv h_2\iv\bc_0\iv
h_2'\bc_{2n}g_2\bc_{2n}\iv h_3\iv\bc_0\iv h_3'\bc_{2n}g_3
...g_{s-1}\bc_{2n}\iv h_s\iv$$ where $g_p$ are non-empty words
over $\sss_{2n+1}$, $p=1,...,s-1$, $h_1,...,h_s,
h_2',...,h_{s-1}'$.

We claim that for every $t=2,...,s-1$, the word $$V\equiv
\Lambda_j(0)\cdot h_1\bc_{2n}g_1\bc_{2n}\iv h_2\iv\bc_0\iv
h_2'\bc_{2n}g_2\bc_{2n}\iv ... g_{t-1}\bc_{2n}\iv h_t\iv\bc_0\iv$$
is equal to $\Lambda_j(0)\cdot T_0(u_t)$ for some word $u_t$ over
$\aaa\cup\aaa\iv$.

Indeed, let $$V'\equiv \Lambda_j(0)\cdot h_1\bc_{2n}g_1\bc_{2n}\iv
h_2\iv\bc_0\iv h_2'\bc_{2n}g_2\bc_{2n}\iv ...
g_{t-1}\bc_{2n}\iv.$$

Then $V\equiv V'\cdot h_t\iv \bc_0\iv$. By Lemma \ref{lm:iji}, the
word $h_t$ is represented in the form $p_1\bc_1p_2\bc_2...p_{2n}$
where $p_i\equiv T_i(v_i)$ for some words $v_i$ over $\sss_i$,
$i=1,...,2n$. Since $V'\cdot h_t\iv$ is in the domain of
$\bc_0\iv$, $V'\cdot h_t\iv$ has the form

$$\tk(j,1)\kappa(j,1)\bk(j,1)\rk(j,1,1)w_1\lk(j,1,1)\rk(j,2,1)w_2...\rk(j,n,1)w_n\lk(j,n,1)\tk(j+1,1).$$

Since $V'\cdot h_t\iv T_1(v_1)$ is in the domain of $\bc_1$, and
$\bc_1$ locks any $\rk(j,i)$-sector, by Lemma \ref{landr}(iii),
each of the words $w_1,...,w_n$ is a copy of $v_1$. Hence $V'\cdot
h_t\iv\bc_0\iv=\Lambda_j(0)\cdot T_0(v_1)$ as required.

Therefore $\Lambda_j(0)$ is in the domain of
$T_0(u_t)\bc_0h_t\bc_{2n}$ and $T_0(u_t)h_t'\bc_{2n}$.

Let $U_{j,t}(2n+1)\equiv \Lambda_j(0)\cdot
T_0(u_t)\bc_0h_t\bc_{2n}$, $U_{j,t}'(2n+1)\equiv \Lambda_j(0)\cdot
T_0(u_t)h'_t\bc_{2n}$. By Lemma \ref{copy2}, $\Lambda_1(0)$ is in
the domain of $T_0(u_t)\bc_0h_t\bc_{2n}$ and
$T_0(u_t)h_t'\bc_{2n}$, and $U_{1,t}(2n+1)=\Lambda_j(0)\cdot
T_0(u_t)\bc_0h_t\bc_{2n}$, $U_{1,t}'(2n+1)=\Lambda_j(0)\cdot
T_0(u_t)h'_t\bc_{2n}$.

By Remark \ref{remark10}, $U_{1,t}(2n+1)$ is a copy of
$U_{j,t}(2n+1)$ and $U_{1,t}'(2n+1)$ is a copy of $U_{j,t}'(2n+1)$
An easy induction on the length of a word $g$ over $\sss_{2n+1}$
gives that $U_{1,t}(2n+1)\cdot g$ (resp. $U_{1,t}'(2n+1)\cdot g$)
is a copy of $U_{j,t}(2n+1)$ (resp. $U_{j,t}'(2n+1)$).

Thus we can conclude that $$
\begin{array}{l}
 (...((\Lambda_1(0)\cdot
h_1\bc_{2n}g_1)\cdot (\bc_{2n}\iv h_2\iv\bc_0\iv
T_0(u_1)\iv))\cdot T_0(u_1)h'_2\bc_{2n}g_2)\cdot ... )\cdot
T_0(u_{s-1})h'_{s-1}\bc_{2n}g_{s-1})\cdot\\ \bc_{2n}\iv h_s\iv
=\Lambda_1(0).\end{array}$$ By Lemma \ref{pseudo}, we obtain
$\Lambda_1(0)\cdot h=\Lambda_1(0)$ as required.
\endproof

\subsection{Stabilizers of arbitrary admissible words}

We call an admissible word $W$ \label{accepted}{\em accepted} if
for some word $h$ over $\sss$, the word $W\cdot h$ is a cyclic
shift of $\Lambda(0)^{\pm 1}$. Clearly $\Lambda(0)$ is an accepted
word and every admissible word of the form $\Lambda(0)\cdot h$ is
accepted.

In the rest of the section, we are going to prove the following
statement.

\begin{prop}\label{propsss}
Let $W$ be an admissible word of the form $z_1w_1z_2...z_sw_sz_1$
(the first and the last $\bkk$-letters are the same) which does
not contain an accepted subword.

(1) Suppose that for some $i=1,...,N$, every sector of $W$ is of
one of the types $[\kappa(i)\iv \kappa(i)]$, $[\kappa(i)\bk(i)]$
or $[\bk(i)\bk(i)\iv]$. Let $h$, $g$, $b$ be words in the alphabet
$\sss$ such that $W\cdot h=W$, $W\cdot ghg\iv=W$, $W\cdot b=W$ and
$b$ is equal to $g\iv hg$ modulo Burnside relations. Then for
every integer $s$ there exists a word $b_s$ which is equal to
$g^shg^{-s}$ modulo Burnside relations and such that $W\cdot
b_s=W$.

(2) Suppose that there is no $i$ from 1 to $N$ such that every
sector of $W$ is of one of the types $[\kappa(i)\iv \kappa(i)]$,
$[\kappa(i)\bk(i)]$ or $[\bk(i)\bk(i)\iv]$. Let $h$, $g$ be words
in the alphabet $\sss$ such that $W\cdot h=W$, $W\cdot ghg\iv=W$.
Then $W\cdot g=W$ or $g$ and $h$ commute modulo Burnside
relations.
\end{prop}

\begin{remark} \label{remark}
{\rm Clearly if $W\cdot g=W$ or $g$ and $h$ commute modulo
Burnside relations then the conclusion of (1) holds, so the
conclusion of (2) is stronger than the conclusion of (1).
}\end{remark}

\medskip


Let us prove condition (1) of the proposition. Thus consider the
case when for some $i=1,...,N$, $W$ contains sectors only of the
form $[\kappa(i)\iv \kappa(i)]$, $[\kappa(i)\bk(i)]$ or
$[\bk(i)\bk(i)\iv]$.

Notice that since $W$ starts and ends with the same $\bkk$-letter,
the definition of admissible words forces $W$ to have all three of
these types of sectors. Indeed, suppose that $W$ contains a sector
of the type $[\kappa(i)\bk(i)]$. Then $W$ must contain sectors of
the types $[z\kappa(i)]$ and $[\bk(i)z']$ for some $z,z'$. By
assumption, each of these sectors can be of one of the three types
$[\kappa(i)\iv \kappa(i)]$, $[\kappa(i)\bk(i)]$ or
$[\bk(i)\bk(i)\iv]$. Hence $z\equiv \kappa(i)\iv$, $z'\equiv
\bk(i)\iv$. Thus sectors of all three types are inside $W$. If $W$
contains a sector of the type $[\kappa(i)\iv \kappa(i)]$ then it
must contain a sector of the type $[\kappa(i)z]$ for some $z$ and
by assumption $z$ cannot be anything but $\bk(i)$. So $W$ contains
a sector of the type $[\kappa(i)\bk(i)]$ which implies that $W$
contains all three types of sectors. Finally if $W$ contains a
sector of the type $[\bk(i)\bk(i)\iv]$ then $W$ contains a sector
of the type $[z\bk(i)]$ for some $z$. This $z$ must be equal to
$\kappa(i)$, so $W$ contains a sector of the type
$[\kappa(i)\bk(i)]$ which, as before, implies that {\em $W$
contains all three types of sectors}.

{\em Therefore by Lemma \ref{lm:badsector} $h$ and $g$ do not
contain rules from $\sss_0\cup\{\bc_0,\bc_0\iv\}$.} We will use
this below several times.

We are going to use lemmas from Section \ref{bgfp}. Consider the
free group freely generated by
$\sss\backslash(\sss_0\cup\{\bc_0,\bc_0\iv\})$ as the free product
of the group $Q$ freely generated by $\sss_{2n+1}\cup\{\bc_{2n}\}$
and the group $P$ freely generated by other generators. Consider
the homomorphism \label{upsilon}$\upsilon$ from $P$ to $\la
A(\kappa(1))\ra$ which takes every rule $\tau$ from $P$ to
$R^a(\tau,\kappa(1))$. This allows us to talk about syllables of
words in $P*Q$ and about $\upsilon$-trivial elements of $P*Q$.

Notice that if $h$ and the reduced form of $ghg\iv$ consist of
rules from $P$, we can apply Lemma \ref{no2n+1} to a sector of the
form $[\kappa(i)\iv \kappa(i)]$ (case (i) of the formulation of
that lemma holds) and conclude that $g, h$ commute modulo Burnside
relations. This implies part (1) of the proposition.

Thus we can assume that $h$ or the reduced form of $ghg\iv$
contain rules from $\sss_{2n+1}\cup \{\bc_{2n}, \bc_{2n}\iv\}$.
Since $W\cdot h=W, W\cdot ghg\iv=W$, $g$ or the reduced form of
$ghg\iv$ must contain a rule from
$\sss_{2n+1}\cup\{\bc_{2n},\bc_{2n}\iv\}$.

Therefore there exists a word $x$ in $\sss$ such that $W'=W\cdot
x$ is in the domain of a rule from
$\sss_{2n+1}\cup\{\bc_{2n},\bc_{2n}\iv\}$. By Lemma \ref{pseudo}
$W'\cdot x\iv hx=W'$, $W'\cdot (x\iv gx)(x\iv hx)(x\iv gx)=W'$,
$W'\cdot x\iv bx=W'$. It is clear that if for every integer $s$
there exists $b'_s$ which is equal to $(x\iv gx)^s(x\iv hx)(x\iv
gx)^{-s}$ modulo Burnside relations and such that, $W'\cdot
b'_s=W'$, then condition (1) holds for $W$, $g, h$. {\em Hence we
can assume that $W=W'$.}

Thus all sectors of the form $[\kappa(i)\bk(i)]$ in $W$ are equal
to $\kappa(i)\bk(i)$ (since every rule from
$\sss_{2n+1}\cup\{\bc_{2n},\bc_{2n}\iv\}$ locks
$\kappa(i)$-sectors).

For every $i=1,...,N$, consider the homomorphism
\label{vtheta}$\vartheta_i:\tau\to L^a(\tau,\kappa(i))$ from $P*Q$
to $\la A(\tk(i))\ra$.

\begin{lm}\label{sectors}
Let $x \in P*Q$. Then the following conditions hold:

(i) $\kappa(1)\bk(1)\cdot x=\kappa(1)\bk(1)$ if and only if $x$ is
$\upsilon$-trivial.

(ii) $\kappa(1)\bk(1)$ is in the domain of $x$ if and only if
every $P$-syllable of $x$ except possibly the last one has trivial
$\upsilon$-image.

(iii) Let $W'\equiv \kappa(i)\iv u\kappa(i)$ be a 1-sector
admissible word, $i=1,...,N$. Then $W'\cdot x=W'$ if and only if
$\vartheta_i(x)$ commutes with $u$ modulo Burnside relations.

(iv) Every 1-sector admissible word of the form $[\kappa(i)\iv
\kappa(i)]$ or $[\bk(i)\bk(i)\iv]$, $i=1,...,N$, is in the domain
of $x$.

(v) For every sector $W'$ of the form $[\bk(i)\bk(i)\iv]$ of the
word $W$, $W'\cdot x=W'$.

(vi) If $i=1$ then the word $W$ is in the domain of $x$ if and
only if all $P$-syllables of $x$ except possibly the last one are
$\upsilon$-trivial. If $i>1$ then the word $W$ is in the domain of
$x$.

(vii) If $i=1$ then $W\cdot x=W$ if and only if $x$ is
$\upsilon$-trivial, and for every sector $W'\equiv \kappa(1)\iv
u\kappa(1)$ of $W$, $\vartheta_1(x)$ commutes with $u$ modulo
Burnside relations. If $i>1$ then $W\cdot x =W$ if and only if for
every sector $W'\equiv \kappa(i)\iv u\kappa(i)$, $\vartheta_i(x)$
commutes with $u$ modulo Burnside relations.
\end{lm}

\proof Since we identify letters $\kappa(i, \omega)$, $\omega\in
\Omega$, with $\kappa(i)$ and $\bk(i,\omega)$, $\omega\in \Omega$,
with $\bk(i)$ ($i=1,...,N$), all rules from $P$ are applicable to
$W$ and to any word of the form $W\cdot x$. Rules from $Q$ are
applicable to such a word if and only if every sector of the form
$[\kappa(1)\bk(1)]$ is equal to $\kappa(1)\bk(1)$ in $\hnka$. This
implies parts (i), (ii) of the lemma.

All rules from $P*Q$ are applicable to any sector of the form
$[\kappa(i)\iv \kappa(i)]$ or $[\bk(i)\bk(i)\iv]$. For every rule
$\tau$ from $P*Q$, $R^a(\tau,\bk(i))=L^a(\tau,\bk(i)_+)=1$. This
implies parts (iii),(iv) and (v) because of Lemma \ref{landr}.

Parts (vi), (vii) follow from the previous parts of the lemma by
Lemma \ref{p2}.
\endproof

\begin{lm} \label{kivk1}
Let $U\equiv \kappa(i)\iv u\kappa(i)$ be a one-sector admissible
word $i\ge 1$. Suppose that $U\cdot h=U, U\cdot ghg\iv=U$. Then
for every positive integer $s$, $U\cdot g^shg^{-s}=U$.
\end{lm}

\proof By Lemma \ref{sectors} (iii) $\vartheta_i(h)$,
$\vartheta_i(ghg\iv)$ commutes with $u$ modulo Burnside relations.
Since $u\ne 1$ modulo Burnside relations (we have agreed to
consider only reduced admissible words), by Lemma
\ref{additional}, then $\vartheta_i(g)$ commutes with
$\vartheta_i(h)$ modulo Burnside relations. Hence
$\vartheta(g^shg^{-s})=\vartheta(h)$ modulo Burnside relations.
Since by Lemma \ref{sectors} (iv) $U$ is in the domain of
$g^shg^{-s}$, by Lemma \ref{burns}(ii), we get $U\cdot
g^shg^{-s}=U\cdot h=U$.\endproof

Now we are ready to prove part (1) of the proposition if $i>1$.

\begin{lm}\label{onlykbk1}
Part (1) of the Proposition \ref{propsss} holds if $i>1$.
\end{lm}

\proof Indeed all rules from $\sss$ lock sectors of the form
[$\kappa(i)\bk(i)$] and are not left active for $\bk(i)$-sectors
and for $\bk(i)\iv$-sectors. Therefore if $W'$ is any of the
sectors of $W$ of the form [$\kappa(i)\bk(i)$] or
[$\bk(i)\bk(i)\iv$] then $W'\cdot g=W'$, so for every integer $s$,
$W'\cdot g^{s}hg^{-s}=W'$. By Lemma \ref{kivk1} the same
conclusion can be drawn for sectors of the form [$\kappa(i)\iv
\kappa(i)$]. Hence by Lemma \ref{p2} $W\cdot g^shg^{-s}=W$ as
required.
\endproof

Now let $i=1$. It is obvious that we can assume that $g,h\ne 1$
modulo Burnside relations.

\begin{lm} \label{x}
If there exists a word $x\in P*Q$ such that $g=xpx\iv$, $h=xqx\iv$
modulo Burnside relations, where $p,q\in P$, then $g$ and $h$
commute modulo Burnside relations.
\end{lm}

\proof Suppose that such a word $x$ exists. By Lemma \ref{sectors}
(i), the word $h$, the reduced form of $ghg\iv$ and the word $b$
from the formulation of part (1) of the proposition are
$\upsilon$-trivial. By Lemma \ref{phitriv}, $g$ is also
$\upsilon$-trivial. Then by Lemma \ref{bernsyl}(iii) (since $g,h
\ne 1$) there exist words $x_1$ and $x_2$ such that
$x_1px_1\iv=g$, $x_2qx_2\iv=h$ modulo Burnside relations,
$x_1px_1\iv$, $x_2qx_2\iv$ are $\upsilon$-trivial. Since $h,g\ne
1$ modulo Burnside relations, $x_1px_1\iv$, $x_2qx_2\iv$ can be
considered reduced. This implies that $\upsilon(p)=\upsilon(q)=1$.

We have that $x\iv x_1$ commutes with $p$ modulo Burnside
relations. Therefore $x\iv x_1$ belongs to the centralizer of $p$
modulo Burnside relations. The centralizer of every non-trivial
element of a free Burnside group is cyclic of order $n$ by Lemma
\ref{VI}. Therefore the centralizer of $p$ in the free Burnside
factor of $P*Q$ and in the free Burnside factor of $P$ are the
same. Thus $x\iv x_1$ belongs to $P$ modulo Burnside relations and
$x=x_1p_1$ for some $p_1\in P$ from the same cyclic subgroup as
$p$. Since $x_1px_1\iv$ is $*$-reduced and $\upsilon$-trivial, all
$P$-syllables of $x_1p_1$ except for, possibly, the last one are
$\upsilon$-trivial. Hence by Lemma \ref{sectors}(vi), $W$ is in
the domain of $x_1p_1$.

Notice that $x_1p_1qp_1\iv x_1\iv = xqx\iv = h$ modulo Burnside
relations. Hence by Lemma \ref{burns}(ii) $W\cdot x_1p_1 qp_1\iv
x_1\iv=W\cdot h=W$. The word $ghg\iv$ is equal to $$x_1px_1\iv
x_1p_1qp_1\iv x_1\iv x_1p\iv x_1\iv= x_1pp_1qp_1\iv p\iv x_1\iv$$
modulo Burnside relations.

Let $W_1=W\cdot x_1$. Then $W_1\cdot p_1qp_1\iv=W_1$ and $W_1\cdot
p(p_1qp_1\iv)p\iv =W_1$. By Lemma \ref{lm:badsector1} applied to a
$[\kappa(i)\iv \kappa(i)]$-sector of $W_1$, $p$ commutes with
$p_1qp_1\iv$ modulo Burnside relations. Conjugating by $x_1$, we
deduce that $g, h$ commute modulo Burnside relations.
\endproof

Now we are ready to {\em prove part (1) of Proposition
\ref{propsss} for $i=1$}.

Recall that we can assume that all sectors of $W$ of the form
$[\kappa(1)\bk(1)]$ are equal to $\kappa(1)\bk(1)$ in $\hnka$.

By Lemma \ref{x}, we can assume that there is no word $x$ such
that $g=xpx\iv$, $h=xqx\iv$ modulo Burnside relations, where
$p,q\in P$.

By Lemma \ref{sectors}(i), $h, ghg\iv$, and $b$ are
$\upsilon$-trivial words. By Lemma \ref{phitriv} then $g$ is an
$\upsilon$-trivial word. Hence by Lemma \ref{sectors}(i),
$W_1\cdot g=W_1$ for every sector of $W$ of the form
$[\kappa(1)\bk(1)]$. This implies that $W_1\cdot g^shg^{-s}=W_1$
for every integer $s$. By Lemma \ref{kivk1}, the same equality is
true for every sector of $W$ of the form $[\kappa(1)\iv
\kappa(1)]$. By Lemma \ref{sectors}(v) the same is true for
sectors of the form $[\bk(1)\bk(1)\iv]$. Thus $W\cdot
g^shg^{-s}=W$ for every integer $s$.

This completes the proof of part (1) of Proposition \ref{propsss}.

\medskip

 Now let us prove part (2) of the proposition. So let us assume that there is
no $i=1,...,N$ such that each sector of $W$ is of one of the forms
$[\kappa(i)\bk(i)]$, $[\kappa(i)\iv \kappa(i)]$,
$[\bk(i)\bk(i)\iv]$. Assume that part (2) of Proposition
\ref{propsss} is not true and that $(|h|,|g|)$ is the smallest in
lexicographic order pair of numbers for all such counterexamples
$(h,g)$ (for any $W$). Clearly the word $h$ is not empty.

The next four lemmas allow us to improve out counterexample
$(h,g)$.

\begin{lm} \label{cyc1}
For every word $x$ over $\sss$ such that $W$ is in the domain of
$x$, the pair $(x\iv hx, x\iv gx)$ is also a counterexample to
statement (2) of Proposition \ref{propsss}.
\end{lm}

\proof Indeed by Lemma \ref{pseudo} $(W\cdot x)\cdot x\iv
hx=W\cdot x$, $$(W\cdot x)\cdot (x\iv gx)(x\iv hx)(x\iv
gx)\iv=W\cdot x.$$ Also since $W$ does not contain accepted
subwords, $W\cdot x$ does not contain accepted subwords.
\endproof

\begin{lm} \label{cyc}
The word $h$ is cyclically reduced.
\end{lm}

\proof Suppose that $h$ is not cyclically reduced. We can
represent $h$ in the form $h_1h_2$ such that the reduced form $h'$
of the word $h_2h_1$ is cyclically reduced and $|h'|<|h|$. By
Lemma \ref{cyc1}, the pair $(h', h_1\iv gh_1)$ is also a
counterexample to part (2) of Proposition \ref{propsss}. This
contradicts the minimality of the counterexample $(h,g)$. Hence
$h$ is cyclically reduced.
\endproof

\begin{lm} \label{lm:square1}
There exists a possibly different counterexample $(\bar h, \bar
g)$ with $|\bar h|=|h|$, $|\bar g|=|g|$ such that $\bar h$ is
cyclically reduced and $\bar g\bar h\bar g\iv$ is reduced.
Moreover for some prefix $x\iv$ of the word $h$, we have $\bar h=x
hx\iv$, $\bar g=xgx\iv$ in the free group.
\end{lm}

\proof By Lemma \ref{cyc}, $h$ is cyclically reduced.

Therefore either $gh$ or $hg\iv$ is a freely reduced word. Let us
assume that $hg\iv$ is freely reduced. The right-left dual of the
argument presented below will work in the case when $gh$ is freely
reduced.

Notice that $g\nee g_1h^{-1}$ for any $g_1$, since otherwise
$g_1hg_1^{-1}$ is freely equal to $ghg^{-1}$, and the
counterexample $(h,g_1)$ would be smaller than $(g,h)$. Thus
$g\equiv g_1f$ where $f^{-1}f'\equiv h$ and $g_1f'$ is a reduced
word. Hence $ghg^{-1}=g_1f'f^{-1}g_1^{-1}$, whence

\begin{equation}\label{formula}
W\cdot g_1f'f\iv g_1\iv=W
\end{equation}
Since $W\cdot f\iv f'=W\cdot h=W$, we have $(W\cdot f\iv)\cdot
f'f\iv = W\cdot f\iv$. Let $W'=W\cdot f\iv$. Then $W'\cdot
f'f\iv=W'$ and $W'\cdot (fg_1)f'f\iv (fg_1)\iv =W'$ by
(\ref{formula}).

Now notice that since $h$ is cyclically reduced, $|f'f\iv|=|h|$.
If $|fg_1|<|g|$ then the minimality of $(|h|,|g|)$ implies that
$W'\cdot (fg_1)=W'$ or $fg_1$ commutes with $f'f\iv$ modulo
Burnside relations. In the first case $W\cdot g_1f=W$. Since
$W\cdot h=W$, we get $W\cdot g=W\cdot (g_1f)= W$, a contradiction.
In the second case conjugating by $f\iv$, we get that $g_1f$
commutes with $f\iv f'=h$ modulo Burnside relations. Therefore
$g\equiv g_1f$ commutes with $h$ modulo Burnside relations, a
contradiction.

Hence the length of the freely reduced form of $fg_1$ is equal to
$|g|$. This implies that  if $g_1$ is not empty then $g$ is
cyclically reduced (since $f$ is not empty by our assumption).

{\bf Case 1.} Suppose that $g_1$ is not empty. Then the word
$fg_1f'f\iv g_1\iv f\iv$ is reduced ($fg_1$ is reduced because
$g\equiv g_1f$, $|fg_1|=|g|$, $g_1f'$ is reduced by the
definition, $f'f\iv$ is reduced because it is a cyclic shift of
the cyclically reduced word $h$, $f\iv g_1\iv$ is reduced because
$g_1f\equiv g$ is reduced). In this case $\bar h\equiv f'f\iv$,
$\bar g\equiv fg_1$ and $x\equiv f$. Therefore $W$ is in the
domain of $x\iv$. By Lemma \ref{cyc1}, the pair $(\bar h, \bar g)$
is a (minimal) counterexample, as required.

{\bf Case 2.} Now let $g_1$ be empty, that is $f\equiv g$. If $g$
is not cyclically reduced, then $g\equiv yg_2y\iv$, $h\equiv
yg_2^{-1}y^{-1}h_2$ for some $g_2, h_2$ and a non-empty $y$.
Applying Lemma \ref{cyc1} for $x\equiv y$, we obtain a
counterexample $(g_2y\iv h_2y, g_2)$. Since $g_2$ is shorter than
$g$, and $g_2y\iv h_2y$ is of the same length as $h$, we get a
contradiction. Hence $g$ is cyclically reduced.

Now let $g^{-s}$ be the maximal power of $g$ which is a prefix of
$h$. Using Lemma \ref{cyc1} with $x\equiv g^{-s}$, we can pass
from the pair $(h,g)$ to another minimal counterexample $(\hat h,
g)= (h_1g^{-s},g)$ where $h_1$ is defined by the formula $h\equiv
g^{-s} h_1$. Notice that since $h$ is cyclically reduced by Lemma
\ref{cyc}, $h_1g^{-s}$ is cyclically reduced. The word $h_1$ is
not trivial because otherwise $g$ and $h$ would commute, so
$(h,g)$ would not be a counterexample. Since $g$ is cyclically
reduced, the product $\hat h g\iv$ is reduced. If $g\hat h$ is
reduced, the product $g\hat hg\iv$ is reduced and we are done.
Otherwise, since $g\iv$ is not a prefix of $\hat h$ (by the choice
of $s$), we can apply the argument of Case 1 and produce the
required counterexample.
\endproof

{\em Lemma \ref{lm:square1} allows us to assume (from now on) that
$ghg\iv$ is reduced.}

\begin{lm} \label{square2}
The words $g$, $h$ are Burnside reduced.
\end{lm}

\proof Let $h$ contain a non-empty subword $c$ which is equal to
$1$ modulo Burnside relations. Then $h\equiv h_1ch_2$ for some
$h_1, h_2$. By Lemma \ref{burns} (iii), $W\cdot h_1h_2=W$. Since
$ghg\iv$ is reduced, we have $W\cdot g\cdot h_1ch_2\cdot g\iv=W$,
so $W\cdot gh_1h_2g\iv=W$. Now it is clear that $(h_1h_2,g)$ is a
counterexample to part (2) of Proposition \ref{propsss}, which
contradicts the minimality of $(h,g)$. Similarly if $g$ is not
Burnside-reduced, we can shorten $g$.
\endproof

\begin{lm} \label{square3}
One cannot represent $h$ as a product $h_1h_2h_3$ where
$h_1=h_3\iv$ modulo Burnside relations and $h_1, h_3$ consist of
rules from the same $\sss_\omega$, $\omega\in \Omega$, and at
least one of the words $h_1$ or $h_3$ is not empty.
\end{lm}

\proof If one of the words $h_1$ or $h_3$ is empty then
$h_1=h_3=1$ modulo Burnside relations, and we get a contradiction
with Lemma \ref{square2}. Hence both $h_1$ and $h_3$ are not
empty. Let $f$ be a shortest word among $h_1, h_3\iv$. By Lemma
\ref{intermediate}, $h_1, h_3\iv, f$ define the same
transformation from $P(\sss)$. Hence we have $$W\cdot h=W\cdot
h_1\cdot h_2\cdot h_3= W\cdot f\cdot h_2\cdot f\iv =W\cdot
fh_2f\iv=W$$ and, since $ghg\iv$ is reduced, $$W\cdot g\cdot
h_1\cdot h_2\cdot h_3\cdot g\iv = W\cdot g(fh_2f\iv)g\iv=W.$$
Since $h=fh_2f\iv$ modulo Burnside relations, the pair $(fh_2f\iv,
g)$ is a counterexample to part (2) of Proposition \ref{propsss}.
Since $|f|\le |h_1|, |h_3|$, $(fh_2f\iv, g)$ is a minimal
counterexample. But $fh_2f\iv$ is not cyclically reduced, a
contradiction with Lemma \ref{cyc}.
\endproof

\begin{lm} \label{twosteps}
Suppose that $ghg\iv$ contains a connecting rule $\bc_\omega^{\pm
1}$ which locks a sector $zwz'$ of $W$, $z\nee \kappa(1),
\bk(1)\iv$. Suppose that all rules in $ghg\iv$ which differ from
$\bc_\omega^{\pm 1}$ are left or right active for $z$-sectors.
Then $ghg\iv$ contains at least two different connecting rules
$\bc_i^{\pm 1}$ and $\bc_{i+1}^{\pm 1}$, $i\in \Omega$.
\end{lm}

\proof By Lemma \ref{lm:badsector}, $z'\equiv z_+$; $z\nee
\kappa(i), \bk(i)\iv$ for any $i$ because all rules lock
$\kappa(i)$-sectors if $i>1$, and $z\nee \kappa(1), \bk(1)\iv$ by
assumption. Hence every connecting rule changes the
$\Omega$-coordinate of $z$-sectors.

Suppose, by contradiction,  that $ghg\iv$ does not contain two
different connecting rules. Since $W\cdot ghg\iv=W$ and $W\cdot
h=W$, and the word $ghg\iv$ is reduced, $W$ and $W\cdot g$ have
the same $\Omega$-coordinate. Hence either $g$ does not contain
$\bc_\omega^{\pm 1}$ or $g$ contains a subword $\bc_\omega^{\pm 1}
f\bc_\omega^{\mp 1}$ where $f$ contains no connecting rules. In
the first case, $h$ contains $\bc_\omega^{\pm 1}$, so (since
$W\cdot h=W$), $h$ contains a subword of the form $\bc_\omega^{\pm
1} f\bc_\omega^{\mp 1}$. In both cases $f$ contains rules from
only one $\sss_i$. All these rules are active for $z$-sectors by
our assumption. Hence by Lemma \ref{lm:iji} (which can be applied
since $z\nee \kappa(1), \bk(1)\iv$), $f=1$ modulo Burnside
relations, a contradiction with Lemma \ref{square2}.
\endproof

\begin{lm} \label{tktk}
$W$ cannot contain a sector of the form $[\rk(i,j)\rk(i,j)\iv]$.
\end{lm}

\proof Suppose that $W$ contains  such a sector $W'\equiv
\rk(i,j,\omega)w\rk(i,j,\omega)\iv$. Notice that since $W$ starts
and ends with the same $\bkk$-letter, it must contain a sector
$W''$ of the form $z\iv w''z$, $z\in \kk(\omega)$, either to the
left or to the right of our sector $W'$. Replacing $W$ by $W\iv$
if necessary, we can assume that such a sector occurs to the left
of $W'$. We also choose $W''$ to be as close to $W'$ on the left
as possible.

By Lemma \ref{lm:badsector}, $ghg\iv$ does not contain occurrences
of $\bc_\omega$ or rules from $\sss_\omega$ which lock
$\rk(i,j)$-sectors. For every $t\in \Omega$ either $\bc_t$ or
$\bc_{t+1}$ locks $\rk(i,j)$-sectors, and $\bc_{2n}$ and every
rule from $\sss_{2n+1}$ lock these sectors. Therefore $ghg\iv$ can
only contain rules from two consecutive $\sss_t$, $\sss_{t+1}$
($\omega=t$ or $\omega=t+1$), and does not contain rules from
$\sss_{2n+1} \cup \{\bc_{2n},\bc_{2n}\iv\}$.

If $ghg\iv$ contains rules only from one $\sss_t$, then by Lemma
\ref{lm:badsector1}(ii), we have that $h$ and $ghg\iv$ are equal
modulo Burnside relations to powers of $\rt_\omega(w)$. By Lemma
\ref{additional}, it follows that $h$ and $g$ commute modulo
Burnside relations, a contradiction.

Thus we can assume that $ghg\iv$ contains rules from two different
sets $\sss_t$ and $\sss_{t+1}$ or $ghg\iv\equiv \bc_{t}^{\pm
1}f\bc_t^{\mp 1}$ for some $f$. In both cases $ghg\iv$ contains an
occurrence of $\bc_t$ and an occurrence of $\bc_t\iv$.

Since $W'$ is not a prefix of $W$, by the definition of admissible
words, $W$ contains a subword of one of the following forms:
$$\begin{array}{l}\rk(i,j,\omega)\iv w_1 \rk(i,j,\omega) w
\rk(i,j,\omega)\iv,\\
\lk(i,j-1,\omega)w_1\rk(i,j,\omega)w\rk(i,j,\omega)\iv \quad
(j>1),\\ \bk(i)w_1\rk(i,1,\omega)w\rk(i,1,\omega)\iv,\end{array}$$
for some word $w_1$ over $A(\rk(i,j)_-)$.

Consider each of these cases.

{\bf Case 1.} $W$ contains $\rk(i,j,\omega)\iv w_1 \rk(i,j,\omega)
w \rk(i,j,\omega)\iv$ for some $w$.

\medskip

{\bf Case 1.1.} Let $j>1$. Then $\rk(i,j)_-=\lk(i,j-1)$. Since for
every $t\ge 0$ if $\bc_t$ is active for $\rk(i,j)$-sectors then it
locks $\lk(i,j-1)$-sectors, $\bc_t$ locks $\lk(i,j-1)$-sectors.
This contradicts Lemma \ref{lm:badsector}.

\medskip

{\bf Case 1.2.} Let $j=1$. Then $\rk(i,j)_-\equiv \bk(i)$. Then by
Lemma \ref{lm:badsector}, $t, t+1\ne 0$. Recall also that
$t,t+1\ne 2n+1$ and that all other (working) rules except for
rules from $\sss_0\cup\{\bc_0\}$ are right active for
$\bk(i)$-sectors. By Lemma \ref{lm:badsector1} (ii) $h$ and
$ghg\iv$ are a powers of $\rt_\omega(w_1)$ modulo Burnside
relations, so by Lemma \ref{additional}, $h$ and $g$ commute
modulo Burnside relations, a contradiction.

\medskip

{\bf Case 2.} $W$ contains
$\lk(i,j-1,\omega)w_1\rk(i,j,\omega)w\rk(i,j,\omega)\iv$, $j>1$.
Then $\bc_t$ locks $\lk(i,j-1)$-sectors (indeed, every connecting
rule which does not lock some $\rk(s,l)$-sectors locks all
$\lk(s,l)$-sectors). Lemma \ref{twosteps} now implies that rules
of one of the sets $\sss_t$ or $\sss_{t+1}$ lock
$\lk(i,j-1)$-sectors. Let $\sss_{t'}$, $t'\in\{t,t+1\}$ consist of
rules which lock $\lk(i,j-1)$-sectors. By definition of $\sss$,
$t'\in \{0, 2j-1, 2j, 2j+1,...,2n+1\}$. Since for $t'\in
\{2j,...,2n+1\}$ the rules of $\sss_{t'}$ lock also
$\rk(i,j)$-sectors, we conclude that $t'\in \{0,2j-1\}$.

{\bf Case 2.1.} Let $t'=2j-1$, $j>1$. Then $t'=t+1$ because rules
from $\sss_{2j}$ lock $\rk(i,j)$-sectors. So $t=2j-2$. Since all
rules from $\sss_{2j-1}$ lock all $\rk(i,1)$-,
$\lk(i,1)$-,...,$\lk(i,j-1)$-sectors, none of the letters
$\rk(i,1,\omega)$, $\lk(i,1,\omega)$, ..., $\lk(i,j-1,\omega)$
belong to the sector $W''$ (by Lemma \ref{lm:badsector}). The
definition of admissible words now implies that in $W$,
$\bkk$-letters to the left of $\lk(i,j-1,\omega)$ are placed in
the natural (decreasing) order from $\lk(i,j-1,\omega)$ to
$\rk(i,1,\omega)$. The $\bkk$-letter in $W$ next to the left of
$\rk(i,1,\omega)$ may be either $\bk(i)$ or $\rk(i,1,\omega)\iv$.
Consider both cases.

\medskip

{\bf Case 2.1.1.} Suppose that $W$ contains a sector of the form
$\bk(i)w'''\rk(i,1,\omega)$. Since neither $t$ nor $t+1$ are equal
to $0$ or $2n+1$ (recall that $t'\in \{t, t+1\}$, $t'=2j-1, j>1$),
all rules from $\sss_t\cup\sss_{t+1}\cup\{\bc_t\}$ are right
active for the $\bk(i)$-sectors and (by Lemma \ref{landr})
$L^a(h,\rk(i,1))$ is a copy of $h$. Since $W\cdot h=W$, we
conclude that $h=1$ modulo Burnside relations, a contradiction.

\medskip

{\bf Case 2.1.2.} Suppose that $W$ contains a sector of the form
$\rk(i,1,\omega)\iv w'''\rk(i,1,\omega)$. Then by Lemma
\ref{lm:badsector1} $h$ and $ghg\iv$ belong to the same cyclic
subgroup generated by $\rt_\omega(w''')$ modulo Burnside
relations. As before, by Lemma \ref{additional}, this implies that
 $h$ and $g$ commute modulo Burnside relations, a contradiction.

\medskip

{\bf Case 2.2.} Let $t'=0$. Then $ghg\iv$ contain rules from
$\sss_0\cup\sss_1\cup\{\bc_0,\bc_0\iv\}$ only.

Notice that $g$ cannot contain a subword of the form $\bc_0
g_1\bc_0\iv$ where $g_1$ is a word over $\sss_1$. Indeed otherwise
by Lemma \ref{lm:iji} applied to the sector
$\lk(i,j-1,\omega)w_1\rk(i,j,\omega)$ of $W$, we get $g_1=1$
modulo Burnside relations (since $\bc_0$ locks this sector, and
all rules from $\sss_1$ are left active for it). This would
contradict Lemma \ref{square2}. Similarly $h$ cannot contain a
subword of the form $\bc_0 h_1\bc_0\iv$ where $h_1$ consists of
rules from $\sss_1$. Hence either $h$ contains rules from only one
$\sss_\omega$ or $h$ has the form $h_1\bc_0\iv h_2\bc_0h_3$ where
$h_1, h_3$ are possibly empty words over $\sss_1$, $h_2$ is a word
over $\sss_0$.

Since $W$ and $W\cdot g$ are in the domain of $h$, $W$ and $W\cdot
g$ have the same $\Omega$-coordinate. Hence the total degree of
$\bc_0$ in $g$ is $0$. Therefore either $g$ does not contain
$\bc_0$ and $\bc_0\iv$ or $g$ has the form $g_1\bc_0\iv g_2\bc_0
g_3$ where $g_1, g_3$ are words over $\sss_1$ (these words can be
empty) and $g_2$ is a non-empty word over $\sss_0$.

Consider these two cases separately.

\medskip

{\bf Case 2.2.1.} Suppose that $g$ does not contain $\bc_0$ and
$\bc_0\iv$. Since $ghg\iv$ must contain $\bc_0$ or $\bc_0\iv$, $h$
has the form $h_1\bc_0\iv h_2\bc_0h_3$ where $h_1, h_3$ are words
over $\sss_1$, $h_2$ is a word over $\sss_0$. Since $W\cdot h=W$,
the $\Omega$-coordinate of $W$ is $1$, and $g$ is a word over
$\sss_1$.

Let us apply Lemma \ref{landr} (i) to the sector
$\lk(i,j-1,\omega)w_1\rk(i,j,\omega)$ of $W$. Since $\bc_0$ locks
this sector, $R^a(gh_1,\lk(i,j-1))=R^a(h_1,\lk(i,j-1))=w_1$ modulo
Burnside relations. By Lemma \ref{landr} (iii), we obtain that
$gh_1=h_1$ modulo Burnside relations, so $g=1$ modulo Burnside
relations, a contradiction.

\medskip

{\bf Case 2.2.2.} Suppose that $g\equiv g_1\bc_0\iv g_2\bc_0 g_3$
where $g_1,g_3$ are words over $\sss_1$, $g_2$ is a word over
$\sss_0$.

If $h$ consists of rules from $\sss_1$, then we can apply Lemma
\ref{lm:iji} to the subword $g_3hg_3\iv$ of $ghg\iv$ and a sector
of the type $[\lk(i,j-1)\rk(i,j)]$ of $W\cdot g_1\bc_0\iv
g_2\bc_0$ to obtain that $g_3hg_3\iv$ is equal to 1 modulo
Burnside relations. Hence $h=1$ modulo Burnside relations, a
contradiction.

Since $h$ does not have subwords of the form $c_0\bar h c_0^{-1}$,
$h$ has the form $h_1\bc_0\iv h_2\bc_0h_3$ where $h_1,h_3$ are
words over $\sss_1$, $h_2$ is a word over $\sss_0$. Then again by
Lemma \ref{lm:iji}, $g_3h_1=1$, $h_3g_3\iv=1$ modulo Burnside
relations. Hence $h_1=h_3\iv$ modulo Burnside relations, which
contradicts Lemma \ref{square3} if $h_1$ or $h_3$ are not empty or
Lemma \ref{cyc} if $h_1$ and $h_3$ are empty.

\begin{remark}{\rm \label{rem90}
Notice that the only properties of the sector
$\lk(i,j-1,\omega)w_1\rk(i,j,\omega)$ used in Case 2.2 are the
properties that all rules from $\sss_0\cup\{\bc_0\}$ lock this
sector and all rules from $\sss_1$ are active for this sector.}
\end{remark}

\medskip

{\bf Case 3.} Let $W$ contain sectors $\bk(i)w_1\rk(i,1,\omega)$
and $\rk(i,1,\omega) w\rk(i,1,\omega)\iv$. Recall that $t, t+1\ne
2n+1$.

\medskip

{\bf Case 3.1.} Suppose that $t\ne 0$. Then all rules from
$\sss_t\cup\sss_{t+1}\cup\{\bc_t\}$ are right active for
$\bk(i)$-sectors. Therefore by Lemma \ref{lm:iji1} $h=1$ modulo
Burnside relations, a contradiction.

\medskip

{\bf Case 3.2.} Let $t=0$. By Remark \ref{rem90}, this case is
similar to Case 2.2 because all rules from $\sss_0\cup\{\bc_0\}$
lock $\bk(i)$-sectors and all rules from $\sss_1$ are (right)
active for these sectors.
\endproof

\begin{lm} \label{tktk1}
$W$ cannot contain a sector of the form $[\lk(i,j)\lk(i,j)\iv]$,
$i=1,...,N, j=1,...,n$.
\end{lm}

\proof Suppose that $W$ contains a sector $W'$ of the form
$\lk(i,j,\omega)w\lk(i,j,\omega)\iv$, $i=1,...,N, j=1,...,n,
\omega\in\Omega$. As in the proof of Lemma \ref{tktk}, we can
assume that $W$ contains a sector $W''$ of the form $[z\iv z]$ to
the left of $W'$, $z\in \bkk$. We also assume that $W''$ is the
closest to $W'$ such sector (in principle, $W''$ can have a common
letter $\lk(i,j,\omega)$ with $W'$).

Since all rules from $\sss_0\cup\sss_{2n+1}\cup\{\bc_{2i},
i=0,...,n\}$ lock $\lk(i,j)$-sectors, none of these rules can
occur in $ghg\iv$ by Lemma \ref{lm:badsector}. Therefore $ghg\iv$
involves rules from at most two sets $\sss_t, \sss_{t+1}$.

If $ghg\iv$ contains rules of only one set $\sss_\omega$ then by
Lemma \ref{no2n+1} $(h,g)$ is not a counterexample to part (2) of
Proposition \ref{propsss}, a contradiction.

Thus $ghg\iv$ contains $\bc_t$ for some $t$. Moreover this $t$
must be odd because connecting rules with even indices lock
$\lk(i,j)$-sectors.

Then $\bc_t$ locks $\rk(k,\ell)$-sectors, $1\le k\le N, 1\le
\ell\le n$. So by Lemma \ref{lm:badsector}, $z\nee \lk(k,\ell)$
for any $k,\ell$.

Hence the first letter $\lk(i,j,\omega)$ of $W'$ does not belong
to $W''$. Therefore $W$ contains a sector of the form
$[\rk(i,j)\lk(i,j)]$. If rules from $\sss_t$ and $\sss_{t+1}$ are
active for $\rk(i,j)$-sectors, and $\bc_t$ locks these sectors, we
get a contradiction with Lemma \ref{twosteps}.

 Hence rules from $\sss_{t+1}$ must lock
$\rk(i,j)$-sectors. Therefore these rules lock all $\rk(i,j')$-
and $\lk(i,j')$-sectors with $j'<j$. Hence $z\nee \rk(i,j')$,
$j'=2,3,...,j$ and $z\nee \lk(i,j')$, $j'=1,2,...,j$. Hence the
$\bkk$-letter preceding $\lk(i,j,\omega)$ in $W$ is
$\rk(i,j,\omega)$, the $\bkk$-letter preceding $\rk(i,j,\omega)$
is $\lk(i,j-1,\omega)$, etc. We conclude that $W$ contains a
sector of one of the forms $[\bk(i)\rk(i,1)]$ or $[\rk(i,1)\iv
\rk(i,1)]$.

In the first case we get a contradiction with Lemma \ref{lm:]tap}.
In the second case by Lemma \ref{lm:badsector1} $h, ghg\iv$ belong
to the same cyclic subgroup modulo Burnside relations, and we get
a contradiction using Lemma \ref{additional}.
\endproof

\begin{lm} \label{notall}
Not all sectors of $W$ have the form $[zz\iv]$ where $z$ is one of
the letters $\tk(i,\omega)$, $\tk(i,\omega)\iv$, $\bk(i)$,
$\bk(i)\iv$.
\end{lm}

\proof Suppose that each sector of $W$ has one of these forms. But
for each $$z\in\{\tk(i,\omega), \tk(i,\omega)\iv, \bk(i),
\bk(i)\iv\},$$ $L^a(\tau,z)=R^a(\tau,z)=1$ for every $\tau\in
\sss$. Notice that $W$ is in the domain of $g$ because $ghg\iv$ is
reduced and $W\cdot ghg\iv=W$.  Now by Lemma \ref{landr} we have
$W\cdot g=W$, a contradiction.
\endproof

\begin{lm}\label{kkiv}
$W$ does not contain sectors of the form
$[\kappa(i)\kappa(i)\iv]$.
\end{lm}

\proof If $i>1$ then the statement immediately follows from Lemma
\ref{lm:badsector} because all rules in $\sss$ lock
$\kappa(i)$-sectors. So let $i=1$ and suppose that $W$ contains a
sector of the form $[\kappa(1)\kappa(1)\iv]$. Then by Lemma
\ref{lm:badsector} $ghg\iv$ does not contain rules from $\sss_0$
and $\sss_{2n+1}$ and connecting rules $\bc_0^{\pm 1},
\bc_{2n}^{\pm 1}$. As before we may assume that there is a sector
$W''$ of the form $z\iv w_1z$ to the left of $W'$ in $W$ (if not
we replace $W$ by $W\iv$). In particular, $W'$ is not a prefix of
$W$. The $\bkk$-letter preceding $\kappa(1)$ of $W'$ in $W$ is
either $\kappa(1)\iv$ or $\tk(1,\omega)$. In the first case we get
a contradiction using Lemma \ref{lm:badsector1} (ii) and Lemma
\ref{sectors}(iii) because by Lemma \ref{landr} (iii)
$\vartheta_1(ghg\iv)$ is a copy of $ghg\iv$ in this case. In the
second case we get a contradiction with Lemma \ref{lm:]tap}.
\endproof

\begin{lm}\label{rivr}
$W$ does not contain sectors of the form $z\iv wz$, where $z\equiv
\rk(i,j,\omega)$ for $j\ne 1$ or $z\equiv \lk(i,j,\omega)$ or
$z\equiv \tk(i,\omega)$.
\end{lm}

\proof Indeed, suppose that $W$ contains such a sector $W '$. Then
by Lemma \ref{lm:badsector}, rules from
$\sss_{2n+1}\cup\{\bc_{2n}, \bc_{2n}\iv\}$ do not occur in
$ghg\iv$ if $z\equiv \rk(i,j,\omega)$ for $j\ne 1$ or $z\equiv
\lk(i,j,\omega)$ or $z\equiv\tk(i,\omega)$.

Since $W$ starts and ends with the same $\bkk$-letter, we can
assume that $W$ contains a sector $W''$ of the form $z_1w_1z_1\iv$
to the right of $W'$.  By Lemmas \ref{tktk} and \ref{tktk1} the
$\bkk$-letter $z_1$ in $W''$ can be of one of the forms
$\tk(s,\omega), \kappa(s), \bk(s)$ for some $s$. By Lemma
\ref{lm:badsector} this implies that rules from $\sss_0$ and
$\bc_0$ cannot occur in $ghg\iv$. As in the proofs of Lemmas
\ref{tktk}, \ref{tktk1}, all rules in $ghg\iv$ are contained in
$\sss_t\cup\sss_{t+1}\cup\{\bc_t, \bc_t\iv\}$, $t\ne 0, 2n$.

\medskip

{\bf Case 1.} Suppose that the $\bkk$-letter in $W'$ is
$\rk(i,j,\omega)$. Then by Lemma \ref{tktk} the $\bkk$-letter next
to the right of $W'$ cannot be $\rk(i,j,\omega)\iv$, so it must be
$\lk(i,j,\omega)$. But rules from $\sss_t\cup\sss_{t+1}$ are
active for $\rk(i,j)$-sectors and the rule $\bc_t$ locks this
sector. Thus we get a contradiction with Lemma \ref{twosteps}.

\medskip

{\bf Case 2.} Suppose now that $z\equiv \lk(i,j,\omega)$. Then by
Lemma \ref{tktk1}, the $\bkk$-letter next to the right of $W'$
cannot be $\lk(i,j,\omega)\iv$, so it must be $\rk(i,j+1,\omega)$
if $j<n$ or $\tk(i+1,\omega)$ if $j=n$. As in Case 1, rules from
$\sss_t\cup\sss_{t+1}$ are active for $\lk(i,j)$-sectors, and the
rule $\bc_t$ locks these sectors, which contradicts Lemma
\ref{twosteps}.

\medskip

{\bf Case 3.} Suppose that the $\bkk$-letter in $W'$ is
$\tk(i,\omega)$. Since by Lemma \ref{notall} not all sectors of
$W$ have the form $zwz\iv$ where $z\in
\{\tk(i,\omega),\tk(i,\omega)\iv\}$, $W$ either has a sector of
the form $[\tk(i)\kappa(i)]$ or it has sectors of the forms
$[\tk(i)\tk(i)\iv]$ and $[\lk(i-1,n)\tk(i)]$.

\medskip

{\bf Case 3.1.} Suppose that $W$ has a sector of the form
$[\tk(i)\kappa(i)]$. Recall that rules from $\sss_0$, $\bc_0,
\bc_{2n}$ do not occur in $ghg\iv$. All other rules are right
active for $\tk(i)$-sectors. This by Lemma \ref{lm:]tap} implies
that $ghg\iv=1$ modulo Burnside relations, a contradiction.

\medskip

{\bf Case 3.2.} Suppose that $W$ has a sector of the form
$[\tk(i)\tk(i)\iv]$ and a sector of the form $[\lk(i-1,n)\tk(i)]$.
Since we have already proved (Cases 1, 2) that $W$ does not
contain sectors of the form $[\lk(p,q)\iv \lk(p,q)]$,
$[\rk(p,q)\iv\rk(p,q)]$, $q\ne 1$, $W$ must contain sectors of the
forms, $[\rk(i-1,n)\lk(i-1,n)]$, $[\lk(i-1,n-1)\rk(i-1,n)]$,
...,$[\rk(i-1,1)\lk(i-1,1)]$. Thus $W$ contains
$\rk(i-1,1,\omega)\iv$. The $\bkk$-letter next to the right of
$\rk(i-1,1,\omega)\iv$ in $W$ is either $\rk(i-1,1,\omega)$ or
$\bk(i-1)\iv$. Since all rules from
$\sss_t\cup\sss_{t+1}\cup\{\bc_t\}$ are left active for
$\rk(i-1,1)\iv$-sectors, in the first case we get a contradiction
by using Lemma \ref{lm:badsector1}(ii) and Lemma \ref{additional}.
In the second case we get a contradiction with Lemma
\ref{lm:]tap}.
\endproof

\begin{lm}\label{lm:ss0}
If $W$ contains a sector $W'$ of the form $[\rk(i,n)\lk(i,n)]$ and
a sector $W''$ of the form $[\lk(i,n)\tk(i+1)]$, and $h$ contains
a subword $h'$ of the form $\bc_s^{\pm 1 }h_1\bc_s^{\mp 1 }$,
where $h_1$ does not contain either connecting rules or rules from
$\sss_{2n+1}$, then $h_1$ contains rules from $\sss_0$.
\end{lm}

\proof Indeed suppose that $h_1$ does not contain rules from
$\sss_0$. Notice that $\bc_s$ locks either $W'$ or $W''$, but all
rules from $\sss\backslash (\sss_0\cup\sss_{2n+1}\cup\{\bc_0,...,
\bc_{2n}\})$ are left active for $\rk(i,n)$-sectors and for
$\lk(i,n)$-sectors. Hence we can apply Lemma \ref{lm:iji} to
either $W'$ or $W''$, $h_1$ is equal to 1 modulo Burnside
relations which contradicts Lemma \ref{square2}.
\endproof

\begin{lm}\label{lm:sss0}
If $W$ contains a sector of the form $[\rk(i,n)\lk(i,n)]$ and a
sector of the form $[\lk(i,n)\tk(i+1)]$ then $h$ contains a rule
from $\sss_0\cup\{\bc_0\}$ or $h$ is a word over $\sss_{2n+1}$.
\end{lm}

\proof Let $\omega$ be the $\Omega$-coordinate of $W$. Suppose $h$
does not contain rules from $\sss_0$ and does not consist of rules
from $\sss_{2n+1}$.

Then by Lemma \ref{lm:ss0}, $h$ has the form $$h_\omega \bc_\omega
h_{\omega+1}...h_{2n}\bc_{2n}h_{2n+1}\bc_{2n}\iv
h'_{2n}...\bc_{\omega}\iv h'_\omega$$ where $h_j$ and $h'_j$ are
words over $\sss_j$, $j=\omega,...,2n+1$.

The connecting rule $\bc_\omega$ locks either $\rk(i,n)$-sectors
or $\lk(i,n)$-sectors. The rule $\bc_{2n}$ locks both
$\rk(i,n)$-sectors and $\lk(i,n)$-sectors. All rules from
$\sss_\omega$, $\omega<2n$, are left active for
$\lk(i,n)$-sectors, and for $\rk(i,n)$-sectors, rules from
$\sss_{2n}$ are left active for $\lk(i,n)$-sectors. Hence by Lemma
\ref{landr}(iii) $h_\omega=(h'_\omega)\iv$ modulo Burnside
relations. By Lemma \ref{square3}, $h_\omega$ and $h_\omega'$ are
empty. Then $h$ is not cyclically reduced (it starts with
$\bc_\omega$ and ends with $\bc_\omega\iv$, a contradiction).
\endproof

\begin{lm} \label{lm:03}
Suppose $W$ contains a sector of the form $[\rk(i,n)\lk(i,n)]$ and
a sector of the form $[\lk(i,n)\tk(i+1)]$. Then $ghg\iv$ contains
rules from $\sss_0$.
\end{lm}

\proof Suppose by contradiction that $ghg\iv$ does not contain
rules from $\sss_0$. By Lemma \ref{lm:sss0} $h$ is a word over
$\sss_{2n+1}$.

Since $W\cdot h=W$, the $\Omega$-coordinate of $W$ is $2n+1$.


Suppose that $g$ consists of rules from $\sss_{2n+1}$. Suppose
that rules from $\sss_{2n+1}$ are active for a sector of the form
$[zz']$ of $W$. Then $z\nee \kappa(1), \bk(1)\iv$. Hence by Lemma
\ref{no2n+1}, either $g$ and $h$ commute modulo Burnside relations
or $W\cdot g=W$, a contradiction. If rules from $\sss_{2n+1}$ are
not active for sectors of $W$ then $W\cdot g=W$ by Lemma
\ref{landr}(i), a contradiction. Hence $g$ contains some rules not
from $\sss_{2n+1}$.

Since $W$ is in the domain of $g$ (indeed $W$ is in the domain of
$ghg\iv$ and $ghg\iv$ is reduced) , $g$ must start with a rule
from $\sss_{2n+1}\cup\{\bc_{2n}\iv\}$. Since $g$ contains some
rules not from $\sss_{2n+1}$, $g$ has a subword of the form
$\bc_t^{\pm 1}f\bc_t^{\mp 1}$ where $f$ is a word over some
$\sss_{t'}$, $t'\ne 0$. By Lemma \ref{lm:iji}, $g$ is not
Burnside-reduced, a contradiction with Lemma \ref{square2}.
\endproof

\begin{lm}\label{norivr}
If $W$ contains an occurrence of $\rk(i,j,\omega)$ or
$\lk(i,j,\omega)$ or $\tk(i,\omega)$ then $ghg\iv$ contains rules
from $\sss_0$.
\end{lm}

\proof Consider the cyclic graph labeled by $\Lambda(\omega)$
where $\omega$ is the $\Omega$-coordinate of $W$. We can consider
this graph as an ``inverse automaton", so that we can read
clockwise and counterclockwise. Of course when we pass an edge
labeled by $z$ in the counterclockwise direction, we read $z\iv$,
when we pass it in the clockwise direction, we read $z$. Then the
word $z_1z_2...z_1$ of $\bkk$-letters from $W$ can be read on this
automaton, and the start edge coincides with the end edge. This
interpretation makes it clear that if $W$ contains
$\rk(i,j,\omega)$, $\lk(i,j,\omega)$ or $\tk(i,\omega)$ and does
not contain sectors of the form $[zz\iv]$ with
$z\in\{\rk(p,q,\omega), \lk(p,q,\omega), \lk(p,q,\omega)\iv,
\tk(p,\omega)\iv\}$ for any $p,q$ (see Lemmas \ref{tktk},
\ref{tktk1}, \ref{rivr}) then $W$ must contain a sector of the
form $[\rk(p,n)\lk(p,n)]$ and a sector of the form
$[\lk(p,n)\tk(p+1,n)]$. Now the statement of the lemma follows
from Lemma \ref{lm:03}.
\endproof

\begin{lm} \label{nobadsectors}
$W$ does not contain sectors of the form $zwz\iv$ where $z\in
\bkk$.
\end{lm}

\proof Suppose that $W$ contains such a sector. By Lemmas
\ref{tktk}, \ref{tktk1}, \ref{kkiv}, $$z\in \{\tk(i,\omega),
\bk(i)\}.$$

\medskip

{\bf Case 1.} Let $z\equiv \tk(i,\omega)$. Since by Lemma
\ref{rivr} $W$ does not contain sectors of the form
$[\tk(i)\iv\tk(i)]$, $W$ must contain a sector of the form
$[\lk(i-1,n)\tk(i)]$. Since by Lemma \ref{rivr} $W$ does not
contain sectors of the form $[\lk(i-1,n)\iv \lk(i-1,n)]$, $W$ must
contain a sector of the form $[\rk(i-1,n)\lk(i-1,n)]$. Now by
Lemma \ref{lm:03}, $ghg\iv$ contains a rule from $\sss_0$. But
this possibility cannot occur because rules from $\sss_0$ lock
$\tk(i)$-sectors.

\medskip

{\bf Case 2.} Let $z\equiv \bk(i)$. Then $ghg\iv$ does not contain
rules from $\sss_0$ by Lemma \ref{lm:badsector}.

Suppose that $W$ contains a $\bkk$-letter of the form
$\rk(i,j,\omega)$, $\lk(i,j,\omega)$ or $\tk(i,\omega)$. Then as
in the proof of Lemma \ref{norivr}, by Lemmas \ref{tktk},
\ref{tktk1}, \ref{rivr} $W$ must contain sectors of the form
$[\rk(i,n)\lk(i,n)]$ and $[\lk(i,n)\tk(i)]$. But this contradicts
Lemma \ref{lm:03}.

Therefore all $\bkk$-letters in $W$ belong to the set
$\{\kappa(i), \kappa(i)\iv, \bk(i), \bk(i)\iv\}$. By Lemma
\ref{notall} not all $\bkk$-letters belong to the set $\{\bk(i),
\bk(i)\iv\}$. By Lemma \ref{kkiv}, not all of them belong to
$\{\kappa(i), \kappa(i)\iv\}$. Therefore for some $i$ from 1 to
$N$ each sector of $W$ has one of the
 four forms $[\kappa(i)\bk(i)]$, $[\kappa(i)\iv \kappa(i)]$,
$[\bk(i)\bk(i)\iv]$, $[\bk(i)\iv \bk(i)]$ (recall that sectors of
the form $[\kappa(i)\kappa(i)\iv]$ cannot appear by Lemma
\ref{kkiv}). Moreover a sector of the form $[\kappa(i)\iv
\kappa(i)]$ must appear among sectors of $W$ (the $\bkk$-letter
preceding $\kappa(i)$ in $W$ cannot be $\tk(i, \omega)$ since $W$
does not contain letters of the form $\tk(j,\omega)$).

If rules from $\sss_{2n+1}$ do not appear in $g, h$ then we can
consider any sector of $W$ of the form [$\kappa(i)\iv \kappa(i)$],
and use Lemma \ref{lm:badsector1} to get a contradiction.
Therefore rules from $\sss_{2n+1}$ do appear in $ghg\iv$. By the
same reason not all rules in $ghg\iv$ are from $\sss_{2n+1}$.
Therefore $ghg\iv$ contains $\bc_{2n}$. Since $\bc_{2n}$ locks
$\kappa(i)$-sectors, sectors of the form [$\bk(i)\iv\bk(i)$] do
not appear in $W$. Therefore for some $i$ all sectors of $W$ have
the form $[\kappa(i)\iv \kappa(i)]$, $[\kappa(i)\bk(i)]$ or
$[\bk(i)\bk(i)\iv]$. But this contradicts the assumption of part
(ii) of Proposition \ref{propsss}.
\endproof

Recall that $W$ is an admissible word which starts and ends with
the same letter and does not contain accepted subwords. By
replacing $W$ with $W\iv$ if necessary we can assume that the
first letter in $W$ belongs to $\kk(\omega)$ for some $\omega\in
\Omega$. Then by Lemma \ref{nobadsectors} all $\bkk$-letters in
$W$ belong to $\kk(\omega)$. Since $W$ is an admissible word, it
must contain an admissible subword $W'$ whose projections onto
$\bkk$-letters is equal to a cyclic shift of $\Lambda(\omega)$.

\begin{lm}\label{final1}
$W'$ is an accepted word.
\end{lm}

\proof By definition, $W$ contains sectors of the form $[zz_+]$
for every $z\in \kk(\omega)$. In particular, $W$ contains sectors
of the form $[\rk(1,n)\lk(1,n)]$ and $[\lk(1,n)\tk(1,n)]$. By
Lemma \ref{lm:03}, $ghg\iv$ contains a rule from $\sss_0$. By
Lemma \ref{no2n+1}, $ghg\iv$ contains rules from at least two
different $\sss_i$.   Hence $ghg^{-1}$ also contains either rules
from ${\bf S}_1$ or ${\bf c}_0^{\pm 1}$.

Let us prove that $ghg^{-1}$ or $gh^{-1}g^{-1}$ contains a subword
${\bf c}_0 h_1 c_1$ where $h_1$ is a word over ${\bf S}_1$.
Suppose that it is not true. Then $ghg\iv$ contains rules from
$\sss_0\cup\{\bc_0, \bc_0\iv\} \cup\sss_1$ only. Hence $ghg\iv$
contains $\bc_0$ and $\bc_0\iv$. Since the $\Omega$-coordinates of
$W$, $W\cdot g$, $W\cdot gh$ are the same, the total exponent of
$\bc_0$ in $g$ is 0  and the total exponent of $\bc_0$  in $h$ is
0. If $g$ or $h$ contains a subword $\bc_0f\bc_0\iv$ where $f$ is
a word over $\sss_1$, then by Lemma \ref{lm:]tap}, applied to any
$\lk(i,j)$-sector of $W$, we get that $f=1$ modulo Burnside
relations, a contradiction with Lemma \ref{square2}. Hence neither
$g$ nor $h$ contain subwords of the form $\bc_0 f\bc_0\iv$. Now
consider two cases.

\medskip

{\bf Case 1.} Suppose that $\omega=0$. Since $W\cdot ghg\iv=W$,
the first occurrence of $\bc_0^{\pm 1}$ in $ghg\iv$ is $\bc_0$.
Since $h$ does not contain subwords of the form $\bc_0f\bc_0\iv$,
this $\bc_0$ occurs in $g$. But since $W\cdot g$ has the same
$\Omega$-coordinate as $W$, $g$ must contain $\bc_0\iv$, hence $g$
contains a subword of the form $\bc_0f\bc_0\iv$, a contradiction.

\medskip

{\bf Case 2.} Suppose that $\omega=1$. Then the first occurrence
of $\bc_0^{\pm 1}$ in $ghg\iv$ is $\bc_0\iv$. Since $W\cdot g=W$,
we can represent $g$ as $g_1g_2$ where $g_2$ is the maximal suffix
of $g$ consisting of rules from $\sss_1$ and $g_1$ is either empty
or ends with $\bc_0$. Since $W\cdot h=W$, we can represent $h$ as
$h_1h_2h_3$ where $h_1$ and $h_3$ consist of rules from $\sss_1$
and either $h_2=h_3$ is empty or $h_2$ starts with $\bc_0\iv$ and
ends with $\bc_0$. Consider any $\lk(1,1)$-sector $U\equiv
\lk(1,1,1)u\rk(1,2,1)$ of $W$.

{\bf Case 2.1.} Suppose that $g_1$ is empty, $g\equiv g_2$
consists of rules from $\sss_1$. Since $ghg\iv$ contains $\bc_0$,
$h$ must contain $\bc_0$, so $h_2$ is not empty.

Since $\bc_0$ locks $U$ and rules from $g$ are left active for
$\lk(1,1)$-sectors, by Lemma \ref{landr} (i), we have
$R^a(gh_1)=u=(R^a(h_3g\iv))\iv$ modulo Burnside relations. By
Lemma \ref{landr} (iii), we have that $gh_1=gh_3\iv$ modulo
Burnside relations. Hence $h_1=h_3\iv$ modulo Burnside relations.
This contradicts Lemma \ref{square3} if $h_1$ or $h_3$ is not
empty or Lemma \ref{cyc} if $h_1$, $h_3$ are empty.

{\bf Case 2.2.} Suppose that $g_1$ is not empty. Then $g_1$ ends
with $\bc_0\iv$. If $h$ does not contain $\bc_0$, that is
$h_2,h_3$ are empty, then by Lemma \ref{lm:]tap} applied to a
$\rk(1,1)$-sector of $W$, $g_2hg_2\iv=1$ modulo Burnside
relations, so $h=1$ modulo Burnside relations, a contradiction.

If $h$ contains $\bc_0$, that is $h_2$ is not empty, then again by
Lemma \ref{lm:]tap} $g_2h_1=1$, $h_3g_2\iv=1$ modulo Burnside
relations. Thus $h_1=h_3\iv$ modulo Burnside relations which as
before contradicts Lemma \ref{square3} or Lemma \ref{cyc}.

\medskip

Thus, since one can always consider $h^{-1}$ instead of $h$, we
can assume that $h$ has the form $e{\bf c}_0h_1{\bf c}_1 f$ for
some words $e,f$. Since $\bc_0$ locks $z$-sectors, $z\in
\{\tk(i),\kappa(i), \bk(i), \lk(i,j)\}$, a cyclic shift of
$W'\cdot e$ has the form

$$\tk(1,0)\kappa(1)\bk(1)\rk(1,1,0)u_{1,1}\lk(1,1,0)\rk(1,2,0)u_{1,2}\lk(1,2,0)...\rk(N,n,0)u_{N,n}\lk(N,n,0).$$
The word $W'\cdot e\bc_0$ has the form

$$\tk(1,1)\kappa(1)\bk(1)\rk(1,1,1)u_{1,1}\lk(1,1,1)\rk(1,2,1)u_{1,2}\lk(1,2,1)...\rk(N,n,1)u_{N,n}\lk(N,n,1)$$
Since $\bc_1$ locks $\rk(i,j)$-sectors, for every $i\in
\{1,...,N\}$ and every $j\in \{1,...,n\}$, we have
$$R^a(h_1,\rk(i,j))=u_{i,j}$$ modulo Burnside relations. Hence all
$u_{i,j}$ are equal modulo Burnside relations, so a cyclic shift
of $W'\cdot e$ is equal to the word
\begin{equation}\label{eq:57}
\tk(1,0)\kappa(1,0)\bk(1,0)\rk(1,1,0)u\lk(1,1,0)\rk(1,2,0)u\lk(1,2,0)...\rk(N,n,0)u\lk(N,n,0)
\end{equation}
for some word $u$ which is equal to $u_{i,j}$ modulo Burnside
relations. Denote by  $\bar u$ the copy of $u$ written in
$\sss_0$. Then by Lemma \ref{landr} $W'\cdot e\bar u=\Lambda(0)$,
so $W'$ is accepted.
\endproof

Lemma \ref{final1} gives us the final contradiction because by our
assumption $W$ cannot contain an accepted subword. Proposition
\ref{propsss} is proved.

\section{Bands, annuli, and trapezia}
\label{bands}

\subsection{Definitions and basic facts}

\label{dbf}

In Section \ref{properties}, we introduced an auxiliary group
$\hnka$ which was the free product of the subgroup $\la\bkk\ra$
freely generated by $\bkk$ and the subgroup $\maa$ generated by
the set $\aaa$ subject to the $R(\kappa(1))$-relations and the
Burnside relations. Recall that the subgroup $\maa$ is the free
product in the variety of Burnside groups of exponent $n$ of its
subgroups $\la A(z)\ra$, $z\in \kk$. From now on let us call all
relations of $\maa$ as \label{aarel}$a$-{\em relations}. All
relations from $Z(\sss)$ involving letters from $\bkk$ will be
called \label{kkrel}$\kappa$-{\em relations} and all relations
involving letters from $\rr$ will be called \label{rrrel}$r$-{\em
relations}. We also will call the commutativity relations from
$Z(\sss)$ involving $\rr$- and $\aaa$-letters
\label{rarel}$ra$-{\em relations}. Accordingly we shall consider
$a$-cells, $\kappa$-cells, $r$-cells and $ra$-cells in \vk
diagrams. Notice that $\kappa$-cells are also $r$-cells. Recall
that the set of relations $Z(\sss,\Lambda)$ consists of $Z(\sss)$
and the hub $\Lambda(0)$. The cells corresponding the hub
relation, will be called {\em hubs} as well.

In this section, we shall consider diagrams over the set
\label{bra}$Z(\sss,a)$ of all $a$-relations and all $r$-relations
(including $\kappa$-relations).

With every word over $\bkk\cup\rr\cup\aaa$ we associate its
$\bkk$-length, its $\rr$-length and its $\aaa$-length in the
natural way. In any diagram over $Z(\sss,a)$ we consider
$\rr$-edges and $\aaa$-edges as $0$-edges and $\bkk$-edges as
non-zero edges. In a diagram is over the presentation consisting
of $a$-relations and $ra$-relations, then we consider $\rr$-edges
as non-zero edges and $\aaa$-edges as zero edges. Thus the metric
on a diagram depends on over what presentation this diagram is.

Let $S$ be one of the sets $\bkk, \rr, \aaa$. Then we can consider
$S$-bands in any diagram over $Z(\sss,a)$ as in Section
\ref{prelim}.

\begin{itemize}
\item \label{bkkband}$\bkk$-bands consist of $\kappa$-cells. A maximal $\bkk$-band
which is not an annulus can start and end on the boundary of a
diagram. All $\bkk$-letters of a $\bkk$-band belong to the same
set $\kk(z)$ for some $z\in \kk\cup\kk\iv$. Thus we can consider
\label{kkband}$\kk(z)$-bands as well.
\item \label{rrband}$\rr$-bands consist of $r$-cells. A maximal $\rr$-band
which is not an annulus must start and end on the boundary of the
diagram. For every $\rr$-band $\ttt$ there exists a rule
$\tau\in\sss$ such that $\rr$-edges in an $\rr$-band have labels
of the form $r(\tau,z)$, $z\in \kk$. Thus we will sometimes call
$\ttt$ a \label{tauband}$\tau$-band.
\end{itemize}

The next lemma immediately follows from the definition of words
$L(.,z)$ and $R(., z)$ and the definition of $\kk$-relations.

\begin{lm} \label{kbands}
Consider a disc or an annular diagram consisting of one
$\kk(z)$-band $\ttt$. Then for some word $h$ over $\sss$ (which is
called the \label{hb}{\em history} of the band), the word written
on the top of $\ttt$ is $R(h,z)$, the word written on the bottom
of $\ttt$ is $L(h,z)$. The word $R(h,z)$ (resp. $L(h,z)$) is
$\rr$-reduced if and only if $\ttt$ does not contain two
consecutive cells that cancel. Moreover $L(h,z)z'=zR(h,z)$ modulo
$\kappa$-relations, where $z,z'$ are the labels of the start and
end edges of the band.
\end{lm}

The next lemma contains easy properties of $\rr$-bands which
follow immediately from the definition of $Z(\sss)$ and Lemma
\ref{lm1}.

\begin{lm} \label{rbands}
Let $$w\equiv u_0z_1u_1z_2...z_su_{s}$$ be the bottom side of a
$\tau$-band $\ttt$, $$w'\equiv
u_0'z_1'u_1'z_2'...z'_{s'}u'_{s'}$$ be the top side of $\ttt$,
where $z_i, z_i'\in \overline{\bkk}\cup(\overline{\bkk})\iv$ $u_i,
u_i'$ are words over $\aaa$. Then $W\equiv z_1u_1z_2...z_s$ and
$W'\equiv z_1'u_1'z_2'...z'_{s'}$ are admissible words and $W\cdot
\tau = W'$.
\end{lm}

\subsection{Forbidden annuli}
\label{nohubs}

Let $\Delta$ be any (disc or annular) diagram over $Z(\sss,a)$. We
call $\Delta$ \label{zsa}{\em $Z(\mathbb{S},a)$-reduced} if every
$\bkk$-band is reduced (that is they do not contain consecutive
$\bkk$-cells that cancel) and every disc subdiagram whose boundary
label is a word over $\aaa$ and is equal to 1 modulo $a$-relations
does not contain $\rr$-edges.

\begin{lm}\label{bandsnohubs}
Let $\Delta$ be a $Z(\mathbb{S},a)$-reduced disc diagram. Then
$\Delta$ does not contain

\begin{enumerate}
\item \label{b1}$\bkk$-annuli,
\item \label{b3}$(\bkk,\rr)$-annuli,
\item \label{b4}$\rr$-annuli,
\item \label{b5} disc subdiagrams with some $r$-cells whose boundary path
contains no $\rr$-edges,

\end{enumerate}
\end{lm}

\proof We are proving these statements by a joint induction on the
number of cells of the inside diagrams of the annulus or the
``bad" subdiagrams. By contradiction take the smallest (with
respect to the number of cells) disc subdiagram $\Delta'$ of
$\Delta$ which either contains one of the annuli forbidden by the
lemma or has boundary path forbidden by the lemma.

\medskip

{\bf Case \ref{b1}.} Assume there is a $\bf K$-annulus $\cal Q$ in
$\Delta$. Then there exists an $\bf R$-edge on a side of $\cal Q$.
Hence $\cal Q$ crosses an $\bf R$-band $\cal T$. Since $\cal Q$
and $\cal T$ form a $(\bf K,\bf R)$-annulus with the inside
diagram smaller than that for $\cal Q$, we obtain a contradiction.

\medskip

{\bf Case \ref{b3}.} Suppose that $\Delta''$ is the inside diagram
of a $(\bkk,\rr)$-annulus ${\cal W}_1\cup {\cal W}_2$ where ${\cal
W}_1$ is a $\bkk$-band and ${\cal W}_2$ is an $\rr$-band.

Then the contour of $\Delta''$ does not contain $\kk$-edges
(otherwise there would be a $(\bkk, \rr)$-annulus with a smaller
inside diagram). Therefore ${\cal W}_1$ consists of two cells.
These two cells have a common $\bkk$-edge and they belong to the
same $r$-band ${\cal W}_2$. Therefore these two $\kappa$-cells
cancel, a contradiction.

\medskip

{\bf Case \ref{b4}.} Suppose that there is a $\kappa$-cell in an
$\bf R$-annulus $\cal Q$. Then we get a contradiction as in Case
1. Therefore the boundary label $W$ of $\Delta''$ is a word in
$\bf A$. Since the counterexample is minimal, $\Delta''$ has no
edges except for $\bf A$-edges. Hence $W=1$ modulo $a$-relations.
Diagram $\Delta'$ has the same boundary label $W$, since all
$(r,a)$-cells in $\cal Q$, are commutativity cells. But this
contradicts the property that $\Delta'$ is $Z(\bf S,a)$-reduced,
because there are $\bf R$-edges in $\cal Q.$

\medskip

{\bf Case \ref{b5}.} Suppose that $\partial(\Delta')$ contains no
$\rr$-edges. Then by part \ref{b4} of the lemma, $\Delta'$
contains no $\rr$-edges, so it cannot contain $r$-cells, a
contradiction.
\endproof

\begin{cy} \label{hnka}
The group $\hnka$ is (naturally) embedded into the group defined
by $Z(\sss,a)$.
\end{cy}

\proof Indeed, let $W$ be a word in $\bkk\cup\aaa$ which is equal
to 1 modulo $Z(\sss,a)$. Then there exists a
$Z(\mathbb{S},a)$-reduced \vk disc diagram $\Delta$ with boundary
label $W$. By Lemma \ref{bandsnohubs}, part \ref{b5}, $\Delta$
contains no $r$-cells, so all cells in $\Delta$ are $a$-cells.
Hence $W=1$ modulo $a$-relations, so $W=1$ in $\hnka$. Thus the
natural map of $\hnka$ into the group given by the relations
$Z(\sss,a)$ is an embedding.
\endproof

\subsection{Trapezia}
\label{trapezia} Here we shall present a tool to translate results
from Section \ref{properties} into the language of \vk diagrams
and vice versa.

\medskip

Let $\Delta$ be a $Z(\mathbb{S},a)$-reduced disc diagram which has
the contour of the form $p_1q_1p_2\iv q_2\iv$ where:\label{TR12}

(TR1) $\Lab(q_1)$, $\Lab(q_2)$ are $\bkk$-reduced words without
$\rr$-letters, starting and ending with $\bkk$-letters, and

(TR2) $p_1$ and $p_2$ are sides of reduced $\bkk$-bands each
containing at least one $\kappa$-cell.

Then $\Delta$ is called a \label{tra}{\em trapezium}. The path
$q_1$ is called the \label{base}{\em bottom base}, the path $q_2$
is called the {\em top base} of the trapezium, the paths $p_1$ and
$p_2$ are called the {\em left and right sides} of the trapezium.
The history of the $\bkk$-band whose side is $p_1$ is called the
\label{ht}{\em history} of the trapezium.

\medskip

By Lemma \ref{lm1}, if a word $h$ over $\sss$ contains just one
rule $\tau$, $W\equiv z_1u_1...z_t$ is a reduced admissible word,
and $W\cdot h=W_1$, then there exists a word $W'$ which is equal
to $W$ modulo $a$-relations and such that $\tau$ is applicable to
$W'$, $W'\cdot h=W_1'$ where $W_1'$ is equal to $W_1$ modulo
$a$-relations. By Lemma \ref{lm1} there exists a trapezium
$\Theta$ consisting of one $\tau$-band with bottom base label
$W'$, top base label $W_1'$. Since $W=W'$, $W_1=W_1'$ modulo
Burnside relations, there exist diagrams $\Delta, \Delta_1$
consisting of $a$-cells with boundary labels $W(W')\iv$ and
$W_1(W_1')\iv$ respectively. Gluing three diagrams $\Delta$,
$\Theta$, $\Delta_1$ in a natural way, we get a trapezium
$T(W,\tau)$ with bottom base label $W$, top base label $W_1$, left
side label $L(\tau,z_1)$, right side label $R(\tau,z_t)$.

If a reduced word $h$ consists of several rules, the corresponding
trapezium can be obtained by induction. Let $h=h'\tau$ where $h'$
is a word of smaller length than $h$. Suppose that we have
constructed a trapezia $T(W,h')$ with the bottom base label $W$
and the top base label $W\cdot h'$, and suppose that $W\cdot h$
exists. Then $W\cdot h'$ is in the domain of $\tau$, so $\tau$ is
applicable to a reduced admissible word $W_1$ which is equal to
$W\cdot h'$ modulo $a$-relations. Then there exists a diagram
$\Delta$ over the $a$-relations with contour $s_1s_2\iv$ where
$\Lab(s_1)=W\cdot h', \Lab(s_2)=W_1$. By Lemma \ref{lm1} there
exists a trapezium $T(W_1,\tau)$ with the bottom base label $W_1$.
Glue $T(W,h')$ with the diagram $\Delta$ identifying the top base
of $T(W,h')$ with $s_1$. Then glue the resulting diagram with
$T(W_1,\tau)$ identifying $s_2$ with the bottom base of
$T(W_1,\tau)$. The $\bkk$-bands in the resulting diagram $T(W,h)$
are reduced because $h$ is a reduced word. Every closed path in
$T(W,h)$ consisting of $\aaa$-edges is contained in one of the
subdiagrams $\Delta$ consisting of $a$-cells. Hence $T(W,h)$ is a
$Z(\mathbb{S},a)$-reduced diagram. It is clear that conditions
(TR1) and (TR2) are satisfied. Hence $T(W,h)$ is a trapezium with
bottom base label $W$, top base label $W\cdot h$, left side label
$L(h,z_1)$ and right side label $R(h,z_t)$.  It has a form of a
sandwich where diagrams consisting of $a$-cells are sandwiched
between $\tau$-bands (see \setcounter{pdtwelve}{\value{ppp}}
Figure \thepdtwelve).

Clearly $T(W,h)$ is not determined uniquely by $W$ and $h$ because
we can choose the $a$-parts of the sandwich in many different
ways. Any trapezium of the form $T(W,h)$ is called a
\label{normal}{\em normal trapezium with history $h$ and base
label $W$}.

\unitlength=1.00mm \special{em:linewidth 0.4pt}
\linethickness{0.4pt}
\begin{picture}(74.22,53.47)
\put(29.78,44.89){\line(1,-6){6.67}}
\put(36.44,4.89){\line(1,0){25.78}}
\put(62.22,4.89){\line(1,5){8.00}}
\put(70.22,44.89){\line(-1,0){40.67}}
\bezier{104}(36.44,4.89)(48.89,8.00)(62.00,4.89)
\bezier{116}(35.33,12.00)(48.89,8.00)(63.78,12.00)
\bezier{116}(35.33,12.00)(52.00,16.00)(63.78,12.00)
\bezier{128}(34.00,20.22)(50.00,17.11)(65.33,20.22)
\bezier{136}(33.78,20.22)(52.00,26.44)(65.33,20.22)
\bezier{140}(32.22,29.78)(50.00,25.78)(67.11,29.78)
\bezier{144}(32.22,29.78)(54.44,34.00)(67.11,29.78)
\bezier{152}(31.33,36.89)(50.67,33.78)(68.67,36.89)
\bezier{152}(31.11,37.11)(54.44,39.33)(68.44,37.11)
\bezier{168}(29.78,44.89)(56.00,39.33)(70.22,44.89)
\put(30.67,44.89){\line(-1,0){4.44}}
\put(26.22,44.89){\line(1,-6){6.67}}
\put(32.89,4.89){\line(1,0){3.78}}
\put(61.56,4.89){\line(1,0){4.44}}
\put(66.22,4.89){\line(1,5){8.00}}
\put(74.22,44.89){\line(-1,0){4.67}}
\put(31.11,36.89){\line(-1,0){3.56}}
\put(28.89,29.56){\line(1,0){3.56}}
\put(30.44,20.22){\line(1,0){3.56}}
\put(31.78,12.00){\line(1,0){3.56}}
\put(68.44,36.89){\line(1,0){4.00}}
\put(67.11,29.56){\line(1,0){4.00}}
\put(65.33,20.22){\line(1,0){4.00}}
\put(63.56,12.00){\line(1,0){4.00}}
\put(47.33,12.00){\circle*{2.22}}
\put(49.11,12.00){\circle*{3.11}}
\put(52.00,12.00){\circle*{2.67}}
\put(50.67,12.00){\circle*{2.67}}
\put(45.56,12.00){\circle*{2.67}}
\put(42.67,12.00){\circle*{1.60}}
\put(44.00,11.56){\circle*{1.78}}
\put(43.78,12.44){\circle*{1.33}}
\put(41.11,11.78){\circle*{1.33}}
\put(41.56,12.00){\circle*{1.78}}
\put(40.00,12.00){\circle*{1.33}}
\put(40.67,12.22){\circle*{1.33}}
\put(38.44,11.78){\circle*{0.89}}
\put(39.11,12.00){\circle*{1.33}}
\put(37.11,12.00){\circle*{0.89}}
\put(37.78,12.00){\circle*{1.33}}
\put(38.22,12.22){\circle*{1.33}}
\put(41.33,12.44){\circle*{1.78}}
\put(39.56,12.67){\circle*{0.89}}
\put(55.20,12.00){\circle*{0.00}}
\put(55.20,12.00){\circle*{3.20}}
\put(53.33,11.20){\circle*{1.89}}
\put(53.47,12.00){\circle*{3.30}}
\put(57.47,12.00){\circle*{2.20}}
\put(59.33,11.87){\circle*{1.89}}
\put(61.20,11.87){\circle*{1.33}}
\put(62.40,12.00){\circle*{0.53}}
\put(61.73,11.87){\circle*{1.19}}
\put(60.27,11.73){\circle*{1.71}}
\put(58.67,12.67){\circle*{0.96}}
\put(56.80,12.80){\circle*{1.33}}
\put(56.27,13.20){\circle*{1.10}}
\put(54.40,13.07){\circle*{1.62}}
\put(51.73,13.20){\circle*{1.36}}
\put(52.67,13.47){\circle*{0.53}}
\put(52.67,13.47){\circle*{1.19}}
\put(46.80,12.93){\circle*{1.33}}
\put(48.00,13.20){\circle*{0.60}}
\put(47.47,13.33){\circle*{0.96}}
\put(48.13,13.20){\circle*{1.51}}
\put(49.07,13.47){\circle*{1.13}}
\put(50.40,13.20){\circle*{1.44}}
\put(49.60,13.47){\circle*{1.51}}
\put(51.07,13.33){\circle*{0.96}}
\put(44.67,12.80){\circle*{1.44}}
\put(44.80,13.20){\circle*{0.75}}
\put(45.47,13.33){\circle*{0.53}}
\put(45.87,13.33){\circle*{1.44}}
\put(46.67,13.60){\circle*{0.75}}
\put(42.53,12.67){\circle*{1.07}}
\put(42.93,12.80){\circle*{1.92}}
\put(43.73,12.93){\circle*{1.55}}
\put(40.40,12.53){\circle*{1.44}}
\put(39.07,12.53){\circle*{1.36}}
\put(47.20,10.27){\circle*{0.84}}
\put(46.53,10.27){\circle*{1.07}}
\put(45.33,10.40){\circle*{0.84}}
\put(45.87,10.27){\circle*{0.96}}
\put(44.67,10.40){\circle*{0.80}}
\put(43.47,10.27){\circle*{1.10}}
\put(42.53,10.67){\circle*{0.80}}
\put(41.47,10.67){\circle*{0.84}}
\put(40.40,10.67){\circle*{0.80}}
\put(39.60,10.93){\circle*{0.53}}
\put(39.20,10.93){\circle*{0.60}}
\put(39.73,10.93){\circle*{0.84}}
\put(40.40,10.93){\circle*{1.13}}
\put(42.40,10.93){\circle*{0.84}}
\put(43.20,10.80){\circle*{0.60}}
\put(47.87,10.53){\circle*{0.80}}
\put(46.93,10.80){\circle*{0.53}}
\put(47.20,10.67){\circle*{0.96}}
\put(46.80,10.80){\circle*{0.53}}
\put(50.67,10.27){\circle*{0.96}}
\put(50.53,10.53){\circle*{1.33}}
\put(48.80,10.13){\circle*{0.80}}
\put(48.27,10.00){\circle*{0.53}}
\put(56.80,10.67){\circle*{0.80}}
\put(58.40,10.93){\circle*{0.75}}
\put(57.60,10.67){\circle*{0.60}}
\put(51.87,10.13){\circle*{1.13}}
\put(52.67,10.27){\circle*{0.60}} \put(51.33,9.87){\circle*{0.60}}
\put(49.60,10.13){\circle*{0.38}}
\put(59.73,12.67){\circle*{0.84}}
\put(59.47,12.93){\circle*{0.60}}
\put(58.13,13.07){\circle*{0.96}}
\put(57.20,13.07){\circle*{1.13}}
\put(61.07,12.53){\circle*{0.53}}
\put(60.53,12.40){\circle*{0.84}}
\put(60.27,12.80){\circle*{0.53}}
\put(62.93,11.87){\circle*{0.53}}
\put(62.40,12.13){\circle*{0.27}}
\put(62.27,12.27){\circle*{0.27}}
\put(42.80,10.93){\circle*{0.96}}
\put(36.13,12.00){\circle*{0.38}}
\put(49.33,20.53){\circle*{4.06}}
\put(53.07,20.13){\circle*{2.53}} \put(52.80,20.80){\circle{4.34}}
\put(52.67,20.67){\circle*{4.14}}
\put(55.73,20.67){\circle*{3.51}}
\put(58.53,20.53){\circle*{2.98}}
\put(60.67,20.67){\circle*{2.20}}
\put(62.40,20.27){\circle*{2.41}}
\put(64.00,20.27){\circle*{0.84}}
\put(46.27,20.27){\circle*{3.94}}
\put(42.80,20.67){\circle*{3.30}}
\put(40.40,20.53){\circle*{2.67}}
\put(38.40,20.53){\circle*{1.92}}
\put(36.93,20.40){\circle*{1.36}}
\put(35.60,20.27){\circle*{1.13}}
\put(34.67,20.27){\circle*{0.60}}
\put(44.40,21.87){\circle*{1.71}}
\put(46.27,22.40){\circle*{0.60}}
\put(47.33,21.73){\circle*{1.71}}
\put(48.13,22.27){\circle*{1.36}}
\put(50.27,22.00){\circle*{1.71}}
\put(48.93,22.53){\circle*{1.19}}
\put(49.73,22.40){\circle*{1.36}}
\put(50.80,22.13){\circle*{1.71}}
\put(51.47,22.67){\circle*{0.96}}
\put(50.27,22.67){\circle*{1.13}}
\put(57.33,21.47){\circle*{1.94}}
\put(56.13,22.00){\circle*{1.44}}
\put(54.93,22.13){\circle*{1.71}}
\put(54.00,22.53){\circle*{1.07}}
\put(53.47,22.80){\circle*{0.53}}
\put(52.00,22.80){\circle*{0.96}}
\put(53.07,22.93){\circle*{0.80}}
\put(52.53,22.93){\circle*{0.80}}
\put(46.80,22.40){\circle*{1.07}}
\put(47.47,22.40){\circle*{1.33}}
\put(45.47,22.00){\circle*{1.60}}
\put(46.00,22.27){\circle*{1.55}}
\put(44.80,22.13){\circle*{1.36}}
\put(59.87,21.20){\circle*{1.36}}
\put(59.20,21.73){\circle*{0.84}}
\put(58.53,22.00){\circle*{1.10}}
\put(57.47,22.13){\circle*{0.96}}
\put(56.93,22.13){\circle*{1.44}}
\put(55.60,22.40){\circle*{1.36}}
\put(56.53,22.67){\circle*{1.10}}
\put(60.27,21.33){\circle*{1.44}}
\put(59.73,21.73){\circle*{1.10}}
\put(54.53,19.33){\circle*{1.07}}
\put(53.60,19.20){\circle*{1.55}}
\put(51.07,19.20){\circle*{0.80}}
\put(44.13,19.20){\circle*{0.84}}
\put(45.20,19.73){\circle*{1.89}}
\put(45.07,20.53){\circle*{1.19}}
\put(39.33,21.07){\circle*{1.44}}
\put(41.33,21.47){\circle*{1.36}}
\put(41.73,21.87){\circle*{1.10}}
\put(43.33,22.00){\circle*{1.79}}
\put(42.27,22.27){\circle*{0.60}}
\put(40.27,21.60){\circle*{0.96}}
\put(40.67,22.00){\circle*{0.53}}
\put(57.20,19.87){\circle*{1.79}}
\put(60.00,19.87){\circle*{1.07}}
\put(61.33,19.87){\circle*{1.44}}
\put(38.93,19.73){\circle*{0.84}}
\put(41.47,19.73){\circle*{1.33}}
\put(47.87,19.20){\circle*{1.36}}
\put(48.80,22.93){\circle*{0.53}}
\put(48.27,22.80){\circle*{0.96}}
\put(46.40,22.67){\circle*{1.10}}
\put(44.67,22.53){\circle*{1.19}}
\put(43.87,22.53){\circle*{0.80}}
\put(54.53,22.67){\circle*{1.07}}
\put(49.20,29.60){\circle*{4.44}}
\put(53.07,29.73){\circle*{4.01}}
\put(56.67,29.73){\circle*{3.48}}
\put(59.33,29.73){\circle*{3.30}}
\put(61.73,29.73){\circle*{2.46}}
\put(63.47,29.73){\circle*{1.36}}
\put(64.67,29.73){\circle*{1.36}}
\put(65.60,29.87){\circle*{0.84}}
\put(66.00,29.87){\circle*{0.60}}
\put(46.00,29.60){\circle*{3.34}}
\put(43.07,29.73){\circle*{3.22}}
\put(40.67,29.60){\circle*{2.46}}
\put(38.53,29.73){\circle*{1.94}}
\put(36.80,29.73){\circle*{1.94}}
\put(35.33,29.73){\circle*{1.13}}
\put(34.27,29.73){\circle*{0.96}}
\put(33.33,29.60){\circle*{0.60}}
\put(51.20,30.53){\circle*{2.53}}
\put(50.67,31.33){\circle*{1.10}}
\put(47.20,30.53){\circle*{2.13}}
\put(44.53,30.53){\circle*{1.60}}
\put(45.60,30.80){\circle*{1.44}}
\put(46.00,31.20){\circle*{0.80}}
\put(46.80,31.07){\circle*{1.44}}
\put(47.33,31.20){\circle*{1.07}}
\put(55.20,30.40){\circle*{2.68}}
\put(54.27,31.07){\circle*{1.62}}
\put(52.67,31.33){\circle*{1.19}}
\put(52.00,31.20){\circle*{1.33}}
\put(41.47,30.67){\circle*{0.80}}
\put(40.40,30.40){\circle*{1.51}}
\put(39.47,29.87){\circle*{1.89}}
\put(37.87,29.87){\circle*{1.62}}
\put(35.73,30.00){\circle*{0.80}}
\put(44.40,29.07){\circle*{1.87}}
\put(47.33,28.40){\circle*{1.36}}
\put(51.07,28.93){\circle*{2.46}}
\put(54.80,28.93){\circle*{2.03}}
\put(62.93,30.00){\circle*{1.87}}
\put(60.40,30.53){\circle*{1.55}}
\put(57.60,30.93){\circle*{1.44}}
\put(56.40,31.20){\circle*{1.10}}
\put(58.27,30.93){\circle*{1.33}}
\put(64.13,30.00){\circle*{1.33}}
\put(39.20,29.33){\circle*{1.71}}
\put(41.87,28.93){\circle*{1.33}}
\put(43.87,28.27){\circle*{0.53}}
\put(44.67,28.27){\circle*{0.84}}
\put(44.53,31.33){\circle*{0.75}}
\put(43.73,31.20){\circle*{0.53}}
\put(42.00,30.67){\circle*{1.44}}
\put(40.93,30.80){\circle*{1.07}}
\put(49.73,36.80){\circle*{2.88}}
\put(52.40,36.93){\circle*{2.53}}
\put(54.93,36.80){\circle*{2.39}}
\put(57.20,36.80){\circle*{2.46}}
\put(59.47,36.53){\circle*{2.46}}
\put(61.47,36.80){\circle*{1.62}}
\put(63.20,36.80){\circle*{1.55}}
\put(64.67,37.07){\circle*{1.19}}
\put(65.73,36.93){\circle*{0.84}}
\put(66.53,36.80){\circle*{1.33}}
\put(64.00,36.40){\circle*{0.53}}
\put(64.53,36.67){\circle*{0.53}}
\put(65.47,36.53){\circle*{0.27}}
\put(65.07,37.07){\circle*{1.33}}
\put(46.67,36.93){\circle*{2.20}}
\put(44.67,36.80){\circle*{2.75}}
\put(42.67,36.53){\circle*{2.08}}
\put(40.53,36.93){\circle*{1.94}}
\put(38.80,36.80){\circle*{1.55}}
\put(37.33,36.93){\circle*{1.79}}
\put(35.73,36.93){\circle*{1.13}}
\put(34.53,36.67){\circle*{1.55}}
\put(33.20,37.07){\circle*{1.10}}
\put(32.27,36.93){\circle*{0.84}}
\put(41.47,36.53){\circle*{1.36}}
\put(41.47,36.53){\circle*{1.94}}
\put(42.67,37.33){\circle*{1.10}}
\put(42.27,37.20){\circle*{1.51}}
\put(41.33,37.33){\circle*{1.07}}
\put(39.47,36.80){\circle*{1.79}}
\put(36.27,36.80){\circle*{1.60}}
\put(47.87,36.27){\circle*{2.03}}
\put(47.73,36.67){\circle*{2.63}}
\put(46.00,37.47){\circle*{1.10}}
\put(47.20,37.60){\circle*{1.07}}
\put(51.20,36.53){\circle*{2.15}}
\put(52.93,36.13){\circle*{1.62}}
\put(51.73,36.00){\circle*{1.69}}
\put(53.87,36.13){\circle*{1.87}}
\put(56.13,36.13){\circle*{1.62}}
\put(57.73,36.27){\circle*{1.79}}
\put(60.67,36.13){\circle*{1.07}}
\put(62.13,36.67){\circle*{1.60}}
\put(64.53,36.53){\circle*{1.07}}
\put(65.47,36.53){\circle*{1.19}}
\put(67.47,36.80){\circle*{0.75}}
\put(43.20,37.20){\circle*{1.60}}
\put(38.53,37.33){\circle*{1.10}}
\put(35.07,36.93){\circle*{1.44}}
\put(53.47,37.33){\circle*{1.60}}
\put(56.13,37.47){\circle*{1.33}}
\put(58.53,37.07){\circle*{2.13}}
\put(60.67,37.33){\circle*{1.33}}
\put(59.87,37.60){\circle*{1.33}}
\put(61.60,37.60){\circle*{0.80}}
\put(62.67,37.20){\circle*{1.62}}
\put(63.73,37.33){\circle*{1.10}}
\put(33.87,36.67){\circle*{1.89}}
\put(51.33,37.60){\circle*{1.60}}
\put(46.27,36.13){\circle*{1.60}}
\put(40.53,36.13){\circle*{1.33}}
\put(61.33,36.00){\circle*{1.13}}
\put(53.87,37.60){\circle*{1.44}}
\put(54.93,37.73){\circle*{0.84}}
\put(31.60,36.93){\circle*{0.38}}
\put(50.53,43.47){\circle*{2.75}}
\put(52.93,43.33){\circle*{2.40}}
\put(54.80,43.20){\circle*{2.72}}
\put(57.33,43.60){\circle*{2.40}}
\put(59.73,43.60){\circle*{2.15}}
\put(62.00,43.73){\circle*{2.13}}
\put(63.87,43.73){\circle*{1.89}}
\put(65.73,44.00){\circle*{1.87}}
\put(47.87,43.47){\circle*{2.67}}
\put(45.60,43.33){\circle*{2.03}}
\put(43.73,43.33){\circle*{2.46}}
\put(41.87,43.87){\circle*{1.92}}
\put(40.40,43.73){\circle*{1.51}}
\put(39.07,43.60){\circle*{1.44}}
\put(37.47,43.87){\circle*{1.71}}
\put(35.73,44.00){\circle*{1.71}}
\put(34.40,44.27){\circle*{1.19}}
\put(33.47,44.40){\circle*{0.80}}
\put(32.53,44.53){\circle*{0.53}}
\put(31.87,44.53){\circle*{0.53}}
\put(31.20,44.53){\circle*{0.53}}
\put(49.20,42.80){\circle*{1.07}}
\put(51.87,42.40){\circle*{0.84}}
\put(56.40,42.80){\circle*{1.33}}
\put(58.80,43.07){\circle*{1.55}}
\put(57.73,42.80){\circle*{1.07}}
\put(60.93,43.20){\circle*{1.62}}
\put(67.07,44.40){\circle*{1.07}}
\put(67.87,44.00){\circle*{1.10}}
\put(68.53,44.40){\circle*{0.53}}
\put(56.00,43.73){\circle*{2.13}}
\put(53.87,44.00){\circle*{1.87}}
\put(51.73,44.00){\circle*{1.79}}
\put(53.07,44.40){\circle*{1.10}}
\put(49.07,43.60){\circle*{2.15}}
\put(49.33,44.40){\circle*{1.19}}
\put(46.40,43.73){\circle*{2.20}}
\put(45.07,44.13){\circle*{1.44}}
\put(44.13,44.13){\circle*{1.62}}
\put(42.80,44.13){\circle*{1.62}}
\put(40.93,44.13){\circle*{1.60}}
\put(39.47,44.00){\circle*{1.94}}
\put(38.27,43.87){\circle*{2.13}}
\put(36.53,44.13){\circle*{1.62}}
\put(34.93,44.40){\circle*{1.07}}
\put(32.93,44.53){\circle*{0.80}}
\put(58.67,44.00){\circle*{1.89}}
\put(60.93,43.87){\circle*{1.89}}
\put(63.07,44.13){\circle*{1.36}}
\put(64.80,44.00){\circle*{1.87}}
\put(63.73,44.40){\circle*{1.07}}
\put(60.13,44.53){\circle*{0.80}} \put(48.80,5.60){\circle*{1.07}}
\put(50.67,5.73){\circle*{0.84}} \put(52.13,5.73){\circle*{1.60}}
\put(53.73,5.73){\circle*{1.62}} \put(55.60,5.60){\circle*{1.36}}
\put(57.20,5.47){\circle*{1.07}} \put(58.40,5.47){\circle*{0.80}}
\put(59.07,5.33){\circle*{0.53}} \put(56.40,5.73){\circle*{1.19}}
\put(54.53,5.73){\circle*{1.69}} \put(52.93,5.73){\circle*{1.89}}
\put(51.20,5.87){\circle*{1.89}} \put(49.87,5.87){\circle*{1.89}}
\put(48.40,5.87){\circle*{1.89}} \put(46.93,5.87){\circle*{1.94}}
\put(45.33,5.60){\circle*{1.55}} \put(44.00,5.73){\circle*{1.60}}
\put(42.67,5.60){\circle*{1.07}} \put(41.73,5.60){\circle*{1.33}}
\put(40.67,5.47){\circle*{1.07}} \put(39.73,5.47){\circle*{1.07}}
\put(38.80,5.33){\circle*{0.80}} \put(38.13,5.20){\circle*{0.53}}
\put(59.60,5.33){\circle*{0.84}} \put(57.73,5.47){\circle*{1.10}}
\put(60.40,5.20){\circle*{0.53}}
\put(33.07,44.40){\line(1,-6){1.20}}
\put(36.27,43.60){\line(1,-6){1.07}}
\put(40.00,43.20){\line(1,-6){0.93}}
\put(44.13,42.53){\line(0,-1){4.67}}
\put(48.80,42.67){\line(0,-1){4.80}}
\put(53.07,42.93){\line(0,-1){5.33}}
\put(53.07,37.60){\line(0,0){0.00}}
\put(57.20,43.07){\line(0,-1){5.87}}
\put(62.67,42.93){\line(0,-1){5.20}}
\put(65.87,43.60){\line(0,-1){6.40}}
\put(34.53,36.13){\line(0,-1){6.40}}
\put(37.20,36.40){\line(0,-1){6.40}}
\put(40.27,36.13){\line(0,-1){6.40}}
\put(43.87,36.00){\line(0,-1){6.53}}
\put(48.93,35.87){\line(0,-1){6.27}}
\put(53.87,35.87){\line(0,-1){6.40}}
\put(58.27,36.13){\line(0,-1){5.47}}
\put(62.67,36.27){\line(0,-1){6.13}}
\put(36.53,29.33){\line(1,-4){2.13}}
\put(40.80,28.67){\line(0,-1){7.07}}
\put(44.80,28.40){\line(0,-1){6.13}}
\put(49.87,28.13){\line(0,-1){6.40}}
\put(54.27,28.67){\line(0,-1){6.93}}
\put(58.80,28.80){\line(0,-1){7.73}}
\put(62.53,29.60){\line(0,-1){9.20}}
\put(39.60,19.60){\line(-1,-6){1.20}}
\put(42.67,19.73){\line(0,-1){6.93}}
\put(42.67,12.80){\line(0,0){0.00}}
\put(46.80,18.93){\line(0,-1){6.13}}
\put(46.80,12.80){\line(0,0){0.00}}
\put(50.80,19.20){\line(0,-1){6.27}}
\put(55.07,19.47){\line(0,-1){6.67}}
\put(59.07,19.60){\line(0,-1){7.07}}
\put(39.33,11.20){\line(1,-3){1.87}}
\put(42.93,10.40){\line(0,-1){4.67}}
\put(47.20,10.40){\line(0,-1){4.40}}
\put(51.47,10.67){\line(0,-1){5.20}}
\put(55.47,11.33){\line(0,-1){5.73}}
\put(59.33,11.60){\line(0,-1){6.27}}
\put(46.40,2.80){\makebox(0,0)[cc]{$W$}}
\put(46.13,47.47){\makebox(0,0)[cc]{$W\cdot h$}}
\put(19.20,22.13){\makebox(0,0)[cc]{$L(h,z_1)$}}
\put(80.53,22.80){\makebox(0,0)[cc]{$R(h,z_t)$}}
\put(72.53,46.93){\makebox(0,0)[cc]{$z_t$}}
\put(63.60,2.80){\makebox(0,0)[cc]{$z_t$}}
\put(34.40,2.80){\makebox(0,0)[cc]{$z_1$}}
\put(27.60,46.93){\makebox(0,0)[cc]{$z_1$}}
\end{picture}

\begin{center}
\nopagebreak[4] Fig. \theppp.

\end{center}
\addtocounter{ppp}{1}

\begin{lm} \label{trap}
Every trapezium is a normal trapezium $T(W,h)$ for some $W$ and
some $h$.
\end{lm}

\proof Consider a trapezium $\Delta$ with bases $q_1, q_2$ and
sides $p_1$, $p_2$. By condition (TR2) $p_1$ and $p_2$ contain
$\rr$-edges. Since by (TR2), $p_1$ and $p_2$ are sides of reduced
$\bkk$-bands,  $\Lab(p_1), \Lab(p_2)$ are $\rr$-reduced. Since
$\Delta$ is $Z(\mathbb{S},a)$-reduced, we can apply Lemma
\ref{bandsnohubs}(\ref{b3}) and conclude that every maximal
$\rr$-band starting on $p_1$ ends on $p_2$ and every maximal
$\rr$-band starting on $p_2$ ends on $p_1$. By Lemma
\ref{bandsnohubs}(\ref{b4}), $\Delta$ does not contain any other
maximal $\rr$-bands.

Since by (TR1) $\Lab(q_1)$ and $\Lab(q_2)$ are $\bkk$-reduced, by
Lemma \ref{uchastok} every maximal $\bkk$-band starting on $q_1$
(resp. $q_2$) ends on $q_2$ (resp. $q_1$) and by Lemma
\ref{bandsnohubs}(\ref{b1}) $\Delta$ does not have any other
maximal $\bkk$-bands.

Let $\ttt, \ttt'$ be the first and the last maximal $\bkk$-bands
in $\Delta$  starting on $q_1$ counting from $p_1$ to $p_2$. Let
$\bb$ and $\bb'$ be the first and the last $\rr$-bands of
$\Delta$, counting from $q_1$ to $q_2$.

Consider the subdiagram $\Delta'$  of $\Delta$ bounded by $q_1$
and the connecting line of $\bb$. By Lemma \ref{rbands} the
boundary label of $\Delta'$ is equal to $\Lab(q_1)W\iv$ where $W$
is an admissible word. Since $\Delta'$ does not contain $r$-cells
(by Lemma \ref{bandsnohubs}, part 3), $W=\Lab(q_1)$ in $\hnka$.
Since $\Lab(q_1)$ is reduced by (TR1), $W$ is a reduced admissible
word and $\Lab(q_1)$ is a reduced admissible word which is equal
to $W$. Let $\Delta''$ be the diagram obtained from $\Delta$ by
removing $\Delta'$ and $\bb$. Then
$\partial(\Delta'')=p_1'q_1'(p_2')\iv (q_2')\iv$ where $p_1'$ is a
side of $\ttt$ without the first cell, $p_2'$ is a side of $\ttt'$
without the first cell, $q_1'$ is a side of $\bb$, $q_2'=q_2$.
Since by Lemma \ref{rbands} $\Lab(q_1')=\Lab(q_1)\cdot \tau$ in
$\hnka$ if $\bb$ is a $\tau$-band, we have that $\Lab(q_1')$ is
$\bkk$-reduced.

Suppose that $p_1'$ does not have $\rr$-edges. Then $p_2'$ does
not have $\rr$-edges either (every maximal $\rr$-band of
$\Delta''$ is a maximal $\rr$-band of $\Delta$), so the boundary
of $\Delta''$ contains no $\rr$-edges. By Lemma \ref{bandsnohubs}
(\ref{b5}) $\Delta''$ does not contain $r$-cells, hence
$\Lab(q_2)=\Lab(q_1')=\Lab(q_1)\cdot \tau$ in $\hnka$, hence
$\Delta=T(\Lab(q_1),\tau)$ as required.

Now if $p_1'$ contains $\rr$-edges then $p_2'$ also contains
$\rr$-edges, so $\Delta''$ satisfies conditions (TR1) and (TR2),
so $\Delta''$ is a trapezium. Since $\Delta''$ contains fewer
maximal $\rr$-bands than $\Delta$, by induction, we can assume
that $\Delta''=T(\Lab(q_1'),h)$ for some $h$. Then
$\Delta=T(\Lab(q_1), \tau h)$ as required
\endproof

\section{The group $\hra$}
\label{arsec}

Consider the group \label{hra}$\hra$ which is the group generated
by $\aaa\cup\rr$ subject to the set of relations consisted of all
$a$-relations and all (commutativity) $ra$-relations from
$Z(\sss,a)$. We consider $\aaa$-letters as zero-letters and
$\rr$-letters as non-zero letters. So in terminology of Section
\ref{segregation}, in this section $Y=\rr$.

We are going to prove that the presentation of $\hra$ satisfies
conditions (Z1), (Z2), (Z3) from Section \ref{segregation}

Locally, in this section, we are going to keep the same notation
as in Section \ref{segregation}. We let ${\cal R}_0={\cal S}_0$ be
the set of $a$-relations and $ra$-relations. The set ${\cal
S}_{1/2}$ is thus empty. Clearly conditions (Z1.1) and (Z1.2)
hold.

Since ${\cal S}_{1/2}$ is empty, conditions (Z2.1), (Z2.2) and
(Z2.3) also hold.

We will use the terminology and results from Section \ref{axioms}
and Section \ref{corollaries}. In particular, we are going to talk
about essential, reduced and cyclically reduced words in $\hra$.

As in Section \ref{segregation}, for every cyclically reduced
essential element $g\in \hra$ we define $\oo(g)=\oo(g)$ as the
maximal subgroup of $\maa$ which contains $g$ in its normalizer.
For an arbitrary essential element $g$ of the form $vuv\iv$ where
$u$ is cyclically reduced, we define $\oo(g)=v\oo(u)v\iv$. By
Lemma \ref{welldef}, $\oo(g)$ is well defined (does not depend on
the presentation $g=vuv\iv$).

The following lemma will imply both (Z3.1) and (Z3.2).

\begin{lm} \label{ra3.2}
Let $g$ be an essential element represented by a cyclically
reduced word $A$, let $x\ne 1$ be a $0$-element such that
$g^{-4}xg^4$ is a $0$-element. Then the following two statements
hold.

(i) There exists a $0$-element $u$ such that $gu^{-1}$ commutes
with every element of $\oo(g)$.

(ii) $x\in \oo(g)$.
\end{lm}

\proof By Lemma \ref{sdvig}, $g^{-2}xg^2$ is $0$-element. Let $B$
be a word in $\aaa$ representing a non-trivial element $x$ from
$\oo(g)$, $C$ be a word in $\aaa$ representing $g^{-2}xg^{2}$.

Let $A\equiv u_1r_1u_2...r_su_{s+1}$ where $r_i\in \rr$,
$i=1,...,s$, $u_i$ are words in $\aaa$, $i=1,...,s+1$. By Lemma
\ref{welldef}  (all of the conditions of that lemma hold), we can
replace $A$ by any cyclic shift of $A$. So we can assume that
$u_1$ is empty.

Then there exists a diagram $\Delta$ of minimal type over the
presentation ${\cal S}_0$ of $\hra$ with boundary $p_1q_1p_2\iv
q_2\iv$ where $\Lab(q_1)\equiv \Lab(q_2)\equiv A^2$,
$\Lab(p_1)\equiv B$, $\Lab(p_2)\equiv C$.

Since $A$ is cyclically reduced and $p_1, p_2$ do not  contain
$\rr$-edges, every maximal $\rr$-band in $\Delta$ starting on
$q_1$ ends on $q_2$ and every maximal $\rr$-band starting on $q_2$
ends on $q_1$. By Lemma \ref{bandsnohubs} (\ref{b4}), $\Delta$
contains no other maximal $\rr$-bands.

Let $\ttt_1,...,\ttt_{2s}$ be all maximal $\rr$-bands starting on
$q_1$, counted from $p_1$ to $p_2$. Notice that every
$ra$-relation is a commutativity relation of the form $ra=ar$
where $r\in \rr(z), a\in A(z)$ for some $z\in \kk$. For every
$i=1,...,2s$, $r_i$ belongs to $\rr(z_i)$ for some $z_i\in \kk$.
Then the labels of the top and the bottom paths of $\ttt_i$ are
equal to a word $v_i$ in the alphabet $A(z_i)\cup A(z_i)\iv$.
Notice that $B$ is conjugate of each of the words $v_i$,
$i=1,...,2s$ in $\hra$. Indeed, the subdiagram of $\Delta$ bounded
by $p_1, q_1,q_2$ and the connecting line of $\ttt_i$ has boundary
label $t\iv Btv_i\iv$ for some word $t$. Since $B\ne 1$ modulo in
$H_{ra}$, none of the words $v_i$ is equal to 1 modulo in
$H_{ra}$.

The connecting lines of the bands $\ttt_1,...,\ttt_{2s}$ divide
$\Delta$ into diagrams $\Delta_0$,...,$\Delta_{2s+1}$ such that
the boundary of each $\Delta_i$ contains no $\rr$-edges. By Lemma
\ref{bandsnohubs} (\ref{b5}), each $\Delta_i$ has no $r$-cells, so
the boundary label of each $\Delta_i$ is equal to 1 in the group
$\maa$.

Since $A$ starts with $r_1$, the subdiagram $\Delta_0$ is bounded
by a side of $\ttt_1$ and $p_1$. Hence $B=\Lab(p_1)=v_1$ modulo
$a$-relations.

The boundary label of $\Delta_1$ is $u_2\iv v_1u_2v_2\iv$. Hence
$u_2\iv v_1u_2=v_2$ in $\maa$. Recall that $\maa$ is the free
product of subgroups $\la A(z)\ra$ in the variety of Burnside
groups of exponent $n$. By Lemma 36.2, 34.11 from \cite{book},
each of these subgroups is malnormal in $\maa$. Therefore $v_2,
u_2\in \la A(z_1)\ra$ modulo $a$-relations. Since $u_2, v_2\ne 1$,
$z_1\equiv z_2$. Now considering subdiagrams
$\Delta_2,...,\Delta_{s+1}$ one by one, we conclude that $u_i,
v_i$ are words in $A(z_1)$, $i=1,...,s+1$, and $z_1=z_2=...=z_s$.
Therefore all $r_i$ belong to $\rr(z_1)$. Since $v_i$ are not
empty, $r_i$ commutes with some letters from $A(z_i)$. By
definition of $Z(\sss)$ then $r_i$ commutes with every letter from
$A(z_1)$ in $\hra$. Hence $r_1,...,r_s$ commute with all letters
from $u_1,...,u_{s+1}$. In particular,
$A=r_1r_2...r_su_1...u_{s+1}$ in $H_{ra}$.

Let $u=u_1...u_s$. This word does not contain letters from $\rr$.
Notice that $\maa$ is a retract of $\hra$, $u$ is the
$\maa$-projection of $g$, so $u$ does not depend in $H_{ra}$ on
the word $A$ representing $g$. Let also $g_1=r_1...r_s$. Since
$\la\rr\ra$ is a retract of $\hra$, $g_1$ does not depend on the
word $A$ representing $g$, $gu^{-1}=g_1$.

We have proved that $B$ is equal modulo $a$-relations to a word in
$A(z_1)$. So $g_1$ commutes with $B=x$. Since $\oo(g)$ is a
$0$-subgroup normalized by $g$, it consists of $0$-elements $y$
such that $g^{-2}yg^2$ is a $0$-element. We have proved that $g_1$
commutes with every such $y$. Hence $gu^{-1}$ commutes with all
elements in $\oo(g)$ which proves (i).

Notice that the subgroup $\la A(z_1)\ra$ is normalized by $g=g_1u$
because $g_1$ commutes with this subgroup and $u$ belongs to it.
Hence $\la A(z_1)\ra\subseteq \oo(g)$. Since $x\in \la A(z_1)\ra$,
$x\in \oo(g)$, which proves (ii).
\endproof

We have checked all properties (Z1), (Z2), (Z3), so by Proposition
\ref{mainprop}, we have the following statement.

\begin{prop} \label{propra}
There exists a graded presentation \label{rra}${\cal R}_{ra}$ of
the factor group \label{hrai}$H_{ra}(\infty)$ of $\hra$ over the
subgroup generated by all $n$-th powers of elements of $\hra$,
such that every g-reduced diagram over ${\cal R}_{ra}$ satisfies
property A from section \ref{conditionA}, and all the lemmas from
Section \ref{segregation} hold.
\end{prop}

As a corollary of Proposition \ref{propra}, we have the following
statements which we are going to use later.

\begin{lm}\label{ct1}
Let $\Delta$ be a g-reduced diagram over ${\cal R}_{ra}$,
$\partial(\Delta)=p_1q_1s\iv q_2\iv$, obtained by gluing two
diagrams $\Gamma$ and $\Pi$ where

(i) $\partial(\Gamma)=p_1q_1p_2\iv q_2\iv$, $\Lab(p_1)$ is a
subword of $R(h,z)$ for some reduced word $h$ over $\sss$ and some
$z\in \bkk\cup\bkk\iv$, $p_1$ starts and ends with $\rr$-edges;
$q_1$ and $q_2$ do not contain $\rr$-edges; $\Gamma$ is a diagram
of rank $0$ (that is over $Z(\sss,a)$);

(ii) $\Pi$ is a cell of corresponding to a relation $v^n=1$ in
${\cal R}_{ra}$ of rank $\ge 1$, $\partial(\Pi)=p_2s\iv$;

(iii) $v$ starts with an $\rr$-letter $r$, $\Lab(p_2)\equiv v^l r$
where $l>1/2\varepsilon n$.

Then there exists a diagram $\Delta'$ over ${\cal R}_{ra}$ which
has the same boundary label as $\Delta$,
$\partial(\Delta')=p_1'q_1'(s')\iv (q_2')\iv$ where
$\Lab(p_1')\equiv \Lab(p_1), \Lab(q_1')\equiv \Lab(q_1)$,
$\Lab(s')\equiv \Lab(s), \Lab(q_2')\equiv \Lab(q_2)$, and

(1) $\Delta'$ contains a subdiagram $\Pi'$ where
$\Lab(\partial(\Pi'))\equiv u^n$ where $u$ is a cyclically
$\rr$-reduced word,

(2) $|\partial(\Pi')\cap p_1'| > 6|u|$,

(3) $\Delta'\backslash \Pi'$ is of rank 0.
\end{lm}

\proof Consider the $\rr$-bands in $\Gamma$ starting on $p_2$.
Since $v$ is a period, it is cyclically $\rr$-reduced (by
Proposition \ref{propra} and the definition of periods from
Section \ref{segregation}. Since $q_1$ and $q_2$ are of
$\rr$-length 0, every $\rr$-band starting on $p_2$ ends on $p_1$.
Since $\Lab(p_1)$ is a subword of an $\rr$- reduced word $R(h,z)$,
every $\rr$-bands starting on $p_1$ ends on $p_2$. Hence the
$\rr$-projections of words $\Lab(p_1)$ and $\Lab(p_2)$ are
identical (the labels of $\rr$-edges in all cells of a $\rr$-band
are the same).

Therefore $\Lab(p_1)\equiv u^l r$ where $|u|_\rr=|v|_\rr$, $u$
starts with $r$. Since $u$ and $u^2$ are subwords of an
$\rr$-reduced word $R(h,z)$ (since $l\ge \varepsilon n>2$), $u$ is
cyclically $\rr$-reduced. In particular, $u$ is an essential word.

Let $\ttt$ and $\ttt'$ be the first and the last maximal
$\rr$-bands in $\Gamma$ starting on $p_1$ (we count the
$\rr$-bands from $(p_1)_-$ to $(p_1)_+$). Consider the subdiagram
$\Gamma_1$ of $\Gamma$ bounded by the top side of $\ttt$, $p_1$,
$p_2$ and the connecting line of $\ttt'$ (so this subdiagram is
equal to $\Gamma$ minus $\ttt'$). Since all cells in $\ttt'$ are
commutativity cells, $x_1=\Lab(\topp(\ttt'))=\Lab(\bott(\ttt'))$.
The subdiagram of $\Gamma$ bounded by a side of $\ttt'$ and $q_1$
does not contain $\rr$-bands, so $\Lab(q_1)=x_1$ modulo
$a$-relations. Similarly $\Lab(q_2)=x_2$ modulo $a$-relations
where $x_2$ is the word written in a side of $\ttt$. Clearly
$|x_1|_\rr=|x_2|_\rr=0$.

Let $\bar p_1$ be the path $p_1$ without the last edge, $\bar p_2$
be $p_2$ without the last edge. Then $\Lab(\bar p_1)\equiv u^l$,
$\Lab(\bar p_2)\equiv v^l$,

\begin{equation} \label{uv}
u^l x_1=x_2 v^l
\end{equation}
modulo $a$-relations and $ra$-relations.

Similarly every subword of the cyclic word $u$ is equal to a
subword of the same $\rr$-length of the cyclic word $v$ up to
multiples of zero length. Hence $v$ is cyclically reduced by Lemma
\ref{reduced}. Therefore $v$ is essential.

Suppose that $u$ is not simple in rank 0. By Lemma \ref{anal18.1}
there exists a word $w$ which is simple in rank 0 (and so $w$ is
simple in rank $1/2$ because ${\cal R}_{ra}$ does not contain
hubs), $|w|_\rr<|u|_\rr$, and a word $t\in \oo(w)$ such that
$u=xw^btx\iv$ in rank 0, for some word $x$ and some integer $b\ne
0$. Since $u$ is cyclically reduced, $|x|_\rr=0$ and
$|w^b|_\rr=|u|_\rr$.

Therefore by $u^l=xw^{bl}t_1x\iv$ for some $t_1\in \oo(w)$ (since
$\oo(w)$ is normalized by $w$). Using (\ref{uv}), we get $$u^l
x_1= xw^{bl}t_1x\iv x_1 = x_2 v^l$$ in rank 0, where both $w$ and
$v$ are simple in rank 0, the $\rr$-lengths of $x, x_1, t_1, x_2$
are equal to 0 (for $x, x_1, x_2$ we have established it above,
for $t_1$, this follows from the definition of $\oo(w)$ since $w$
is cyclically $\rr$-reduced). Now we can apply Lemma
\ref{anal18.8} and conclude that $v=yw^{\pm 1}t_2y\iv$ in rank 0
for some word $y$ and some $t_2\in \oo(w)$. This contradicts the
fact that $v$ is simple in rank 0 (since $|w|_\rr <
|u|_\rr=|v|_\rr$).

Thus $u$ is simple in rank $0$ (and so it is simple in rank
$1/2$). Then again by Lemma \ref{anal18.8}, we can conclude that

\begin{equation} \label{uv1}
x_2v^{\pm 1}t_2x_2\iv=u
\end{equation}
in rank 0 for some $t_2\in \oo(v)$. Suppose that the exponent of
$v$ in (\ref{uv1}) $-1$, that is $x_2v\iv t_2x_2\iv =u$. Therefore
$$u^l=x_2(v\iv t_2)^l x_2\iv = x_2v^{-l}t'x_2\iv$$ for some $t'\in
\oo(v)$ (since $\oo(v)$ is normalized by $v$). By (\ref{uv}),
$$x_2\iv u^l= v^{l}x_1\iv = v^{-l}t'x_2\iv$$ whence
$$v^{2l}=t'x_2\iv x_1.$$ Since $t'$ is of $\rr$-length 0 (this
follows from the definition of $\oo(v)$ since $v$ is cyclically
$\rr$-reduced), $x_1$ and $x_2$ are also of $\rr$-length $0$, we
have that $v^{2l}$ is of $\rr$-length 0 which contradicts the fact
that $v$ is cyclically $\rr$-reduced. This contradiction shows
that in fact the exponent of $v$ in (\ref{uv1}) is 1.

By Lemma \ref{star}, $(vt_2)^n=v^n$ in rank 0. Hence $u^n =
x_2v^nx_2\iv$ in rank 0.

Recall that $\Lab(\bar p_1)\equiv u^l$. Let $\bar\Delta$ be a
diagram with contour $\bar p_1(p_1')\iv$ where $\Lab(p_1')\equiv
u^l$, which is a result of gluing two mirror copies $\Pi'$ and
$\bar\Pi$ with boundary labels $u^{\pm n}$ over ${\cal R}_{ra}$
where $p_1'$ is a part of the boundary of $\Pi'$, $\bar p_1$ is a
part of the boundary of $\Pi$. Consider the union of $\Delta$ and
$\bar\Delta$. Notice that $|p_1'|\ge l|u|\ge
\frac{1}{2}\varepsilon n|u|>6|u|$. It remains to observe that in
the union of $\Delta$ and $\bar\Delta$, the $\bar\Pi$ is connected
with $\Pi$ by a path with label $x_2$, hence this pair of
subdiagrams and the connecting path can be replaced by a diagram
of rank 0 corresponding to the equality $u^n=x_2v^nx_2\iv$.  The
lemma is proved.
\endproof

\section{The group $\hkra$}

\label{hkrasec}  Recall that in Section \ref{bands} we considered
the set of defining relations $Z(\sss,a)$ which is the union of
the set of all $a$-relations and the set $Z(\sss)$ of relations
corresponding to the $S$-machine $\sss$.

Now consider the union
\label{tza}$\tza=Z(\sss,a)\cup\{\Lambda(0)\}\cup {\cal R}_{ra}$
and the group \label{hkra}$\hkra$ generated by $\bkk\cup\rr\cup
\aaa$ subject to the relations from $\tza$. Each relation of
$\tza$ is either the hub, an $a$-relation, a (commutativity)
$ra$-relation, or a $\kappa$-relation or belongs to ${\cal
R}_{ra}$. Burnside relations from ${\cal R}_{ra}$ will be called
\label{bbra}{\em Burnside $ra$-relations}. The corresponding
cells in diagrams will be called {\em Burnside $ra$-cells}. The
set of all relations from $\tza$ except the hub will be denoted
by \label{tz}$\tz$.


As in Section \ref{segregation}, we divide the set of generators
of $\hkra$ into the set of $0$-letters and the set of
non-$0$-letters.

Our main goal is to prove that the presentation $\tza$ satisfies
properties (Z1), (Z2), (Z3) if we consider $\bkk$ as the set of
non-zero letters and the $\rr\cup\aaa$ as the set of $0$-letters.

But at first we show how to get rid of Burnside $ra$-cells in
diagrams over $\tz$.

As in Section \ref{bands}, we can consider $\bkk$- and
$\rr$--bands in diagrams over $\tz$. The definition of these bands
is the same as in Section \ref{bands} (Burnside cells do not
belong to $\rr$-bands). The only difference is that now
$\rr$-bands can start and end on Burnside $ra$-cells.

\subsection{Getting rid of Burnside $ra$-cells}
\label{rrmain}

Consider (temporarily) $\rr$ as the set of non-$0$-letters and all
other letters as $0$-letters. We can consider $\rr$-bands
connecting two Burnside $ra$-cells (a Burnside $ra$-cell and a
part of the boundary of a diagram) as $0$-bonds, and define bonds
and \ct subdiagrams of all ranks. (We preserve the grading of the
set of Burnside $ra$-relations ${\cal R}_{ra}$.)

Suppose that a disc \vk diagram $\Delta$ is obtained by gluing
together two disc diagrams $\Theta$ and $\Pi$ over $\tz$ along
pieces of their boundaries, where

(1) $\Theta$ is a trapezium $T(W,h)$ with contour $p_1q_1p_2\iv
q_2\iv$ ($q_1, q_2$ are the bottom and the top bases);

(2) $\Lab(\Pi)\equiv u^n$ for some cyclically $\rr$-reduced word
$u$ over $\rr\cup\aaa$, the diagram $\Pi$ does not contain
$\bkk$-edges;

(3) $\Pi$ and $\Theta$ are glued along a subpath $p$ of $p_2$;
$|p|\ge 3|u|$.

Then we shall say that $\Delta$ has the form \label{tsp}$\Theta
*\Pi$.

Suppose that a disc \vk diagram $\Delta$ is obtained by gluing
together two diagrams $\Pi$ and $\Theta$ where

(1) $\Theta$ is a trapezium $T(W,h)$ with contour $p_1q_1p_2\iv
q_2\iv$ ($q_1, q_2$ are the bases);

(2) $\Lab(\Pi)\equiv u^n$ for some cyclically $\rr$-reduced word
$u$ over $\rr\cup\aaa$ containing $\rr$-letters; $\Pi$ does not
contain $\bkk$-edges;

(3) $\Pi$ and $\Theta$ are glued along a subpath $p$ of $p_1$.

Then we shall say that $\Delta$ has the form \label{pct}$\Pi\circ
\Theta$. It is easy to see that $$T(W,h) * \Pi=\Pi \circ T(W\cdot
h,h\iv)$$ (this means that if the right hand side of the equality
exists then both sides exist and are equal).

\begin{lm}\label{starcirc}
For every \vk disc diagram $\Delta$ of the form $\Theta *\Pi$
there exists a \vk disc diagram $\Delta'$ of the form $\Pi'\circ
\Theta'$ with the same boundary label and such that:

(i) the label of the bottom (top) base of $\Theta$ is equal to the
label of the bottom (top) base of $\Theta'$;

(ii) The histories of the trapezia $\Theta$ and $\Theta'$ are
equal modulo Burnside relations (see
\setcounter{pdeleven}{\value{ppp}} Figure \thepdeleven).
\end{lm}

\unitlength=1.00mm \special{em:linewidth 0.4pt}
\linethickness{0.4pt}
\begin{picture}(105.34,59.50)
\put(20.83,45.00){\line(0,-1){37.67}}
\put(20.83,7.33){\line(1,0){20.33}}
\put(41.17,7.33){\line(0,1){37.50}}
\put(41.17,44.83){\line(-1,0){20.33}}
\put(46.33,28.33){\oval(10.33,18.67)[]}
\put(85.17,57.00){\line(0,-1){50.00}}
\put(85.17,7.00){\line(1,0){20.17}}
\put(105.34,7.00){\line(0,1){50.00}}
\put(105.34,57.00){\line(-1,0){20.33}}
\put(78.84,31.25){\oval(12.33,25.17)[]}
\put(62.83,29.00){\makebox(0,0)[cc]{$\longrightarrow$}}
\put(31.50,25.00){\makebox(0,0)[cc]{$\Theta$}}
\put(18.33,27.34){\makebox(0,0)[cc]{$p_1$}}
\put(30.83,9.00){\makebox(0,0)[cc]{$q_1$}}
\put(43.83,16.50){\makebox(0,0)[cc]{$p_2$}}
\put(30.50,42.50){\makebox(0,0)[cc]{$q_2$}}
\put(46.33,28.00){\makebox(0,0)[cc]{$\Pi$}}
\put(39.00,28.00){\makebox(0,0)[cc]{$p$}}
\put(30.33,5.00){\makebox(0,0)[cc]{$W$}}
\put(30.67,47.33){\makebox(0,0)[cc]{$W\cdot h$}}
\put(95.00,4.50){\makebox(0,0)[cc]{$W$}}
\put(94.00,59.50){\makebox(0,0)[cc]{$W\cdot h$}}
\put(94.50,31.50){\makebox(0,0)[cc]{$\Theta'$}}
\put(78.84,30.50){\makebox(0,0)[cc]{$\Pi'$}}
\end{picture}

\begin{center}
\nopagebreak[4] Fig. \theppp.

\end{center}
\addtocounter{ppp}{1}

\proof Let $\Theta=T(W,h)$. Then $\Lab(p_2)\equiv R(h,z)\iv$ where
$z\in\bkk\cup\bkk\iv$ is the last letter of $W$. Let $p$ be the
common subpath of the contour of $\Pi$ and the right side of the
trapezium $\Theta$, $\Lab(\partial(\Pi))\equiv u^n$ for a
cyclically $\rr$-reduced word $u$. Then $\Lab(p)\equiv R(h,z)\iv$
contains a cube of a cyclic shift of $u$ as a subword. Without
loss of generality we can assume that this cyclic shift is $u$
itself. The word $R(h,z)\iv$ is a product of blocks of the form
$R(\tau,z)$, $ \tau\in\sss$, each of which is either a one-letter
word $r\in \rr\cup\rr\iv$ or a two-letter word $(ra)^{\pm 1}$
where $r\in \rr$, $a\in \aaa\cup\aaa\iv$. Every block is
completely determined by its $\rr$-letter. Consider the occurrence
of $u^3$ in $R(h,z)\iv$.

The word $u^3$ contains subword $v^2$ where $v$ is a cyclic shift
of $u$ and the first occurrence of $v$ (in $v^2$) starts with the
beginning of a block. Suppose that $v$ does not end with the end
of a block. Then the second occurrence of $v$ in $v^2$ starts in
the middle of a block, so the first letter of $v$ is both a
beginning of a block and an end of a block. Therefore this letter
belongs to $\aaa^{\pm 1}$. Hence $v\equiv ar_1v'r_2$ where $a\in
\aaa\cup\aaa\iv$, $r_1\in \rr\iv$, $r_2\in \rr$, $ar_1$ and $r_2a$
are blocks. But then the second letter of the second occurrence of
$v$ must be $r_1$, and it must be the beginning of a block. Thus
the same $\rr$-letter $r_1$ is a beginning and the end of a block.
This contradicts the fact that every block is completely
determined by its $\rr$-letter.

Hence $v$ starts with a beginning of a block and ends with the end
of a block. Therefore $v\equiv R(h_1,z)$ for some cyclically
reduced word $h_1$ over $\sss$. Thus $h_1^2$ is a subword of $h$.
Let $h\equiv h'h_1^2h''$ for some words $h', h''$. Then $W\cdot
h'$ is in the domain of $h_1^2$. By Lemma \ref{burns}(iv), then
$W\cdot h'$ is in the domain of $h_1^s$ for every integer $s$. In
particular, $W\cdot h'$ is in the domain of $h_1^{2-n}$. By Lemma
\ref{burns}(ii),

$$W\cdot h=W\cdot h'h_1^2h'=W\cdot h'h_1^{2-n}h''.$$

Since $v\equiv R(h_1,z)$ is a cyclic shift of $u$, one can read
$v^n$ on the boundary of $\Pi$, and the boundary label of $\Theta
*\Pi$ is freely equal as a cyclic word to

\begin{equation}\label{trap1}
\begin{array}{l}
WR(h'',z)\iv v^{n-2}R(h',z)\iv (W\cdot h)\iv L(h,z')=\\
WR(h'',z)\iv R(h_1^{2-n},z)\iv R(h',z)\iv (W\cdot h)\iv L(h,z')=\\
WR(h'h_1^{2-n}h'',z)\iv (W\cdot h)\iv L(h'h_1^2h'',z')
\end{array}
\end{equation}
where $z'$ is the first letter of $W$.

Consider now the trapezium $\Theta'=T(W,h'h_1^{2-n}h'')$ which
exists by Lemma \ref{trap}(i) since $W$ is in the domain of
$h'h_1^{2-n}h''$. Since $W\cdot h'h_1^{2-n}h''=W\cdot h$, we can
assume that the top base label of $\Theta'$ is the same as the top
base label of $\Theta$, that is $W\cdot h$.

The boundary label of $\Theta'$ is thus freely equal (as a cyclic
word) to

$$WR(h'h_1^{2-n}h'',z)\iv (W\cdot h)\iv L(h'h_1^{2-n}h'',z').$$

Notice that by Proposition \ref{propra} there exists a \vk disc
diagram $\Pi'$ containing only $a$-cells, $ra$-cells and Burnside
$ra$-cells, that has boundary label $L(h_1,z')^{n}$. We can glue
$\Pi'$ and $\Theta'$ along the path labeled by $L(h_1,z')^{2-n}$,
we get a diagram $\Delta'$ of the form $\Pi'\circ \Theta'$. The
boundary label of this diagram will be freely equal to

$$WR(h'h_1^{2-n}h'',z)\iv (W\cdot h)\iv L(h'h_1^2h'',z').$$

Comparing this with (\ref{trap1}) we conclude that the boundaries
of $\Delta$ and $\Delta'$ are freely equal. The labels of the
bottom (top) bases of $\Theta$ and $\Theta'$ are equal by
construction, the history $h\equiv h'h_1^2h''$ of $\Theta$ is
equal to the history $h'h^{2-n}h''$ of $\Theta'$ modulo Burnside
relations, as required.
\endproof

Lemma \ref{starcirc} allows us to move subdiagrams with boundary
labels of the form $u^n$ from one side of a trapezium to another
side, preserving the history of the trapezium modulo Burnside
relations. The next lemma is the first application of this trick.

\begin{lm}\label{kkolca}
For every disc diagram $\Delta$ over $\tz$ there exists another
disc diagram $\Delta'$ with the same boundary label and without
$\bkk$-annuli.
\end{lm}

\proof We can assume that $\Delta$ contains exactly one
$\bkk$-annulus $\bb$ and $\partial(\Delta)$ is equal to
$t=\topp(\bb)$. Then the subdiagram $\Delta'$ bounded by
$\bott(\bb)$ is a diagram over ${\cal R}_{ra}$. We can assume that
this subdiagram is $g$-reduced, and that the $\bkk$-annulus $\bb$
is reduced.

By Proposition \ref{propra} $\Delta'$ is an $A$-map. Suppose that
$\Delta'$ contains a Burnside $ra$-cell. By Lemma \ref{Anal16.2}
there exists a Burnside $ra$-cell $\Pi$ in $\Delta'$ and its \ct
subdiagram $\Gamma$ to $\bott(\bb)$ with \ct degree at least
$\varepsilon$. By Lemma \ref{ct1}, we can replace $\Pi\cup\Gamma$
by a new subdiagram $\Psi$ with the same boundary label which
contains a subdiagram $\Pi'$ such that $\Lab(\partial(\Pi'))\equiv
u^n$ for some word $u$ and $|\partial(\Pi')\cap \bott(\bb)|>6|u|$,
$\Psi\backslash \Pi'$ is a diagram of rank 0. Now using Lemma
\ref{starcirc}, we can move $\Pi'$ through $\bb$ without producing
new $\bkk$-annuli.

After a number of such transformations we replace $\Delta$ by a
new disc diagram $\bar\Delta$ with the same boundary label
containing exactly one $\bkk$-annulus which bounds a subdiagram
without Burnside $ra$-cells. By Lemma \ref{bandsnohubs}, part 1,
applied to the subdiagram bounded by that $\bkk$-band, we can
replace that subdiagram by a subdiagram without $\bkk$-edges.
Since this operation reduces the number of $\bkk$-annuli, the
statement of the lemma follows.
\endproof

Let $\Delta$ be an annular diagram over $\tz$ with inner contour
$p_1$ and outer contour $p_2$, and let $q$ be a path in $\Delta$.
Suppose that $\Delta, q$ satisfy the following conditions:\label{T123}

(T1) $p_1$ and $p_2$ are sides of $\bkk$-annuli $\ttt_1$ and
$\ttt_2$ respectively; $\Lab(p_1)\ne 1$ modulo $\tz$;

(T2) $q$ connects $(p_1)_-=(p_1)_+$ with $(p_2)_-=(p_2)_+$; the
first and the last edges of $q$, $e_1, e_2$, are $\bkk$-edges;

(T3) Every maximal $\bkk$-band in $\Delta$ is an annulus
surrounding the hole of $\Delta$.

Maximal $\bkk$-bands in $\Delta$ divide $\Delta$ into annular
subdiagrams (without $\kappa$-cells) which will be called
\label{chamber}{\em chambers}. We count the chambers of $\Delta$
from $p_1$ to $p_2$.

\begin{lm} \label{chambers}
Suppose that all chambers of $\Delta$ except the last one,
$\bar\Delta$, do not contain Burnside $ra$-cells and $\Lab(q)$ is
$\bkk$-reduced.

Then there exists another annular diagram $\Delta'$ over
$Z(\sss,a)$ (that is without Burnside $ra$-cells) with contours
$p_1'$ and $p_2'$, and a path $q'$ in $\Delta'$ satisfying the
conditions (T1),(T2), (T3) such that $\Lab(p_1')=\Lab(p_1)$,
$\Lab(p_2')=\Lab(p_2)$ modulo Burnside relations,
$\Lab(q')=\Lab(q)$ modulo $\tz$, $\Lab(q')$ is $\bkk$-reduced.
\end{lm}

\proof Among all pairs $(\Delta, q)$ satisfying the conditions of
the lemma and having the same boundary labels of contours $p_1,
p_2$ modulo Burnside relations and the same labels of $q$ modulo
$\tz$ let us choose a pair where $\Delta$ has the smallest
possible number of chambers and the last chamber $\bar\Delta$ has
the smallest type in the sense of Section \ref{segregation}.

First of all let us reduce all chambers of $\Delta$. We can get
rid of $j$-pairs of Burnside $ra$-cells by first replacing the
path $q$ by a path avoiding the $j$-pair and its connecting path
and then reducing the $j$-pair. Similarly, if $\Delta$ contains a
disc subdiagram $\Psi$ over $Z(\sss,a)$ containing $r$-cells,
whose boundary label is a word over $\aaa$ and is equal to 1
modulo $a$-relations, then we can replace this subdiagram by a
subdiagram with the same boundary label containing no $r$-cells.
The resulting diagram will satisfy conditions (T1), (T2), (T3).
These operations will not increase the number of chambers in
$\Delta$ or the type of $\bar\Delta$ or the $\bkk$-length of
$\Lab(q)$. Thus we can assume that  $\Delta$ does not contain disc
subdiagrams whose boundary labels are words over $\aaa$ equal 1
modulo $a$-relations. In particular, all chambers except for the
last one are $Z(\mathbb{S},a)$-reduced in the sense of Section
\ref{nohubs}.

We can also reduce all $\bkk$-annuli in $\Delta$ and also the
bands $\ttt_1$ and $\ttt_2$ without changing the $\bkk$-length of
$q$. Notice that we do not assume that  $\ttt_1$ and $\ttt_2$ are
{\bf cyclically} reduced because it may be impossible to reduce
the annuli $\ttt_1$ and $\ttt_2$ without changing $q$.

Hence we can assume that $\bar\Delta$ is g-reduced in the sense of
Section \ref{segregation} as a diagram over $\tz$ and
$\Delta\backslash(\bar\Delta\cup \ttt_2\cup\ttt_1)$ which is a
diagram over $Z(\sss,a)$ is $Z(\mathbb{S},a)$-reduced.

Since $\Lab(q)$ is $\bkk$-reduced in $\hkra$, $q$ intersects every
$\bkk$-annulus in $\Delta$ only once. Therefore the $\bkk$-length
of $q$ is equal to the number of $\bkk$-annuli in $\Delta$.

Let $\bb$ be a maximal $\rr$-band connecting the contours of
$\Delta''=\Delta\backslash(\bar\Delta\cup\ttt_2)$. Such a band
exists because if every $\rr$-band starting on the boundary of
$\ttt_1$ ends on the boundary of $\ttt_1$ then the label of a side
of $\ttt_1$ is equal to 1 modulo $\tz$, a contradiction with (T1).

Since the label of a side of $\bb$ is equal to the label of the
portion of $q$ in that subdiagram multiplied by words of
$\bkk$-length 0 (read along the boundaries of the subdiagram), the
label of any side of $\bb$ is $\bkk$-reduced.

By Proposition \ref{propra}, chamber $\bar\Delta$ is an $A$-map.
Let $o=(e_2)_-$.  Let $s$ be the inner contour of $\bar\Delta$,
and let $t$ be the outer contour of $\bar\Delta$. We consider $o$
as the beginning and the end of $t$.

Suppose that $\bar\Delta$ contains Burnside $ra$-cells. Then by
Lemma \ref{th16.2}, part c),  there exists a Burnside $ra$-cell
$\Pi$ in $\bar\Delta$ and its \ct subdiagram $\Gamma$ of rank 0 to
one of the contours $s$ or $t$ of \ct degree $\ge \varepsilon$.

Replacing $\Gamma$ (if necessary) by a smaller \ct diagram, we can
assume that  conditions (i)-(iv) of Lemma \ref{ct1} hold. By this
lemma, we can replace $\Gamma\cup \Pi$ by a subdiagram $\Psi$ with
the same boundary label, which has a subdiagram $\Pi'$ such that
$\partial(\Pi')\equiv u^n$ where $u$ is a cyclically $\rr$-reduced
word and $|\partial(\Pi')\cap s|>6|u|$ or $|\partial(\Pi')\cap
t|>6|u|$ and $\Psi\backslash \Pi'$ has rank 0.

Suppose first that $|\partial(\Pi')\cap t|>6|u|$. By Lemma
\ref{starcirc} we can move $\Pi'$ through $\ttt_2$ replacing
$\ttt_2$ by a $\bkk$-band with a history which is equal to the
history of $\ttt_2$ modulo Burnside relations. The path $q$ can be
replaced in the new annular diagram by a path with the same label
modulo $\tz$. This operation decreases the type of the last
chamber of $\Delta$.

Suppose now that $|\partial(\Pi')\cap s|>6|u|$. Let
$x=\partial(\Pi')\cap s$. Consider the annular subdiagram
$\Delta''$ bounded by $s$ and $p_1$. The diagram $\Delta''$ is a
$Z(\mathbb{S},a)$-reduced diagram. By Lemma \ref{bandsnohubs}
every maximal $\rr$-band starting on $x$ must end on $p_1$. Indeed
otherwise the $\rr$-band would have two intersections with a
reduced $\bkk$-band in a disc subdiagram of $\Delta''$ which
contradicts Lemma \ref{bandsnohubs}, part \ref{b3}. Since such an
$\rr$-band cannot intersect twice a $\bkk$-annulus, and $\phi(q)$
is a $\bkk$-reduced word, the label of a side of any of these
$\rr$-bands is a $\bkk$-reduced word. There exists a disc
subdiagram $\Theta$ of $\Delta''$ containing at least half of
these $\rr$-bands, bounded by a subpath $y$ of $x$, a subpath of
$p_1$ and by two $r$-bands. Now   $\Theta$ is a trapezium because
it satisfies conditions (TR1) and (TR2). By Lemma \ref{trap},
$\Theta$ is a normal trapezium. The union of $\Theta$ and $\Pi'$
is, by definition, equal to $\Theta * \Pi'$. As before we can move
$\Pi'$ through $\Theta$ replacing $\Theta$ by a trapezium with the
same history modulo Burnside relations. Again we decrease the type
of $\Delta$.

 The
new diagram satisfies conditions (T1), (T2), (T3) and  has a
smaller type than $\Delta$.

Thus we can complete the proof by induction on the type of
$\Delta$.
\endproof

Lemma \ref{chambers} implies the following stronger statement.

\begin{lm} \label{chambers1}
The conclusion of Lemma \ref{chambers} holds if one removes the
restriction that only the last chamber contains Burnside
$ra$-cells.
\end{lm}

\proof Since $\Lab(q)$ is reduced, we can assume that it crosses
the connecting line of every $\bkk$-annulus only once.

Let us use induction on the number of chambers in $\Delta$. If the
number of chambers is 1, we can apply Lemma \ref{chambers}.
Consider the diagram
$\Delta_1=\Delta\backslash(\bar\Delta\cup\ttt_2)$ where
$\bar\Delta$ is the last chamber. Let $q_1$ be the part of $q$ in
$\Delta_1$. The number of chambers in $\Delta_1$ is smaller than
the number of chambers in $\Delta$. Therefore we can apply the
induction hypothesis and conclude that there exists a diagram
$\Delta_2$ and a path $q_2$ which satisfy the conclusion of Lemma
\ref{chambers}.

Since the label of the outer contour of $\Delta_2$ is equal to the
label of the outer contour of $\Delta_1$ modulo Burnside
relations, we can glue to the outer contour of $\Delta_2$ a
diagram over ${\cal R}_{ra}$ such that the resulting diagram
$\Delta_3$ has the same outer boundary label as $\Delta_1$.

Now replace the subdiagram $\Delta_1$ in $\Delta$ by $\Delta_3$
and the subpath $q_1$ of $q$ by $q_2$. Let $(\Delta_4, q_4)$ be
the resulting pair of an annular diagram and a cutting path.
Reducing if necessary, the last $\bkk$-annulus of subdiagram
$\Delta_3$ of $\Delta_4$, we obtain a pair $(\Delta_5,q_5)$
satisfying conditions of Lemma \ref{chambers}. Applying that
lemma, we complete the proof. \endproof

The following lemma is an algebraic consequence of Lemma
\ref{chambers1}.

\begin{lm} \label{chambers2}
Let $P_1, P_2$ be words over $\aaa\cup\rr$ which are not equal to
1 modulo $\tz$. Let $A\iv P_1A=P_2$ in $\hkra$ for some word $A$
which is $\bkk$-reduced in $\hkra$ and starts and ends with a
letter $z\in \bkk\cup\bkk\iv$. Then $A=R_1BR_2$ in $\hkra$ where
$R_1=L(h_1,z)$, $R_2=R(h_2,z)$, $z$ is in the domain of
$h_1,h_2\iv$, $B$ is a reduced admissible word starting and ending
with $z_1$ such that $R_1z_1=zR(h_1,z)$, $z_1R_2=L(h_2,z_1)z$
modulo $Z(\sss,a)$.
\end{lm}

\proof Equality $P_2=A\iv P_1A$ gives us an annular diagram
$\Delta$ with inner contour $p_1$ labeled by $P_1$, outer contour
$p_2$ labeled by $P_2$ and a cutting path $q$ labeled by $A$.

The connecting line of any $\bkk$-annulus in $\Delta$ not
surrounding the hole of $\Delta$ cannot cross $q$ because
otherwise by Lemma \ref{kkolca} we would be able to decrease the
$\bkk$-length of $q$ by going along the boundary of such a
$\bkk$-annulus instead of crossing it. Hence by Lemma \ref{kkolca}
we can assume that $\Delta$ does not contain $\bkk$-annuli not
surrounding the hole of $\Delta$.

Let $e$ and $f$ be the first and the last edges of $q$. Since $e$
and $f$ are $\bkk$-edges, they belong to the first and the last
$\bkk$-annulus of $\Delta$ counting from $p_1$ to $p_2$. The label
of $p_1$ is equal modulo $\tz$ to the label of a side $p_1'$ of
the first $\bkk$-annulus in $\Delta$ and the label of $p_2$ is
equal in $\tz$ to the label of a side of the last $\bkk$-annulus
in $\Delta$. Hence we can assume that $p_1$ and $p_2$ are the
sides of these annuli.

Thus $\Delta$ satisfies conditions (T1), (T2), (T3). Since
$\Lab(q)$ is $\bkk$-reduced, we can apply Lemma \ref{chambers1}
and conclude that there exists another annular diagram $\Delta'$
over $Z(\sss,a)$ with contours $p_1'$ and $p_2'$, and a path $q'$
in $\Delta'$ satisfying the conditions (T1),(T2), (T3) such that
$\Lab(p_1')=\Lab(p_1)$, $\Lab(p_2')=\Lab(p_2)$ modulo Burnside
relations, $\Lab(q')=\Lab(q)=A$ modulo $\tz$, $\Lab(q')$ is
$\bkk$-reduced.

Assume that every maximal $\bf R$-band starting on the first $\bf
K$-annulus of $\Delta'$ does not end on $p'_2$. Then the label of
$p'_1$ is equal to $1$ modulo $Z(\sss, ra),$ a contradiction.
Hence $\Delta'$ contains an $\bf R$-band $\cal B$, starting on
$p'_1$ and ending on $p'_2$.

Let $t$ be  a side of $\bb$. Then $q$ is homotopic to $(p_1')^s
(p_1'') tp_2''$ for some integer $s$, where $p_1''$ is a subpath
of $p_1'$ and $p_2''$ is a subpath of $p_2'$. By Lemma
\ref{rbands} $\Lab(t)\equiv B$ is an admissible word starting and
ending with a $\bkk$-letter $z'$. The $\bkk$-length of $B$ is the
same as the $\bkk$-length of $A$ (since $\bb$ intersects each of
the $\bkk$-annuli exactly once), whence $B$ is reduced in $\hkra$.

The word $\Lab((p_1')^s(p_1'')$ is written on the top side of
$z'$-band, so it is equal to $R_1\equiv L(h_1,z)$ for some $h_1$.
Similarly $\Lab(p_2'')\equiv R(h_2,z)$ for some $h_2$. The last
equalities from the lemma follow from Lemma \ref{kbands}.
\endproof

\begin{lm}\label{commutativity}
Let $A$ be a $\bkk$-reduced word without $\rr$-letters, starting
and ending with letters $z_1,z_2$ from $\bkk\cup\bkk\iv$, and
suppose that for some words $P_1$ and $P_2$ of $\bkk$-length $0$,
$P_1\ne 1$ modulo $\tz$, we have $A\iv P_1A=P_2$ modulo $\tz$.
Then:

(1) $P_1=L(h,z_1)$, $P_2=R(h,z_2)$ modulo $\tz$ for some word $h$
over $\sss$;

(2) $A$ is a reduced admissible word;

(3) $A\cdot h=A$.
\end{lm}

\proof Consider a \vk disc diagram $\Delta$ over $\tz$ for the
equality $A\iv P_1A=P_2$, $\partial(\Delta)=p_1q_1p_2\iv q_2\iv$
where $\Lab(p_i)\equiv P_i\iv, \Lab(q_i)\equiv A$, $i=1,2$. By
Lemma \ref{kkolca} we can assume that $\Delta$ does not contain
$\bkk$-annuli. We can also assume that all $\bkk$-bands in
$\Delta$ are reduced. Hence every maximal $\bkk$-band in $\Delta$
connects $q_1$ and $q_2$. These maximal $\bkk$-bands divide
$\Delta$ into diagrams over $\tz$ without $\bkk$-edges, which we
shall call {\em chambers} as in the annular case before. As in the
proof of Lemma \ref{chambers1}, one can use induction on the
number of chambers, Proposition \ref{propra}, Lemma \ref{ct1} and
Lemma \ref{starcirc}, to move all Burnside $ra$-cells out of
$\Delta$. Since only two sides of $\Delta$, namely $p_1$ and $p_2$
may contain $\rr$-edges, only these two sides will change during
this process. Their labels will stay the same modulo $\tz$. Hence
there exists a \vk diagram $\Delta_1$ over $Z(\sss,a)$ (i.e. over
$\tz$ but without Burnside $ra$-cells) with contour
$p_1'q_1(p_2')\iv q_2\iv$ where $\Lab(p_i')=\Lab(p_i)$, modulo
Burnside relations $i=1,2$.

Let $\ttt_1$ and $\ttt_2$ be the first and the last $\bkk$-bands
in $\Delta_1$ starting on $q_1$ and ending on $q_2$ (we count the
$\bkk$-bands starting at $p_1'$). Consider the subdiagram
$\Delta_2$ of $\Delta_1$ bounded by $\ttt_1$, $\ttt_2$, $q_1$,
$q_2$ (and containing $\ttt_1, \ttt_2$). Then
$\partial(\Delta_2)=p_1''q_1(p_2'')\iv q_2\iv$ where $p_1''$ is a
side of $\ttt_1$, $p_2''$ is a side of $\ttt_2$.

Considering the subdiagrams of $\Delta_1$ bounded by
$p_1'(p_1'')\iv$ and $p_2'(p_2'')\iv$, we conclude that
$\Lab(p_1')=\Lab(p_1'')$, $\Lab(p_2')=\Lab(p_2'')$ modulo $\tz$.
This gives part (1) of the lemma.

The subdiagram $\Delta_2$ is a trapezium (by definition). Hence it
is a normal trapezium $T(A,h)$ for some word $h$ (by Lemma
\ref{trap}). In particular, $A$ is an admissible word which proves
part (2) of the lemma. Moreover,
$\Lab(p_1)=\Lab(p_1')=\Lab(p_1'')=L(h,z_1)\iv$,
$\Lab(p_2)=\Lab(p_2')=\Lab(p_2'')=R(h,z_2)\iv$ modulo $\tz$ for
some word $h$ over $\sss$. Since the top base of $\Delta_2$ has
label $A$, we have $A\cdot h=A$. This gives part (3) of the lemma.
\endproof

\begin{lm}\label{commutativity1}
Let $B$ be a $\bkk$-reduced word, starting with a letter $z\in
\bkk\cup\bkk\iv$ and ending with a letter $z'\in \bkk\cup\bkk\iv$
which does not contain $\rr$-letters . Let $U$ be a word over
$\aaa$, let $P,R$ be words of $\bkk$-length $0$  such that
$BURBUz$ is $\bkk$-reduced, $P\ne 1$ in $Z(\sss,ra)$. Suppose that
$(BURBUz)\iv P(BURBUz)$ is a $0$-word $Q$ modulo $\tz$. Then $BUz$
is an admissible word, and there exist $0$-words $P', R'$ and $Q'$
which are equal modulo $\tz$ respectively to $P, R,Q$, such that
$(BUR'BUz)\iv P'(BUR'BUz)=Q'$ modulo $Z(\sss,a)$.
\end{lm}

\proof Consider a disc diagram $\Delta$ over $\tz$ corresponding
to the equality $$(BURBUz)\iv P(BURBUz)=Q,$$
$\partial(\Delta)=p_1q_1p_2\iv q_2\iv$ where $\Lab(p_1)\equiv P$,
$\Lab(p_2)\equiv Q$, $\Lab(q_1)\equiv \Lab(q_2)\equiv BURBUz$. Let
$q_1=s_1x_1t_1s_1'$, $q_2=s_2x_2t_2s_2'$ where $\Lab(s_i)\equiv
B$, $\Lab(x_i)\equiv U$, $\Lab(t_i)\equiv R$, $\Lab(s_i')\equiv
BUz$, $i=1,2$. Let $\ttt_1, \ttt_2$ be the first and the last
$\bkk$-bands starting on $s_1$, $\ttt_3, \ttt_4$ be the first and
the last $\bkk$-bands starting on $s_1'$. Since $BUR$ is
cyclically reduced, $\ttt_1,\ttt_2$ end on the first and the last
edges of $s_2$, $\ttt_3,\ttt_4$ end on the first and the last
edges of $s_2'$. Then we can assume that $p_1$ is a side of
$\ttt_1$, $p_2$ is a side of $\ttt_4$. Let $\bar p_1$ be the side
of $\ttt_2$ which connects $(s_2)_+$ with $(s_1)_+$, let $\bar
p_2$ be a side of $\ttt_3$ which connects $(s_2')_-$ with
$(s_1')_-$.

By Lemma \ref{commutativity} $\Lab(p_1)=L(h_1,z)$,
$\Lab(p_2)=R(h_2,z)$ modulo $\tz$, $B$ and $BUz$ are admissible
words, and $B\cdot h_1=B$, $BUz\cdot h_2=BUz$. Thus as in Lemma
\ref{chambers1} we can move all Burnside cells out of the
subdiagrams $\Delta_1$ and $\Delta_2$ of $\Delta$ bounded by
$p_1s_1(\bar p_1)\iv s_2\iv$ and by $\bar p_2s_1'p_2\iv
(s_2')\iv$, and replace them by trapezia $T(B,h_1)$ and $T(B,h_2)$
in the expense of possible increasing the number of cells in the
subdiagram $\Gamma$ bounded by $\bar p_1x_1t_1 \bar p_2\iv
(x_2t_2)\iv$. Hence we can assume that $\Delta_1=T(B,h_1)$,
$\Delta_2=T(B,h_2)$.

Thus all Burnside $ra$-cells in $\Delta$ are in $\Gamma$. Consider
the diagram $\Gamma$. By Lemma \ref{kkolca} we can assume that
$\Gamma$ does not have $\bkk$-annuli. Since the boundary of
$\Gamma$ does not contain $\bkk$-edges, the entire diagram
$\Gamma$ does not have $\bkk$-edges.

Since $\Lab(t_1)\equiv \Lab(t_2)$, we can identify $t_1$ and $t_2$
to obtain an annular diagram $\bar\Gamma$ and a path $t$
connecting the contours $x_2\iv\bar p_1x_1$ and $\bar p_2$ of
$\Gamma$. We can eliminate all $j$-pairs of $\bar\Gamma$,
replacing $t$ by a path with the same label modulo $\tz$. Hence we
can assume that $\bar\Gamma$ does not have $j$-pairs. Hence
$\bar\Gamma$ is an $A$-map by Proposition \ref{propra}.

If $\bar\Gamma$ does not have Burnside $ra$-cells, we are done,
because we can again cut $\bar\Gamma$ along the path $t$. Suppose
that $\bar\Gamma$ contains Burnside $ra$-cells. Consider
$\bar\Gamma$ as a map satisfying condition A from Section
\ref{conditionA} with boundary divided into four parts:
$x_1,x_2,\bar p_1, \bar p_2$. Recall that paths $x_1$, $x_2$ do
not contain $\rr$-edges. By Lemma \ref{th16.2}, there exists a
Burnside $ra$-cell $\Pi$ and a \ct subdiagram $\Gamma'$ of rank
$0$ of $\Pi$ to $\bar p_1$ or $\bar p_2$ with \ct degree at least
$\varepsilon$. We can assume that $t$ avoids $\Gamma'$. Hence
$\Pi$ and $\Gamma'$ are subdiagrams of the (disc) diagram
$\Gamma$. As in the proofs of Lemmas \ref{chambers},
\ref{chambers1}, we can use Lemmas \ref{ct1} and \ref{starcirc}
and move $\Pi$ through either the trapezium $\Delta_1$ or the
trapezium $\Delta_2$ leaving $\Lab(p_1)$ and $\Lab(p_2)$ the same
modulo $\tz$, leaving paths $s_i, x_i$ untouched and decreasing
the number of Burnside $ra$-cells in $\Gamma$. The proof now can
be completed by induction on the number of Burnside $ra$-cells in
$\Gamma$.
\endproof

\subsection{Conditions (Z1), (Z2), (Z3)}

\label{z1z4}

Now we are ready to check that the presentation $\tza$ satisfies
conditions (Z1), (Z2), (Z3) from Section \ref{segregation}.

Here we are going to use notation from Section \ref{segregation}
again, as in Section \ref{arsec}.

Let ${\cal R}_0={\cal S}_0$ be the set $\tz$ of all relations in
$\tza$ except for the hub, and let ${\cal
S}_{1/2}=\{\Lambda(0)\}$, so $\tza={\cal R}_{1/2}$. We let ${\cal
Y}=\bkk$, so all letters from $\rr\cup\aaa$ are 0-letters.

Condition (Z1.1) holds by an easy inspection. Condition (Z1.2)
holds by Proposition \ref{propra}.

The length of $\Lambda$ is $(2n+3)N>n$, so (Z2.1) holds. The
condition (Z2.2) is obviously true.

The next lemma gives us (Z2.3).

\begin{lm} \label{hubs}
Assume that $v_1w_1$ and $v_2w_2$ are cyclic permutations of
$\Lambda(0)^{\pm 1}$ and $|v_1|\ge\varepsilon|\Lambda(0)|.$ Then
if $u_1v_1=v_2u_2$ holds modulo ${\cal S}_0$ for some 0-words
$u_1$, $u_2$, then we have $u_2w_1=w_2u_1$ modulo ${\cal S}_0$.
\end{lm}

\proof  Suppose that $u_1v_1=v_2u_2$ modulo ${\cal S}_0$. Then
there exists a \vk diagram $\Gamma$ over ${\cal S}_0$ with
boundary $p_1q_1p_2\iv q_2\iv$ such that $\Lab(p_1)\equiv u_1$,
$\Lab(p_2)\equiv u_2$, $\Lab(q_1)\equiv v_1$, $\Lab(q_2)\equiv
v_2$ (see Figure 6). The $\bkk$-bands of $\Gamma$ starting on
$q_1$ end on $q_2$ since $\Lab(q_1)$ is a linear word and
$u_1,u_2$ do not contain letters from $\bkk\cup\bkk\iv$. The
labels of the start and end edges of each of these bands are
equal. Hence $v_2\equiv \Lab(q_2)\equiv \Lab(q_1)\equiv v_1\equiv
v$. Since $|v|>0$, and $\Lambda(0)$ and $\Lambda(0)\iv$ do not
have common letters, $v$ cannot be a subword of both $\Lambda(0)$
and $\Lambda(0)\iv$. Without loss of generality assume that $v$ is
a subword of $\Lambda(0)$.

Recall that for every $j=1,...,N$, we denote $\Lambda_j(0)$ to be
the subword of $\Lambda(0)$ starting with $\tk(j)$ and ending with
$\tk(j+1)$ ( as usual ``$N+1$" here is 1).

By Lemma \ref{commutativity} there exists a word $h$ over $\sss$
such that $v\cdot h=v$, $u_1=L(h,z_1)$, $u_2=R(h,z_2)$ modulo
$\tz$ for some $z_1, z_2$. Hence we can assume that
$u_1=L(h,z_1)$, $u_2=R(h,z_2)$ in the free group and that
$\Gamma=T(v,h)$. Since $T(v,h)$ is a diagram over $Z(\sss,a)$,
$\Gamma$ does not contain Burnside $ra$-cells.

Since $v\equiv \Lab(q_1)\equiv \Lab(q_2)$ is a subword of
$\Lambda(0)$, we can attach to $\Gamma$ two hubs $\Pi_1$ and
$\Pi_2$ along the paths $q_1$ and $q_2$, and consider the
resulting diagram $\Delta$ with boundary label $w_2u_1w_1\iv
u_2\iv$. We need to show that $\Delta$ can be replaced by a
diagram over ${\cal S}_0$ with the same boundary label.

Since $|v_1|\ge \varepsilon|\Lambda(0)|$ and $N>n>3/\varepsilon$,
we conclude that $v_1$ contains $\Lambda_j(0)$ for some $j>1$.

Let $q_1'$ be the subpath of $q_1$ such that $\Lab(q_1')\equiv
\Lambda_j(0)$. Since $\Lambda$ is linear, every $\bkk$-band
starting on $q_1$ ends on $q_2$. Let $\ttt_1$ be the $\tk(j)$-band
and $\ttt_2$ be the $\tk(j+1)$-band starting on $q_1'$. Then
$\ttt_1$ ends on the edge of $q_2$ labeled by $\tk(j,0)$ and
$\ttt_2$ ends on the edge of $q_2$ labeled by $\tk(j+1,0)$.

The bands $\ttt_1$, $\ttt_2$, the path $q_1'$ and the path $q_2$
bound a subdiagram $\Gamma'$ of $\Gamma$ which is a trapezium. By
Lemma \ref{trap}, $\Gamma'=T(\Lambda_j(0),h)$ and
$\Lambda_j(0)\cdot h=\Lambda_j(0)$.  Let $t$ be the subpath of
$\partial(\Gamma')$ which is a side of $\ttt_1$.

Consider the subdiagram $\Delta'$ of $\Delta$ consisting of
$\Pi_1$, $\Pi_2$ and $t$. Since $\Delta\backslash \Delta'$
consists of $0$-cells, it suffices to show that the boundary label
of $\Delta'$ is equal to 1 modulo ${\cal S}_0$.

Thus we need to show that $\Lab(t)=L(h,\tk(j))$ commutes with

\begin{equation}\label{eqj}
\Lambda_j(0)\tk(j+1,0)\iv\Lambda_{j+1}(0)\tk(j+2,0)\iv...\Lambda_{j-1}(0)\tk(j,0)\iv.
\end{equation}
By Proposition \ref{copy1}, for every $i=1,...,N$ we have
$\Lambda_i(0)\cdot h=\Lambda_i(0)$. The trapezium
$T(\Lambda_i(0),h)$ gives us the equality

$$\Lambda_i(0)\iv
L(h,\tk(i))\Lambda_i(0)=R(h,\tk(i+1))=\tk(i+1,0)\iv
(L(h,\tk(i+1))\tk(i+1,0), i=1,...,N$$ modulo ${\cal S}_0$. Hence

$$ \tk(i+1,0)\Lambda_i(0)\iv L(h,\tk(i))\Lambda_i(0)\tk(i+1,0)\iv=
L(h,\tk(i+1)), i=1,...,N. $$ Hence if we conjugate $L(h,\tk(i))$
by the word (\ref{eqj}), we get $L(h,\tk(i))$ as
required.\endproof

As in Section \ref{segregation}, for every cyclically reduced
essential element $g\in \hkra$ we define $\oo(g)$ as the maximal
subgroup of $\hra$ which contains $g$ in its normalizer. For an
arbitrary essential element $g$ of the form $vuv\iv$ where $u$ is
cyclically reduced, we define $\oo(g)=v\oo(u)v\iv$. By Lemma
\ref{welldef}, $\oo(g)$ is well defined (does not depend on the
presentation $g=vuv\iv$).

The following lemma will imply both (Z3.1) and (Z3.2).

\begin{lm} \label{ra3.21}
Let $w$ be an essential element represented by a cyclically
minimal in rank $1/2$ word $W$, let $x\ne 1$ be a $0$-element in
$\hkra$ such that $w^{-4}xw^4$ is a $0$-element. Then the
following two statements hold.

(i) There exists a $0$-element $u$ such that $wu^{-1}$ commutes
with every element of $\oo(w)$.

(ii) $x\in \oo(w)$.
\end{lm}

\proof We divide the proof into several steps.

1. Let $P_1$ be a $0$-word representing $x$. By Lemma
\ref{cyclically} every cyclic shift of $W$ is a cyclically
$\bkk$-reduced word. Therefore we can assume that $W$ starts with
a letter $z\in \bkk\cup\bkk\iv$ (by Lemma \ref{welldef}).

By Lemma \ref{sdvig}, $(Wz)\iv P_1 Wz$ is equal in $\hkra$ to a
$0$-word $P$. By Lemma \ref{chambers2} $Wz=CBz'D$ in $\hkra$ where
$B$ does not contain $\rr$-letters and starts with a $\bkk$-letter
$z'$ such that $Cz'=zC'$, $z'D=D'z$ modulo $Z(\sss,a)$ where
$D'\equiv L(f,z)$, $C\equiv L(f',z)$ for some words $f,f'$ over
$\sss$.

In particular since $Wz=CBz'D=CBD'z$, we have $W=CBD'$ in $\hkra$.
Since $D'Cz'=z'DC'$, passing if necessary to a cyclic shift $BD'C$
of $W$, we can assume that $C$ is empty, so $z'=z$,  and replace
$W$ by $BD'$ where $D'z=zD$ modulo $Z(\sss,a)$. The word $B$ is
equal to $zB_1U$ where $B_1$ is either empty or ends with a letter
$z'\in \bkk\cup\bkk\iv$, $U$ is a word over $\aaa$.

2. By Lemma \ref{sdvig}, we have the following equalities in
$\hkra$:

\begin{equation} \label{11}
(zB_1UD'zB_1Uz)\iv P_1 (zB_1UD'zB_1Uz) = P_2
\end{equation}

\begin{equation} \label{12}
(zB_1UD'zB_1UD'zB_1Uz)\iv P_1 (zB_1UD'zB_1UD'zB_1Uz) = P_3
\end{equation}
where $P_2$ and $P_3$ are $0$-words.

By Lemma \ref{commutativity1} (which can be applied because
$zB_1UD'zB_1Uz$ is reduced since this word is equal to a subword
of the same length of the reduced word $W^3$) we can assume that
equality (\ref{11}) holds modulo $Z(\sss,a)$ (replacing $P_1$ and
$P_2$ by equal in $\hkra$ words if necessary), $Bz$ is a reduced
admissible word. By Lemma \ref{sdvig}, since $D'z=zD$, we have
that $(zB_1Uz)\iv P_1(zB_1Uz)$ is a $0$-word modulo $Z(\sss,a)$.
Hence by Lemma \ref{commutativity} $P_1=L(h,z)$ modulo $Z(\sss,a)$
and

\begin{equation}\label{13}
Bz\cdot h=Bz.
\end{equation}

3. Now consider a diagram $\Delta$ over $Z(\sss,a)$ corresponding
to equality (\ref{11}), $\partial(\Delta)=p_1q_1p_2\iv q_2\iv$
where $\Lab(p_1)\equiv P_1$, $\Lab(p_2)\equiv P_2$,
$\Lab(q_1)\equiv BD'Bz\equiv \Lab(q_2)$. Then $q_1=s_1t_1s_1'$,
where $\Lab(s_1)\equiv B$, $\Lab(t_1)\equiv D'$, $\Lab(s_1')\equiv
Bz$. Consider the maximal $z$-band $\ttt_1$ in $\Delta$ starting
on the first edge of $s_1'$. As usual, we can assume that $p_2$ is
a side of a $z$-band $\ttt_2$. The bands $\ttt_1$ and $\ttt_2$
bound a subdiagram $\Delta_1$ of $\Delta$ which is a trapezium. By
Lemma \ref{trap}, $\Delta_1=T(Bz,h_1)$ where $h_1$ is the history
of $\ttt_1$.

Consider the natural homomorphism $\psi$ from $\hkra$ to the free
group freely generated by $\sss$. This homomorphism takes $\bkk$
and $\aaa$ to 1, and takes each letter from $\rr$ to the
corresponding rule from $\sss$. It is clear that this homomorphism
is correctly defined because killing all letters from
$\bkk\cup\aaa$ in every relation from $Z(\sss,a)$ and replacing
letters from $\rr$ by the corresponding rules from $\sss$ gives us
a trivial relation $1=1$. Since $D'=L(f,z)$, $\psi(D')=f$.
Applying this homomorphism to the boundary label of the subdiagram
of $\Delta$ bounded by $p_1$ and a side of $\ttt_1$, we get
$h_1=f\iv hf$. Hence we have

\begin{equation}\label{14}
Bz\cdot f\iv hf=Bz.
\end{equation}

4. Now let $\Delta_1$ be a \vk diagram corresponding to equality
(\ref{12}). Applying Lemma \ref{commutativity} to the subdiagram
$\Delta_2$ of $\Delta_1$ bounded by two $z$-bands $\ttt_1$,
$\ttt_2$ corresponding to the two last distinguished occurrences
of $z$ in $zB_1UD'zB_1UD'zB_1Uz$, we obtain $zB_1Uz\cdot
h_2=zB_1Uz$, i.e. $Bz\cdot h_2=Bz$ where $h_2$ is the history of
the $z$-band $\ttt_1$. Applying the homomorphism $\psi$ to the
boundary label of the subdiagram of
$\Delta_2'=\Delta_1\backslash\Delta_2$ we obtain the equality
$\psi(D')^{-2}h\psi(D')^2h_2\iv=f^{-2}hf^2h_2\iv$ modulo Burnside
relations. Hence $f^{-2}hf^2h_2\iv=1$ modulo Burnside relations.
We conclude that

\begin{equation}\label{15}
Bz\cdot h_2=Bz
\end{equation}
where $h_2=f^{-2}hf^2$ modulo Burnside relations.

Let $W_1=Bz$, $h_1=f\iv hf$, $g=f$, $b=h_2$. Then $W_1\cdot
h_1=W_1$ by (\ref{14}), $W_1\cdot gh_1g\iv =W_1$ by (\ref{13}) and
$W_1\cdot b=W_1$ by (\ref{15}) where $b=g\iv h_1g$ modulo Burnside
relations.

Suppose that $Bz$ contains an accepted subword $B'$. Then $B'\cdot
h'$ is a cyclic shift $C$ of $\Lambda(0)$ for some word $h'$ over
$\sss$. Then for some admissible subword $B''$ of $B$ and $B'$,
$B''\cdot h'=C'$ where $C'C''\equiv C$, $|C''|\le 1$. By Lemma
\ref{trap} $B''=UC'V$ in $\hkra$ where $U$ and $V$ are $0$-words.
Hence $B''=U(C'')\iv V$ in rank $1/2$. Since $|U(C'')\iv V
|<|UC'V|$, and $B''$ is a subword of $W$, we get a contradiction
with the assumption that $W$ is cyclically minimal in rank 1/2.

Thus all conditions of Proposition \ref{propsss} hold. By that
proposition and by Remark \ref{remark}, we can conclude that for
every integer $s$ there exists an word $b_{s+1}$ over $\sss$ which
is equal to $g^sh_1g^{-s}=f^{-s-1}hf^{s+1}$ modulo Burnside
relations and such that

$$Bz\cdot b_{s+1}=Bz. $$

5. For arbitrary integer $s$ consider the trapezia
$T_s=T(Bz,b_s)$. Let $\partial(T_s)=p_1q_1p_2\iv q_2\iv$, where
$\Lab(p_1)=L(b_s,z)$, $\Lab(p_2)=R(b_s,z)$, $\Lab(q_i)=Bz$,
$i=1,2$.

Considering the subdiagram of $T_s$ obtained by removing the last
$z$-band (whose side is $p_2$), we get the equality

$$ L(b_s,z)BL(b_s,z)\iv =B $$
 modulo $Z(\sss,a)$. Thus $B$
commutes with $L(b_s,z)$ modulo $Z(\sss,a)$. Notice that $b_0=h$.
Hence $B$ commutes with $L(f,z)^{-s}L(h,z)L(f,z)^s$ modulo $\tz$.
Recall that $L(f,z)\equiv D'$, $L(h,z)\equiv P_1\equiv L(b_0,z)$.
Hence $(D')^{-s}x(D')^s$ commutes with $B$ for every integer $s$.

Therefore for every non-negative integer $s$, we have the
following equalities in $\hkra$:

\begin{equation}\label{18}
\begin{array}{l} W^{-s}P_1W^s=\\(BD')^{-s}P_1(BD')^s=(D')\iv
B\iv...(D')\iv B\iv(D')\iv  \underline{(B\iv P_1 B)}D'BD'...BD'=\\
(D')\iv B\iv...(D')\iv \underline{B\iv((D')\iv P_1 D')B}D'...BD'
=\\ (D')\iv B\iv ... \underline{B\iv
((D')^{-2}P_1(D')^2)B}...BD'=...\\ (D')^{-s}P_1(D')^s\end{array}
\end{equation}
(we underline the parts of the words which are changing). Similar
equalities hold for negative $s$. Therefore for every integer $s$,
$w^{-s}xw^s$ is a $0$-element in $\hkra$. Hence $x\in \oo(w)$
(since the subgroup generated by $w^{-s}xw^s$, $s=0, \pm 1,...$,
consists of $0$-elements and is normalized by $w$).

6. Let us prove that $B$ commutes with every element $y\in \oo(w)$
modulo $\tz$. Let $P_y$ be a word in $\rr\cup\aaa$ representing
$y$. Since $y\in \oo(w)$, $W^{-2}P_yW^2$ is equal modulo $\tz$ to
a $0$-word. By Lemma \ref{sdvig}, $(BD'z)\iv P_y(BD'z)$ is a
$0$-word (modulo $\tz$). Since $D'z=zD$ modulo $\tz$, we have that
$(BzD)\iv P_y(BzD)$ is a $0$-word modulo $\tz$. By Lemma
\ref{sdvig}, $(Bz)\iv P_y(Bz)$ is a $0$-word $P_y'$ modulo $\tz$.
Considering a \vk diagram for this equality, we deduce (by a usual
argument) that $P_y=L(h_y,z)$, $P_y'=R(h_y',z)$ for some words
$h_y, h_y'$ over $\sss$. Applying the homomorphism $\psi$ to the
equality $$ (Bz)\iv L(h_y,z)(Bz)=R(h_y',z)$$ we get that
$h_y=h_y'$ modulo Burnside relations. The last equality implies

\begin{equation}\label{vynos}B\iv P_yB=B\iv L(h_y,z)B =
zR(h_y',z)z\iv=L(h_y',z)=P_y\end{equation} modulo $\tz$ (by Lemma
\ref{kbands}). Hence $B$ commutes with $P_y$ modulo $\tz$, as
required.

7. Now let $y\in \oo(w)$ be represented by a word $P_y$. Let us
prove by induction on $l$ that $(D')^l P_y(D')^{-l}$ belongs to
$\oo(w)$ for every $l\ge 0$. This is clearly true for $l=0$.
Suppose that it is true for some $l\ge 0$. Then

\begin{equation}\label{222}
(D')^{l+1} P_y(D')^{-l-1}=(D') (D')^l P_y (D')^{-l} (D')\iv= B\iv
(BD')(D')^l P_y (D')^{-l} (BD')\iv B.
\end{equation}
 Since $(D')^l P_y (D')^{-l}$ belongs to $\oo(w)=\oo(BD')$, we have
that $(BD')(D')^l P_y (D')^{-l} (BD')\iv$ belongs to $\oo(w)$ (by
the definition of $\oo(w)$). Since $B$ commutes with every element
of $\oo(w)$, (\ref{222}) implies that $$(D')^{l+1} P_y(D')^{-l-1}=
(BD')(D')^l P_y (D')^{-l} (BD')\iv \in \oo(w).$$ Similarly one can
prove that $(D')^l P_y(D')^{-l}$ belongs to $\oo(w)$ for all
$l<0$.

8. Let $u$ be the element of $\hkra$ represented by $D'$. Then
clearly $u$ is a $0$-element. Let $w$ be the element represented
by $W(D')^{-1}$. We need to prove that $w$ commutes with $y$, that
is, $W(D')^{-1}$ commutes with $P_y$ modulo $\tz$, or,
equivalently,
$$(BD')^{-1}P_y(BD')=(D')^{-1}P_y(D')$$ modulo $\tz$. This has
been proved already in (\ref{vynos}).
\endproof

This completes verification of properties (Z1), (Z2), (Z3). By
Proposition \ref{mainprop}, we have the following statement.

\begin{prop} \label{propkra}
There exists a graded presentation \label{rkra}${\cal R}_{kra}$
containing $\tz$ of the factor group \label{hkrai}$\hkra(\infty)$
of $\hkra$ over the subgroup generated by all $n$-th powers of
elements of $\hkra$, such that every relation of ${\cal R}_{kra}$
of rank $\ge 1$ has the form $u^n=1$ for some word $u$, and every
minimal diagram over ${\cal R}_{kra}$ satisfies property A from
Section \ref{conditionA}, and all the lemmas from Section
\ref{segregation}.
\end{prop}

\section{Proof of Theorem \ref{CEP}}
\label{proof}

\begin{lm} \label{final}
For every word $w(x_1,...,x_m)$, if
$w(a_1(\kappa(1)),...,a_m(\kappa(1)))=1$ in $\hkra(\infty)$, then
$w(a_1,...,a_m)=1$ in $\la A(\kappa(1))\ra$.
\end{lm}

\proof Let $\Delta$ be a g-reduced diagram over ${\cal R}_{kra}$
with boundary label $w(a_1(\kappa(1)),...,a_m(\kappa(1)))$.
Suppose that $\Delta$ contains cells of rank $> 0$.  By
Proposition \ref{propkra} $\Delta$ is an $A$-map. By Lemma
\ref{th16.2} $\Delta$ contains a cell $\Pi$ of rank $> 0$ and a
\ct subdiagram of $\Pi$ to $\partial(\Delta)$. Since non-0 letters
in presentation ${\cal R}_{kra}$ are  letters from $\bkk$, by
definition of \ct subdiagrams from Section \ref{amaps},
$\partial(\Delta)$ contains a $\bkk$-edge, a contradiction.

Hence $\Delta$ does not contain cells of rank $>0$. Therefore
$\Delta$ is a diagram over $\hkra$. So we can assume that $\Delta$
is a g-reduced diagram over $\hkra$. By Lemma \ref{kkolca} we can
assume that $\Delta$ does not contain  $\bkk$-annuli. Hence
$\Delta$ does not contain $\bkk$-edges. Therefore $\Delta$ is a
diagram over ${\cal R}_{ra}$ where $\rr$ is the set of non-0
letters.

Suppose that $\Delta$ contains a cell of rank $>0$. Then by
Proposition \ref{propra} and Lemma \ref{Anal16.2} $\Delta$
contains a cell $\Pi$ of rank $>0$ and a \ct subdiagram of $\Pi$
to $\partial(\Delta)$. Hence $\partial(\Delta)$ contains
$\rr$-edges, a contradiction. Therefore $\Delta$ does not contain
cells of rank $>0$. Hence $\Delta$ is a diagram over $\hra$. By
Lemma \ref{bandsnohubs}, part (\ref{b4}), $\Delta$ does not
contain $\rr$-annuli. Hence $\Delta$ does not contain $\rr$-edges.
This means that all cells of $\Delta$ are $a$-cells. Therefore
$\Delta$ is a diagram over $\maa$.

Since $\la A(\kappa(1))\ra$ is a retract of $\maa$, we obtain that
$w(a_1(\kappa(1)),...,a_m(\kappa(1)))$ is equal to 1 modulo
$R(\kappa(1))$-relations and Burnside relations, as required.
\endproof

Given the Burnside group $B(m,n)$ generated by $a_1,\dots,a_m$, we
have, by Lemma \ref{lmbk}, a homomorphism $\phi: B(m,n)\to {\cal
H}$ such that $\phi(a_i)=a_i(\kappa(1))$, $i=1,\dots,m$.

Assume that $\phi(w(a_1,\dots,a_m))\in {\cal H}^{n}$. Then
$w(a_1(\kappa(1)),\dots a_m(\kappa(1)))=1$ in
$H_{kra}(\infty)={\cal H}(n)={\cal H}/{\cal H}^n$ by Proposition
\ref{propkra}. Setting $R(\kappa(1))=\varnothing$, we have from
Lemma \ref{final} that $w=1$ in $\la A(\kappa(1))\ra
=B(a_1,\dots,a_m; n)$. Hence $\phi$ is an embedding, and moreover,
the intersection of $\phi(B(m,n))\cap {\cal H}^n$ is trivial. This
proves Statement (1) of Theorem \ref{CEP}.

Let now $L$ be an arbitrary normal subgroup of $B(m,n)$ (the
latter can now be identified with its image in ${\cal H}/{\cal
H}^n$). To prove Statement (2) of Theorem \ref{CEP}, we have to
show that $L=L^{{\cal H}/{\cal H}^n}\cap B(m,n)=L$. Let us define
$R(\kappa(1))$ to be the set of all words in
$a_1(\kappa(1)),\dots, a_m(\kappa(1))$ representing elements of
$L$. The new group $H_{kra}(\infty)$ (whose construction depends
now on the new choice of $R(\kappa(1))$) is a homomorphic image of
${\cal H}/{\cal H}^n$ by Proposition \ref{propkra}. Therefore it
suffices to prove that every word
$w(a_1(\kappa(1)),\dots,a_m(\kappa(1)))$ that is equal to 1 in
$H_{kra}(\infty)$, is also trivial in $B(m,n)/L\cong\la
A(\kappa(1))\ra$. But again, this claim is a part of Lemma
\ref{final}.

The proof of Theorem \ref{CEP} is complete.

\twocolumn

\noindent {\bf Subject index} \label{sind}

\medskip

\noindent accepted admissible word \pr{accepted}\\ adjacent edges
\pr{adj}\\ admissible word \pr{adm}, \pr{adm1.5}, \pr{adm2}\\
annulus \pr{annulus} \\ $(S,T)$-annulus \pr{stanl}\\ $c$-aperiodic
word \pr{caper}\\ $\bkk$-band \pr{bkkband}\\ $\kk$-band
\pr{kkband}\\ $\rr$-band \pr{rrband}\\ $S$-band \pageref{sband}\\
$\tau$-band \pr{tauband}\\ bases of a trapezium \pr{base}\\ bond
\pr{princ}\\ $0$-bond \pr{obond}\\ bottom side of a band
\pr{bottom}\\ Burnside cells \pr{bercel}\\ Burnside-reduced word
\pr{brw}\\ Burnside $ra$-relations \pr{bbra}\\ $0$-cell
\pageref{ocell}  \\
chamber \pr{chamber}\\ cleaning rule \pr{cler2}\\ compatibility of
cells \pr{comp1}, \pr{compq2}\\ congruence extension property
(CEP) \pr{CEPr}\\
 connecting line of a band \pageref{cl}\\
connecting line of a bond \pr{clb}\\ connecting line of a
contiguity submap \pr{cls}\\ connecting rule \pr{cor1},
\pr{cor2}\\  contiguity arc \pr{car}\\ contiguity degree
\pr{cde}\\ contiguity submap \pr{csm}\\ converting $S$-rules into
sets of relations \pr{crr}\\ $\Omega$-coordinate of a word
\pr{oca}\\ copy of a word \pr{copy}\\ crossing bands \pr{cross}\\
cyclically $Y$-reduced element \pr{cyre}\\ cyclically $Y$-reduced
word \pr{crw}\\ diagram divided by a connecting line (band)
\pr{divides}\\ domain of a word in $\sss$ \pr{domain}\\ $0$-edge
\pageref{oedge}\\ $0$-element \pr{0el}\\ essential element
\pr{ess}\\ graded map \pageref{graded}\\ graded presentation
\pr{gp}\\ history of a computation \pr{history}\\ history of a
$\bkk$-band \pr{hb}\\ history of a trapezium \pr{ht}\\ hoop
\pr{hoop}\\ homomorphism $\vartheta_i$ \pr{vtheta}\\ homomorphism
$\upsilon$ \pr{upsilon}\\ inside diagram of an annulus
\pr{inside}\\ inverse rule \pr{invr}\\ inverse semigroup $P(\sss)$
\pr{pss}\\ Lowest Parameter Principle \pr{lpp}\\ length of a path
\pageref{lop}\\ $Y$-length \pr{yl}\\ $Y$-letter \pr{yt}\\
$0$-letter \pr{0let}\\ map \pageref{mapp}\\ A-map \pr{amp}\\
narrow bond \pr{nb}\\ normal trapezium \pr{normal}\\ $j$-pair
\pr{jp1}\\ $1/2$-pair \pr{12p}\\
period of rank $i$ \pr{period1}, \pr{period2}\\ $A$-periodic word
\pr{per}\\ phase decomposition \pr{phase}\\ positive edge (cell)
\pageref{pedge}\\ positive and negative $S$-rules \pr{pnr}\\
principal cell of a bond \pr{princ}\\ properties A1, A2, A3
\pr{A123}\\ properties (M1),(M2), (M3) \pr{m123}\\ properties
(T1), (T2), (T3) \pr{T123}\\ properties (TR1), (TR2) \pr{TR12}\\
properties (Z1),(Z2), (Z3) \pr{z123}\\ rank of a cell
\pageref{rank}\\ rank of a map \pageref{rankmap}\\ reduced
admissible word \pr{radm}\\ reduced band \pr{redd}\\ reduced path
\pr{redpath}\\ g-reduced diagram \pr{gred}\\ $Z(\sss,a)$-reduced
diagram \pr{zsa}\\ *-reduced word \pr{str}\\ $Y$-reduced word
\pr{yr}\\ $T$-relator \pr{trel}\\ $a$-relations \pr{aarel}\\
$k$-relations \pr{kkrel}\\ $r$-relations \pr{rrrel}\\
$ra$-relations \pr{rarel}\\  $S$-rule \pr{srule}\\ $S$-rule
applicable to an admissible word \pr{app}\\ $S$-rule locking a
sector \pr{lock}\\ $S$-rule left (right) active for a sector
\pr{lra}\\ section of a contour of a diagram \pr{sect}\\
$\varepsilon$-section \pr{epss}\\ $H$-section \pr{hsec}\\ sector
\pr{sector}\\ sector of the type (form) $[zz']$ \pr{z1z2}\\ simple
in rank $i$ word \pr{simple1}, \pr{simple2}\\ smooth section
\pr{smooth}\\ $\alpha$-series \pr{alpha}\\ side arc of a \ct
subdiagram \pr{sar}\\ sides of a trapezium \pr{base}\\ start and
end edges of a band \pr{start}\\
$0$-subgroup \pr{0sub}\\ top side of a band \pr{bottom}\\
trapezium \pr{tra}\\ $\psi$-trivial word \pr{ptriv}\\ type of a
map (diagram) \pageref{type}\\ \vk diagram \pageref{vk}\\ \vk
lemma \pageref{vkl}\\ vertices in the same phase \pr{sp}\\
$0$-word \pr{0word}\\ $H$-word \pr{hword}\\ working rule
\pr{wor1}, \pr{wor2}\\

\noindent $\aaa$ \pr{aa}\\ ${\cal A}$  \pr{aaaa}\\
 $|A|_Y$ \pr{ay}\\
  $A(z)$ \pr{az}\\
$\alpha, \beta, \gamma, \delta, \varepsilon, \zeta, \iota, h, n,
\bar\alpha, \bar\beta, \bar\gamma$ \pr{param}\\ $B(m,n)$
\pr{bmn}\\ $\bott(\bb)$ \pr{bottom}\\
$\partial(\Pi_1,\Gamma,\Pi_2)$ \pr{par}\\ $\partial(\Pi)$,
$|\partial(\Pi)|$ \pageref{perimeter}\\ $e_-$, $e_+$
\pageref{emin}\\ $\equiv$ \pr{equiv}\\
 $\Gamma\bigwedge\Pi$
\pr{wedge}\\ ${\bf G}(i), {\bf G}(\infty)$ \pr{gi}\\ ${\cal H}$
 \pr{ash}, \pr{ash1} \\ $H_a$ \pr{ha}\\ $\hnka$ \pr{hka}\\
 $\hra$ \pr{hra}\\
$H_{ra}(\infty)$ \pr{hrai}\\
 $H_{kra}$\pr{hkra}\\
$H_{kra}(\infty)$ \pr{hkrai}\\
 $\kk$ \pr{kk}\\
  $\kk(z)$ \pr{kz}\\
$\bkk$ \pr{bkk}\\
 $\Lambda$ \pr{lambda}\\ $\Lambda_i(0)$
\pr{lami}\\
 $L(\tau,z)$ \pr{lrtz}\\
 $L^a(\tau, z)$ \pr{lara}\\
$\oo(g)$ \pr{0g}\\
$\Omega$ \pr{omeg}\\
$(\Pi_1,\Gamma,\Pi_2)$ \pr{cde1}\\
 $\Pi \circ\Theta$ \pr{pct}\\
$R(\kappa(1))$ \pr{rkappa1}\\
$R(\tau,z)$ \pr{lrtz}\\
$R^a(\tau,z)$ \pr{lara}\\
$\rr$ \pr{RR}\\
$\rr(z)$ \pr{rrz}\\ ${\cal R}(\infty)$
\pr{rr}\\ ${\cal R}_i$ \pr{ri}\\
 ${\cal R}_{kra}$ \pr{rkra}\\
 ${\cal R}_{ra}$ \pr{rra}\\
$r(\Pi)$ \pageref{rank}\\
$r(q)$ \pr{rq}\\
$r(\tau,z)$ \pr{tia}\\
 ${\cal S}_i$ \pr{si}\\
  $\sss^+$ \pr{s+1},\pr{s+2}\\
$\sss_\omega^+$ \pr{s+o1}\\
 $\sim_i$ \pr{sim}\\
$\Theta *\Pi$ \pr{tsp}\\
$\tool$ \pr{tool}\\
 $\bt{z}$ \pr{btz}\\
 $\topp(\bb)$ \pr{bottom}\\
$U(\omega)$ \pr{uo}\\ $W\cdot h$ \pr{wh}\\ $\xxx$ \pr{xxx}\\
${\cal X}_i$ \pr{xi}\\ $z_-, z_+$ \pr{zmp}\\ $Z(\sss)$ \pr{zss}\\
$Z(\sss,\Lambda)$ \pr{zsl}\\ $Z(\sss,a)$ \pr{bra}\\ $\tz$
\pr{tz}\\ $\tza$ \pr{tza} \onecolumn
\begin{minipage}[t]{3 in}
\noindent Alexander Yu. Ol'shanskii\\ Department of Mathematics\\
Vanderbilt University \\ olsh@math.vanderbilt.edu\\
http://www.math.vanderbilt.edu/$\sim$olsh\\ and\\ Department of
Higher Algebra, MEHMAT\\
 Moscow State University\\
olshan@shabol.math.msu.su\\
\end{minipage}
\begin{minipage}[t]{3 in}
\noindent Mark V. Sapir\\ Department of Mathematics\\
Vanderbilt University\\
http://www.math.vanderbilt.edu/$\sim$msapir\\
\end{minipage}

\addtocontents{toc}{\contentsline {section}{\numberline {
}References \hbox {}}{\pageref{bibbb}}}

\addtocontents{toc}{\contentsline {section}{\numberline { }Subject
index \hbox {}}{\pageref{sind}}}

\end{document}